\documentclass[preprint,12pt]{elsarticle}
\usepackage{amsfonts,amsmath,amssymb}
\usepackage{amsthm}
\usepackage{tikz}
\usepackage[colorlinks=true, citecolor=black, linkcolor=black, urlcolor=black]{hyperref}

\newtheorem{thm}{Theorem}[section]
\newtheorem{lem}[thm]{Lemma}

\newtheorem{cor}[thm]{Corollary}
\newtheorem{prb}[thm]{Problem}

\newdefinition{df}{Definition}[section]
\newdefinition{rem}{Remark}[section]
\newdefinition{ex}{Example}[section]
\newdefinition{coex}{Counterexample}[section]
\newproof{pf}{Proof}
\newproof{pot}{Proof of Theorem}
\numberwithin{equation}{section}

\journal{
Journal of Functional Analysis
}

\begin{document}
\begin{frontmatter}


\title{Irreducibility of the Koopman representations for the group ${\rm GL}_0(2\infty,{\mathbb R})$
acting on  three infinite rows
\footnote{\bf To all fearless Ukrainians defending not only their country, but the whole civilization 
against putin's 
\href{https://en.wikipedia.org/wiki/Ruscism}{rashism}
}
}
%

\author{A.V.~Kosyak\corref{cor1}}
\ead{kosyak02@gmail.com}
\address{Institute of Mathematics, Ukrainian National Academy of Sciences,\\
3 Tereshchenkivs'ka Str., Kyiv, 01601, Ukraine}
\address{
London Institute for Mathematical Sciences,\\
21 Albemarle St, London W1S 4BS, UK
}

\cortext[cor1]{Corresponding author}

\author{P.~Moree}
\address{Max-Planck-Institut f\"ur Mathematik, Vivatsgasse 7, D-53111 Bonn, Germany}
\begin{abstract}
In [25] the first author started with the development of harmonic analysis on infinite-dimensional groups.  
In this article,  following these ideas, we construct an analogue of quasi-regular representations, when the group $G$ acts on
a $G-$space $X$ equipped with a quasi-invariant measure. For the group $G$ we
take
the inductive limit of the general linear
groups  ${\rm
GL}_0(2\infty,{\mathbb R})$ $= \varinjlim_{n}{\rm
GL}(2n\!-\!1,{\mathbb R})$,
acting on the space $X_m$  of $m$ rows, infinite in both
directions,   with Gaussian measure.  This measure is the infinite tensor product of
one-dimensional arbitrary Gaussian non-centered measures. 
We prove an irreducibility criterion for $m\!=\!3$. In 2019, the first
author [28] established a criterion for $m\!\le\! 2$. Our proof is in the same spirit,
but the details are far more involved.
\end{abstract}


\end{frontmatter}
\tableofcontents

\newpage
\section{Representations of the inductive limit of the general linear groups
${\rm GL}_0(2\infty,{\mathbb R})$
}
\subsection{Finite-dimensional case}
Consider the space 
$$
X_{m,n}=\Big\{x=\sum_{1\leq k\leq m}\,\,\sum_{-n\leq r\leq n} x_{kr}E_{kr},\,\,x_{kr}\in
{\mathbb R}\Big\},
$$
where $E_{kn},\,\,k,n\in {\mathbb Z}$ are infinite matrix unities, with the measure   (see \eqref{mu(b,a)})
$$
\mu_{(b,a)}^{m,n}(x)=\otimes_{k=1}^m\otimes_{r=-n}^{n}
\mu_{(b_{kr},a_{kr})}(x_{kr}).
$$
 Two groups act on the space $X_{m,n}$, namely ${\rm GL}(m,\mathbb R)$ from the left, and $\text{GL}(2n+1,\mathbb R)$ from the right, and their actions commute. 
Therefore, two von Neumann algebras ${\mathfrak A}_{1,n}$  and ${\mathfrak A}_{2,n} $ in the Hilbert space $L^2(X_{m,n},\mu_{(b,a)}^{m,n})$  generated 
respectively by the left and the right actions of the corresponding groups have the property that ${\mathfrak A}_{1,n}'\subseteq {\mathfrak A}_{2,n}$, where 
${\mathfrak A}'$ is a commutant of a von Neuman algebra ${\mathfrak A}$.
We study what happens as $n\!\to\infty$. In the limit we obtain some unitary representation $T^{R,\mu,m}$ (see (\ref{T^(R,mu,m)})) of the group  
$G:={\rm GL}_0(2\infty,{\mathbb R})\!=\!\varinjlim_{n,i^s}{\rm GL}(2n+1,{\mathbb R})$ acting from the right on $X_m$. In the generic case, the representation 
$T^{R,\mu,m}$ is reducible. Indeed, if there exists a non-trivial element $s\in \text{GL}(m,\mathbb R)$ such the  left action is {\it admissible} for 
the measure $\mu_{(b,a)}^m$, i.e., $(\mu_{(b,a)}^m)^{L_s}\!\sim\! \mu_{(b,a)}^m$  the operator $T^{L,\mu,m}_s$ naturally associated with the left action,
is well defined 
and $[T^{R,\mu,m}_t, T^{L,\mu,m}_s]\!=\!0$ for all $t\!\in\! G,\,s\!\in\! {\rm GL}(m,\mathbb R)$. 
We use notation $\mu^f(\Delta)=\mu\big(f^{-1}(\Delta)\big)$ for $f:X\to X$, where $\Delta$ is some measurable set in $X$.

{\it The main result} of the article is the following.
The representation  $T^{R,\mu,m}$ is irreducible, see Theorem~\ref{5.t.irr}
if and only if {\it no  left actions are admissible}, i.e.,
when $(\mu_{(b,a)}^m)^{L_s}\perp \mu_{(b,a)}^m$ for all $s\in \text{GL}(m,\mathbb R)\backslash \{e\}$.
This is again a manifestation of the {\it Ismagilov conjecture}, see \cite{Kos_B_09}.

Here, as in the case of the regular \cite{Kos90,Kos92} and   quasiregular \cite{Kos02.2,Kos02.3} representations of the group $B^{\mathbb N}_0$,
which is an inductive limit of upper-triangular real matrices,
we obtain the remarkable result that  the {\it irreducible representations} can be obtained as the {\it inductive limit of reducible representations\,}!
\subsection{Infinite-dimensional case}
Let us denote by ${\rm Mat}(2\infty,{\mathbb R})$ the space of all
real  matrices that are infinite in both directions:
\begin{equation}
\label{5.Mat(2inf,R)} {\rm Mat}(2\infty,{\mathbb R})=
\Big\{x=\sum_{ k,n\in{\mathbb Z}}x_{kn}E_{kn},\,\,x_{kn}\in
{\mathbb R}\Big\}.
\end{equation}

%
%
The group ${\rm GL}_0(2\infty,{\mathbb R})=\varinjlim_{n,i^s}{\rm
GL}(2n+1,{\mathbb R})$ is defined as the inductive limit of the
general linear groups $G_n={\rm GL}(2n+1,{\mathbb R})$ with
respect to the {\it symmetric embedding} $i^s$:
\begin{equation}
\label{N.i^s} 
G_n\ni x\mapsto
i^s_{n+1}(x)=x+E_{-(n+1),-(n+1)}+E_{n+1,n+1}\in %
G_{n+1}.
\end{equation}
For a fixed natural number $m$, consider a $G$-space $X_m$  as the
following subspace of the space ${\rm Mat}(2\infty,{\mathbb R})$:
\begin{equation}
\label{5.X_m} X_m=\Big\{x\in {\rm Mat}(2\infty,{\mathbb
R})\mid x=\sum_{k=1}^m\sum_{n\in{\mathbb
Z}}x_{kn}E_{kn}\Big\}.
\end{equation}
The  group ${\rm GL}_0(2\infty,{\mathbb R})$ acts from the right
on the space $X_m.$ Namely, the right action of the group ${\rm
GL}_0(2\infty,{\mathbb R})$ is correctly defined on the space
$X_m$ by the formula $R_t(x)=xt^{-1},\,\,t\in G,\,\,x\in X_m$. We
define a Gaussian non-centered product measure
$\mu:=\mu^m:=\mu_{(b,a)}^m$
on the space $X_m:$
\begin{equation}
\label{5.mu^m}
\mu_{(b,a)}^m(x)=\otimes_{k=1}^m\otimes_{n\in{\mathbb
Z}}\mu_{(b_{kn},a_{kn})}(x_{kn}),
\end{equation}
 where
\begin{equation}
\label{mu(b,a)} 
d\mu_{(b_{kn},a_{kn})}(x_{kn})=\sqrt{\frac{ b_{kn}}{\pi}}e^{-b_{kn}(x_{kn}-a_{kn})^2}
dx_{kn}
\end{equation}
and $b=(b_{kn})_{k,n},\,\,b_{kn}>0,\,a=(a_{kn})_{k,n},\,a_{kn}\in
{\mathbb R},\,1\leq k\leq m,\,n\in {\mathbb Z}.$ Define the unitary 
representation $T^{R,\mu,m}$ of the group ${\rm
GL}_0(2\infty,{\mathbb R})$ on the space $L^2(X_m,\mu^m_{(b,a)})$
by the formula:
\begin{equation}
\label{T^(R,mu,m)}
(T^{R,\mu,m}_tf)(x)=\big(d\mu_{(b,a)}^m(xt)/d\mu_{(b,a)}^m(x)\big)^{1/2}f(xt),\,\,f\in L^2(X_m,\mu^m_{(b,a)}).
\end{equation}
Obviously, the {\it centralizer} $Z_{{\rm Aut}(X_m)}(R(G))\subset
{\rm Aut}(X_m)$ contains  the group $L({\rm GL}(m,{\mathbb R}))$,
i.e., the image of the group ${\rm GL}(m,{\mathbb R})$ with
respect to the left action $L:{\rm GL}(m,{\mathbb R})\rightarrow
{\rm Aut}(X_m),\,L_s(x)\!=\!sx,\,s\in {\rm GL}(m,{\mathbb
R}),\,x\in X_m.$  We prove the following theorem.  
\begin{thm}
\label{5.t.irr} The representation $T^{R,\mu,m}\!:\!{\rm
GL}_0(2\infty,{\mathbb R})\!\rightarrow\! U\Big(L^2(X_m,\mu^m_{(b,a)})\Big)$ is
irreducible, for $m=3$, if and only if  
\begin{eqnarray*}
(i)&&(\mu^m_{(b,a)})^{L_{s}}\perp
\mu^m_{(b,a)}\quad\text {for all}\quad s\in {\rm GL}(m,{\mathbb
R})\backslash\{e\};\\
(ii)&&\text{the measure} \quad\mu^m_{(b,a)}\quad \text{is $G$-{\rm ergodic}}.
\end{eqnarray*}
\end{thm}
In {\rm \cite{Kos_B_09,KosJFA17}}
this result was proved for  $m\leq 2$ .
Note that conditions (i) and (ii) are necessary conditions for irreducibility.
\begin{rem}
Any Gaussian product-measure $\mu_{(b,a)}^m$ on $X_m$ is ${\rm
GL}_0(2\infty,{\mathbb R})$-right-ergodic \cite[\S 3, Corollary 1]{ShFDT67}, see Definition~\ref{d.erg}. For non-product-measures
this is not true in general.
\end{rem}
In order to study the condition $(\mu^m_{(b,a)})^{L_{t}}\perp
\mu^m_{(b,a)}$ for $t\in {\rm GL}(m,{\mathbb R})\setminus \{e\}$ set
\begin{equation}
\label{5.X_n(t)}
t=(t_{rs})_{r,s=1}^m\in{\rm GL}(m,{\mathbb R}), \,\,\,B_n={\rm
diag}(b_{1n},b_{2n},...,b_{mn}),\,\,\,X_n(t)=B_n^{1/2}tB_n^{-1/2}.
\end{equation}
Let  $M^{i_1i_2...i_r}_{j_1j_2...j_r}(t)$ be the {\it minors} of
the matrix $t$ with $i_1,i_2,...,i_r$ rows and $j_1,j_2,...,j_r$
columns, $1\leq r\leq m.$ Let $\delta_{rs}$ be the Kronecker
symbols.
\begin{lem}[\cite{Kos_B_09}, Lemma~10.2.3; \cite{KosJFA17}, Lemma~2.2]
\label{perp1} For the measures $\mu_{(b,a)}^m$, with  $m$ a natural number,
the relation
\begin{equation*}
 (\mu_{(b,a)}^m)^{L_t}\perp\mu_{(b,a)}^m\,\,\,\text{for all}\,\,\, t\in
 {\rm GL}(m,{\mathbb R})\backslash\{e\}\,\,
\end{equation*}
 holds if and only if
 \begin{eqnarray*}
 &&
\prod_{n\in {\mathbb Z}}\frac{1}{2^m\vert{\rm det}\,\,t\vert}\, {\rm
det}\left(I\!+\!X_n^*(t)X_n(t)\right)\!+\! \sum_{n\in {\mathbb
Z}}\sum_{r=1}^m
b_{rn}\Big(\!\sum_{s=1}^m(t_{rs}\!-\!\delta_{rs})a_{sn}\!\Big)^2
\!=\!\infty,\\
&&
\text{where}\\
&&
{\rm det}\Big(I+X_n^*(t)X_n(t)\Big)= 1+\sum_{r=1}^{m} \sum_{1\leq i_1<i_2<...<i_r\leq m\atop 1\leq
j_1<j_2<...<j_r\leq m}
\left(M^{i_1i_2...i_r}_{j_1j_2...j_r}(X_n(t))\right)^2.
\end{eqnarray*}
\end{lem}
For the convenience of the reader this lemma is proved in Section~\ref{sec.Gauss-meas}.
\begin{rem}
\label{r.irr-idea-G} (The idea of the proof of
irreducibility.) Let us denote by ${\mathfrak A}^m$ the  {\it von
Neumann algebra} generated by the representation $T^{R,\mu,m}$, i.e., $
{\mathfrak A^m}=(T^{R,\mu,m}_t\mid t\in G)''. $ For
$\alpha\!=\!(\alpha_k)\!\in\!\{0,1\}^m$ define the von Neumann algebra
$L^\infty_\alpha(X_m,\mu^m)$ as follows:
$$
L^\infty_\alpha(X_m,\mu^m)\!=\!\Big(\exp(itB^\alpha_{kn})\mid
1\leq k\leq m,\,\,t\in {\mathbb R},\,\,n\in {\mathbb
Z}\Big)'',\,\,
$$
where $
B^\alpha_{kn}\!=\!\left\{\begin{array}{rcc}
 x_{kn},&\text{if}&\alpha_k=0\\
 i^{-1}D_{kn},&\text{if}&\alpha_k=1
\end{array}\right.$ and $D_{kn}=\partial/\partial x_{kn}-b_{kn}(x_{kn}-a_{kn})$.
\vskip 0.3 cm
{\bf The proof of the  irreducibility is based on  four
facts}:
\par 1) we can approximate by generators
$A_{kn}=A_{kn}^{R,m}=\frac{d}{dt}T^{R,\mu,m}_{I+tE_{kn}}\vert_{t=0}$
the set of operators $(B^\alpha_{kn})_{k=1}^m,\,n\!\in\!{\mathbb
Z}$ {\it for some} $\alpha\!\in\!\{0,1\}^m$ depending on the measure
$\mu^m$  using the orthogonality condition
$(\mu^m)^{L_{s}}\perp \mu^m$ for all $s\in {\rm GL}(m,{\mathbb
R})\backslash{\{e\}}$,
\par 2)
it is sufficient to  verify the approximation  only on the  cyclic vector
${\bf 1}(x)\!\equiv\! 1$, since the representation $T^{R,\mu,m}$ is {\it cyclic}, 
\par 3) the subalgebra $L^\infty_\alpha(X_m,\mu^m)$
is a {\it maximal abelian subalgebra} in ${\mathfrak A}^m$,
\par 4)  the measure $\mu^m$ is $G$-ergodic. 

Here the 
{\it generators}
$A_{kn}$ are given by the formulas:
$$
A_{kn}\!=\!\sum_{r=1}^{m}x_{rk}D_{rn},\quad k,n \in
{\mathbb Z},\quad\text{where}\quad D_{kn}=\partial/\partial x_{kn}-b_{kn}(x_{kn}-a_{kn}).
$$
\end{rem}
\begin{rem}
\label{r.shem} 
{\it Scheme of the proof.} We prove the irreducibility as follows
%
\begin{eqnarray}
\label{shem}
&&\left(\mu^{L_s}\perp \mu\,\,\, \text{for all} \,\,\, 
s\in {\rm GL}(3,\mathbb R)\setminus \{e\} \right)\Leftrightarrow
\left(\begin{smallmatrix}
\text{criteria}\\
\text{of}\\
\text{orthogonality}
\end{smallmatrix}\right)
\& \\
\nonumber
&&
\left(\begin{smallmatrix}
\text{Lemma~\ref{l.min=proj.3}}\\
\text{about}\\
\text{three vectors}\, f,g,h\not\in l_2
\end{smallmatrix}\right)
\Rightarrow
\left(\begin{smallmatrix}
\text{some of}        &\Delta^{(1)},&\Delta_1\\
\text{the expressions}&\Delta^{(2)},&\Delta_2\\
\text{are divergent:} &\Delta^{(3)},&\Delta_3
\end{smallmatrix}
\right)
\Rightarrow
\text{irreducibility}, \\
\label{Si,si}
&&
\text{where}\quad \Delta^{(i)}:=\Delta(Y_i^{(i)},Y_j^{(i)},Y_k^{(i)}),\quad 
\Delta_i:=\Delta(Y_i,Y_j,Y_k), 
\end{eqnarray}
$\Delta(f,g,h)$ is defined by \eqref{Delta(f,g,h)}, 
and $\{i,j,k\}$ is a cyclic permutation 
of $\{1,2,3\}$, see for details
Lemmas~\ref{d1.3}--\ref{d3.3}, Lemmas~\ref{x1x1.3}--\ref{x3x3.3} and Lemmas~\ref{exp(x)A3}-\ref{l.Re-Im-exp.3}.
\end{rem}
\begin{rem}
\label{r.key-lemma2}
The fact that the conditions $(\mu^3)^{L_{t}}\!\perp\! \mu^3$ for all $t\in {\rm GL}(3,{\mathbb
R})\backslash{\{e\}}$ imply the possibility of the approximation of $x_{kn}$ and $D_{kn}$ 
by combinations of generators
is based on some
completely i{\it ndependent statement} about three infinite vectors $f,g,h\not \in l_2$
such that 
\begin{equation}
 \label{not-in-l_2}
C_1 f+C_2 g+C_3 h\not \in l_2\quad \text{for arbitrary}\quad (C_1,C_2,C_3)\in {\mathbb R}^3\setminus\{0\},
\end{equation}
see Lemma~\ref{l.min=proj} for $m=2$ \cite{KosJFA17}
and Lemma~\ref{l.min=proj.3} for $m=3$. 
A similar result for general $m$ is studied in \cite{Kos-m-Arx23}.
These  lemmas are the {\it key ingredients}
of the proof of the irreducibility of the representation.
\end{rem}
\begin{rem}
\label{r.key-lemma}
Note that in the case of the ``nilpotent group'' $B_0^{\mathbb N}$ and the infinite  product of {\it arbitrary}  Gaussian measures on ${\mathbb R}^m$ 
(see \cite{KosAlb06J}) the  proof of the {\it irreducibility} is also  based on another  completely {\it independent statement} namely,  
{\it the Hadamard -- Fischer inequality}, see Lemma~\ref{l.det(det)1}.
\begin{lem}[Hadamard -- Fischer inequality
\cite{HornJon89,HornJon91}]
\label{l.det(det)1}  For any
positive definite matrix $C\in{\rm Mat}(m,{\mathbb
R}),\,\,m\in{\mathbb N}$ and  any two subsets $\alpha$ and $\beta$
with $\emptyset\subseteq\alpha,\,\beta\subseteq\{1,...,m\}$ the
following inequality holds:
\begin{equation}
\label{det(det)1}
 \left|\begin{array}{cc}
M(\alpha)&M(\alpha\bigcap\beta)\\
M(\alpha\bigcup\beta)&M(\beta)
\end{array}\right|=
 \left|\begin{array}{cc}
A(\hat{\alpha})&A(\hat{\alpha}\bigcup\hat{\beta})\\
A(\hat{\alpha}\bigcap\hat{\beta})&A(\hat{\beta})
\end{array}\right|
\geq 0,
\end{equation}
where $M(\alpha)=M^\alpha_\alpha(C),\,\,
A(\alpha)=A^\alpha_\alpha(C)$ and
$\hat{\alpha}=\{1,...,m\}\setminus\alpha.$
\end{lem}
For the details see \cite[p.573]{HornJon89} and \cite[Chapter
2.5, problem 36]{HornJon91}.
In \cite{KosAlb06J} the conditions of orthogonality $\mu^{L_{t}}\perp \mu$ with respect to the left action of the group $B(m,{\mathbb R})$ on $X^m$ were expressed as the divergence of some series, $S^L_{kn}(\mu)=\infty,\,\,1\leq k<n\leq m$.
Conditions on the measure $\mu$ for the variables $x_{kn}$ to be approximated
by combinations of generators $A_{pq}$ were expressed in terms of the divergence of another series $\Sigma_{kn}$. 
The proof of the fact that conditions $S^L_{kn}(\mu)=\infty,\,\,1\leq\! k\!<\!n\!\leq m$ imply the conditions
$\Sigma_{kn}=\infty,\,\,1\leq k<n\leq m$ was based on the
Hadamard -- Fischer inequality.
\end{rem}
\section{Some orthogonality problem in  measure theory}
\label{sec.orth.m=3}
\subsection{General setting}
Our aim now is to find the minimal generating set of conditions of the orthogonality
$(\mu^m_{(b,a)})^{L_{t}}\!\perp\! \mu^m_{(b,a)}$ for all
$t\in {\rm GL}(m,{\mathbb R})\setminus\{e\}$. To be more precise, consider the following more general situation.
Let $\alpha:G\rightarrow {\rm
Aut}(X)$  be a {\it measurable action} of a group $G$  on a
measurable space $(X,
\mu)$ with the following property:
 $\mu^{\alpha_t}\perp\mu$ for all $t\in G\setminus\{e\}$. 
Define a {\it generating subset} $G^\perp(\mu)$ in the group $G$ 
as follows:
\begin{equation}
\label{min-perp}
\text{if}\quad \mu^{\alpha_t}\perp\mu\,\,\,\text{for all}\,\,\, t\in G^\perp(\mu),
\quad\text{then}\quad \mu^{\alpha_t}\perp\mu\,\,\,\text{for all}\,\,\,
t\in G\setminus\{e\}.
\end{equation}
\begin{prb}
\label{pr.min-perp}
Find a {\rm minimal generating subset} $G^\perp_0(\mu)$ 
satisfying {\rm (\ref{min-perp})}.
\end{prb}
\begin{df}
\label{d.erg}
Recall that the probability measure $\mu$ on some $G$-space $X$ is called {\it ergodic} if any function
$f\in L^1(X,\mu)$ with property $f(\alpha_t(x))=f(x)\,\,{\rm mod}\,\mu$ is constant.
\end{df}

\subsection{Orthogonality criteria  $\mu^{L_t}\perp \mu$ for $t\in {\rm GL}(2,\mathbb R)\setminus\{e\}$}
\begin{rem}
\label{r.m=2.crit-orthog} 
By Lemma~4.1 proved in  \cite{KosJFA17} or Lemma~10.4.1 in \cite{Kos_B_09} for $m=2$ we conclude that
the minimal generating set  $G^\perp_0(\mu)={\rm GL}(2,\mathbb R)^\perp_0(\mu)$  (see Problem~\eqref{pr.min-perp})
is reduced to the following subgroups, families and elements:
\begin{equation}
\label{m=2.crit-orthog.1} 
\exp(tE_{12})=I+tE_{12}\!=\! \left(\begin{array}{cc}
 1&t\\
 0&1
\end{array}\right),\,\,
\exp(tE_{21})\!=\!I+tE_{21}=\left(\begin{array}{cc}
 1&0\\
 t&1
\end{array}\right),
\end{equation}
\begin{equation}
\label{m=2.crit-orthog.2} 
\exp(tE_{12})P_1= \left(\begin{array}{cc}
 -1&t\\
  0&1
\end{array}\right),\quad
\exp(tE_{21})P_2= \left(\begin{array}{cc}
 1&0\\
 t&-1
\end{array}\right),
\end{equation}
\begin{equation}
\label{m=2.crit-orthog.3}
\tau_{-}(\phi,s)= \left(
\begin{array}{cc}
\cos\phi&s^2\sin\phi\\
s^{-2}\sin\phi&-\cos\phi
\end{array}
\right)=
D_2(s)
\left(\begin{array}{cc}
\cos\phi&-\sin\phi\\
\sin\phi&\cos\phi
\end{array}\right)
D^{-1}_2(s)P_2. 
\end{equation}
The families \eqref{m=2.crit-orthog.1}, \eqref{m=2.crit-orthog.2} are one-parameter, 
the family \eqref{m=2.crit-orthog.3} is two-parameter.
All elements are of order $2$ except the elements in  subgroups given in \eqref{m=2.crit-orthog.1}!
It suffices to verify the conditions  \eqref{m=2.crit-orthog.1}  only for some $t\in \mathbb R\setminus \{0\}.$ 
The family $\tau_{-}(\phi,s)$,  
actually,  coincides with
$D_2(s){\rm O}(2)D^{-1}_2(s)P_2$, where $ D_2(s)={\rm diag}(s,s^{-1})$.
All points $t$ in \eqref{m=2.crit-orthog.2} and all points $(\phi,s)$ in \eqref{m=2.crit-orthog.3} are essential, i.e., we can not remove any single point.
\end{rem}
\begin{rem}
\label{r.5.5}
We note \cite[Chapter V, \S 8 Problems, 2, p. 147]{Knap86} that every element of   ${\rm SL}(2,{\mathbb R})$ is
conjugate to at least one matrix of the form
$$
\left(\begin{array}{cc}
 a&0\\
  0&a^{-1}
\end{array}\right),\,\,a\not=0,\,\,\left(\begin{array}{cc}
 1&t\\
  0&1
\end{array}\right),\,\,\left(\begin{array}{cc}
 -1&t\\
  0&-1
\end{array}\right),\,\,
\left(\begin{array}{cc}
\cos\phi&\sin\phi\\
-\sin\phi&\cos\phi
\end{array}
\right).
$$
\end{rem}
\begin{rem}
\label{r.a^2=1} 
Some elements  $a=gP_r$ in the set $G^\perp_0(\mu)={\rm GL}(2,\mathbb R)^\perp_0(\mu)$ are of order $2$ 
(see Remark~\ref{r.m=2.crit-orthog}): 
\begin{equation}
\label{a^2=1} 
a^2=(gP_r)^2=1.
\end{equation}
This follows from the relation
\begin{equation}
\label{PgP=g^{-1}} 
P_rgP_r=g^{-1}.
\end{equation}
\end{rem}
To see this we note that 
if \eqref{PgP=g^{-1}} holds, then we get \eqref{a^2=1}:
\begin{equation*}
a^2=(gP_r)^2=gP_rgP_r=gg^{-1}=1. 
\end{equation*}
For example, for  $g=\exp(tE_{12})$ we get $P_1gP_1=g^{-1}$, for  $g=\exp(tE_{21})$ we get $P_2gP_2=g^{-1}$
and for $g=\tau(\phi,s)$ we get $P_2gP_2=g^{-1}$.
See \eqref{m=2.crit-orthog.2}  and \eqref{m=2.crit-orthog.3} for details, where
\begin{equation*}
\tau_{-}(\phi,s)= \left(
\begin{array}{cc}
\cos\phi&s^2\sin\phi\\
s^{-2}\sin\phi&-\cos\phi
\end{array}
\right),\quad
\tau(\phi,s)= \left(
\begin{array}{cc}
\cos\phi&-s^2\sin\phi\\
s^{-2}\sin\phi&\cos\phi
\end{array}
\right).
\end{equation*}
We recall some useful lemmas from \cite{KosJFA17}.
\begin{lem}
\label{l7.5}
 For $t=\left(\begin{smallmatrix}
t_{11}& t_{12}\\
t_{21}& t_{22}
  \end{smallmatrix}
  \right)
 \in {\rm GL}(2,{\mathbb R})\setminus\{e\}$ we have, if ${\rm det}\,t>0$,
$$
(\mu_{(b,0)}^2)^{L_t}\perp\mu_{(b,0)}^2\quad\Leftrightarrow\quad
$$
\begin{equation}
\sum_{n\in{\mathbb Z}}\Big[ (1-\mid{\rm
det}\,t\mid)^2+(t_{11}-t_{22})^2+\Big(t_{12}\sqrt{\frac{b_{1n}}{b_{2n}}}+
t_{21}\sqrt{\frac{b_{2n}}{b_{1n}}}\Big)^2\Big]=\infty.
\end{equation}
If ${\rm det}\,t<0$ we have
$$
(\mu_{(b,0)}^2)^{L_t}\perp\mu_{(b,0)}^2\quad\Leftrightarrow\quad
$$
\begin{equation}
\sum_{n\in{\mathbb Z}} \Big[ (1-\mid{\rm det}\,t\mid)^2 +
(t_{11}+t_{22})^2+\Big(t_{12}\sqrt{\frac{b_{1n}}{b_{2n}}}-
t_{21}\sqrt{\frac{b_{2n}}{b_{1n}}}\Big)^2\Big]=\infty.
\end{equation}
\end{lem}
\begin{lem}
\label{l7.6}
 For $t\in {\rm GL}(2,{\mathbb R})\setminus\{e\}$ we have
\begin{equation*}
 (\mu_{(b,a)}^2)^{L_t}\perp\mu_{(b,a)}^2\quad{\it
if}\quad\vert {\rm det}\,t\vert\not=1.
\end{equation*}
If ${\rm det}\,t=1$, we have
$$
(\mu_{(b,a)}^2)^{L_t}\perp\mu_{(b,a)}^2\quad\Leftrightarrow\quad\Sigma^+(t):=\Sigma_1^+(t)+
\Sigma_2(t)=\infty.
$$
If ${\rm det}\,t=-1$, we have
$$
(\mu_{(b,a)}^2)^{L_t}\perp\mu_{(b,a)}^2\quad\Leftrightarrow\quad\Sigma^{-}(t):=\Sigma_1^{-}(t)+
\Sigma_2(t)=\infty,
$$
where
$$
\Sigma_1^+(t)=\sum_{n\in{\mathbb
Z}}\Big[(t_{11}-t_{22})^2+\Big(t_{12}\sqrt{\frac{b_{1n}}{b_{2n}}}+
t_{21}\sqrt{\frac{b_{2n}}{b_{1n}}}\Big)^2\Big],
$$
$$
\Sigma_1^{-}(t)=\sum_{n\in{\mathbb
Z}}\Big[(t_{11}+t_{22})^2+\Big(t_{12}\sqrt{\frac{b_{1n}}{b_{2n}}}-
t_{21}\sqrt{\frac{b_{2n}}{b_{1n}}}\Big)^2\Big],
$$
\begin{equation}
\label{Sigma_2}
\Sigma_2(t^{-1})=\sum_{n\in{\mathbb
Z}}\Big(b_{1n}\big[(t_{11}-1)a_{1n}+t_{12}a_{2n}\big]^2+b_{2n}\big[t_{21}a_{1n}+(t_{22}-1)a_{2n}\big]^2\Big).
\end{equation}
\end{lem}
\begin{rem}
\label{r.pm-SL(2)}
By Lemma~\ref{l7.6} we have
\begin{equation*}
 (\mu_{(b,a)}^2)^{L_t}\perp\mu_{(b,a)}^2\quad
\text{ for}\quad t\in {\rm GL}(2,{\mathbb R})\setminus\{e\}
\end{equation*}
if and only if this holds  for two subsets  
of the group 
$\pm{\rm SL}(2,\mathbb R)$
defined as follows:
\begin{eqnarray}
\label{G^+_2.d} 
G^+_2 
&=&\{t\in {\rm SL}(2,{\mathbb R})\mid t_{11}=A^1_1(t)\},\\
\label{G^-_2.d}
G^-_2
&=&\{t\in -{\rm SL}(2,{\mathbb R})\mid t_{11}=-A^1_1(t)\}.
\end{eqnarray}
The set $G^+_2$
is reduced to  two  families of one-parameter subgroups \eqref{m=2.crit-orthog.1}.
The set  $G^-_2$
is reduced to  the one-parameter family \eqref{m=2.crit-orthog.2}, the reflections of 
\eqref{m=2.crit-orthog.1} by $P_2$, and two parameter family \eqref{m=2.crit-orthog.3} of elements 
from $D_2(s){\rm O}(2)D^{-1}_2(s)P_2$.
\end{rem}

\begin{lem}
If $t\in  G^+_2$
we have
$$
(\mu_{(b,a)}^2)^{L_t}\perp\mu_{(b,a)}^2\quad\Leftrightarrow\quad\Sigma^+(t):=\Sigma_1^+(t)+
\Sigma_2(t)=\infty.
$$
If $t\in G^-_2$
we have
$$
(\mu_{(b,a)}^2)^{L_t}\perp\mu_{(b,a)}^2\quad\Leftrightarrow\quad\Sigma^{-}(t):=\Sigma_1^{-}(t)+
\Sigma_2(t)=\infty,
$$
where  $\Sigma_2(t^{-1})$ is defined by \eqref{Sigma_2} and
\begin{eqnarray}
\label{Sigma_1^+(t).m=2}
 &&\Sigma_1^+(t)=\sum_{n\in{\mathbb
Z}}\Big(t_{12}\sqrt{\frac{b_{1n}}{b_{2n}}}+
t_{21}\sqrt{\frac{b_{2n}}{b_{1n}}}\Big)^2,\\
\label{Sigma_1^{-}(t).m=2}
&&\Sigma_1^{-}(t)=\sum_{n\in{\mathbb
Z}}\Big(t_{12}\sqrt{\frac{b_{1n}}{b_{2n}}}-
t_{21}\sqrt{\frac{b_{2n}}{b_{1n}}}\Big)^2.
\end{eqnarray}
\end{lem}
The conditions of 
orthogonality with respect to elements defined by
\eqref{m=2.crit-orthog.1}--\eqref{m=2.crit-orthog.3} 
are transformed in the divergence of the following series:
\begin{equation}
\label{S^L_(kn)}
S^L_{12}(\mu)=\sum_{n\in {\mathbb Z}}
\frac{b_{1n}}{2}\left(\frac{1}{2b_{2n}}+a_{2n}^2\right),\,\,
S^L_{21}(\mu)=\sum_{n\in {\mathbb Z}}
\frac{b_{2n}}{2}\left(\frac{1}{2b_{1n}}+a_{1n}^2\right),
\end{equation}
\begin{equation}
\label{S^(L,-)_(kn)}
S^{L,-}_{kn}(\mu,t)=\frac{t^2}{4}\sum_{m\in{\mathbb
Z}}\frac{b_{km}}{b_{nm}}+ \sum_{m\in{\mathbb
Z}}\frac{b_{km}}{2}(-2a_{km}+ta_{nm})^2,\,\,t\in\mathbb R,
\end{equation}
\begin{equation}
 \label{(phi,s)}
\Sigma^-_{12}\big(\tau_{-}(\phi,s)\big)=\sin^2\phi\,\Sigma_1(s)+\Sigma^-_2\big(\tau_{-}(\phi,s)\big),\,\,\phi\in[0,2\pi),\,\,s>0,
\end{equation}
\begin{equation}
\label{sigma1(s)}
\text{where}\quad\Sigma_1(s):=\!\sum_{n\in {\mathbb Z}}
\Big(s^2\sqrt{\frac{b_{1n}}{b_{2n}}}\!-\!s^{-2}\sqrt{\frac{b_{2n}}{b_{1n}}}\Big)^2,
\end{equation}
{\small
\begin{equation}
\label{sigma2(s)} \Sigma^-_2\big(\tau_{-}(\phi,s)\big)\!:=\!\sum_{n\in
{\mathbb Z}}
\Big(4b_{1n}\sin^2\frac{\phi}{2}\!+\!4s^{-4}b_{2n}\cos^2\frac{\phi}{2}\Big)
\Big(a_{1n}\sin\frac{\phi}{2}\!-\!s^2a_{2n}\cos\frac{\phi}{2}\Big)^2.
\end{equation}
}

Recall Remark 4.3 from \cite{KosJFA17}.
\begin{rem}
\label{perp2-1} The following three conditions are equivalent:
\begin{eqnarray*}
(i)&\quad\mu^{L_{\tau_-(\phi,s)}}\perp\mu,&\phi\in[0,2\pi),\,\,s>0,\\
(ii)&\Sigma_{12}^-\big(\tau_{-}(\phi,\!s)\big)\!=\!\sin^2\phi\,\Sigma_1(s)\!+\!\Sigma^-_2\big(\tau_{-}(\phi,s)\big)\!=\!
\infty,&\phi\in[0,2\pi),\,\,s>0,\\
(iii)&\quad\Sigma_1(s)+\Sigma_2(C_1,C_2)=\infty,&s\!>\!0,\,(C_1,C_2)\!\in\!{\mathbb R}^2\!\setminus\!\{0\},
\end{eqnarray*}
where $\Sigma_1(s)$ is defined by (\ref{sigma1(s)}) and
\begin{equation}
\label{sigma_2(C,C)} 
\Sigma_2(C_1,C_2):=
\sum_{n\in {\mathbb Z}}(C_1^2b_{1n}+C_2^2b_{2n})(C_1a_{1n}+C_2a_{2n})^2.
\end{equation}
\end{rem}
\subsection{Equivalent series and equivalent sequences}
There is an extensive 
theory of convergent and divergent series. In our context we are only interested in when a series with positive coefficients is divergent or convergent.  
\begin{df}
\label{df.two-sim-ser} We say  that two series $\sum_{n\in\mathbb N}a_n$ and  $\sum_{n\in\mathbb N}b_n$ with positive $a_n,\,b_n$ are {\it equivalent} if  they are divergent or convergent simultaneously. We will write
$\sum_{n\in\mathbb N}a_n\sim 
\sum_{n\in\mathbb N}b_n$. 
We say that two sequences $(a_n)_{n\in\mathbb N}$ and $(b_n)_{n\in\mathbb N}$ are {\it equivalent} if for some 
$C_1,C_2>0$ we have $C_1b_n\leq a_n\leq C_2b_n$ for all $n\in \mathbb N$. We will use the same notation $a_n\sim b_n$.
\end{df}
\begin{lem}

\label{l.two-ser-alpha}
Let $1+c_n>0$   for all $n\in \mathbb Z$. Then two series 
are equivalent:
\begin{equation}
\label{two-ser-alpha} 
\Sigma_1
:=\sum_{n\in \mathbb Z}
\frac{ c_n^2}{1+c_n},\quad \Sigma_2
:=\sum_{n\in \mathbb Z}c_n^2.
\end{equation}
\end{lem}
\begin{pf}
Fix some $\varepsilon\in (0,1)$ and a big $N$. We have three cases:\\
(a) $1+c_n\in (\varepsilon,N)$,\\
(b) for infinite subset $\mathbb Z_1$ we have $\lim_{n\in \mathbb Z_1}c_n=\infty$,\\
(c)   for infinite subset $\mathbb Z_1$ we have $\lim_{n\in \mathbb Z_1}(1+c_n)=0$.

In the case (a) we have
\begin{equation}
\label{(a)}
\frac{1}{N}\sum_{n\in \mathbb Z} c_n^2<
 \sum_{n\in \mathbb Z}
\frac{ c_n^2}{1+c_n}<\frac{1}{\varepsilon}\sum_{n\in \mathbb Z}c_n^2,
\end{equation}
In the case (b) and (c) both series are divergent.
\qed\end{pf}

We will make systematic use of the following statement.
\begin{rem}[\cite{Kos_B_09}]
Let $a_n,b_n>0$ for all $n\in \mathbb N$. The following two series  are equivalent
\label{r.two-sim-ser}
\begin{equation}
 \label{two-sim-ser}
 \sum_{n\in\mathbb N}\frac{a_n}{a_n+b_n}\sim \sum_{n\in\mathbb N}\frac{a_n}{b_n}.
\end{equation}
\end{rem}

\subsection{Orthogonality criteria  $\mu^{L_t}\perp \mu$ for $t\in {\rm GL}(3,\mathbb R)\setminus\{e\}$}
For $m=2$ and ${\rm det}\,t>0$ we have,
here $H_{m,n}(t)$ is defined by \eqref{Hel1-sl}
{\small
\begin{eqnarray*}
&& 
2^2\mid{\rm det}\,t\mid\big(H_{2,n}^{-2}(t)-1\big)\!=\!\Big[(1\!-\!\mid{\rm
det}\,t\mid)^2\!+\!
(t_{11}\!-\!t_{22})^2\!+\!\Big(t_{12}\sqrt{\frac{b_{1n}}{b_{2n}}}\!+\!
t_{21}\sqrt{\frac{b_{2n}}{b_{1n}}}\Big)^2\Big]\!=\\
&&
\left[\big(M^\emptyset_\emptyset(X(t))\!-\!A^\emptyset_\emptyset(X(t))\big)^2+
\big(M^1_1(X(t))\!-\!A^1_1(X(t))\big)^2+\big(M^1_2(X(t))\!-\!A^1_2(X(t))\big)^2
\right].
\end{eqnarray*}
}
For $m=3$  using (\ref{5.X_n(t)}) we have $X(t)=B^{1/2}tB^{-1/2}$, hence
$$
X(t)=
 \left(\begin{array}{ccc}
 b_{1n}&0&0\\
 0&b_{2n}&0\\
 0&0&b_{3n}
\end{array}\right)^{1/2}\!\!
 \left(\begin{array}{ccc}
 t_{11}&t_{12}&t_{13}\\
 t_{21}&t_{22}&t_{23}\\
 t_{31}&t_{32}&t_{33}\\
\end{array}\right)
 \left(\begin{array}{ccc}
b_{1n}&0&0\\
 0&b_{2n}&0\\
 0&0&b_{3n}
\end{array}\right)^{-1/2}=\!\!
$$
$$
 \left(\begin{array}{ccc}
 t_{11}&\sqrt{ \frac{b_{1n}}{ b_{2n}}}t_{12}&\sqrt{ \frac{b_{1n}}{b_{3n}}}t_{13}\\
 \sqrt{ \frac{b_{2n}}{ b_{1n}}}t_{21}&t_{22}&\sqrt{ \frac{b_{2n}}{ b_{3n}}}t_{23}\\
\sqrt{ \frac{b_{3n}}{ b_{1n}}}t_{31}& \sqrt{ \frac{b_{3n}}{ b_{2n}}}t_{32}&t_{33}
\end{array}\right).
$$
Therefore, using  
\eqref{G.detC-LI.2}
and the fact that $X=X^*(t)X(t)$ we obtain
\begin{eqnarray*}
 &&2^3\mid{\rm det}\,t\mid H_{3,n}^{-2}(t)=\Big(1+\mid{\rm
det}\,t\mid^2+t_{11}^2
+\frac{b_{1n}}{b_{2n}}t_{12}^2+\frac{b_{1n}}{b_{3n}}t_{13}^2+\frac{b_{2n}}{b_{1n}}t_{21}^2\\
&&+t_{22}^2+\frac{b_{2n}}{b_{3n}}t_{23}^2+
\frac{b_{3n}}{b_{1n}}t_{31}^2+\frac{b_{3n}}{b_{2n}}t_{32}^2+t_{33}^2+ 
(M^{12}_{12}(t))^2+\frac{b_{2n}}{b_{3n}}(M^{12}_{13}(t))^2+
\end{eqnarray*}
\begin{eqnarray*}
&&
\frac{b_{1n}}{b_{3n}}(M^{12}_{23}(t))^2+\frac{b_{3n}}{b_{2n}}(M^{13}_{12}(t))^2+(M^{13}_{13}(t))^2+\frac{b_{1n}}{b_{2n}}(M^{13}_{23}(t))^2\\
&&
+\frac{b_{3n}}{b_{1n}}(M^{23}_{12}(t))^2+\frac{b_{2n}}{b_{2n}}(M^{23}_{13}(t))^2+(M^{23}_{23}(t))^2
\Big)\\
&&=
1+\mid{\rm det}\,t\mid^2+ \sum_{1\leq i\leq j\leq 3}\Big[\left(t^i_j\sqrt{\frac{b_{in}}{b_{jn}}}\right)^2 +
\left(A^i_j\sqrt{\frac{b_{jn}}{b_{in}}}\right)^2\Big]\\
&&=1+\mid{\rm det}\,t\mid^2+
\sum_{1\leq i\leq j\leq 3}\left(\vert M^i_j(X(t))\vert^2+\vert A^i_j(X(t))\vert^2\right).
\end{eqnarray*}

Using the notation $t^i_j=t_{ij}$ and the fact
$$
{\rm det}\,t=t^k_1A^k_1+t^k_2A^k_2+t^k_3A^k_3,\quad k=1,2,3,
$$
we get
$$
2^3\mid{\rm det}\,t\mid\big(H_{3,n}^{-2}(t)-1\big)=
(1-\mid{\rm det}\,t\mid)^2+
\sum_{1\leq i,j\leq 3}\Big(M^i_j(X(t))-A^i_j(X(t))\Big)^2
$$
\begin{equation}
 \label{(x-y)^2=.3}
=(1-\mid{\rm det}\,t\mid)^2+ \sum_{1\leq i\leq j\leq 3}\Big(t^i_j\sqrt{\frac{b_{in}}{b_{jn}}} -A^i_j(t)\sqrt{\frac{b_{jn}}{b_{in}}}\Big)^2.
\end{equation}
Similar to \cite[Lemmas~2.22]{KosJFA17} in the case $m=2$, or \cite[Lemma~10.4.30]{Kos_B_09} we get the following lemma, for $m=3$.
\begin{lem}
\label{l.mu^3(b,0)^t-perp}
 For $t\in {\rm GL}(3,{\mathbb R})\setminus\{e\}$ we have, if ${\rm det}\,t>0$,
$$
(\mu_{(b,0)}^3)^{L_t}\perp\mu_{(b,0)}^3\quad\Leftrightarrow\quad
$$
\begin{equation}
\label{mu^3(b,0)^t-perp.>0} 
\sum_{n\in{\mathbb Z}}\Big[ 
(1\!-\!\mid{\rm det}\,t\mid)^2+\! \sum_{1\leq i\leq 3}\Big(t^i_i -A^i_i(t)\Big)^2
+\!\!\!\!\sum_{1\leq i<j\leq 3}\Big(t^i_j\sqrt{\frac{b_{in}}{b_{jn}}} -A^i_j(t)\sqrt{\frac{b_{jn}}{b_{in}}}\Big)^2
\Big]=\infty.
\end{equation}
If ${\rm det}\,t<0$ we have
$$
(\mu_{(b,0)}^3)^{L_t}\perp\mu_{(b,0)}^3\quad\Leftrightarrow\quad
$$
\begin{equation}
\label{mu^3(b,0)^t-perp.<0}
\sum_{n\in{\mathbb Z}} \Big[ 
(1\!-\!\mid{\rm det}\,t\mid)^2+\! \sum_{1\leq i\leq 3}\Big(t^i_i +A^i_i(t)\Big)^2+
\!\!\!\!\sum_{1\leq i<j\leq 3}\Big(t^i_j\sqrt{\frac{b_{in}}{b_{jn}}} +A^i_j(t)\sqrt{\frac{b_{jn}}{b_{in}}}\Big)^2
\Big]=\infty.
\end{equation}
\end{lem}

By Lemma \ref{l7.3}
and (\ref{Hel2-sl}) the following lemma holds true.
\begin{lem}
\label{l.mu^3(b,a)^t-perp}
 For $t\in {\rm GL}(3,{\mathbb R})\setminus\{e\}$ we have
\begin{equation*}
 (\mu_{(b,a)}^3)^{L_t}\perp\mu_{(b,a)}^3\quad{\it
if}\quad\mid{\rm det}\,t\mid\not=1.
\end{equation*}
If ${\rm det}\,t=1$, we have
$$
(\mu_{(b,a)}^3)^{L_t}\perp\mu_{(b,a)}^3\quad\Leftrightarrow\quad\Sigma^+(t):=\Sigma_1^+(t)+
\Sigma_2(t)=\infty.
$$
If ${\rm det}\,t=-1$, we have
$$
(\mu_{(b,a)}^3)^{L_t}\perp\mu_{(b,a)}^3\quad\Leftrightarrow\quad\Sigma^{-}(t):=\Sigma_1^{-}(t)+
\Sigma_2(t)=\infty,
$$ 
where  
\begin{equation}
\label{Sigma_1^+(t)} 
\Sigma_1^+(t)\!=\!\sum_{n\in{\mathbb
Z}}\Big[
\sum_{k=1}^3(t_{kk}\!-\!A_k^k(t))^2\!+\!\!\!\!
 \sum_{1\leq i<j\leq 3}\Big(t^i_j\sqrt{\frac{b_{in}}{b_{jn}}}\! -\!A^i_j(t)\sqrt{\frac{b_{jn}}{b_{in}}}\Big)^2
\Big],
\end{equation}
\begin{equation} 
\label{Sigma_1^{-}(t)}
\Sigma_1^{-}(t)\!=\!\sum_{n\in{\mathbb
Z}}\Big[
\sum_{k=1}^3(t_{kk}\!+\!A_k^k(t))^2\!+\!\!\!\!\!
 \sum_{1\leq i<j\leq 3}\Big(t^i_j\sqrt{\frac{b_{in}}{b_{jn}}}\!+\!A^i_j(t)\sqrt{\frac{b_{jn}}{b_{in}}}\Big)^2
\Big],
\end{equation}
\begin{eqnarray}
\label{Sigma_3}
&&\Sigma_2(t^{-1})=\sum_{n\in{\mathbb
Z}}\Big[b_{1n}\big((t_{11}-1)a_{1n}+t_{12}a_{2n}+t_{13}a_{3n}\big)^2+\\
\nonumber
&&
b_{2n}\big(t_{21}a_{1n}+(t_{22}\!-\!1)a_{2n}+t_{23}a_{3n}\big)^2
+b_{3n}\big(t_{31}a_{1n}+t_{32}a_{2n}+(t_{33}\!-\!1)a_{3n}\big)^2\Big].
\end{eqnarray}
\end{lem}

\begin{rem}
\label{r.pm-SL(3)}
By Lemma~\ref{l.mu^3(b,a)^t-perp}, it suffices to verify, the condition  of orthogonality 
\begin{equation*}
 (\mu_{(b,a)}^3)^{L_t}\perp\mu_{(b,a)}^3\quad
\text{ for}\quad t\in {\rm GL}(3,{\mathbb R})\setminus\{e\}
\end{equation*}
for the following two subsets 
of the group 
$\pm{\rm SL}(3,\mathbb R)$:
\begin{eqnarray}
 \label{G_3^+}
&&G_3^+:=\{t\in {\rm SL}(3,\mathbb R)\mid t_{kk}= A^k_k(t),\,\,1\leq k\leq 3\},\\
 \label{G_3^-}
&&G_3^-:=\{t\in -{\rm SL}(3,\mathbb R)\mid t_{kk}= -A^k_k(t),\,\,1\leq k\leq 3\}.
\end{eqnarray}
\end{rem}
\begin{lem}
\label{l.SL(3)-perp}
If $t\in G_3^\pm $, we have respectively 
\begin{eqnarray}
\nonumber
&&(\mu_{(b,a)}^3)^{L_t}\perp\mu_{(b,a)}^3\quad\Leftrightarrow\quad\Sigma^\pm(t)=\Sigma_1^\pm(t)+
\Sigma_2(t)=\infty,\\
\label{Sigma_1^+(t).m=3}
 &&\Sigma_1^+(t)=\!\!\sum_{1\leq i< j\leq 3}\sum_{n\in \mathbb Z}\left(t^i_j\sqrt{\frac{b_{in}}{b_{jn}}} -A^i_j(t)\sqrt{\frac{b_{jn}}{b_{in}}}\right)^2\!\!=\!\!\sum_{1\leq i<j\leq 3}\Sigma^+_{ij}(t),\\
 \label{Sigma_1^{-}(t).m=3}
&&\Sigma_1^{-}(t)=\!\! \sum_{1\leq i< j\leq 3}\sum_{n\in \mathbb Z}\left(t^i_j\sqrt{\frac{b_{in}}{b_{jn}}} +A^i_j(t)\sqrt{\frac{b_{jn}}{b_{in}}}\right)^2\!\!=\!\!\sum_{1\leq i< j\leq 3}\Sigma^{-}_{ij}(t),\\
\label{Sigma(pm)(ij)}
&&
\Sigma^\pm_{ij}(t)\!=\!\sum_{n\in \mathbb Z}\left(t^i_j\sqrt{\frac{b_{in}}{b_{jn}}}\! \mp\!A^i_j(t)\sqrt{\frac{b_{jn}}{b_{in}}}\right)^2,
\end{eqnarray}
where $\Sigma_2(t)$ is defined by \eqref{Sigma_3}.
\end{lem}
\begin{rem}
\label{r.conj.m=3}
We note that the {\it Iwasawa decomposition} holds for ${\rm SL}(3,{\mathbb R})$, i.e.,
${\rm SL}(3,{\mathbb R})=KAN$,  where $K={\rm O}(3),$
{\small
\begin{equation}
\label{A,N.m=3} 
A\!=\!\left\{\!D_3(s)\!=\!\left(\begin{array}{ccc}
 s_1&0&0\\
  0&s_2&0\\
  0&0&s_3
\end{array}\right)\!,\,\,{\rm det}D_3(s)\!=\!1\right\},\,N\!=\!\left\{\!\left(\begin{array}{ccc}
 1&x&z\\
  0&1&y\\
  0&0&1
\end{array}\right)\!\mid x,y,z\in \mathbb R\right\}.
\end{equation}
}
\end{rem}
Next we will show that the set $G_3^+$ can be reduced to  the 
six family of one-parameter subgroups $\exp(tE_{kr}),$ $\,\,1\leq k\not=r\leq 3$, see
\eqref{m=3.subrg-orth}, or three families of two-parameter subgroups, 
see \eqref{2-par.m=3}.
The set $G_3^-$  can be reduced to  the three two-parameter family  \eqref{2-par-P_r.m=3} 
reflections of \eqref{2-par.m=3}  by $P_r$. The remaining part is reduced to the sets $D_3(s){\rm O}(3)D^{-1}_3(s)P_r$
or five parameter family 
of elements $\tau_r(t,s)\!=\!D_3(s)tD^{-1}_3(s)P_r$, see \eqref{tau_r(t,s)}. 
\begin{lem}
\label{l.m=3.crit-orthog} 
In case $m=3$ the {\rm minimal generating set} ${\rm GL}(3,{\mathbb R})^\perp_0(\mu)$ is defined as follows
(compare with Remark~\ref{r.m=2.crit-orthog}) :
\begin{eqnarray}
\nonumber 
&&
 {\rm GL}(3,{\mathbb R})^\perp_0(\mu)=\{e_r(t,s),e_r(t,s)P_r,\mid 1\leq r\leq 3,\,\,(t,s)\in {\mathbb R}^2\}\\
\label{G_0-min-m=3} 
 &&
\cup\{
{\rm O}_r^A(3),\,\,1\leq r\leq 3\} ,\quad\text{where}
\end{eqnarray}
\begin{equation}
 \label{m=3.subrg-orth}
e_{kn}(t):=\exp(tE_{kn})=I+tE_{kn},\,\,1\leq k\not=n\leq 3,\,\,t\in \mathbb R,
\end{equation}
\begin{equation}
\label{2-par.m=3}
 e_1(t,s)=
 \left(\begin{smallmatrix}
 1&t&s\\
  0&1&0\\
  0&0&1
\end{smallmatrix}\right),\,\,
 e_2(t,s)= 
\left(\begin{smallmatrix}
1&0&0\\
t&1&s\\
0&0&1
\end{smallmatrix}\right),\,\,
 e_3(t,s)=   
 \left(\begin{smallmatrix}
1&0&0\\
  0&1&0\\
t&s&1
\end{smallmatrix}\right),
\end{equation}
\begin{equation}
\label{2-par-P_r.m=3}
 e_1(t,s)P_1\!=\!\left(\begin{smallmatrix}
 -1&t&s\\
  0&1&0\\
  0&0&1
\end{smallmatrix}\right),\,
 e_2(t,s)P_2\!=\! \left(\begin{smallmatrix}
1&0&0\\
t&-1&s\\
0&0&1
\end{smallmatrix}\right),\,
e_3(t,s)P_3\!=\! \left(\begin{smallmatrix}
1&0&0\\
  0&1&0\\
t&s&-1
\end{smallmatrix}\right),
\end{equation}
\begin{equation}
\label{m=3.P_k-orth} 
P_1=\left(\begin{smallmatrix}
-1&0&0\\
  0&1&0\\
0&0&1
\end{smallmatrix}\right),\quad
P_2=\left(\begin{smallmatrix}
1&0&0\\
0&-1&0\\
0&0&1
\end{smallmatrix}\right),\quad
P_3=\left(\begin{smallmatrix}
1&0&0\\
0&1&0\\
0&0&-1
\end{smallmatrix}\right),
\end{equation}
\begin{equation}
 \label{O^A(3)}
{\rm O}^A(3):=\{D_3(s){\rm O}(3)D^{-1}_3(s)\mid D_3(s)\in A\},
\end{equation}
\begin{equation}
 \label{O_r^A,-(3)}
{\rm O}_r^A(3):=\{D_3(s){\rm O}(3)D^{-1}_3(s)P_r\mid D_3(s)\in A\}
,\,\,1\leq r\leq 3,
\end{equation}
\begin{equation}
 \label{tau_r(t,s)}
\tau_r(t,s):=D_3(s)tD_3^{-1}(s)P_r,\,\,
\,\,\,t\in O(3),\,\,D_3(s)={\rm diag}(s_1,s_2,s_3)\in A,
\end{equation}
and $A$ is defined by \eqref{A,N.m=3}.
\end{lem}

The families  \eqref{m=3.subrg-orth} give us respectively the divergence of the following series:
\begin{equation}
\label{S^L(kn)}
S^L_{kr}(\mu)=\sum_{n\in \mathbb Z}\frac{b_{kn}}{2}\Big(\frac{1}{2b_{rn}}+a_{rn}^2\Big),\quad1\leq k,r\leq 3,\,\,k\not=r.
\end{equation}
The families  \eqref{2-par.m=3} give us, respectively, the  divergence of the following series:
{\small
\begin{eqnarray}
&& 
\label{S^{L}_{1,23}}
S^{L}_{1,23}(\mu,t,s)=\sum_{n\in{\mathbb
Z}}\left[\frac{t^2}{4}\frac{b_{1n}}{b_{2n}}+\frac{s^2}{4}\frac{b_{1n}}{b_{3n}}+
\frac{b_{1n}}{2}\big(\!-\!2a_{1n}\!+\!ta_{2n}\!+\!sa_{3n}\big)^2\right],\,
\\
&& 
\label{S^{L}_{2,31}}
S^{L}_{2,13}(\mu,t,s)=\sum_{n\in{\mathbb
Z}}\left[\frac{t^2}{4}\frac{b_{2n}}{b_{1n}}+ \frac{s^2}{4}\frac{b_{2n}}{b_{3n}}+
\frac{b_{2n}}{2}\big(ta_{1n}-2a_{2n}+sa_{3n}\big)^2\right],\,\,\,\,\\
&& 
\label{S^{L}_{3,12}}
S^{L}_{3,12}(\mu,t,s)=\sum_{n\in{\mathbb
Z}}\left[\frac{t^2}{4}\frac{b_{3n}}{b_{1n}}+ \frac{s^2}{4}\frac{b_{3n}}{b_{2n}}+
\frac{b_{3n}}{2}\big(ta_{1n}+sa_{2n}-2a_{3n}\big)^2\right].
\end{eqnarray}
}

The families  \eqref{tau_r(t,s)} give us the  conditions 
\eqref{O(3).perp}, see Lemma~\ref{l.O(3).perp} below.
\begin{pf}
Consider the subset  ${\rm GL}(3,{\mathbb R})^\perp_{0}(\mu)$  of ${\rm GL}(3,{\mathbb R})$ described by \eqref{G_0-min-m=3}. 
The fact that this set is {\it minimal generating} will follow from Lemma~\ref{l.approx-(x,D).m=3}, more precisely, 
from the following implications:
\begin{eqnarray}
&& 
\label{1=0}
\Big(\mu^{L_t}\perp\mu\,\,\,\text{for all}\,\,\, t\in {\rm GL}(3,{\mathbb R})^\perp_0(\mu)\Big)\Rightarrow
\Big(\text{irreducibility}\Big)\\
\nonumber
&&
\Rightarrow
\Big(\mu^{L_t}\perp \mu\,\,\,\text{for all}\,\,\,t\in {\rm GL}(3,{\mathbb R})\setminus\{e\}\Big).
\end{eqnarray}
The first implication is just Lemma~\ref{l.approx-(x,D).m=3}.
The second implication follows from irreducibility. 
Indeed, suppose that ${\rm GL}(3,{\mathbb R})^\perp_0(\mu)$ is not a minimal generating set. Then   we can find  an $s\in {\rm GL}(3,{\mathbb R})\setminus\{e\}$ having the property
$$
\big(\mu_{(b,a)}^3\big)^{L_s}\sim \mu_{(b,a)}^3.
$$
Hence the non-trivial operator $T^{L,\mu,3}_s$
can be defined by
\begin{equation}
\label{T^(L,mu,3)}
(T^{L,\mu,3}_sf)(x)=\big(d\mu_{(b,a)}^3(s^{-1}x)/d\mu_{(b,a)}^3(x)\big)^{1/2}f(s^{-1}x),\,\,f\in L^2(X_3,\mu^3_{(b,a)}).
\end{equation}
This operator  commutes  with  the representations $T^{R,\mu,3}$: 
$$
[T^{R,\mu,3}_t,T^{L,\mu,3}_s]=0\quad\text{for all}\quad t\in G,
$$
contradicting the irreducibility.

The relations \eqref{S^L(kn)}--\eqref{S^{L}_{3,12}}
follows from \eqref{Sigma_1^+(t)}--\eqref{Sigma_3}.
The relation \eqref{S^{L}_{2,31}}, for example, follows from \eqref{Sigma_1^{-}(t)} and \eqref{Sigma_3}.
The relation \eqref{O(3).perp}
we obtain from \eqref{Sigma_1^+(t)} for $\tau_r(t,s),\,\,t\!\in \!{\rm O}(3),\,\,s\!\in \!\big({\mathbb R^*}\big)^3$ 
defined by 
\begin{equation}
\label{tau(t,s)}
\tau_r(t,s)=D_3(s)tD_3^{-1}(s)P_r,\quad\text{where}\quad D_3(s)={\rm diag}(s_1,s_2,s_3).
\end{equation}
\qed\end{pf}
\begin{lem}
\label{l.O(3).perp}
Set  
$$\tau(s,t):=D_3(s)tD_3^{-1}(s),\quad \tau_r(s,t):
=\tau(s,t)P_r
$$
for $t\!\in\! \pm {\rm O}(3),\,\,D_3(s)\!=\! {\rm diag}(s_1,s_2,s_3),\,\,s\!=\!(s_1,s_2,s_3)\!\in\! (\mathbb R^*)^3$ 
and $1\leq r\leq 3$. Then 
\begin{equation}
 \label{O(3).perp}
 \big(\mu_{(b,a)}^3\big)^{L_{\tau_r(s,t)}}\perp \mu_{(b,a)}^3\Leftrightarrow
\Sigma_1^\pm\big(\tau_r(s,t)\big)
+\Sigma_2\big(
\tau_r(s,t)\big)
=\infty, 
\end{equation}
where $\Sigma_1^\pm(t)$
are defined by \eqref{Sigma_1^+(t).m=3},\,\,\eqref{Sigma_1^{-}(t).m=3} and 
$\Sigma_2(t)$ is defined by \eqref{Sigma_3}.

In particular, if we denote $ s_{ij}=s_is_j^{-1}$ we get
\begin{equation}
\label{Sigma-123(T)}
\Sigma_{1}^+\big(\tau(t,s)\big)=\Sigma_{1}^+(t,s)=
t^2_{12}\Sigma_{12}(s_{12}^{1/2})+
t^2_{13}\Sigma_{13}(s_{13}^{1/2})+
t^2_{23}\Sigma_{23}(s_{23}^{1/2}).
\end{equation}
\end{lem} 
\begin{pf}
For $T:=\tau(s,t)$ and $T(3):=\tau_3(s,t)$ we have respectively:
\begin{eqnarray}
 \label{T}
&&T=
D_3(s)tD_3^{-1}(s)\!=\!
\left(\!
\begin{smallmatrix}
t_{11}&\frac{s_1}{s_2}t_{12}&\frac{s_1}{s_3}t_{13}\\
\frac{s_2}{s_1}t_{21}&t_{22}&\frac{s_2}{s_3}t_{23}\\
\frac{s_3}{s_1}t_{31}&\frac{s_3}{s_2}t_{32}&t_{33}
\end{smallmatrix}
\!\right),\\
\label{T(3)}
&& 
\left(
\begin{smallmatrix}
t_{11}&\frac{s_1}{s_2}t_{12}&-\frac{s_1}{s_3}t_{13}\\
\frac{s_2}{s_1}t_{21}&t_{22}&-\frac{s_2}{s_3}t_{23}\\
\frac{s_3}{s_1}t_{31}&\frac{s_3}{s_2}t_{32}&-t_{33}
\end{smallmatrix}
\right)=
D_3(s)tD_3^{-1}(s)P_3=:T(3).
\end{eqnarray}
By Lemma~\ref{l.O(3)-M^i_j} we have for $t\in {\rm O}(3)$ 
\begin{equation}
\label{A(t)}
t_{kr}= A^k_r(t),\,\, 1\leq k,r\leq 3. 
\end{equation}
Therefore, for $T$ and $T(3)$ we have for $1\!\leq\! k,r\!\leq\! 3$:   
\begin{eqnarray}
\label{A(T)}
&&
M^k_r(T)=T_{kr}\!=\!\frac{s_k}{s_r}t_{kr},\,\,A^k_r(T)\!=\!\frac{s_r}{s_k}A^k_r(t)\stackrel{\eqref{O(3)-M^i_j}}{=}
\frac{s_r}{s_k}t_{kr},\quad
M^k_r(T(3))
\quad \quad
\\
\label{A(T(3))}
&&
\!=\!(-1)^{\delta_{3,r}}\frac{s_k}{s_r}t_{kr},\,
\, A^k_r(T(3))\!=(-1)^{\delta_{3,r}}\frac{s_r}{s_k}A^k_r(t)=
(-1)^{\delta_{3,r}}
\frac{s_r}{s_k}t_{kr}.
\end{eqnarray}
Finally, we get 
\begin{eqnarray*}
&& 
\Sigma_{1}^+(T)=\Sigma_{1}^+\big(\tau(s,t)\big)
=\sum_{n\in{\mathbb
Z}}\Big[
 \sum_{1\leq i<j\leq 3}\Big(M^i_j(T)\sqrt{\frac{b_{in}}{b_{jn}}} -A^i_j(T)\sqrt{\frac{b_{jn}}{b_{in}}}\Big)^2
\Big]\\
&&=\!\sum_{n\in {\mathbb Z}}\Big[
t_{12}^2\Big(
s_{12}\sqrt{\frac{b_{1n}}{b_{2n}}}\!-\!
s_{12}^{-1}\sqrt{\frac{b_{2n}}{b_{1n}}}\Big)^2
+
t_{13}^2\Big(
s_{13}\sqrt{\frac{b_{1n}}{b_{3n}}}\!-\!
s_{13}^{-1}\sqrt{\frac{b_{3n}}{b_{1n}}}\Big)^2+\\
&&t_{23}^2\Big(
s_{23}\sqrt{\frac{b_{2n}}{b_{3n}}}\!-\!
s_{23}^{-1}\sqrt{\frac{b_{3n}}{b_{2n}}}\Big)^2
\Big]=t^2_{12}\Sigma_{12}(s_{12}^{1/2})+
t^2_{13}\Sigma_{13}(s_{13}^{1/2})+
t^2_{23}\Sigma_{23}(s_{23}^{1/2}).
\end{eqnarray*}
 $\Sigma_{1}^-(T(3))=t^2_{12}\Sigma_{12}(s_{12}^{1/2})+
t^2_{13}\Sigma_{13}(s_{13}^{1/2})+
t^2_{23}\Sigma_{23}(s_{23}^{1/2})$
\qed\end{pf}
%
\begin{lem}
 \label{l.O(3)-M^i_j}
For an  arbitrary  orthogonal matrix $t\in \pm{\rm O}(3)$  
we have 
\begin{equation}
 \label{O(3)-M^i_j} 
 t_{kn}=\pm A^k_n(t),\,\, 1\leq k,n\leq 3,  \quad\text{where}\quad
t=\left(\begin{smallmatrix}
t_{11}&t_{12}&t_{13}\\
t_{21}&t_{22}&t_{23}\\
t_{31}&t_{32}&t_{33}
\end{smallmatrix}
\right).
\end{equation}
\end{lem}
\begin{pf}
Denote the three rows of the matrix $t$ 
by, respectively, $t_1,t_2,t_3\in \mathbb R^3$. Since $t\in \pm{\rm O}(3)$ we get
\begin{equation}
\label{O(3)}
\Vert t_1\Vert^2=\Vert t_2\Vert^2=\Vert t_3\Vert^2=1 \quad\text{and} \quad t_l\perp t_r,\,\,l\not=r.
\end{equation}
Moreover, since $t_1$ is orthogonal to the hyperplane $V_{23}$ generated by the vectors $t_2$ and $t_3$ 
and $t\in \pm{\rm O}(3)$  we get respectively $t_l=\pm[t_r,t_s]$, where
$[x,y]$ is the {\it vector product} or {\it cross product}  
of two vectors $x,y\in \mathbb R^3$ and 
the triple  $\{l,r,s\}$ 
denotes any cyclic permutations of $\{1,2,3\}$.
For $t\in {\rm O}(3)$ and $l=1$ we get 
\begin{equation}
\label{[t_2,t_3].3} 
t_1=[t_2,t_3]=
\left|\begin{smallmatrix}
i&j&k\\
t_{21}&t_{22}&t_{23}\\
t_{31}&t_{32}&t_{33}
\end{smallmatrix}\right|=
i
\left|\begin{smallmatrix}
t_{22}&t_{23}\\
t_{32}&t_{33}
\end{smallmatrix}\right|-j\left|\begin{smallmatrix}
t_{21}&t_{23}\\
t_{31}&t_{33}
\end{smallmatrix}\right|+k
\left|\begin{smallmatrix}
t_{21}&t_{22}\\
t_{31}&t_{32}
\end{smallmatrix}\right|,
\end{equation}
where $i,j,k$ is the standard orthonormal basis in $\mathbb R^3$, i.e., 
$$i=(1,0,0),\,\,\,j=(0,1,0),\,\,\,k=(0,0,1).$$  Define  $X$ {\it formally} as the matrix:  
$$
X=\left(\begin{smallmatrix}
i&j&k\\
x_1&x_2&x_3\\
y_1&y_2&y_3
\end{smallmatrix}\right).
$$
Then
$$
t_1=(t_{11},t_{12},t_{13})=\left(A^1_1(X),A^1_2(X),A^1_3(X)\right),
$$
thus proving \eqref{O(3)-M^i_j} for $k=1$. For other rows the proof is similar.
\qed\end{pf}
\begin{rem}
\label{O(n)}
For $t\in \pm{\rm O}(n)$ we can prove a similar statement.
\end{rem}
\section{Irreducibility, the cases $m=1$ and $m=2$}
For convenience  of the reader, we recall  the previous results (see details
in \cite{KosJFA17}). 

\subsection{Case $m=1$}
Let us denote by $\langle f_n\mid n\in {\mathbb
N}\rangle$ the {\it closure of the linear space}
generated by the set of
vectors $(f_n)_{n\in{\mathbb N}}$ in a Hilbert space $H.$ Consider the measure $\mu_{(b,a)}$
on 
the space $X_1\sim \mathbb R^{\mathbb Z}=\otimes_{n\in\mathbb Z}\mathbb R$,
the infinite product of the real lines:
\begin{equation*}
 d\mu_{(b,a)}(x)=\otimes_{n\in \mathbb Z}\sqrt{\frac{b_{1n}}{\pi}}
e^{-b_{1n}(x_{1n}-a_{1n})^2}dx_{1n},
\end{equation*}
for $b=(b_{1n})_{n\in \mathbb Z}$ and $a=(a_{1n})_{n\in \mathbb Z}$ with $b_{1n}>0,\,\,a_{1n}\in \mathbb R$ where $x=(x_{1n})_{n\in \mathbb Z}$.  
In the case $m=1$ the generators $A_{kn}:=A_{kn}^{R,1}$ have the form
$$
A_{kn}=x_{1k}D_{1n},\quad k,n \in {\mathbb Z},
$$
where $D_{kn}=\frac{\partial }{\partial x_{kn}}-b_{kn}(x_{kn}-a_{kn}).$ The following lemmas are proved in \cite{KosAlb05J}
\begin{lem}
\label{5.1} The following  three conditions are equivalent:\par
(i)\,$\,\,(\mu_{(b,a)})^{L_t}\perp\mu_{(b,a)}$ for all $t\in
{\rm GL}(1,{\mathbb R})\setminus\{e\}$, \par
(ii)\,\,\,$(\mu_{(b,a)})^{L_{-E_{11}}}\perp\mu_{(b,a)}$,\par
(iii)\,\,\,\,$S^L_{11}(\mu)=4\sum_{n\in{\mathbb Z}}{b_{1n}}a_{1n}^2=\infty. $
\end{lem}
\begin{lem} 
\label{5.2}
 For $k,m\in {\mathbb Z}$ we have
$$
x_{1k}x_{1m}{\bf 1}\in\langle A_{kn}A_{mn}{\bf
1}=x_{1k}x_{1m}D_{1n}^2{\bf 1} \mid n\in {\mathbb Z}\rangle.
$$
\end{lem}
\begin{lem}
\label{5.3}
 For any $k\in {\mathbb Z}$ we have
$$
x_{1k}{\bf 1}\in\langle x_{1k}x_{1n}{\bf 1}\mid n\in {\mathbb Z}\rangle
\Leftrightarrow S^L_{11}(\mu)=\infty.
$$
\end{lem}
So, the operators $x_{1k},\,\,k\in {\mathbb Z}$  are {\it affiliated}
with the von Neumann algebra ${\mathfrak A}^1$ generated by the representation, which completes
the proof of the  irreducibility for $m=1.$
\begin{df}
\label{A-eta-M} Recall (see \cite{Dix69W}) that, a  not
necessarily bounded self-adjoint operator $A$ in a Hilbert space
$H$, is said to be {\it affiliated} with a von Neumann algebra $M$
of operators in this Hilbert space $H$ if $e^{itA}\in M$ for all
$t\in{\mathbb R}$. One writes $A\,\,\eta\,\,M$.
\end{df}

\subsection{Case $m=2$, approximation of $x_{kn}$ and $D_{kn}$}
In this case the generators
$A_{kn}:=A_{kn}^{R,2}:=\frac{d}{dt}T^{R,\mu,2}_{I+tE_{kn}}\mid_{t=0}$
have the form:
$$
A_{kn}=x_{1k}D_{1n}+x_{2k}D_{2n},\quad k,n \in {\mathbb Z}.
$$
We will formulate several  
useful lemmas for
approximation of the operators of multiplication by the independent variables $x_{kn}$ and operators
$D_{kn}$ by combinations of the generators $A_{kn}$. 
In what follows we use the following notation for $f,g\in{\mathbb R}^m$ (see Remark~\ref{r.two-sim-ser} for notations $\sim$)
{\small
\begin{equation}
\label{D(f,g)}
\Delta(f,g):=\frac{\Gamma(f)+\Gamma(f,g)}{\Gamma(g)+1}\sim
\frac{I+\Gamma(f)+\Gamma(g)+\Gamma(f,g)}{I+\Gamma(g)}=\frac{{\rm det}(I+\gamma(f,g))}{{\rm det}(I+\gamma(g))},
\end{equation}
}
%
where $\Gamma(x_1,\dots,x_n)$ is the Gram determinant 
of vectors $x_1,x_2,..., x_n\in {\mathbb R}^m$ defined  by (\ref{Gram-det}).
To make the notations of the article \cite{KosJFA17} compatible with the notations in the case $m=3$ (see \eqref{Y_r^{(s)}}), we denote
%
{\small
\begin{eqnarray}
\nonumber
&& \Vert Y_1^{(1)} \Vert^2:=\!\Vert f^1 \Vert^2=\sum_{k\in \mathbb Z}\frac{b_{1k}^2}{b_{1k}^2+2b_{1k}b_{2k}},\quad
\Vert Y_2^{(1)} \Vert^2:=\Vert g^1\Vert^2=\sum_{k\in \mathbb Z}\frac{b_{2k}^2}{b_{1k}^2+2b_{1k}b_{2k}},\\
\nonumber
&&\Vert Y_1^{(2)} \Vert^2:=\!\Vert g^2 \Vert^2=\sum_{k\in \mathbb Z}\frac{b_{1k}^2}{b_{2k}^2+2b_{1k}b_{2k}},\quad
\Vert Y_2^{(2)} \Vert^2:=\!\Vert f^2\Vert^2=\sum_{k\in \mathbb Z}\frac{b_{2k}^2}{b_{2k}^2+2b_{1k}b_{2k}},\\
\label{Y_r^{(s)}.m=2}
&&\Vert Y_1 \Vert^2:=\Vert f \Vert^2=\sum_{k\in \mathbb Z}\frac{a_{1k}^2}{\frac{1}{2b_{1k}}+\frac{1}{2b_{2k}}},\quad
\Vert Y_2 \Vert^2:=\Vert g\Vert^2=\sum_{k\in \mathbb Z}\frac{a_{2k}^2}{\frac{1}{2b_{1k}}+\frac{1}{2b_{2k}}}.
\end{eqnarray}
}
\begin{lem}
\label{x1x1}
  For any  $k,t\in {\mathbb Z}$ one has
$$
x_{1k}x_{1t}\in\langle A_{kn}A_{tn}{\bf 1}\mid n\in{\mathbb
Z}\rangle \,\,\Leftrightarrow\,\,\Delta(Y_1^{(1)},Y_2^{(1)}).
$$

\end{lem}
\begin{lem}
\label{x2x2} For any  $k,t\in {\mathbb Z}$ we have
$$
x_{2k}x_{2t}\in\langle A_{kn}A_{tn}{\bf 1}\mid n\in {\mathbb
Z}\rangle\,\,\Leftrightarrow\,\,\Delta(Y_2^{(2)},Y_1^{(2)})=\infty.)
$$
\end{lem}
\begin{lem}
\label{x1,dA}
Set $\Sigma^{rs}=\sum_{n\in{\mathbb Z}}\frac{b_{rn}}{b_{sn}},\,\,1\leq r,s\leq 2$. For any  $k\in {\mathbb Z}$ we get
 $$
 x_{1k}{\bf 1}\in\langle D_{1n}A_{kn}{\bf 1}\mid n\in{\mathbb Z}\rangle
 \quad\Leftrightarrow\quad \Sigma^{12}=
\infty.
$$
\end{lem}
\begin{lem}
\label{x2,dA}
 For any  $k\in {\mathbb Z}$ we have
 $$
 x_{2k}{\bf 1}\in\langle D_{2n}A_{kn}{\bf 1}\mid n\in{\mathbb Z}\rangle
 \quad\Leftrightarrow\quad  \Sigma^{21}=\infty.
$$
\end{lem}

\begin{lem}
\label{d1}
 For any  $n\in {\mathbb Z}$ we have
$$
D_{1n}{\bf 1}\in\langle A_{kn}{\bf 1}\mid k\in {\mathbb Z}\rangle
\quad\Leftrightarrow\quad 
\Delta(Y_1,Y_2)=\infty.
$$
\end{lem}

\begin{lem}
\label{d2} For any  $n\in {\mathbb Z}$ we have
$$
D_{2n}{\bf 1}\in\langle A_{kn}{\bf 1}\mid k\in {\mathbb Z}\rangle
\quad\Leftrightarrow\quad 
\Delta(Y_2,Y_1)=\infty.
$$
\end{lem}
\subsubsection{Technical part of the proof of irreducibility}
\label{sec.3.2}
\begin{lem}
\label{approx-(x,D).2} If $\mu^{L_t}\perp\mu$ for all $t\in
{\rm GL}(2,{\mathbb R})\setminus\{e\}$, we can approximate 
by combinations of generators 
at least one of
the following  pairs of operators:
$(x_{1n},x_{2n}),\,\,(x_{1n},D_{2n}),\,\,(D_{1n},x_{2n})\,\,$ or
$(D_{1n},D_{2n})$.
\end{lem}
\begin{pf} 
Recall the orthogonality conditions for the case $m=2$ 
\begin{eqnarray*}
&&
S^L_{kr}(\mu)=\sum_{n\in {\mathbb Z}}
\frac{b_{kn}}{2}\left(\frac{1}{2b_{rn}}+a_{rn}^2\right),\quad 1\leq k,r\leq 2,\,\, k\not=r,\\
&&
S^{L,-}_{kr}(\mu,t)=\frac{t^2}{4}\sum_{n\in{\mathbb
Z}}\frac{b_{kn}}{b_{rn}}+ \sum_{n\in{\mathbb
Z}}\frac{b_{kn}}{2}(-2a_{kn}+ta_{rn})^2,\quad 1\leq k\not=r\leq 2,\\
&&\Sigma^-_{12}\big(\tau_{-}(\phi,s)\big)=\sin^2\phi\,\Sigma_1(s)+\Sigma^-_2\big(\tau_{-}(\phi,s)\big),\quad\text{where}\\
&&
\Sigma_1(s):=\!\sum_{n\in {\mathbb Z}}
\Big(s^2\sqrt{\frac{b_{1n}}{b_{2n}}}\!-\!s^{-2}\sqrt{\frac{b_{2n}}{b_{1n}}}\Big)^2,\\ 
&&
\Sigma^-_2\big(\tau_{-}(\phi,s)\big)\!=\!\sum_{n\in
{\mathbb
Z}}\Big(4b_{1n}\sin^2\frac{\phi}{2}\!+\!4s^{-4}b_{2n}\cos^2\frac{\phi}{2}\Big)
\Big(a_{1n}\sin\frac{\phi}{2}\!-\!s^2a_{2n}\cos\frac{\phi}{2}\Big)^2.
\end{eqnarray*}
Let ${\mathfrak A}^2$ be the von Neumann algebra generated by our representation. 
In order to approximate
operators  $x_{kn}$ or $D_{kn}$  by the corresponding generators,
by Lemmas~\ref{x1x1}--\ref{x2x2}, Lemmas~\ref{5.2}--~\ref{5.2} and Lemmas~\ref{d1}--~\ref{d2} we have:
\begin{eqnarray*}  
&& x_{1n}x_{1t}\,\,\eta\,\, {\mathfrak A}^2\Leftrightarrow 
\Delta(Y_1^{(1)},Y_2^{(1)})=\infty,\quad
x_{2n}x_{2t}\,\,\eta\,\, {\mathfrak A}^2\Leftrightarrow 
\Delta(Y_2^{(2)},Y_1^{(2)})=\infty,\\
&& D_{1n}\,\,\eta\,\, {\mathfrak A} ^2\Leftrightarrow \Delta(Y_1,Y_2)=\infty,\quad 
D_{2n}\,\,\eta\,\, {\mathfrak A}^2\Leftrightarrow \Delta(Y_2,Y_1)=\infty,
\end{eqnarray*}
where $Y_r^{(s)}$ and $Y_r$ for $1\leq r,s\leq 2$ are defined by \eqref{Y_r^{(s)}.m=2}.
\qed\end{pf}

\subsubsection{Scheme of the proof for two lines}

There are two different cases: 
\begin{itemize}
	\item[${\rm I.}$] Approximation of $x_{rk}x_{rt}$ for $1\leq r\leq 2$ by $A_{kn}A_{tn}$,
    \item[${\rm II.}$]  Approximation of $D_{rk}$ for $1\leq r\leq 2$  by $A_{kn}$.   
\end{itemize}
In the case $m=2$ the analysis of the divergence
$$
\Delta(Y_1^{(1)},Y_2^{(1)})=\infty\quad \text{and} \quad \Delta(Y_2^{(2)},Y_1^{(2)})=\infty
$$ 
is governed by the 
convergence or divergence
of 
$\Sigma^{12}$ and $\Sigma^{21}$, since 
\begin{eqnarray*}
&&\Vert Y_1^{(1)} \Vert^2=\sum_{k\in \mathbb Z}\frac{b_{1k}^2}{b_{1k}^2+2b_{1k}b_{2k}}\sim \frac{1}{2}\sum_{k\in \mathbb Z}\frac{b_{1k}^2}{b_{1k}b_{2k}}=
\frac{1}{2}\Sigma^{12},\\
&&\Vert Y_2^{(2)}\Vert^2=\sum_{k\in \mathbb Z}\frac{b_{2k}^2}{b_{2k}^2+2b_{1k}b_{2k}}\sim \frac{1}{2}\sum_{k\in \mathbb Z}\frac{b_{2k}^2}{b_{1k}b_{2k}}
=\frac{1}{2}\Sigma^{21}.
\end{eqnarray*}
%
%
\begin{rem}
\label{S_r(2)}
To guess the right generalisation for the case $m=3$ we note that $\Sigma^{12}=S_1(2)$ and $\Sigma^{21}=S_2(2)$, where  we denote
 %
%
\begin{eqnarray*}
&& S_1(2):=\sum_{k\in \mathbb Z}\frac{b_{1k}^2}{b_{1k}b_{2k}},\,\,\,
S_2(2):=\sum_{k\in \mathbb Z}\frac{b_{2k}^2}{b_{1k}b_{2k}},\,\,\,\Sigma^{12}=\sum_{k\in \mathbb Z}\frac{b_{1k}}{b_{2k}},\,\,\,
\Sigma^{21}=\sum_{k\in \mathbb Z}\frac{b_{2k}}{b_{1k}}.
 \end{eqnarray*}
 We also observe that
\begin{equation}
  \label{Y_r(s)-S_r(3).m=2} 
\Vert Y_1^{(1)} \Vert^2\sim S_1(2),\quad  \Vert Y_2^{(2)} \Vert^2\sim S_2(2).
\end{equation}
Hence, in the case  $m=3$ it is natural to replace  $\Sigma^{12}=S_1(2)$ and $\Sigma^{21}=S_2(2)$ by   
$S_r(3)$, which is  defined as follows:
\begin{equation}
\label{S_r(3)}
 S_r(3)=\sum_{n\in \mathbb Z}\frac{b_{rn}^2}{b_{1n}b_{2n}+b_{1n}b_{3n}+b_{2n}b_{3n}},\quad 1\leq r\leq 3.
\end{equation}
\end{rem}
\subsubsection{I, Approximation of $x_{1n}$ and $x_{2n}$}
\begin{lem}
\label{l.S1+S2.m=2}
We have 
 \begin{eqnarray}
 \label{S1+S2}
 &&S_1(2)+S_2(2)=\infty,\\
\label{Yr-sim-Sr.m=2}
&& \Vert Y_r^{(r)}\Vert^2\sim S_r(2)\quad{for\,\,\, all}\quad 1\leq r\leq 2,\\ 
\label{Yr^s-<-Sr.m=2}
&& \Vert Y_1^{(2)}\Vert^2< \frac{1}{2}S_1(2),\quad
 \Vert Y_2^{(1)}\Vert^2< \frac{1}{2}S_2(2),\\
 \label{i_1,i_2.m=2} 
&&\Vert Y_1^{(i)}\Vert^2+\Vert Y_2^{(j)} \Vert^2=\infty,\quad
i,j\in \{1,2\}.
\end{eqnarray}
\end{lem}
\begin{pf} Since $a^2+b^2\geq 2ab$ we get
 \begin{equation*}
S_1(2)+S_2(2)=\sum_{k\in \mathbb Z}\frac{b_{1k}^2+b_{2k}^2}{b_{1k}b_{2k}}\geq \sum_{k\in \mathbb Z}2=\infty.
\end{equation*}
Further, for $1\leq r\leq 2$
\begin{eqnarray*}
&&
\Vert Y_r^{(r)} \Vert^2=\sum_{k\in \mathbb Z}\frac{b_{rk}^2}{b_{rk}^2+2b_{1k}b_{2k}}\sim  
\sum_{k\in \mathbb Z}\frac{b_{rk}^2}{2b_{1k}b_{2k}}=\frac{1}{2}S_r(2),\\
 &&\Vert Y_1^{(2)} \Vert^2=\sum_{k\in \mathbb Z}\frac{b_{1k}^2}{b_{2k}^2+2b_{1k}b_{2k}}<
\sum_{k\in \mathbb Z}\frac{b_{1k}^2}{2b_{1k}b_{2k}}= \frac{1}{2}S_1(2),\\
 &&\Vert Y_2^{(1)} \Vert^2=\sum_{k\in \mathbb Z}\frac{b_{2k}^2}{b_{1k}^2+2b_{1k}b_{2k}}<
\sum_{k\in \mathbb Z}\frac{b_{2k}^2}{2b_{1k}b_{2k}}= \frac{1}{2}S_2(2).
\end{eqnarray*}
In addition,  we have
\begin{eqnarray*}
&& \Vert Y_1^{(i)}\Vert^2+\Vert Y_2^{(j)} \Vert^2=
\sum_{k\in \mathbb Z}\frac{b_{1k}^2}{b_{ik}^2+2b_{1k}b_{2k}}+\sum_{k\in \mathbb Z}\frac{b_{2k}^2}{b_{jk}^2+2b_{1k}b_{2k}}\\
&&>\sum_{k\in \mathbb Z}\frac{b_{1k}^2+b_{2k}^2}{(b_{1k}+b_{2k})^2}=\infty.
\hskip 8.3cm \Box
\end{eqnarray*}
\end{pf}
{\bf
\begin{rem}
\label{r.(0,1)} 
In what follows  if some expression $<\infty$ (resp. $=\infty$)
we denote this case  by $0$ (respectively,  by $1$).
\end{rem}
}
Set $S=(S_1(2),S_2(2))$, since $S_1(2)+S_2(2)=\infty$ we have two cases:\\
I(1) $S=(0,1)$,  
if $S=(1,0)$ we interchange $(b_{1n},a_{1n})_n$ with  $(b_{2n},a_{2n})_n$,\\
I(2) $S=(1,1)$,
i.e., $S_1(2)=\infty,\,\,S_2(2)=\infty$.
\subsubsection{Case  $S=(0,1)$}
\begin{lem}
\label{(0,1).m=2} 
In the case  $S=(0,1)$ the  representation is irreducible, moreover we can approximate: \par
{\rm (1)} $x_{2k}x_{2t}$   by $A_{kn}A_{tn}$, since $\Delta(Y_2^{(2)},Y_1^{(2)})=\infty$, \par
{\rm (2)} $D_{1n},\,D_{2n}$ by $A_{kn}$, since   $\Delta(Y_1,Y_2)=\infty$ and $\Delta(Y_2,Y_1)=\infty$.
\end{lem}
\begin{pf}
(1) Set (compare with \eqref{y^(k).m=3} in the case $m=3$)
\begin{eqnarray}
 \label{y^(k).m=2}
&&y^{(k)}=(y_1^{(k)},y_2^{(k)})=(\Vert Y_1^{(k)}\Vert^2,\Vert Y_2^{(k)}\Vert^2),\,\,1\leq k\leq 2,\\
\label{y-mat.m=2}
&&\,\,\text{or}\,\,y=
\left(\begin{array}{c}
y^{(1)}\\
y^{(2)}
\end{array}\right)=\left(\begin{array}{cc}
y_1^{(1)}&y_2^{(1)}\\
y_1^{(2)}&y_2^{(2)}
\end{array}\right).
\end{eqnarray}
In the case $S=(0,1)$ we have by Lemma~\ref{l.S1+S2.m=2} (see Remark~\ref{r.(0,1)})
\begin{equation}
\label{S=(0,1)} 
y^{(1)}=(0,1),\,\,y^{(2)}=(0,1)\,\,\text{or}\,\,
\left(\begin{array}{c}
y^{(1)}\\
y^{(2)}
\end{array}\right)=\left(\begin{array}{cc}
0&1\\
0&1
\end{array}\right).
\end{equation}
Indeed,   $\Vert Y_1^{(1)} \Vert^2\sim S_1(2)<\infty$, hence $\Vert Y_2^{(1)} \Vert^2=\infty$ by \eqref{i_1,i_2.m=2}.
Further, $\Gamma(Y_2^{(2)})=\Vert Y_2^{(2)} \Vert^2\sim S_2(2)=\infty$, and $\Gamma(Y_1^{(2)})<\infty$, since 
$S_1(2)<\infty$. Therefore,
$$
\Delta(Y_2^{(2)},Y_1^{(2)})=\frac{\Gamma(Y_2^{(2)})+\Gamma(Y_2^{(2)},Y_1^{(2)})}{1+\Gamma(Y_1^{(2)})}>
\frac{\Gamma(Y_2^{(2)})}{1+\Gamma(Y_1^{(2)})}=\infty,
$$
and we conclude that $x_{2n}\,\,\eta\,\,{\mathfrak A}^2$. 

(2) Further, since 
\begin{equation}
\label{Y_1,Y_2,Y_1+Y_2}
\Vert Y_1 \Vert^2= \Vert Y_2 \Vert^2= \Vert Y_1 -sY_2 \Vert^2=\infty
\end{equation}
by Lemma~\ref{l.min=proj}, we conclude that $\Delta(Y_1,Y_2)=\Delta(Y_2,Y_1)=\infty$,
so $D_{1n},\,D_{2n}\,\,\eta\,\,{\mathfrak A}^2 $ by Lemmas~\ref{d1} and \ref{d2}.
Finally, $x_{2n},\,D_{1n},\,D_{2n}\,\eta\,\,{\mathfrak A}^2$. Now we get
$$
A_{kn}-x_{2k}D_{2n}=x_{1k}D_{1n},\quad k,n\in {\mathbb Z},
$$ 
and the proof is complete since
we are in the case $m=1$.

Relations \eqref{Y_1,Y_2,Y_1+Y_2} follows from 
$S_1(2)\!=\!\Sigma^{12}\!=\!\sum_{k\in \mathbb Z}\frac{b_{1k}}{b_{2k}}<\infty$. Indeed, we have
$$
\Vert Y_1\Vert^2=\sum_{k\in{\mathbb Z}}\frac{a_{1k}^2}{\frac{1}{2b_{1k}}+\frac{1}{2b_{2k}}}=
\sum_{k\in{\mathbb Z}}\frac{b_{1k}a_{1k}^2}{\frac{1}{2}+\frac{b_{1k}}{2b_{2k}}}\sim 2\sum_{k\in{\mathbb Z}}b_{1k}a_{1k}^2=
S^L_{11}(\mu)=\infty,
$$
$$
\Vert Y_2\Vert^2=
\sum_{k\in{\mathbb Z}}\frac{b_{1k}a_{2k}^2}{\frac{1}{2}+\frac{b_{1k}}{2b_{2k}}}\sim
\sum_{k\in{\mathbb Z}}b_{1k}a_{2k}^2\sim
 \sum_{k\in{\mathbb Z}}\frac{b_{1k}}{2}\Big(\frac{1}{2b_{2k}}+a_{2k}^2\Big)=S^L_{12}(\mu)=\infty,
$$
$$
\Vert Y_1-sY_2\Vert^2\!=\!\!\sum_{k\in{\mathbb Z}}\frac{b_{1k}(a_{1k}\!-\!s a_{2k})^2}{\frac{1}{2}+\frac{b_{1k}}{2b_{2k}}}\!\sim\!
\sum_{k\in{\mathbb Z}}b_{1k}(a_{1k}\!-\!s a_{2k})^2\!=\!\frac{1}{4}\sum_{k\in{\mathbb Z}}b_{1k}(-2a_{1k}+2sa_{2k})^2
$$
$$
\sim
\frac{1}{2}\Big[\frac{(2s)^2}{4}\sum_{k\in{\mathbb
Z}}\frac{b_{1k}}{b_{2k}}+ \sum_{k\in{\mathbb
Z}}\frac{b_{1k}}{2}(-2a_{1k}+2sa_{2k})^2\Big]
= \frac{1}{2}S^{L,-}_{12}(\mu,t)=\infty,
$$
for $t=2s$.
\qed\end{pf}
\subsubsection{Case  $S=(1,1)$}
In the case $S=(1,1)$ we have three possibilities (see Remark~\ref{r.(0,1)} )
\begin{equation}
\label{S=(1,1)}
{\rm I(2a)}\,\,y=\left(\begin{array}{cc}
1&0\\
1&1
\end{array}\right),\,\,
{\rm I(2b)}\,\,y=\left(\begin{array}{cc}
1&1\\
0&1
\end{array}\right),\,\,
{\rm I(2c)}\,\,y=\left(\begin{array}{cc}
1&1\\
1&1
\end{array}\right),
\end{equation}
i.e., {\rm I(2a)}\,\, $y_2^{(1)}<\infty$ and $y_1^{(2)}=\infty$,
{\rm I(2b)}\,\,$y_2^{(1)}=\infty$ and $y_1^{(2)}<\infty$,
{\rm I(2c)}\,\, $y_2^{(1)}=y_1^{(2)}=\infty$,
since $\Vert Y_1^{(1)}\Vert^2\sim S_1(2)=\infty$ and $\Vert Y_2^{(2)}\Vert^2\sim S_2(2)=\infty$
and  $y_2^{(1)}+y_1^{(2)}=\infty$ by \eqref{i_1,i_2.m=2}.

In the first case I(2a) or in the second case I(2b), i.e., if $\Vert Y_2^{(1)}\Vert^2\!<\!\infty$  or $\Vert Y_1^{(2)}\Vert^2\!<\!\infty$, we conclude  
respectively that
$\Delta(Y_1^{(1)},Y_2^{(1)})\!=\!\infty$ or $\Delta(Y_2^{(2)},Y_1^{(2)})=\infty$ hence, $x_{1n}x_{1t}\,\,\eta\,\, {\mathfrak A} ^2$ or 
$x_{2n}x_{2t}\,\,\eta\,\, {\mathfrak A} ^2$. So, we will get respectively  $x_{1n}$ or $x_{2n}$.
We shall come back to these cases later.

It remains to consider the case I(2c): $\Vert Y_2^{(1)}\Vert^2=\Vert Y_1^{(2)}\Vert^2=\infty$.
In the case $S=(1,1)$ set $c_n=\frac{b_{2n}}{b_{1n}},\,\,n\in{\mathbb Z}$. 
Then  $\sum_{k\in{\mathbb Z}}\frac{b_{1k}}{b_{2k}}\!=\!\sum_{k\in{\mathbb Z}}\frac{b_{2k}}{b_{1k}}\!=\!\infty.$
We get
\begin{eqnarray}
\label{Y11,Y12}
&&\Vert Y_1^{(1)}\Vert^2=\sum_{n\in \mathbb Z}\frac{1}{1+2c_n}\sim\sum_{n\in \mathbb Z}\frac{1}{c_n},\quad
\Vert Y_2^{(1)}\Vert^2=\sum_{n\in \mathbb Z}\frac{c_n^2}{1+2c_n},\\
\label{Y21,Y22}
&&\Vert Y_1^{(2)}\Vert^2=\sum_{n\in \mathbb Z}\frac{1}{c_n^2+2c_n},\quad
\Vert Y_2^{(2)}\Vert^2=\sum_{n\in \mathbb Z}\frac{c_n^2}{c_n^2+2c_n}\sim\sum_{n\in \mathbb Z}c_n.
\end{eqnarray}

\subsubsection{II, Approximation of $D_{1n}$ and $D_{2n}$}
Set
\begin{eqnarray}
\label{y=.m=2}
&&y_{12}=(y_1,y_2)=(\Vert Y_1\Vert^2,\Vert Y_2\Vert^2),\\
\label{Y1,Y2.m=2}
&&\text{where}\quad\Vert Y_1\Vert^2=\sum_{n\in {\mathbb Z}}\frac{a_{1n}^2}{\frac{1}{2b_{1n}}+\frac{1}{2b_{2n}}},\quad
\Vert Y_2\Vert^2=\sum_{n\in {\mathbb Z}}\frac{a_{2n}^2}{\frac{1}{2b_{1n}}+\frac{1}{2b_{2n}}}.
\end{eqnarray}

The case ${\rm II}$ splits into four subcases
\begin{equation}
\label{m=2.II-4-cas}
(1)\,(y_1,y_2)\!=\!(1,0),\,\,(2)\,(y_1,y_2)\!=\!(0,1),\,\,(3)\,(y_1,y_2)\!=\!(1,1),\,\,(4)\,(y_1,y_2)\!=\!(0,0).
\end{equation}
We have $4=2^2$ possibilities for $y=(y_1,y_2)\in \{0,1\}^2$:
$$
\begin{array}{cccccccc}
   &(1)&(2)&(3a)&(3b)&(3c)&
   (4)\\
y_1&1&0&1&1&1&0\\
y_2&0&1&1&1&1&0\\
\alpha&&&1&0&C_1\leq \alpha_m\leq C_2&
\end{array}
$$
where $\alpha=\lim_{m\to\infty}\alpha_m,$ with $\alpha_m\!=\!\frac{\Vert Y_1(m)\Vert^2}{\Vert Y_2(m)\Vert^2}$, and
$$
\Vert Y_1(m) \Vert^2=\sum_{n=-m}^{m}\frac{a_{1n}^2}{\frac{1}{2b_{1n}}+\frac{1}{2b_{2n}}},\quad
\Vert Y_2(m) \Vert^2=\sum_{n=-m}^{m}\frac{a_{2n}^2}{\frac{1}{2b_{1n}}+\frac{1}{2b_{2n}}}.
$$
All the different cases are presented in the following table:
%
 \vskip 0.1cm
\begin{tabular}{|p{2.34cm}|p{1.4cm}|p{1.4cm}|p{1.4cm}|p{1.4cm}|p{1.5cm}|p{0.9cm}|
} \hline
\text{Table}\,\,II &(1)&(2)&(3a)&(3b)&(3c)&(4)
\\   \hline
$\Vert Y_1 \Vert^2$&$\infty$&$<\infty$&$\infty$&$\infty$&$\infty$&$<\infty$
\\   \hline
$\Vert Y_2 \Vert^2$&$<\infty$&$\infty$&$\infty$&$\infty$&$\infty$&$<\infty$
\\   \hline
$
\alpha_m\!=\!\frac{\Vert Y_1(m)\Vert^2}{\Vert Y_2(m)\Vert^2}
$&&&$\to\infty$&$\to 0$&
$C_1\!\!\leq\!\!\alpha_m
\!\!\leq\!\!C_2$&
\\   \hline
Lemma &\ref{d1} &\ref{d2}&\ref{d1}   &\ref{d2} &\ref{d1} ,\,\,\ref{d2}   &
\\
&\ref{x1,dA} &\ref{x2,dA}&\ref{x1,dA}&\ref{x2,dA}&\ref{the-last},\,\,\ref{l.min=proj}  &   \\
\hline
&$D_{1n},\,\,x_{1n}$&$D_{2n},\,\,x_{2n}$&$D_{1n},\,\,x_{1n}$&$D_{2n},\,\,x_{2n}$&$D_{1n},\,\,D_{2n}$&
\\   \hline
\end{tabular}
 \vskip 0.3cm
\begin{rem}
\label{r.II.1} 
We show that if $\Vert Y_2 \Vert^2<\infty$ and
$S^L_{12}(\mu)=\infty$, then $\sum_n\frac{b_{1n}}{b_{2n}}=\infty$.
Indeed, let us suppose that $\sum_n\frac{b_{1n}}{b_{2n}}<\infty$,
then
\begin{equation}
\label{y2=inf}
\Vert Y_2 \Vert^2=\sum_{n\in{\mathbb Z}}\frac{a_{2n}^2}{\frac{1}{2b_{1n}}+\frac{1}{2b_{2n}}}
\sim
\sum_{n\in{\mathbb Z}}b_{1n}a_{2n}^2
\sim
\sum_{n\in{\mathbb Z}}\frac{b_{1n}}{2}\Big(\frac{1}{2b_{2n}}+a_{2n}^2\Big)=S_{12}^L(\mu)=\infty.
\end{equation}
\end{rem}

We explain  Table ${\rm II}$ in detail. 
 The first two case (1) and (2) are independent
of the case {\rm I(2)}, i.e.,
$S=(1,1)$.

{\bf (1)} If $\Vert Y_2 \Vert^2<\infty$ and $\Vert Y_1 \Vert^2=\infty$, then
 $D_{1k}\,\,\eta\,\,{\mathfrak A}^2$ by Lemma \ref{d1}. The
condition $\Vert Y_2 \Vert^2<\infty$ implies $\sum_{k\in{\mathbb
Z}}\frac{b_{1k}}{b_{2k}}=\infty$, by  Remark~\ref{r.II.1}
therefore, $x_{1k}\,\,\eta\,\,{\mathfrak A}^2$,  by Lemma
\ref{x1,dA}. Further, $A_{kn}-x_{1k}D_{1n}=x_{2k}D_{2n},\,\,k,n\in
{\mathbb Z}$ and the proof is complete since we are reduced to the
case $m=1$.

{\bf (2)} If $\Vert Y_2 \Vert^2=\infty$ and $\Vert Y_1 \Vert^2<\infty$, then
  $D_{2k}\,\,\eta\,\,{\mathfrak A}^2$ by Lemma~\ref{d2}. 
Reasoning as in  
 Remark~\ref{r.II.1}, we conclude that $\sum_{k\in{\mathbb Z}}\frac{b_{2k}}{b_{1k}}=\infty$ and therefore,
 $x_{2k}\,\,\eta\,\,{\mathfrak A}^2$ by Lemma~\ref{x2,dA} and $A_{kn}-x_{2k}D_{2n}=x_{1k}D_{1n},\,\,k,n\in {\mathbb Z}$, case $m=1$.

{\bf (3)}  Consider now the case  I(2). 
Suppose that both series are divergent: $\Vert Y_2 \Vert^2\!=\!\infty$ and
$\Vert Y_1 \Vert^2\!=\!\infty$.  We show that in the case $(B)$   {\rm (see (\ref{A,B}))} holds $$\Vert Y_1+sY_2\Vert^2\!=\!\infty\quad\text{ for all}\quad s\in {\mathbb R}$$
by Lemma~\ref{the-last} therefore, 
by Lemma~\ref{l.min=proj}, we can approximate $D_{1n}$ and $D_{2n}$.
To be more precise  consider three possibilities:

{\bf (3a)} Let $\frac{\Vert Y_1(m)\Vert^2}{\Vert Y_2(m)\Vert^2}\to\infty$, then
$D_{1k}\,\,\eta\,\,{\mathfrak A}^2$. Since $\sum_n\frac{b_{1n}}{b_{2n}}=\infty$, case I(2), we have
$x_{1n}\,\,\eta\,\,{\mathfrak A}^2$ by  Lemma \ref{x1,dA} and
finally, $x_{1n},\,D_{1n}\,\,\eta\,\,{\mathfrak A}^2,\,\,n\in
{\mathbb Z}$. We are reduced to the case $m=1$.

{\bf (3b)}  let $\frac{\Vert Y_1(m)\Vert^2}{\Vert Y_2(m)\Vert^2}\to 0$, then
$D_{2k}\,\,\eta\,\,{\mathfrak A}^2$.  Since
$\sum_n\frac{b_{2n}}{b_{1n}}=\infty$, case I(2), we get
$x_{2n}\,\,\eta\,\,{\mathfrak A}$, by  Lemma \ref{x2,dA} and
finally, $x_{2n},\,D_{2n}\,\,\eta\,\,{\mathfrak A}^2,\,\,n\in
{\mathbb Z}$. We are reduced to the case $m=1$.

{\bf (3c)} The case when  $\Vert Y_1 \Vert^2\!=\!\Vert Y_2 \Vert^2\!=\!\infty$ and $C_1\!\leq
\frac{\Vert Y_1(m)\Vert^2}{\Vert Y_2(m)\Vert^2}\leq C_2,\,\,\,
m\in \mathbb N$.

{\bf (4)} The case when $\Vert Y_1 \Vert^2+\Vert Y_2 \Vert^2<\infty.$

{\it To complete the proof} of the lemma it remains to consider I(2), i.e., $S=(1,1)$
and the last two  cases in the table ${\rm
II}$, i.e., ${\rm II(3c)}$ and  ${\rm II(4)}$, where:
\begin{equation}
\label{I(3c3)}
{\rm I(2)}\quad
\sum_{k\in{\mathbb Z}}\frac{b_{1k}}{b_{2k}}=\sum_{k\in{\mathbb Z}}\frac{b_{2k}}{b_{1k}}=\infty,
\end{equation}
\begin{equation}
{\rm II(3c)}\quad \sum_{k\in{\mathbb Z}}a_{1k}^2\Big(\frac{1}{2b_{1k}}+\frac{1}{2b_{2k}}\Big)^{-1}=
\sum_{k\in{\mathbb Z}}a_{2k}^2\Big(\frac{1}{2b_{1k}}+\frac{1}{2b_{2k}}\Big)^{-1}=\infty,
\end{equation}
\begin{equation}
\label{I+II(d)}
{\rm II(4)}\quad
\sum_{k\in{\mathbb Z}}\Big(a_{1k}^2+a_{2k}^2\Big)\Big(\frac{1}{2b_{1k}}+\frac{1}{2b_{2k}}\Big)^{-1}<\infty.
\end{equation}

We come back to the condition $\mu^{L_t}\perp\mu$. By Remark~\ref{perp2-1} we have
$$
\mu^{L_{\tau_-(\phi,s)}}\perp\mu,\,\,\,\phi\in[0,2\pi),\,\,s>0\Leftrightarrow
\Sigma_1(s)+\Sigma_2(C_1,C_2)\!=\!\infty,\,\,s>0,
$$
for $(C_1,C_2)\in{\mathbb R}^2\setminus\{0\}$.
To make the notation consistent for the case $m=3$ we replace everywhere $\Sigma_1(s)$ (defined by  \eqref{sigma1(s)}) with $\Sigma_{12}(s)$ and 
$\Sigma_2(C_1,C_2)$  defined by \eqref{sigma2(s)} with
$\Sigma_{12}(C_1,C_2)$ for $
(C_1,C_2)\in {\mathbb R}^2\setminus\{0\}
$:
\begin{eqnarray}
\label{sigma12(s)}
&&
\Sigma_{12}(s)
=\sum_{n\in {\mathbb Z}}
\Big(s^2\sqrt{\frac{b_{1n}}{b_{2n}}}-s^{-2}\!\sqrt{\frac{b_{2n}}{b_{1n}}}\Big)^2,\quad s\in \mathbb R\setminus\{0\}, \\
\label{Sigma12(s)}
&&
\Sigma_{12}(C_1,C_2)=\sum_{n\in {\mathbb Z}}(C^2_1b_{1n}+C^2_2b_{2n})(C_1a_{1n}+C_2a_{2n})^2.
\end{eqnarray}
The condition $\Sigma_{12}(s)+\Sigma_{12}(C_1,C_2)\!=\!\infty$, splits into two cases:
\begin{equation}
\label{A,B}
\begin{array}{cccc}
(A) & \Sigma_{12}(s)=\infty,&&\\
(B) &\Sigma_{12}(s)<\infty&\text{and}& \Sigma_{12}(C_1,C_2)=\infty.
\end{array}
\end{equation}

\vskip 0.5cm 
 Finally, we need to consider the  following 12 possibilities:
 \vskip 0.2 cm
\begin{tabular}{|c|c|c|c|} 
 \hline
A&${\rm I}(3a)$&${\rm I}(3b)$&${\rm I}(3c)$\\
 \hline
${\rm II}(3c)$&&&\\ 
 \hline
 $ {\rm II}(4)$&&&\\ 
 \hline
\end{tabular}\,.
\quad
\begin{tabular}{|c|c|c|c|} 
 \hline
B&${\rm I}(3a)$&${\rm I}(3b)$&${\rm I}(3c)$\\
 \hline
${\rm II}(3c)$&&&\\ 
 \hline
 $ {\rm II}(4)$&&&\\ 
 \hline
\end{tabular}\,.

 \vskip 0.3 cm
{\it Briefly}:\\
 ${\bf (A)\&I(2)}$.  In this case
independently of the conditions ${\rm II(3c)}$ and ${\rm II(4)}$ we can approximate
$x_{1n}$ and $x_{2n}$ using Lemma \ref{x1x1} and \ref{x2x2}.\\
${\bf (B)\&II(3c)}$ In this case
we can  approximate $D_{1n}$ and $D_{2n}$ using Lemmas~\ref{d1} and
\ref{d2} respectively. More precisely, to be able to use Lemma~\ref{l.min=proj} we show that conditions
(\ref{norm=infty}) are satisfied for the two vectors $Y_1$ and $Y_2$ defined by \eqref{Y_r^{(s)}.m=2} (see Lemma~\ref{the-last}).\\
${\bf (B)\&II(4)}$  Since $\Sigma_{12}(C_1,C_2)=\infty$, this  case
(see (\ref{I+II(d)})) cannot be realized 
\vskip 0.1 cm
{\it More details}:\\
{\bf Case} ${\rm (A)\&I(2)}$. Using Lemma~\ref{l.min=proj} we  conclude that
 \begin{equation}
\label{xx3}
\Delta(Y_1^{(1)},Y_2^{(1)})=\infty\quad \text{and} \quad \Delta(Y_2^{(2)},Y_1^{(2)})=\infty.
\end{equation}
To use Lemma~\ref{l.min=proj}, it is sufficient to show that in the case $(A)$ relations (\ref{norm=infty}) hold for
$Y_1^{(1)},Y_2^{(1)}$ and $Y_2^{(2)},Y_1^{(2)}$, i.e., for all
$s\in {\mathbb R}\setminus\{0\}$ we have (see Lemma~\ref{l.park})
\begin{eqnarray}
 \label{norm=infty(1,2)}
&&\Vert Y_1^{(1)}\Vert^2\!=\!\Vert Y_2^{(1)}\Vert^2\!=\!\Vert Y_1^{(1)}+s Y_2^{(1)}\Vert^2=\infty,\\
\nonumber
&&\Vert Y_2^{(2)}\Vert^2\!=\!\Vert Y_1^{(2)}\Vert^2\!=\!\Vert Y_2^{(2)}+s Y_1^{(2)}\Vert^2=\infty.
\end{eqnarray}
Consider the following three possibilities in the case ${\rm I(2)}$:\\
%

I(2a) If $\Vert Y_2^{(1)}\Vert<\infty$, then $\Vert Y_1^{(1)}\Vert=\infty$ by \eqref{i_1,i_2.m=2} therefore, 
$\Delta(Y_1^{(1)},Y_2^{(1)})\\= \infty$
so, $x_{1n}\,\,\eta\,\,{\mathfrak A}^2$ by Lemma~\ref{l.min=proj0} (a). In the case $(A)$ by Lemma~\ref{l.park} holds
$$\Vert Y_2^{(2)} \Vert^2=\Vert Y_1^{(2)} \Vert^2=\Vert Y_2^{(2)} +s Y_1^{(2)} \Vert^2=\infty,$$ therefore, 
$x_{2n}\,\,\eta\,\,{\mathfrak A}^2$ by Lemma~\ref{l.min=proj}.

I(2b) If $\Vert  Y_1^{(2)} \Vert<\infty$, then $\Vert  Y_2^{(2)} \Vert=\infty$ by \eqref{i_1,i_2.m=2} therefore,
$\Delta(Y_2^{(2)},Y_1^{(2)})\\=\infty$ so,
$x_{2n}\,\,\eta\,\,{\mathfrak A}^2$ by Lemma~\ref{l.min=proj0} (a). In the case $(A)$ by Lemma~\ref{l.park} we have
$$\Vert Y_2^{(1)} \Vert^2=\Vert Y_1^{(1)} \Vert^2=\Vert Y_1^{(1)} +s Y_2^{(2)} \Vert^2=\infty,$$
and therefore, $x_{1n}\,\,\eta\,\,{\mathfrak A}^2$ by Lemma~\ref{l.min=proj}.

I(2c) If $\Vert Y_2^{(1)} \Vert=\Vert Y_1^{(2)}\Vert=\infty$, then by Lemma~\ref{l.park} all relations (\ref{norm=infty(1,2)}) hold in the case $(A)$ and
therefore, $x_{1n},\,\,x_{2n}\,\,\eta\,\,{\mathfrak A}^2$.
To prove (\ref{norm=infty(1,2)}) we need Lemma~\ref{l.light.123}.
\begin{lem}
\label{l.park} If $\Sigma_{12}(s)=\infty$ for any $s>0$, then
$$
\Vert Y_1^{(1)}-CY_2^{(1)}\Vert^2=\infty\quad\text{and}\quad \Vert Y_2^{(2)}-CY_1^{(2)}\Vert^2=\infty,\quad\text{for any}\quad 
C\in \mathbb R\setminus\{0\}.
$$
\end{lem}
So, in the case ${\rm (A)\&I(2)}$  we can approximate $x_{1n}$ and $x_{2n}$.

{\bf Case} ${\rm (B)\&II(3c)}$.
\begin{lem}
\label{the-last} When 
$\Sigma_{12}(s)<\infty$ and
$\Sigma_{12}(C_1,C_2)=\infty$,
we get
\begin{equation}
\label{C_1f+C_2g=}
\sigma(C_1,C_2):=
\Vert C_1Y_1+C_2Y_2\Vert^2\!=\!\sum_{n\in {\mathbb Z}}\frac{(C_1a_{1n}+C_2a_{2n})^2}{\frac{1}{2b_{1n}}+\frac{1}{2b_{2n}}}\!=\!\infty,\,\,
(C_1,C_2)\!\in\!{\mathbb R}^2\setminus\{0\}.
\end{equation}
\end{lem}
Finally, we can approximate $D_{1n}$ and $D_{2n}$ in the case ${\rm (B)\&II(3c)}$.

{\bf Case} ${\rm (B)\&II(4)}$. The last case  ${\rm (B)\&II(4)}$ (see (\ref{I+II(d)})) can not be realized.
Indeed, in this case
$\Sigma_{12}(s)<\infty$ and $\Sigma_{12}(C_1,C_2)=\infty$. Therefore, by Lemma~\ref{l.light.123} 
we have $s^4\lim_{n\to\infty}\frac{b_{1n}}{b_{2n}}=1$ and hence,  
$$ \sigma(C_1,C_2)\sim \Sigma_{12}(C_1,C_2)=\infty
$$
(see also the proof of  Lemma~\ref{the-last}  in \cite{KosJFA14}). 
This contradicts (\ref{I+II(d)}):
 $$
 \sigma(1,1)=\sum_{k\in{\mathbb Z}}\big(a_{1k}^2+a_{2k}^2\big)
\Big(\frac{1}{2b_{1k}}+\frac{1}{2b_{2k}}\Big)^{-1}<\infty.
$$ Thus   the proof of Lemma~\ref{approx-(x,D).2} for $m=2$ is completed.
%

The proof of  the irreducibility for $m=2$ follows  from
Remark~\ref{r.irr-idea-G}. Depending on the measure, we can
approximate four different families of commuting operators
$B^\alpha=(B^\alpha_{1n},B^\alpha_{2n})_{n\in{\mathbb Z}}$ for
$\alpha\in \{0,1\}^2$:
$$
B^{(0,0)}\!\!=\!\!(x_{1n},x_{2n})_{n},\,\,B^{(0,1)}\!\!=\!\!(x_{1n},D_{2n})_{n},\,\,
B^{(1,0)}\!\!=\!\!(D_{1n},x_{2n})_{n},\,\,B^{(0,0)}\!\!=\!\!(D_{1n},D_{2n})_{n}.
$$
The von Neumann algebra $L^\infty_\alpha(X_2,\mu^2)$ consists of
all essentially bounded functions $f(B^\alpha)$ in the commuting
family of operators $B^\alpha$ (see, e.g., \cite{Ber86}) as, in
particular, $L^\infty_{(0,0)}(X_2,\mu^2)=L^\infty(X_2,\mu^2)$.
Since the von Neumann algebras $L^\infty_\alpha(X_2,\mu^2)$ are
maximal abelian, the commutant $\big({\mathfrak A}^2\big)'$ of the
von Neumann algebra ${\mathfrak A}^2$ generated by the
representation is contained in  $L^\infty_\alpha(X_2,\mu^2)$.
Hence, the bounded operator $A\in \big({\mathfrak A}^2\big)'$ will
be some function $A=a(B^\alpha)\in L^\infty_\alpha(X_2,\mu^2)$.
The commutation relation $[A,T^{R,\mu,2}_t]=0$ gives us the
following relations:
$a(B^\alpha)^{R_t}=a(B^\alpha)$ for all $t\in {\rm
GL}_0(2\infty,{\mathbb R}).$
Set $B^\alpha_r=(B^\alpha_{rn})_n,\,\,x_r=(x_{rn})_n,\,\,D_r=(x_{rn})_n,\,\,r=1,2,\,\,n\in{\mathbb
Z}$ and set as before, $E_{kn}(t):=I+tE_{kn},\,\,t\in {\mathbb R},\,\,k,n\in
{\mathbb Z},\,\,k\not=n$. Then the action $(B^\alpha)^{R_t}$ is
defined as follows:
%
\begin{eqnarray*}
(B^\alpha_1,B^\alpha_2)^{R_t}=((B^\alpha_1)^{R_t},(B^\alpha_2)^{R_t}),\quad (x_r)^{R_t}=x_rt,\quad (D_r)^{R_t}=D_rt^T,\\
a(\dots, x_{rk},\dots ,x_{rn},\dots)^{R_{E_{kn}(t)}}:=a(\dots, x_{rk},\dots ,x_{rn}+tx_{rk},\dots),\\
a(\dots, D_{rk},\dots ,D_{rn},\dots)^{R_{E_{kn}(t)}}:=a(\dots, D_{rk}+tD_{rn},\dots ,D_{rn},\dots),\,\,\,t\in {\mathbb R}.
\end{eqnarray*}
In all the cases, by ergodicity of the measure $\mu^2$, we conclude that $a$ is constant.
\section{Irreducibility, the case $m=3$}
\label{sec.Irr.m=3}
\subsection{Technical part of the proof of irreducibility}
\label{sec.4.1}
\begin{lem}
\label{l.approx-(x,D).m=3}
If $\mu^{L_t}\perp\mu$ for all $t\in {\rm GL}(3,{\mathbb R})\setminus \{e\}
$, we can approximate at least  one of  the following  eight triplets of operators:
\begin{eqnarray*}
 &&(x_{1n},x_{2n},x_{3n}),\,\,(x_{1n},x_{2n},D_{3n}),\,\,(x_{1n},D_{2n},x_{3n}),\,\,(D_{1n},x_{2n},x_{3n}),\\
 &&(x_{1n},D_{2n},D_{3n}), (D_{1n},x_{2n},D_{3n}),\,\,(D_{1n},D_{2n},x_{3n}),\,\, (D_{1n},D_{2n},D_{3n}).
\end{eqnarray*}
\end{lem}
\begin{pf}
By Lemma~\ref{l.SL(3)-perp}, 
the condition of orthogonality $(\mu_{(b,a)}^3)^{L_t}\perp\mu_{(b,a)}^3$ 
for   $t\in \pm{\rm SL}(3,\mathbb R)\setminus\{e\}$ are,
\begin{equation}
 \label{III}
\Sigma^{\pm}(t)=\Sigma_1^{\pm}(t)+\Sigma_2(t)=\infty,
\end{equation}
where    $\Sigma_2(t)$  is defined by \eqref{Sigma_3} and 
$\Sigma_1^+(t),\,\,\Sigma_1^{-}(t)$ are defined by \eqref{Sigma_1^+(t).m=3}, \eqref{Sigma_1^{-}(t).m=3}.
Let ${\mathfrak A}^3$ be the von Neumann algebra generated by the representation. 
We write compactly:
\begin{eqnarray}
 \label{x,D-eta-A^3}
&&
x_{kn}\,\,\eta\,\, {\mathfrak A}^3\Leftrightarrow \Delta^{(k)}=\infty,\quad
D_{kn}\,\,\eta\,\, {\mathfrak A}^3\Leftrightarrow \Delta_k=\infty,\\
\label{Delta^k(krs)}
&&\text{where}\,\,\Delta^{(k)}:=\Delta(Y_k^{(k)},Y_r^{(k)},Y_s^{(k)}),\quad
\Delta_k:=\Delta(Y_k,Y_r,Y_s),
\end{eqnarray}
and $\{k,r,s\}$ is a cyclic permutation of $\{1,2,3\}$.

{\bf Case I.} {\it  Approximation of $x_{rk}x_{rt}$ for $1\leq r\leq 3$ by $A_{kn}A_{tn}$.}\\
Set $B_{3k}=b_{1k}+b_{2k}+b_{3k}$. To approximate the operators  $x_{kn}$  by the corresponding operators,
by  Lemmas~\ref{x1x1.3}--\ref{x3x3.3} we get:
\begin{equation} 
\label{I} 
x_{1n}x_{1t}\,\eta\, {\mathfrak A}^3\Leftrightarrow \,\Delta^{(1)}=\infty,\,\,
x_{2n}x_{2t}\,\,\eta\,\,{\mathfrak A}^3\Leftrightarrow \Delta^{(2)}=\infty,\,\,
x_{3n}x_{3t}\,\eta\, {\mathfrak A}^3\Leftrightarrow \Delta^{(3)}=\infty,
\end{equation}
%
where
\begin{equation}
\label{Y_r^{(s)}}
\Vert Y_s^{(r)} \Vert^2=\sum_{k\in \mathbb Z}
\frac{b_{rk}^2}{B_{3k}^2-(b_{1k}^2+b_{2k}^2+b_{3k}^2-b_{sk}^2)},\quad 1\leq r,s\leq 3. 
\end{equation}

{\bf Case II.} {\it  Approximation of $D_{rn}$ by $A_{kn}$.}\\

By Lemmas~\ref{d1.3}--\ref{d2.3} we have for $1\leq r\leq 3$ (see \eqref{Delta^k(krs)}):
\begin{equation}
\label{II}
D_{rn}\,\,\eta\,\,{\mathfrak A}^3\Leftrightarrow \Delta_r=\infty,\quad\text{where}\quad \Vert Y_r\Vert^2=\sum_{k\in \mathbb Z}\frac{a_{rk}^2}{\frac{1}{2b_{1k}}+\frac{1}{2b_{2k}}+\frac{1}{2b_{3k}}}.
\end{equation}

{\bf Case III.} {\it  Approximation of $x_{rk}$ by $D_{rn}A_{kn}$.}
\begin{eqnarray*}
&&D_{1n}A_{kn}=x_{1k}D_{1n}^2+x_{2k}D_{1n}D_{2n}+x_{3k}D_{1n}D_{3n},\\
&&D_{2n}A_{kn}=x_{1k}D_{1n}D_{2n}+x_{2k}D_{2n}^2  +x_{3k}D_{2n}D_{3n},\\
&&D_{3n}A_{kn}=x_{1k}D_{1n}D_{3n}+x_{2k}D_{2n}D_{3n}+x_{3k}D_{3n}^2.
\end{eqnarray*}
By Lemmas~\ref{x1,dA.m=3}--\ref{x3,dA.m=3}  we have 
\begin{eqnarray*}
&& 
x_{1n}{\bf 1}\in\langle D_{1n}A_{kn}{\bf 1}\mid k\in {\mathbb Z}\rangle
\quad\Leftrightarrow\quad 
\Sigma_1=\infty,\\
&& 
x_{2n}{\bf 1}\in\langle D_{2n}A_{kn}{\bf 1}\mid k\in {\mathbb Z}\rangle
\quad\Leftrightarrow\quad 
\Sigma_2=\infty,\\
&& 
x_{3n}{\bf 1}\in\langle D_{3n}A_{kn}{\bf 1}\mid k\in {\mathbb Z}\rangle
\quad\Leftrightarrow\quad 
\Sigma_3=\infty,
\end{eqnarray*}
where $\Sigma_r=\sum_{k\in \mathbb Z}\frac{b_{rk}}{b_{1k}+b_{2k}+b_{3k}}.$

{\bf Case IV.} {\it  Approximation of $D_{rn}$ by $x_{rk}A_{kn}$.}
\begin{eqnarray*}
&&x_{1k}A_{kn}=x^2_{1k}D_{1n}+x_{1k}x_{2k}D_{2n}+x_{1k}x_{3k}D_{3n},\\
&&x_{2k}A_{kn}=x_{1k}x_{2k}D_{1n}+x_{2k}^2D_{2n}+x_{2k}x_{3k}D_{3n},\\
&&x_{3k}A_{kn}=x_{1k}x_{3k}D_{1n}+x_{2k}x_{3k}D_{2n}+x_{3k}^2D_{3n}.
\end{eqnarray*}
By Lemmas~\ref{dx1.3}--\ref{dx3.3}  we have 
\begin{eqnarray*}
&& 
D_{1n}{\bf 1}\in\langle x_{1k}A_{kn}{\bf 1}\mid k\in {\mathbb Z}\rangle
\quad\Leftrightarrow\quad 
\Delta(Y_{11},Y_{12},Y_{13})=\infty,\\
&& 
D_{2n}{\bf 1}\in\langle x_{2k}A_{kn}{\bf 1}\mid k\in {\mathbb Z}\rangle
\quad\Leftrightarrow\quad 
\Delta(Y_{22},Y_{23},Y_{21})=\infty,\\
&& 
D_{3n}{\bf 1}\in\langle x_{3k}A_{kn}{\bf 1}\mid k\in {\mathbb Z}\rangle
\quad\Leftrightarrow\quad 
\Delta(Y_{33},Y_{31},Y_{32})=\infty,
\end{eqnarray*}
where $Y_{kr}$ for $1\leq k,r\leq 3$ are defined by 
\eqref{y_1(123)}--\eqref{y_3(123)}.

{\bf Case V.}  By Lemma~\ref{l.SL(3)-perp} we have two conditions 
\begin{eqnarray}
\nonumber
&&(A)\quad \Sigma_1^+(t)=\infty,\,\,\,{\rm det}\,t=1\quad \text{or} \quad \Sigma_1^-(t)=\infty,\,\,\,{\rm det}\,t=-1,\\
\label{(A)(B).m=3}
&&(B)\quad \Sigma_1^+(t)<\infty \quad \text{or}\quad  \Sigma_1^-(t)<\infty\quad \text{but}\quad\Sigma_2(t)=\infty,
\end{eqnarray}
where $\Sigma_1^+(t),\,\,\Sigma_1^{-}(t),\quad\Sigma_2(t)$ are defined respectively by 
\eqref{Sigma_1^+(t).m=3}, \eqref{Sigma_1^{-}(t).m=3}  
and \eqref{Sigma_3}.
The rest
of this section will be devoted to the proof of Lemma~\ref{l.approx-(x,D).m=3}.
\qed\end{pf}

\subsubsection{Notations and the change of the variables}

In what follows we will systematically use the following notations:
\begin{eqnarray}
&& 
\label{S_r(3).1}
S_r(3)=\sum_{n\in \mathbb Z}\frac{b_{rn}^2}{b_{1n}b_{2n}+b_{1n}b_{3n}+b_{2n}b_{3n}},\quad 1\leq r\leq 3,\\
\label{Sigma_r}
&&\Sigma_r:=\sum_{n\in \mathbb Z}\frac{b_{rn}}{b_{1n}+b_{2n}+b_{3n}},\quad  1\leq r\leq 3,\\
\label{Sigma^(rs)}
&&\Sigma^{rs}:=\sum_{k\in{\mathbb Z}}\frac{b_{rk}}{b_{sk}},\quad 1\leq r\not=s\leq 3,
\quad
C_k=\frac{1}{2b_{1k}}+\frac{1}{2b_{2k}}+\frac{1}{2b_{3k}},
\\
\label{y(123)}
&&y_{123}=(y_1,y_2,y_3),\quad\text{where}\quad y_r:=\Vert Y_r\Vert^2,\\ 
 \label{y^(k).m=3}
&&y^{(k)}\!=\!(y_1^{(k)},y_2^{(k)},y_3^{(k)})\!=\!(\Vert Y_1^{(k)}\Vert^2,\Vert Y_2^{(k)}\Vert^2,\Vert Y_3^{(k)}\Vert^2),\,1\leq k\leq 3,\\
\label{y-mat.m=3}
&&y=\!
\left(\begin{array}{c}
y^{(1)}\\
y^{(2)}\\
y^{(3)}
\end{array}\right)=\left(\begin{array}{ccc}
y_1^{(1)}&y_2^{(1)}&y_3^{(1)}\\
y_1^{(2)}&y_2^{(2)}&y_3^{(2)}\\
y_1^{(3)}&y_2^{(3)}&y_3^{(3)}\\
\end{array}\right),\,\,\text{where}\,\, y_s^{(r)}:=\Vert Y_s^{(r)}\Vert^2,\,\,\,\\
\label{Sigma(s)} 
&&
\Sigma_{123}(s)
=(\Sigma_{12}(s_{12}),\Sigma_{23}(s_{23}),\Sigma_{13}(s_{13})),\quad s=(s_{12},s_{23},s_{13}).
\end{eqnarray}
The expressions $S_r(3)$ in the case $m=3$ can be generalized for an arbitrary  $m\in \mathbb N$ as follows:
\begin{equation}
\label{S_k(m)}
S_k(m)=\sum_{n\in \mathbb Z}\frac{b_{kn}^2}{\sum_{1\leq r<s\leq m}b_{rn}b_{sn}},\quad 1\leq k\leq m.
\end{equation}
\begin{lem}
\label{l.S1+S2+S3.m=3}
We have 
\begin{eqnarray}
 \label{S1+S2+S3}
&& S_1(3)+S_2(3)+S_3(3)=\infty,\\
\label{Yr-sim-Sr}
&& \Vert Y_r^{(r)}\Vert^2\sim S_r(3)\quad{for\,\,\, all}\quad 1\leq r\leq 3,\\
\label{Yr^s-<-Sr}
&& \Vert Y_r^{(s)}\Vert^2< \frac{1}{2}S_r(3)\quad{for\,\,\, all}\quad 1\leq r\not=s\leq 3,\\
\label{i_r.m=3} 
&&
\Vert Y_1^{(i_1)}\Vert^2+\Vert Y_2^{(i_2)}\Vert^2+\Vert Y_3^{(i_3)}\Vert^2=\infty,\quad 
i_1,i_2,i_3\in \{1,2,3\}.
\end{eqnarray}
\end{lem}
\begin{pf} Since $3(a^2+b^2+c^2)\geq 2(ab+ac+bc)$ we get 
 \begin{equation*}
  S_1(3)+S_2(3)+S_3(3)=\sum_{n\in \mathbb Z}\frac{b_{1n}^2+b_{2n}^2+b_{3n}^2}{b_{1n}b_{2n}+b_{1n}b_{3n}+b_{2n}b_{3n}}\geq\sum_{k\in \mathbb Z}2/3=\infty.
\end{equation*}
Further by \eqref{Y_r^{(s)}}
\begin{eqnarray*}
&& \Vert Y_r^{(r)}\Vert^2 =\sum_{n\in \mathbb Z}\frac{b_{rn}^2}{b_{rn}^2+2(b_{1n}b_{2n}+b_{1n}b_{3n}+b_{2n}b_{3n})}\sim S_r(3),\\
&& \Vert Y_r^{(s)}\Vert^2 =\sum_{n\in \mathbb Z}\frac{b_{rn}^2}{b_{sn}^2+2(b_{1n}b_{2n}+b_{1n}b_{3n}+b_{2n}b_{3n})}<\frac{1}{2} S_r(3),\,\,s\not=r.
\end{eqnarray*}
To prove \eqref{i_r.m=3} we get by \eqref{Y_r^{(s)}} 
\begin{eqnarray*}
&&\Vert Y_1^{(i_1)}\Vert^2+\Vert Y_2^{(i_2)}\Vert^2+\Vert Y_3^{(i_3)}\Vert^2=\\
&&
\sum_{r=1}^3\sum_{n\in \mathbb Z}\frac{b_{rn}^2}
{b_{i_rn}^2+2(b_{1n}b_{2n}+b_{1n}b_{3n}+b_{2n}b_{3n})}
>
\sum_{n\in \mathbb Z}\frac{\sum_{r=1}^3b_{rn}^2}{(\sum_{r=1}^3b_{rn})^2}=\infty.\hskip 2.3cm\Box
\end{eqnarray*}
\end{pf}
%
We make the following change of the variables:
\begin{equation}
 \label{(b,a)-(b',a')}
\left(\begin{smallmatrix}
b_{1n}&b_{2n}&b_{3n}\\
a_{1n}&a_{2n}&a_{3n}
\end{smallmatrix}
\right)\to
\left(\begin{smallmatrix}
b_{1n}'&b_{2n}'&b_{3n}'\\
a_{1n}'&a_{2n}'&a_{3n}'
\end{smallmatrix}
\right)=
\left(\begin{smallmatrix}
1&d_{2n}:=\frac{b_{2n}}{b_{1n}}&d_{3n}:=\frac{b_{3n}}{b_{1n}}\\
a_{1n}\sqrt{b_{1n}}&a_{2n}\sqrt{b_{1n}}&a_{3n}\sqrt{b_{1n}}
\end{smallmatrix}
\right),
\end{equation}
motivated by
the following formulas:
\begin{eqnarray*}
&&d\mu_{(b,a)}(x)\!=\!\sqrt{\frac{b}{\pi}}
\exp(-b(x\!-\!a)^2)dx\!=\!\sqrt{\frac{1}{\pi}}
\exp(-(x'\!-\!a')^2)dx'\!=\!d\mu_{(b',a')}(x'),\\
&&d\mu_{(b_2,a_2)}(x)=\sqrt{\frac{b_2}{\pi}}
\exp(-b_2(x\!-\!a_2)^2)dx=\sqrt{\frac{b_2}{b_1\pi}}
\exp\Big(-\frac{b_2}{b_1}(x'\!-\!a_2')^2\Big)dx'\\
&&=d\mu_{(b_2',a_2')}(x'),\quad (b',a')=(1,a\sqrt{b}),\quad (b_2',a_2')=\Big(\frac{b_2}{b_1},a_2\sqrt{b_1}\Big).
\end{eqnarray*}
\begin{rem}
\label{r.(b,a)-(b',a')} 
All the expressions, given in  the list \eqref{Sigma_1^+(t)} \eqref{Sigma_1^{-}(t)},  \eqref{Sigma_3} and \eqref{III} 
are invariant under the transformations \eqref{(b,a)-(b',a')}
$$
S^L_{kr}(\mu)=\sum_{n\in {\mathbb Z}}
\frac{b_{kn}}{2}\left(\frac{1}{2b_{rn}}+a_{rn}^2\right),\quad 
Y_r=\Big(a_{rk}\Big(\frac{1}{2b_{1k}}+\frac{1}{2b_{2k}}+\frac{1}{2b_{3k}}\Big)^{-1/2}\Big)_{k\in{\mathbb Z}},
$$
etc., and $S_r(3)$ that are defined by \eqref{S_r(3)}.
\end{rem}

\subsection{Approximation scheme}
 %
%
{\bf Case I.}
{\it  Approximation of $x_{kt}$ by $A_{kn}A_{tn}$.}
Recall that  we will write $1$ if some values  $=\infty$ and $0$ in the case $<\infty$ (see Remark~\ref{r.(0,1)}).
We use the following notation $S:=(S_1(3), S_2(3),S_3(3))$.
By Lemma~\ref{l.S1+S2+S3.m=3} we get $$\sum_{r=1}^3S_r(3)=\infty.$$ 
Therefore, without loss of generality, it suffices to consider the following
three cases:
\begin{equation}
\label{S=1,2,3} 
(1)\,\,S=(0,0,1),\quad (2)\,\,S=(0,1,1),\quad (3)\,\,S=(1,1,1).
\end{equation}
By Lemma~\ref{l.SL(3)-perp}, 
the condition of orthogonality $(\mu_{(b,a)}^3)^{L_t}\perp\mu_{(b,a)}^3$ 
for   $t\in \pm{\rm SL}(3,\mathbb R)\setminus\{e\}$, i.e.,
$
\Sigma^{\pm}(t)=\Sigma_1^{\pm}(t)+\Sigma_2(t)=\infty,
$
splits into two cases:
\begin{equation}
\label{A,B.3}
\begin{array}{cccc}
(A) & \Sigma_1^{\pm}(t)=\infty,&& \Sigma_1^{\pm}(t)
=\sum_{1\leq i< j\leq 3}\Sigma^{\pm}_{ij}(t),\\
(B) &\Sigma_1^{\pm}(t)<\infty&\text{ but}& \Sigma_2(t)=\infty,
\end{array}
\end{equation}
where $\Sigma_1^\pm(t),\,$ $\Sigma^{\pm}_{ij}(t)$ and $\Sigma_2(t)$ are defined by \eqref{Sigma_1^+(t).m=3},
\eqref{Sigma_1^{-}(t).m=3} , \eqref{Sigma(pm)(ij)} and \eqref{Sigma_3}.
\subsection{Case $S=(0,0,1)$}
\begin{lem}
\label{l.(001)} 
The case $S=(0,0,1)$ is equivalent with
\begin{equation}
 \label{S(31)}
\Sigma^{13}+\Sigma^{23}<\infty,
\quad S_3(3)\sim \sum_n\frac{b_{3n}^2}{b_{1n}b_{2n}}=\infty.
\end{equation}
\end{lem}
%
\begin{pf}
To prove the first part of \eqref{S(31)}
we set $c_n=\frac{b_{3n}}{b_{1n}+b_{2n}}$ and note that
\begin{eqnarray*}
 &&\infty>S_1(3)+S_2(3)=\sum_{n\in \mathbb Z}\frac{b^2_{1n}+b^2_{2n}}{b_{1n}b_{2n}+b_{1n}b_{3n}+b_{2n}b_{3n}}\stackrel{\eqref{two-sim-ser}}{\sim}\\
&&
\sum_{n\in \mathbb Z}\frac{b_{1n}^2+b_{2n}^2}{(b_{1n}+b_{2n}+b_{3n})^2-b_{3n}^2}\sim 
\sum_{n\in \mathbb Z}\frac{(b_{1n}+b_{2n})^2}{(b_{1n}+b_{2n}+b_{3n})^2-b_{3n}^2}=\\
&&
\sum_{n\in \mathbb Z}\frac{1}{(1+c_n)^2-c_n^2}
\!=\!\sum_{n\in \mathbb Z}\frac{1}{1+2c_n}\!\stackrel{\eqref{two-sim-ser}}{\sim}\!  \sum_{n\in \mathbb Z}\frac{1}{c_n}=\sum_{n\in \mathbb Z}\frac{b_{1n}+b_{2n}}{b_{3n}}=\Sigma^{13}+\Sigma^{23}.
\end{eqnarray*}
To prove the second part of \eqref{S(31)} we have by the first part of \eqref{S(31)}
\begin{equation*}
 S_3(3)\!=\!\sum_{n\in \mathbb Z}\frac{b^2_{3n}}{b_{1n}b_{2n}\!+\!b_{1n}b_{3n}\!+\!b_{2n}b_{3n}}
\!=\!\sum_{n\in \mathbb Z}\frac{1}{\frac{b_{1n}b_{2n}}{b^2_{3n}}+\frac{b_{1n}}{b_{3n}}+\frac{b_{2n}}{b_{3n}}}\!
\sim\!
\sum_{n\in \mathbb Z}\frac{b^2_{3n}}{b_{1n}b_{2n}}.\,\,\,\Box
\end{equation*}
\end{pf}
%
In the case $S=(0,0,1)$ we have
$$
\Delta(Y_3^{(3)}, Y_1^{(3)},Y_2^{(3)})\sim \Delta(Y_3^{(3)})\sim \Vert Y_3^{(3)}\Vert^2=\infty,
$$
so we can approximate $x_{3n}x_{3t}$ using Lemma~\ref{x3x3.3} and after that we can approximate $x_{3n}$ by an analogue of Lemma~\ref{5.3}. 
{\it From now on we will say that we can approximate $x_{3n}$ using Lemma~\ref{x3x3.3}, without mentioning   Lemma~\ref{5.3}}. 

We can not approximate  $x_{1n}$ and $x_{2n}$ using  Lemma~\ref{x1x1.3}-\ref{x2x2.3}, since we have
$$
\Delta(Y_1^{(1)}, Y_2^{(1)},Y_3^{(1)})+\Delta(Y_2^{(2)}, Y_3^{(2)},Y_1^{(2)})<\infty.
$$
We can try to approximate some of $D_{rn}$ for $1\leq r\leq 3$ using Lemmas~\ref{d1.3}--\ref{d3.3}, see Section~\ref{s.444} for details. 
We have for $1\leq k\leq 3$ (see \eqref{Delta^k(krs)}):
\begin{equation*}
D_{kn}\,\,\eta\,\,{\mathfrak A}^3\Leftrightarrow \Delta_k=\infty,\quad\text{where}\quad 
\Delta_k:=\Delta(Y_k,Y_r,Y_s),
\end{equation*}
and $\{k,r,s\}$ is a cyclic permutation of $\{1,2,3\}$.
Recall that by $\Sigma^{12}+\Sigma^{13}<\infty$ we get  (see \eqref{II} for the expressions of $\Vert Y_r \Vert^2,\,\,1\leq r\leq 3$)
\begin{equation}
\label{Y(r)-(001)}
\Vert Y_1 \Vert^2\sim
\sum_{n\in \mathbb Z}b_{1n}a_{1n}^2,\quad
\Vert Y_2 \Vert^2\sim
\sum_{n\in \mathbb Z}b_{1n}a_{2n}^2,\quad
\Vert Y_3\Vert^2\sim
\sum_{n\in \mathbb Z}b_{1n}a_{3n}^2. 
\end{equation}
By \eqref{S(31)}
we have $
\Sigma^{13}+\Sigma^{23}<\infty$. We distinguish two cases:\\
(1) $\Sigma^{12}<\infty$,\\
(2) $\Sigma^{12}=\infty$.

In the case (1), since $
\Sigma^{12}+\Sigma^{13}<\infty$ we have
\begin{eqnarray*}
&&
S^{L}_{1,23}(\mu,t,s)\stackrel{\eqref{S^{L}_{1,23}} }{=}\sum_{n\in{\mathbb
Z}}\left[\frac{t^2}{4}\frac{b_{1n}}{b_{2n}}+\frac{s^2}{4}\frac{b_{1n}}{b_{3n}}+
\frac{b_{1n}}{2}\big(\!-\!2a_{1n}\!+\!ta_{2n}\!+\!sa_{3n}\big)^2\right]\sim\\
&&
\sum_{n\in{\mathbb
Z}}
\frac{b_{1n}}{2}\big(\!-\!2a_{1n}\!+\!ta_{2n}\!+\!sa_{3n}\big)^2\stackrel{\eqref{Y(r)-(001)}}{\sim} 
\Vert C_1Y_1+\!C_2Y_2\!+\!C_3Y_3\Vert^2.
\end{eqnarray*}
Finally, in the case (1) we can approximate all $D_{rn},\,1\leq r\leq 3$  using Lemmas~\ref{d1.3}--\ref{d3.3} and 
Lemma~\ref{l.min=proj.3}
and the proof is finished.

The case (2) 
can be divided into three
cases (if necessary, we can chose an appropriate subsequence of $\Big(\frac{b_{1n}}{b_{2n}}\Big)_n$\big):
\begin{equation}
\label{2.(001)}
\lim_n\frac{b_{1n}}{b_{2n}}=\!\left\{\begin{array}{cl}
 (a)&0\\
 (b)&b>0\\
 (c)&\infty
\end{array}\right.. 
\end{equation}
The case (c) is reduced to the case (a) by exchanging $(b_{2n},a_{2n})$ with $(b_{1n},a_{1n})$.

This transfomation does not change the first condition in \eqref{S(31)}.
In the case  (2.a-b), by \eqref{II}  we obtain the following  expressions for  $\Vert Y_r \Vert^2,\,\,1\leq r\leq 3$: 
\begin{eqnarray*}
 &&\Vert Y_1\Vert^2=\sum_{n\in \mathbb Z}\frac{a_{1n}^2}{\frac{1}{2b_{1n}}+\frac{1}{2b_{2n}}+\frac{1}{2b_{3n}}}= 
\sum_{k\in \mathbb Z}\frac{2b_{1n}a_{1n}^2}{1+\frac{b_{1n}}{b_{2n}}+\frac{b_{1n}}{b_{3n}}
}\stackrel{\Sigma^{13}<\infty}{\sim}
\sum_{n\in \mathbb Z}2b_{1n}a_{1n}^2,\\
&&\Vert Y_2 \Vert^2\sim
\sum_{n\in \mathbb Z}2b_{1n}a_{2n}^2,\quad
\Vert Y_3\Vert^2=
\sum_{n\in \mathbb Z}b_{1n}a_{3n}^2.
\end{eqnarray*}
Since
$$
\Vert Y_1 \Vert^2\sim
\sum_{n\in \mathbb Z}b_{1n}a_{1n}^2\sim S^L_{11}(\mu)=\infty,
$$ 
we have four possibilities 
for $y_{23}:=(y_2,y_3)\in \{0,1\}^2$ as in \eqref{(y_1,y_2,y_3)(011)}, see Section~\ref{s.444}:
\begin{equation*}
\begin{array}{ccccc}
   &(1.0)&(1.1)&(1.2)&(1.3)\\
y_1&1&1&1&1\\
y_2&0&1&0&1\\
y_3&0&0&1&1\\
\end{array}
\end{equation*}
We just follow the instructions given in Remark~\ref{r.(011)-(1-3)(0-1)}.
We note that the cases (1.0) and (1.1) can not occur since the following conditions are contradictory:
$$
S^L_{13}(\mu)\stackrel{\eqref{S^L(kn)} }{=}\sum_{n\in{\mathbb
Z}}\frac{b_{1n}}{2}\Big(\frac{1}{2b_{3n}}+a_{3n}^2\Big)\!=\!\infty, \,\,\,\Vert Y_3\Vert^2\sim \sum_{n\in \mathbb Z}b_{1n}a_{3n}^2<\infty,\,\,\, \Sigma^{13}\stackrel{\eqref{S(31)}}{<}\infty.
$$
We have two cases  (1.2.1) and (1.3.1) according to whether respectively the expressions in
\eqref{nu(12)} 
 or \eqref{nu(123)} are 
divergent. We can approximate in these cases respectively  
$D_{1n}$ and $D_{3n}$ in \eqref{(1.2)}  and all $D_{1n},\,D_{2n},\,D_{3n}$ in \eqref{(1.3)}. The proof of irreducibilty is finished in both cases
because we have  
$x_{3n},\,D_{3n}\,\,\eta\,\, {\mathfrak A}^3$ 
%
and the problem is reduced to the case $m=2$ \cite{KosJFA17}, since
$
A_{kn}\!=\!\sum_{r=1}^{3}x_{rk}D_{rn}-x_{3k}D_{3n}=
\sum_{r=1}^{2}x_{rk}D_{rn}.
$

If the opposite holds, 
we have two different cases (1.2.0) and (1.3.0).
We try
to approximate $D_{3n}$  using Lemma~\ref{l.Re-Im-exp.3}. If one of the expressions  $\Sigma_3(D,s)$ or $\Sigma^{\vee}_3(D,s)$ is divergent
for some sequence $s=(s_k)_{k\in \mathbb Z}$,
we can approximate $D_{3k}$ and the proof is finished,  since we have  $x_{3n},\,D_{3n}\,\,\eta\,\, {\mathfrak A}^3$ 
and the problem is reduced to the case $m=2$. 
Let us suppose,  as in Remark~\ref{r.Re-Im-(011)}, that for every sequence $s=(s_k)_{k\in \mathbb Z}$ we have
$$
\Sigma_3(D,s)+\Sigma^{\vee}_3(D,s)<\infty.
$$
Then, in particular, we have for $s^{(3)}=(s_{k})_{k\in \mathbb Z}$
with $\frac{s^2_{k}}{b_{3k}}\equiv 1$
\begin{eqnarray*}
&& 
 \infty> \Sigma_3(D,s^{(3)})+\Sigma^{\vee}_3(D,s^{(3)})\sim \Sigma_3(D)+\Sigma^{\vee}_3(D)=\\
&&
\sum_k\frac{\frac{1}{2b_{3k}}+a^2_{3k}}{C_k+a^2_{1k}+a^2_{2k}+a^2_{3k}}\stackrel{\eqref{two-sim-ser}}{\sim}
\sum_k\frac{
\frac{1}{2b_{3k}}+a^2_{3k}}{\frac{1}{2b_{1k}}+a^2_{1k}+
\frac{1}{2b_{2k}}+a^2_{2k}}=\\
&&
 \sum_k\frac{\frac{b_{1k}}{b_{3k}}+2b_{1k}a^2_{3k}}{1+2b_{1k}a^2_{1k}+\frac{b_{1k}}{b_{2k}}+2b_{1k}a^2_{2k}}\stackrel{\eqref{S(31)}}{\sim}
\sum_k\frac{2b_{1k}a^2_{3k}
}
{1+2b_{1k}a^2_{1k}+2b_{1k}a^2_{2k}}=:\Sigma^+_3(D).
\end{eqnarray*}
\begin{rem}
\label{r.(001)} 
Finally, we have $\Sigma^+_3(D)\!\sim\! \sum_k\frac{2a^2_{3k}}
{1+2a^2_{1k}+2a^2_{2k}}$, we take $b_{1n}\!\equiv\! 1$ by  \eqref{(b,a)-(b',a')}. 
In the case (1.2.0) we have
$\Vert Y_2 \Vert^2\sim
\sum_{n\in \mathbb Z}b_{1n}a_{2n}^2<\infty$, and therefore
$\Sigma^+_3(D)\sim \sum_k\frac{2a^2_{3k}}
{1+2a^2_{1k}}$,  and hence 
$\Sigma^+_3(D)=\infty$ by Lemma~\ref{l.(011)-fin}.
In the case (1.3.0) we have $a_3=\pm a_1\pm a_2+h$ or $a_3-h=\pm a_1\pm a_2$,
see the proof of  Lemma~\ref{l.(011)-fin.3}.
Therefore,
\begin{eqnarray}
\nonumber
&&
\infty>\Sigma^+_3(D)\sim 
\sum_k\frac{a^2_{3k}}
{1+a^2_{1k}+a^2_{2k}}\geq
\sum_k\frac{a^2_{3k}}
{1+a^2_{1k}+2\vert a_{1k}\vert \vert a_{2k}\vert +a^2_{2k}}=
\\
\label{(001).1}
&&
\sum_k\frac{a^2_{3k}}
{1+\big(\vert a_{1k} \vert +\vert a_{2k}\vert \big)^2},\quad \infty>\Sigma^+_3(D)\sim 
\sum_k\frac{a^2_{3k}}
{1+a^2_{1k}+a^2_{2k}}\geq\\
\nonumber
&&
\sum_k\frac{a^2_{3k}}
{1+a^2_{1k}+a^2_{2k} +(\vert a_{1k}\vert-\vert a_{2k}\vert)^2
}\sim\sum_k\frac{a^2_{3k}}
{1+2a^2_{1k}-2\vert a_{1k}\vert \vert a_{2k}\vert +2a^2_{2k}}\\
\label{(001).2}
&&
\sim
\sum_k\frac{a^2_{3k}}
{1+a^2_{1k}-2\vert a_{1k}\vert \vert a_{2k}\vert +a^2_{2k}}\sim
\sum_k\frac{a^2_{3k}}
{1+\big(\vert a_{1k} \vert -\vert a_{2k}\vert \big)^2}.
\end{eqnarray}
Hence, we have by \eqref{(001).1} and \eqref{(001).2}  
\begin{equation}
\infty>\Sigma^+_3(D)\geq
\sum_k\frac{a^2_{3k}}
{1+\big(\pm a_{1k}\pm a_{2k}\big)^2}
=\sum_k\frac{a^2_{3k}}
{1+\big( a_{3k}-h_{k}\big)^2}=\infty 
\end{equation}
by Lemma~\ref{l.(011)-fin}, contradiction.
Therefore, in both cases we can approximate $D_{3n}$ and the proof is finished.
\end{rem}

\subsection{Case $S=(0,1,1)$}
\begin{lem}
\label{l.(011)} 
In the case $S=(0,1,1)$
we have 
\begin{equation}
\label{(011)} 
\lim_nd_{2n}=\lim_nd_{3n}=\infty.
\end{equation}
\end{lem}
\begin{pf} Setting  as before $d_{rn}\!=\!\frac{b_{rn}}{b_{1n}}$,  we obtain
by  \eqref{S_r(3)}  and \eqref{two-sim-ser}
\begin{eqnarray}
\label{S_1(3)<}
 &&S_1(3)=
 \sum_{n\in \mathbb Z}\frac{1}{d_{2n}\!+\!d_{3n}\!+\!d_{2n}d_{3n}}\!\sim\! 
  \sum_{n\in \mathbb Z}\frac{1}{(1\!+\!d_{2n})(1\!+\!d_{3n})}
 \!<\!\infty,\\
\label{S_2(3)=inf}
&&
S_2(3)\!=\!\!
 \sum_{n\in \mathbb Z}\frac{d_{2n}^2}{d_{2n}\!+\!d_{3n}\!+\!d_{2n}d_{3n}}\!\!
 \stackrel{\eqref{two-sim-ser}}{\sim} \!
  \sum_{n\in \mathbb Z}\frac{d_{2n}^2}{(1\!+\!d_{2n})(d_{2n}\!+\!d_{3n})}\!=\!\infty,\,\,\\
  \label{S_(3)=inf}
  &&
  S_3(3)\!=\!\!\
 \sum_{n\in \mathbb Z}\frac{d_{3n}^2}{d_{2n}\!+\!d_{3n}\!+\!d_{2n}d_{3n}}\!\!\stackrel{\eqref{two-sim-ser}}{\sim} \! \!
  \sum_{n\in \mathbb Z}\frac{d_{3n}^2}{(1\!+\!d_{3n})(d_{2n}\!+\!d_{3n})}\!=\!\infty.\quad
\end{eqnarray}
Suppose that $d_{2n}\leq C$ for all $n\in \mathbb Z$. Then by \eqref{S_1(3)<} and \eqref{S_2(3)=inf} we conclude
\begin{eqnarray*}
 &&S_1(3) \sim\sum_{n\in \mathbb Z}\frac{1}{(1\!+\!d_{2n})(1\!+\!d_{3n})}\!\sim\! 
 \sum_{n\in \mathbb Z}\frac{1}{1\!+\!d_{3n}}\sim
  \sum_{n\in \mathbb Z}\frac{1}{d_{3n}}
 \!<\!\infty,\quad\infty=S_2(3)\\
&&
\sim\!
  \sum_{n\in \mathbb Z}\frac{d_{2n}^2}{(1\!+\!d_{2n})(d_{2n}\!+\!d_{3n})}\sim
  \sum_{n\in \mathbb Z}\frac{d_{2n}^2}{d_{2n}\!+\!d_{3n}}\!\leq\! 
\sum_{n\in \mathbb Z}\frac{C^2}{C\!+\!d_{3n}}\!\stackrel{\eqref{two-sim-ser}}{\sim}\! \sum_{n\in \mathbb Z}\frac{1}{d_{3n}}
\!<\!\infty,
\end{eqnarray*}
a contradiction. 
We use the fact that for any fixed $D>0$  the function
$$
f_D(x)=\frac{x^2}{x+D}
$$  is  strictly increasing when $x>0$.
Similarly, if we suppose that $d_{3n}\leq C$ for all $n\in \mathbb Z$ we will obtain a contradiction too.
\qed\end{pf}
\begin{lem}
\label{l.S_(123)} 
The case  $S=(0,1,1)$ is equivalent with
\begin{equation}
\label{S_(123)}
S_1(3)\!\sim\!\sum_n\frac{b_{1n}^2}{b_{2n}b_{3n}}\!<\!\infty,\,\,
S_2(3)\sim\sum_n\frac{1}{d_n}=\infty,
\,\,
S_3(3)\sim\sum_nd_n=\infty.
\end{equation}
\end{lem}
\begin{pf} 
Recall that $d_{n}\!=\!\frac{d_{3n}}{d_{2n}}$.
Denote $D_n:=1+d^{-1}_{2n}+d^{-1}_{3n}$.
By Lemma~\ref{l.(011)} we have
\begin{equation}
\label{alpha}
\quad 1\leq D_n=1+d^{-1}_{2n}+d^{-1}_{3n}\leq C,\quad \text{for all}\quad n\in\mathbb Z.
\end{equation}
Therefore, we get
\begin{eqnarray*}
 && S_1(3)\!=\!\sum_{n\in \mathbb Z}\frac{1}{d_{2n}+d_{3n}\!+\!d_{2n}d_{3n}}=
\sum_{n\in \mathbb Z}\frac{1}{D_nd_{2n}d_{3n}}\sim 
\sum_{n\in \mathbb Z}\frac{1}{d_{2n}d_{3n}}=\sum_n\frac{b_{1n}^2}{b_{2n}b_{3n}},
\\ 
 &&
S_2(3)\!=\!\sum_{n\in \mathbb Z}\frac{d_{2n}^2}{d_{2n}\!+\!d_{3n}\!+\!d_{2n}d_{3n}}=
\sum_{n\in \mathbb Z}\frac{d_{2n}^2}{D_nd_{2n}d_{3n}}\sim\sum_{n\in \mathbb Z}\frac{1}{d_n}, \\
 && S_3(3)=\sum_{n\in \mathbb Z}\frac{d_{3n}^2}{d_{2n}\!+\!d_{3n}\!+\!d_{2n}d_{3n}}
=\sum_{n\in \mathbb Z}\frac{d_{3n}^2}{D_nd_{2n}d_{3n}}\sim \sum_{n\in \mathbb Z}d_n.
\hskip 2.9cm\Box
 \end{eqnarray*}
 \end{pf}

By Lemma~\ref{l.S1+S2+S3.m=3}, \eqref{Yr^s-<-Sr} we get $\Vert Y_1^{(r)}\Vert^2<\infty,\,\,1\leq r\leq 3$ therefore,
we get
\begin{lem}
\label{l.S=(0,1,1)} 
In the case  $S=(0,1,1)$ we have 
\begin{eqnarray}
\nonumber
&&\Delta(Y_1^{(1)},Y_2^{(1)},Y_3^{(1)})<\infty,
\quad
\Delta(Y_2^{(2)},Y_3^{(2)},Y_1^{(2)})\sim  \Delta(Y_2^{(2)},Y_3^{(2)}),\\
\label{(011)-(11)}
&&\Delta(Y_3^{(3)}, Y_1^{(3)},Y_2^{(3)})\sim \Delta(Y_3^{(3)},Y_2^{(3)}).
\end{eqnarray}
\end{lem}
\begin{pf} 
Set $(f_1,f_2,f_3)=(Y_3^{(3)}, Y_1^{(3)},Y_2^{(3)})$. Then
\begin{eqnarray*}
&&\Delta(f_1,f_2,f_3)\overset{\eqref{Delta(f,g,h)}}{=}
\frac{\Gamma(f_1)+\Gamma(f_1,f_2)+\Gamma(f_1,f_3)+\Gamma(f_1,f_2,f_3)}{1+\Gamma(f_2)+\Gamma(f_3)+\Gamma(f_2,f_3)}>\\
&&
\frac{\Gamma(f_1)+\Gamma(f_1,f_3)}{1+\Gamma(f_2)+\Gamma(f_3)+\Gamma(f_2,f_3)}\stackrel{\eqref{G(f,g)-leq()()}}{\geq}
\frac{\Gamma(f_1)+\Gamma(f_1,f_3)}{(1+\Gamma(f_2))(1+\Gamma(f_3))}\sim  \Delta(f_1,f_3),
\end{eqnarray*} 
since $f_2\!\in\! l_2(\mathbb Z)$.  Indeed, for $f,g\!\in\! l_2(\mathbb Z)$ and  $f\!\in\! l_2(\mathbb Z), g\!\not\in\! l_2(\mathbb Z)$ we have respectively 
\begin{eqnarray}
\label{G(f,g)-leq()()}
&&
\Gamma(f,g)\leq \Gamma(f)\Gamma(g)<\infty,\,\,\,
\Gamma(f,g)\leq \Gamma(f)\Gamma(g),\,\,\,\text{where}\quad\Gamma(f,g),
\\
\nonumber
&&
\Gamma(g)\quad \text{are defined by}\quad \Gamma(f,g):=\lim_n \Gamma(f_{(n)},g_{(n)})\quad\Gamma(g):=\lim_n\Gamma(g_{(n)}),
\end{eqnarray}
and $g_{(n)}\!:=\!(g_k)_{k=-n}^n\!\in\!\mathbb R^{2n+1}$. 
Similarly, set $(f_1,f_2,f_3)\!=\!(Y_2^{(2)},Y_3^{(2)},Y_1^{(2)})$,
then
\begin{eqnarray*}
&&\Delta(f_1,f_2,f_3)\overset{\eqref{Delta(f,g,h)}}{=}
\frac{\Gamma(f_1)+\Gamma(f_1,f_2)+\Gamma(f_1,f_3)+\Gamma(f_1,f_2,f_3)}{1+\Gamma(f_2)+\Gamma(f_3)+\Gamma(f_2,f_3)}>\\
&&
\frac{\Gamma(f_1)+\Gamma(f_1,f_2)}{1+\Gamma(f_2)+\Gamma(f_3)+\Gamma(f_2,f_3)}\stackrel{\eqref{G(f,g)-leq()()}}{\geq}
\frac{\Gamma(f_1)+\Gamma(f_1,f_2)}{(1+\Gamma(f_2))(1+\Gamma(f_3))}
\sim  \Delta(f_1,f_2),
\end{eqnarray*}
since $f_3\in  l_2(\mathbb Z)$.
Finally, we  
derive both equivalences
in \eqref{(011)-(11)}.
To prove that \\$\Delta(Y_1^{(1)},Y_2^{(1)},Y_3^{(1)})<\infty$ we set $(f_1,f_2,f_3)=( Y_1^{(1)},Y_2^{(1)},Y_3^{(1)})$, and note that
\begin{eqnarray*}
&&\Delta(f_1,f_2,f_3)\overset{\eqref{Delta(f,g,h)}}{=}
\frac{\Gamma(f_1)+\Gamma(f_1,f_2)+\Gamma(f_1,f_3)+\Gamma(f_1,f_2,f_3)}{1+\Gamma(f_2)+\Gamma(f_3)+\Gamma(f_2,f_3)}\leq \\
&&\frac{\Gamma(f_1)\Big(1+\Gamma(f_2)+\Gamma(f_3)+\Gamma(f_2,f_3)\Big)}
{1+\Gamma(f_2)+\Gamma(f_3)+\Gamma(f_2,f_3)}=\Gamma(f_1)<\infty.
\hskip 3cm \Box
\end{eqnarray*}
\end{pf} 
In order
to approximate $x_{2n}$ or $x_{3n}$, 
it remains to study when  
\begin{equation}
\label{(011).m=3}
\Delta(Y_2^{(2)},Y_3^{(2)})=\infty,\quad \Delta(Y_3^{(3)},Y_2^{(3)})=\infty,
\end{equation}
 where $
\Delta(f_1,f_2)=\frac{\Gamma(f_1)+\Gamma(f_1,f_2)}{1+\Gamma(f_2)}$. For $2\leq r\leq 3$, 
denote 
\begin{eqnarray}
 \label{rho_r(C,C)}
 &&\rho_r(C_2,C_3):=\Vert C_2Y_2^{(r)}+C_3Y_3^{(r)}\Vert^2,\quad (C_2,C_3)\in \mathbb R^2,\\
 \label{nu(C,C,C)}
 &&\nu( C_1, C_2, C_3):=\Vert C_1Y_1+ C_2Y_2+C_3Y_3\Vert^2,
 \,\, (C_1,C_2,C_3)\in \mathbb R^3.\,\,\,
\end{eqnarray}
\begin{lem}
\label{rho,nu,(011)} 
In the case $S=(0,1,1)$   
we have 
\begin{eqnarray}
  \label{rho_r-al-(011)}
&&\rho_2(C_2,C_3)\!\sim\!\sum_{n\in \mathbb Z}\frac{\big(C_2\!+\!C_3d_{n}\big)^2}{1+2d_n},\,\,
\rho_3(C_2,C_3)\!\sim\!\sum_{n\in \mathbb Z}\frac{\big(C_2\!+\!C_3d_{n}\big)^2}{d^2_{n}+2d_n}\,\,\,\,\,\\
\label{nu(CCC)} 
&&
=\sum_{n\in \mathbb Z}\frac{\big(C_2l_n\!+\!C_3\big)^2}{1+2l_n},\,\,\,
\,\,\,
\nu( C_1, C_2, C_3)\sim\sum_{n\in \mathbb Z}b_{1n}\Big(\sum_{r=1}^3C_ra_{rn}\Big)^2.
\end{eqnarray}
\end{lem}
\begin{pf}
Set as before $d_n=\frac{d_{3n}}{d_{2n}}$.
By \eqref{Y_r^{(s)}} and \eqref{II} we get 
{\small
\begin{eqnarray*}
\nonumber
 && \Vert Y_2^{(2)}\Vert^2\!=\!\sum_{n\in \mathbb Z}\frac{d_{2n}^2}{d_{2n}^2+\!2(d_{2n}+d_{3n}\!+\!d_{2n}d_{3n})}\!=\!
\sum_{n\in \mathbb Z}\frac{d_{2n}^2}{d_{2n}^2\!+\!2D_nd_{2n}d_{3n}}\!\sim\! \sum_{n\in \mathbb Z}\frac{1}{D_nd_n},\\ 
\nonumber
 &&
 \Vert Y_3^{(2)} \Vert^2\!=\!\sum_{n\in \mathbb Z}\frac{d_{3n}^2}{d_{2n}^2\!+\!2(d_{2n}\!+\!d_{3n}\!+\!d_{2n}d_{3n})}=
\sum_{n\in \mathbb Z}\frac{d_{3n}^2}{d_{2n}^2\!+\!2D_nd_{2n}d_{3n}}=  \\
 \nonumber
 &&
\sum_{n\in \mathbb Z}\frac{d_n^2}{1+2D_nd_n},\quad  
 \Vert Y_2^{(3)} \Vert^2=\sum_{n\in \mathbb Z}\frac{d_{2n}^2}{d_{3n}^2+2(d_{2n}\!+\!d_{3n}\!+\!d_{2n}d_{3n})}
=\\
&&
\sum_{n\in \mathbb Z}\frac{d_{2n}^2}{d_{3n}^2\!+\!2D_nd_{2n}d_{3n}}= \sum_{n\in \mathbb Z}\frac{1}{d_n^2+2D_nd_n},
 \end{eqnarray*}
 \begin{eqnarray}
\nonumber
 &&\Vert Y_3^{(3)}\Vert^2=\sum_{n\in \mathbb Z}\frac{d_{3n}^2}{d_{3n}^2+2(d_{2n}\!+\!d_{3n}\!+\!d_{2n}d_{3n})}
  =\sum_{n\in \mathbb Z}\frac{d_{3n}^2}{d_{3n}^2\!+\!2D_nd_{2n}d_{3n}}\sim \sum_{n\in \mathbb Z}\frac{d_n}{D_n},\\
 \nonumber 
 &&\Vert Y_1\Vert^2=\sum_{n\in \mathbb Z}\frac{a_{1n}^2}{\frac{1}{2b_{1n}}+\frac{1}{2b_{2n}}+\frac{1}{2b_{3n}}}= 
\sum_{k\in \mathbb Z}\frac{2b_{1n}a_{1n}^2}{1+d_{2n}^{-1}+d_{3n}^{-1}}=
\sum_{n\in \mathbb Z}\frac{2b_{1n}a_{1n}^2}{D_n},\\
\label{Y_r(011)}
&&\Vert Y_2 \Vert^2=
\sum_{n\in \mathbb Z}\frac{2b_{1n}a_{2n}^2}{D_n},\quad
\Vert Y_3\Vert^2=
\sum_{n\in \mathbb Z}\frac{b_{1n}a_{3n}^2}{D_n}.
\end{eqnarray}
}
Recall that
$d_{rn}\!=\!\frac{b_{rn}}{b_{1n}}$. By \eqref{alpha}, we obtain
\begin{eqnarray}
\label{S=(011)Y^r_s(alpha)}
&& \Vert Y_2^{(2)}\Vert^2\sim \sum_{n\in \mathbb Z}\frac{1}{1+2d_n}
\sim \sum_{n\in \mathbb Z}\frac{1}{d_n},\quad
 \Vert Y_3^{(2)} \Vert^2\sim \sum_{n\in \mathbb Z}\frac{d_n^2}{1+2d_n},\\
 \nonumber
&&\Vert Y_2^{(3)} \Vert^2\sim \sum_{n\in \mathbb Z}\frac{1}{d_n^2+2d_n},\quad
\Vert Y_3^{(3)}\Vert^2 \sim \sum_{n\in \mathbb Z}\frac{d_n^2}{d_n^2+2d_n}
\sim \sum_{n\in \mathbb Z}d_n,
\end{eqnarray}
\begin{eqnarray}
 \nonumber
&&\Vert Y_1\Vert^2\sim \sum_{n\in \mathbb Z}b_{1n}a_{1n}^2,\,\,
\Vert Y_2\Vert^2\sim \sum_{n\in \mathbb Z}b_{1n}a_{2n}^2,\,\,
\Vert Y_3\Vert^2\sim \sum_{n\in \mathbb Z}b_{1n}a_{3n}^2,\\
\label{S=(011)Y_r(alpha)}
&&
\Vert C_1Y_1+C_2Y_2+C_3 Y_3\Vert^2\stackrel{\eqref{alpha}}{\sim}\sum_{n\in \mathbb Z}b_{1n}\big(C_1a_{1n}\!+\!C_2a_{2n}\!+\!C_3a_{3n}\big)^2.\quad
\end{eqnarray}
By \eqref{S=(011)Y^r_s(alpha)} and \eqref{S=(011)Y_r(alpha)}
the proof is finished.
\qed\end{pf}

\subsubsection{Approximation of $x_{2n},\,x_{3n}$}
To approximate $x_{2n},\,x_{3n}$, we need several lemmas.
Denote $l_n=d_n^{-1}$.
\begin{lem}
\label{l.light.123} 
The following five series are equivalent:
\begin{eqnarray}
\label{rho(2)}
&&(i-ii)\quad \sum_{n\in \mathbb Z}\frac{\big(C_2\!-\!C_3d_{n}\big)^2}{1+2d_n}\sim \sum_{n\in \mathbb Z}c_n^2,\\
\label{rho(3)}
&&(iii-iv)\quad\sum_{n\in \mathbb Z}\frac{\big(C_2l_n\!-\!C_3\big)^2}{1+2l_n}\sim \sum_{n\in \mathbb Z}e_n^2,\\
\label{sigma(23)(s)}
&&(v)\,\,\Sigma_{23}(s)\!=
\!
\sum_{n\in \mathbb Z}\Big(
s^2\sqrt{\frac{b_{2n}}{b_{3n}}}\!-\!
s^{-2}\sqrt{\frac{b_{3n}}{b_{2n}}}\Big)^2\!=
\!\sum_{n\in \mathbb Z}\Big(
\frac{s^2}{\sqrt{d_n}}\!-\!
\frac{\sqrt{d_n}}{s^2}\Big)^2
,\quad
\end{eqnarray}
where 
\begin{equation}
\label{c_n,e_n}
d_n\!=\!
C_2C_3^{-1}(1+c_n),\,\,\,l_n\!=\!C_3C_2^{-1}(1+e_n),\,\,s^4\!=\!C_2C_3^{-1}>0,\quad l_n\!=\!d_n^{-1}.
\end{equation}
\end{lem}
\begin{pf}
To prove  \eqref{rho(2)}  and \eqref{rho(3)} 
we get by Lemma~\ref{l.two-ser-alpha} using \eqref{c_n,e_n}
\begin{eqnarray*}
&&
\sum_{n\in \mathbb Z}\frac{\big(C_2\!-\!C_3d_{n}\big)^2}{1+2d_n}=
\sum_{n\in \mathbb Z}\frac{C_2^2c_n^2}{1+2C_2C_3^{-1}(1+c_n)}
\sim \sum_{n\in \mathbb Z}c_n^2,\\
&&
\sum_{n\in \mathbb Z}\frac{\big(C_2l_n\!-\!C_3\big)^2}{1+2l_n}=
\sum_{n\in \mathbb Z}\frac{C_3^2e_n^2}{1+2C_3C_2^{-1}(1+e_n)}\sim \sum_{n\in \mathbb Z}e_n^2.
\end{eqnarray*}
To finish the proof we make use of the following lemma
\qed\end{pf}
\begin{lem}
\label{l.2-eq-ser}
Let $(c_n)_{n\in \mathbb Z}$ be a sequence of real numbers 
with $1\!+c_n\!>\!0$ and $(1+c_n)(1+e_n)=1$.
Then the following three  series are equivalent:
$$\sum_{n\in \mathbb Z}
\left((1+e_n)^{1/2}-(1+e_n)^{-1/2} \right)^2,\quad\sum_{n\in \mathbb Z}c_n^2\quad\text{and}\quad \sum_{n\in \mathbb Z}e_n^2.
$$  
\end{lem}
\begin{pf}
 Set $s^4=C_2C_3^{-1}$, replacing $1+c_n$ by $(1+e_n)^{-1}$
in Lemma~\ref{l.eq-prod}  gives
\begin{equation*}
 \Sigma_{23}(s)=
\sum_{n\in \mathbb Z}
\left((1+c_n)^{-1/2}-(1+c_n)^{1/2} \right)^2=
\sum_{n\in \mathbb Z}
\left((1+e_n)^{1/2}-(1+e_n)^{-1/2} \right)^2.
\end{equation*}
Therefore, $\sum_{n\in \mathbb Z}\frac{c_n^2}{1\!+\!c_n}\!=\!\sum_{n\in \mathbb Z}\frac{e_n^2}{1\!+\!e_n}$ and hence,  by Lemma~\ref{l.two-ser-alpha}, the two series are equivalent: $ \sum_{n\in \mathbb Z}c_n^2\!\sim\! \sum_{n\in \mathbb Z}e_n^2$.
\qed\end{pf}

\subsubsection{Two remaining possibilities}
\label{ss.2-cases}
By Lemma~\ref{l.light.123}  there are only two cases:\\
(1) when $\rho_2(C_2,C_3)=\rho_3(C_2,C_3)=\infty$ for all $(C_2,C_3)\in \mathbb R^2\setminus\{0\}$,\\
(2) when both $\rho_2(C_2,C_3)$ and $\rho_3(C_2,C_3)$ are finite and hence, $\Sigma_{23}(s)<\infty$.
To illustrate this we start with the following example
\begin{ex}
 \label{ex.(011)(a)}
Set $d_n=n^\alpha$ for $n\in \mathbb N$ with $\alpha\in \mathbb R$. We have 
\end{ex}
\begin{equation}
\label{d_n-ex}
\lim_nd_n=\!\left\{\begin{array}{rcc}
 \infty&\text{if}&\alpha>0\\
 1&\text{if}&\alpha=0\\
 0&\text{if}&\alpha<0
\end{array}\right.. 
\end{equation}
For the general sequence $(d_n)_{n\in \mathbb Z}$ 
we have four cases (if necessary, we can chose an appropriate subsequence):
\begin{equation}
\label{d_n-gen}
\lim_nd_n=\!\left\{\begin{array}{clcc}
(a)& \infty&&\\
(b)&d>0&\text{with}&\sum_nc_n^2=\infty\\
(c)&d>0&\text{with}&\sum_nc_n^2<\infty\\
(d)&d=0&
\end{array}\right., 
\end{equation}
where $d_n=d(1+c_n)$ and $\lim_nc_n=0$.

\subsubsection{Cases (a), (b), (d)}
\begin{rem}
\label{r.(011)} 
In the cases (a) 
we see by \eqref{rho_r-al-(011)} that $$\qquad \rho_2(C_2,C_3)=\rho_3(C_2,C_3)=\infty\quad\text{
for all}\quad (C_2,C_3)\in \mathbb R^2\setminus\{0\}.$$
The case (d) is reduced to the case (a) by exchanging $(b_{2n},a_{2n})$ with $(b_{3n},a_{3n})$.
In the cases (b) by  
Lemma~\ref{l.light.123}  and \eqref{d_n-gen}  we  conclude that 
$$
\rho_2(C_2,-C_3)=\rho_3(C_2,-C_3)=\infty\quad \text{ for}\quad C_2C_3^{-1}>0
$$
hence, $\rho_2(C_2,C_3)=\rho_3(C_2,C_3)=\infty$
for all $(C_2,C_3)\in \mathbb R^2\setminus\{0\}$. Therefore,   in cases (a), (b) and (d) we get
$x_{2n},\,x_{3n}\,\eta\,\, {\mathfrak A}^3$.

To finish the proof in these cases, it is sufficient to approximate one of operators $D_{rn},\,\,1\!\leq\! r\!\leq\! 3$ 
by opertors $(A_{kn})_{k\in \mathbb Z}$ using Lemmas~\ref{d1.3}--\ref{d3.3}, see Section~\ref{s.444}. Alternatively we can try to approximate $D_{3n},\,D_{2n}$ using Lemma~\ref{l.x(3)A} and its analogue,
see Section~\ref{s.445}, or to approximate $D_{3n},\,D_{2n}$  using Lemma~\ref{exp(x)A3} and its analogue, see Section~\ref{s.446}.
\end{rem}
Note that by Lemma~\ref{l.(011)} we have $\lim_nb_{2n}=\lim_nb_{3n}=\infty$.
In the cases  (a) and (b) the conditions \eqref{S_(123)} are expressed by  \eqref{d_n-gen} as follows:
\begin{eqnarray}
\label{(011)(a)}
&&b=(1,b_{2n},d_nb_{2n}),\,\,\, \sum_n\frac{1}{b_{2n}^2d_n}<\infty, \,\,\,\,\sum_n\frac{1}{d_n}=\infty,\,\,\,\,\lim_nd_n=\infty,\qquad
\\
\label{(011)(b)}
 &&
b=\big(1,b_{2n},db_{2n}(1+c_n)\big),\quad
\sum_n\frac{1}{b_{2n}^2}<\infty,\quad \sum_nc_n^2=\infty.\qquad
\end{eqnarray}
Indeed, to get \eqref{(011)(a)}  we observe that \eqref{S_(123)} are expressed as follows:
\begin{eqnarray*}
&& 
S_1(3)\!\sim\!\sum_n\frac{1}{b_{2n}b_{3n}}=\sum_n\frac{1}{b_{2n}^2d_n}<\infty,\quad 
S_2(3)\sim\sum_n\frac{1}{d_n}=\infty.
\end{eqnarray*}
Condition $S_3(3)\sim\sum_nd_n=\infty$ holds by $\lim_nd_n=\infty$.

In order to get \eqref{(011)(b)}, we express the conditions \eqref{S_(123)} as follows:
\begin{eqnarray*}
 &&
S_1(3)\!\sim\!\sum_n\frac{1}{b_{2n}db_{2n}(1+c_n)}\sim \sum_n\frac{1}{b_{2n}^2}
\!<\!\infty,\\
&&
S_2(3)\sim\sum_n\frac{1}{d_n}=
\sum_n\frac{1}{1+c_n}=\infty,\quad S_3(3)\sim\sum_nd_n\!=\!\sum_n(1+c_n)\!=\!\infty.
\end{eqnarray*}
The  conditions $S_2(3)=\infty$  holds by $\lim_nc_n=0$.

\subsubsection{Approximation of $D_{rn},\,\,1\leq r\leq 3$, 1}
\label{s.444}
By Lemmas~\ref{d1.3}--\ref{d3.3} we have for $1\leq k\leq 3$ (see \eqref{Delta^k(krs)}):
\begin{equation*}
D_{kn}\,\,\eta\,\,{\mathfrak A}^3\Leftrightarrow \Delta_k=\infty,\quad\text{where}\quad 
\Delta_k:=\Delta(Y_k,Y_r,Y_s),
\end{equation*}
and $\{k,r,s\}$ is a cyclic permutation of $\{1,2,3\}$.

Recall that by \eqref{Y_r(011)} we have 
\begin{equation}
\label{Y(r)=}
\Vert Y_1 \Vert^2\sim
\sum_{n\in \mathbb Z}b_{1n}a_{1n}^2,\quad
\Vert Y_2 \Vert^2\sim
\sum_{n\in \mathbb Z}b_{1n}a_{2n}^2,\quad
\Vert Y_3\Vert^2\sim
\sum_{n\in \mathbb Z}b_{1n}a_{3n}^2. 
\end{equation}
Since $\Vert Y_1\Vert^2\sim\sum_{n\in \mathbb Z}
b_{1n}a^2_{1n}\sim S^L_{11}(\mu)
=\infty$, we have  four 
possibilities for $y_{23}:=(y_2,y_3)\in \{0,1\}^2$:
\begin{equation}
 \label{(y_1,y_2,y_3)(011)}
\begin{array}{ccccc}
   &(1.0)&(1.1)&(1.2)&(1.3)\\
y_1&1&1&1&1\\
y_2&0&1&0&1\\
y_3&0&0&1&1\\
\end{array}
\end{equation}
In the case (1.0) we have $\Delta(Y_1,Y_2,Y_3)\sim \Vert Y_1\Vert^2=\infty$, so we can approximate $D_{1n}$ using Lemma~\ref{dx1.3} and the proof is finished.
We should consider the three following cases:
$$
\begin{array}{cccc}
   &(1.1)&(1.2)&(1.3)
\end{array}
$$
In the cases (1.1), (1.2) and (1.3) we have respectively (see the proof of Lemma~\ref{l.S=(0,1,1)})
\begin{eqnarray}
\label{(1.1)}
&& 
\Delta(Y_1,Y_2,Y_3)\sim \Delta(Y_1,Y_2),\quad
\Delta(Y_2,Y_3,Y_1)\sim  \Delta(Y_2,Y_1),\\
\label{(1.2)}
&& 
\Delta(Y_1,Y_2,Y_3)\sim \Delta(Y_1,Y_3),\quad 
\Delta(Y_3,Y_1,Y_2)\sim  \Delta(Y_3,Y_1),\\
\label{(1.3)}
&& 
\Delta(Y_1,Y_2,Y_3),\quad \Delta(Y_2,Y_3,Y_1),\quad
\Delta(Y_3,Y_1,Y_2).
\end{eqnarray}
By \eqref{Y_r(011)} and Lemma~\ref{l.(011)} 
we have respectively in the cases (1.1)--(1.3):
\begin{eqnarray}
\label{nu(12)}
&& 
\nu_{12}( C_1, C_2):=\Vert C_1Y_1+C_2Y_2\Vert^2\sim \sum_{n\in \mathbb Z}b_{1n}\Big( C_1a_{1n}+C_2a_{2n}\Big)^2,\\
\label{nu(13)}
&& 
\nu_{13}( C_1, C_3):=\Vert C_1Y_1+C_3Y_3\Vert^2\sim \sum_{n\in \mathbb Z}b_{1n}\Big( C_1a_{1n}+C_3a_{3n}\Big)^2,\\
\nonumber
&& 
\nu( C_1, C_2, C_3)\!=\!\Vert C_1Y_1+\!C_2Y_2\!+\!C_3Y_3\Vert^2\!\sim\!\\
\label{nu(123)}
&&
\sum_{n\in \mathbb Z}b_{1n}\Big( C_1a_{1n}\!+\!C_2a_{2n}
\!+\!C_3a_{3n}\Big)^2.
\end{eqnarray}
\begin{rem}
\label{r.(011)-(1-3)(0-1)} 
We have three cases (1.1.1), (1.2.1) and (1.3.1) according to whether respectively the expressions in
\eqref{nu(12)}, \eqref{nu(13)} or \eqref{nu(123)} are 
divergent. We can approximate in these cases respectively   $D_{1n}$ and $ D_{2n}$ in \eqref{(1.1)},  $D_{1n}$ and $D_{3n}$ in \eqref{(1.2)} all $D_{1n},\,D_{2n},\,D_{3n}$ in \eqref{(1.3)}. The proof of irreducibilty is finished in these cases
because we have  $D_{rn},\, x_{2n},\,x_{3n}\,\eta\,\, {\mathfrak A}^3$ for some $1\leq r\leq 3$.

If the opposite holds, 
we have three different cases:
\begin{eqnarray*}
&&
(1.1.0)\quad\Vert C_1Y_1+C_2Y_2\Vert<\infty\quad  \text{for some}\quad
(C_1,C_2)\in {\mathbb R}^2\setminus\{0\},\\
&&
(1.2.0)\quad\Vert C_1Y_1+C_3Y_3\Vert<\infty\quad  \text{for some}\quad
(C_1,C_3)\in {\mathbb R}^2\setminus\{0\},\\
&&
(1.3.0)\quad\nu(C_1,C_2,C_3)<\infty
\quad \text{for some}\quad
(C_1,C_2,C_3)\in {\mathbb R}^3\setminus\{0\}.
\end{eqnarray*}
\end{rem}
Recall that by \eqref{S^{L}_{1,23}} we have
$$
S^{L}_{1,23}(\mu,t,s)=\sum_{n\in{\mathbb
Z}}\left[\frac{t^2}{4}\frac{b_{1n}}{b_{2n}}+\frac{s^2}{4}\frac{b_{1n}}{b_{3n}}+
\frac{b_{1n}}{2}\big(\!-\!2a_{1n}\!+\!ta_{2n}\!+\!sa_{3n}\big)^2\right].
$$
\begin{rem}
\label{r.(11)-(13).(0)}
In the case (1.1.0) we have $
\Sigma^{12}=\sum_{n\in \mathbb Z}\frac{b_{1n}}{b_{2n}}=
\infty$, since 
$$S^{L}_{1,23}(\mu,t,0)=\infty,\quad\text{ but} 
\quad\nu_{12}( C_1, C_2)<\infty,$$
and $\Sigma^{13}=\infty$, since 
$
S^L_{13}(\mu)=\sum_{n\in{\mathbb
Z}}\frac{b_{1n}}{2}\Big(\frac{1}{2b_{3n}}+a_{3n}^2\Big)=\infty,\quad\text{but}  \quad\Vert Y_3\Vert^2\sim \sum_{n\in \mathbb Z}b_{1n}a_{3n}^2<\infty$, see \eqref{S^L(kn)} for definition of $S^L_{kr}(\mu)$.

In the case (1.2.0) we conclude that $\Sigma^{13}
=\infty$, since $S^{L}_{1,23}(\mu,0,s)=\infty$, but $\nu_{13}( C_1, C_3)<\infty$, and $\Sigma^{12}=\infty$, since  $S^L_{12}(\mu)=\sum_{n\in{\mathbb
Z}}\frac{b_{1n}}{2}\Big(\frac{1}{2b_{2n}}+a_{2n}^2\Big)=\infty$, but $\Vert Y_2\Vert^2\sim\sum_{n\in \mathbb Z}b_{1n}a_{2n}^2<\infty$.

In the case (1.3.0) we have $\Sigma^{12}=\Sigma^{13}=\infty$, since $S_{1,23}^{L}(\mu,t,s)=\infty,$ but $\nu(C_1,C_2,C_3)\sim\! \sum_{n\in \mathbb Z}b_{1n}\Big( C_1a_{1n}\!+\!C_2a_{2n}
\!+\!C_3a_{3n}\Big)^2<\infty$. 
\end{rem}
So, it remains to consider only the three following cases, when $\Sigma^{12}\!=\!\Sigma^{13}\!=\!\infty$:
$$
\begin{array}{cccc}
   &(1.1.0)&(1.2.0)&(1.3.0)
\end{array}
$$
\subsubsection{Approximation of $D_{2n}$ and $D_{3n}$, 2}
\label{s.445}
 Recall that by Lemma~\ref{l.x(3)A} we have $D_{3n}\,\eta\,\, {\mathfrak A}^3\Leftrightarrow \Sigma_3(\mu)=\infty$ where $\Sigma_3(\mu)$ is defined by  \eqref{x(3)A}
$$
\Sigma_3(\mu):=
\sum_{k\in \mathbb Z}\frac{\frac{1}{2b_{3k}}}{C_k+a_{1k}^2+a_{2k}^2+a_{3k}^2},\quad
C_k=\frac{1}{2b_{1k}}+\frac{1}{2b_{2k}}+\frac{1}{2b_{3k}}.
$$
Similarly, by analogue of Lemma~\ref{l.x(3)A} we can prove that 
$D_{2n}\,\eta\,\, {\mathfrak A}^3\Leftrightarrow \Sigma_2(\mu)=\infty$, where $\Sigma_2(\mu)$ is defined as follows
$$
\Sigma_2(\mu):=\sum_{k\in \mathbb Z}\frac{\frac{1}{2b_{2k}}}{C_k
+a_{1k}^2+a_{2k}^2+a_{3k}^2}.
$$
If one of $\Sigma_2(\mu)$ or $\Sigma_3(\mu)$ is infinite, we can approximate $D_{2n}$ or $D_{3n}$ and the proof is finished. 
If  $\Sigma_2(\mu)+\Sigma_3(\mu)<\infty$, we conclude that
\begin{eqnarray}
\label{sigma(23)(mu)}
&& 
\Sigma_{23}(\mu):=\sum_{k\in \mathbb Z}\frac{\frac{1}{2b_{2k}}
+\frac{1}{2b_{3k}}
}{C_k
+a_{1k}^2+a_{2k}^2+a_{3k}^2}<\infty.
\end{eqnarray}
\subsubsection{Approximation of $D_{2n}$ and $D_{3n}$, 3}
\label{s.446}
By Lemmas~\ref{dx2.3}--\ref{dx3.3} we have 
\begin{eqnarray*}
&& 
D_{2l}{\bf 1}\in\langle x_{2k}A_{kl}{\bf 1}\mid k\in {\mathbb Z}\rangle
\quad\Leftrightarrow\quad 
\Delta(Y_{22},Y_{23},Y_{21})=\infty,\\
&&
D_{3l}{\bf 1}\in\langle x_{3k}A_{kl}{\bf 1}\mid k\in {\mathbb Z}\rangle
\quad\Leftrightarrow\quad 
\Delta(Y_{33},Y_{31},Y_{32})=\infty,
\end{eqnarray*}
where vectors $Y_{rs}$ for $2\leq r\leq3,\,\, 1\leq s\leq 3 $ are defined by \eqref{lambda(2)}--\eqref{y_3(123)}.
We can not prove that $ \Delta(Y_{22},Y_{23},Y_{21})=\infty$ or $\Delta(Y_{33},Y_{31},Y_{32})=\infty$.
We can try to approximate $D_{3n}$  using Lemma~\ref{exp(x)A3}  or
to approximate $D_{2n}$ using an analogue of Lemma~\ref{exp(x)A3}, but it does not work.
Therefore, to approximate  $D_{3n}$ we are forced to prove Lemma~\ref{l.Re-Im-exp.3} the refinement of Lemma~\ref{exp(x)A3} and its analogue for $D_{2n}$, see Remark~\ref{r.Re-Im-(011)} below.
\subsubsection{Two technical lemmas}
\begin{lem}
\label{l.(011)-fin}
Let $a_1,a_2\not\in l_2(\mathbb Z)$ and  $C_1a_1+C_2a_2\in l_2$ for some $(C_1,C_2)\in \mathbb R^2\setminus\{0\},\,\,C_2\not=0$, where $a_r=(a_{rk})_{k\in \mathbb Z},\,\,1\leq r\leq 2$.
Then we have 
\begin{equation}
\label{(011)-fin} 
\sum_{k\in \mathbb Z}\frac{a^2_{1k}}{
1+a_{2k}^2}=\infty.
\end{equation}
\end{lem}
\begin{pf}
We set $Y_r=a_r$, in the case (1.1.0) when $C_1Y_1+C_2Y_2=h\in l_2(\mathbb Z)$ with $C_1C_2>0$ 
(we have $C_1C_2\!\not=\!0$) 
we should take $a_2=-a_1+h$, in the case when  $C_1C_2<0$ we  take $a_2=a_1+h$.  The series  $\sum_{k\in \mathbb Z}\frac{a^2_{1k}}{
1+a_{2k}^2}$ will remain equivalent  with the initial one, if we replace
$(C_1,C_2)$ with $(\pm 1,1)$ in the expression for $h$. 
Fix a small $\varepsilon\!>\!0$ and a large $N
\!\in \!\mathbb N$.
Since $ \vert\pm a+b\vert\leq \vert a\vert +\vert b\vert$,
we get
\begin{eqnarray*}
&& 
\sum_{k\in \mathbb Z}\frac{a^2_{1k}}{
1+a_{2k}^2}=\sum_{k\in \mathbb Z}\frac{a^2_{1k}}{
1+\big(\pm a_{1k}+h_k\big)^2}\geq 
\sum_{k\in \mathbb Z}\frac{a^2_{1k}}{
1+a^2_{1k}+2\vert a_{1k}\vert \vert h_k\vert +h^2_k}\stackrel{\eqref{two-sim-ser}}{\sim}\\
&&\sum_{k\in \mathbb Z}\frac{a^2_{1k}}{
1+2\vert a_{1k}\vert \vert h_k\vert +h^2_k}\stackrel{(*)}{>}
\sum_{k\in \mathbb Z_N}\frac{a^2_{1k}}{
1+2\vert a_{1k}\vert \varepsilon +\varepsilon^2}\stackrel{\eqref{two-ser-alpha} }{\sim} \sum_{k\in \mathbb Z_N}a^2_{1k}=\infty,
\end{eqnarray*}
where $\mathbb Z_N:=\{n\in \mathbb Z\mid \vert n\vert>N\}$. The inequality (*) holds,
since $h\in l_2(\mathbb Z)$ and  we have $\sum_{k\in \mathbb Z_N}h^2_k<\varepsilon^2$ for sufficiently large $N\in \mathbb N$.
\qed\end{pf}

\begin{lem}
\label{l.(011)-fin.3}
Let  $a_1,a_2,a_3\not\in l_2(\mathbb Z)$
and 
$C_1a_1+C_2a_2+C_3a_3\in l_2(\mathbb Z)$ for some $(C_1,C_2,C_3)\in \mathbb R^3\setminus\{0\},\,\,C_3\not=0$, where 
$a_r=(a_{rk})_{k\in \mathbb Z}$ for $1\leq r \leq 3$.  
Then we have 
\begin{equation}
\label{(011)-fin.3} 
\sum_{k\in \mathbb Z}\frac{a^2_{1k}+a^2_{2k}}{
1+a_{3k}^2}=\infty.
\end{equation}
\end{lem}
\begin{pf}
We set $Y_r=a_r$,
in the case (1.3.0), we have $C_1a_1+C_2a_2+C_3a_3=h\in l_2(\mathbb Z)$ for some $(C_1,C_2,C_3)\in \mathbb R^3$, see Remark~\ref{r.(011)-(1-3)(0-1)}. We can take $C_3=1$, then $a_3=-C_1a_1-C_2a_2+h$. 
When $C_1=0$ or $C_2=0$ lemma is reduced to Lemma~\ref{l.(011)-fin}.
Suppose $C_1C_2\not=0$.
The series  $\sum_{k\in \mathbb Z}\frac{a^2_{1k}+a^2_{2k}}{
1+a_{3k}^2}$ will remain equivalent  with the initial one, if we replace
$(C_1,C_2,C_3)$ with $(\pm 1,\pm 1,1)$ in the expression for $h$.
Fix a small $\varepsilon>0$ and a large $N
\in \mathbb N$. Suppose the opposite, i.e.,
\begin{eqnarray*}
&& 
\infty\!>\!
\sum_{k\in \mathbb Z}\frac{a^2_{1k}+ a^2_{2k}}{
1+\big(\pm a_{1k}\pm a_{2k}+h_k\big)^2},\,\,\,
\text{then}
\,\,\,
\infty\!>\!\sum_{k\in \mathbb Z}\frac{(\vert a_{1k}\vert+\vert a_{2k}\vert)^2}{
1+\big(\pm a_{1k}\pm a_{2k}+h_k\big)^2}
\\
&&
\geq\sum_{k\in \mathbb Z}\frac{(\vert a_{1k}\vert+\vert a_{2k}\vert)^2}{
1+a^2_{1k}+a^2_{2k}+2\vert a_{1k}\vert \vert a_{2k}\vert+2\vert a_{1k}\vert \vert h_k\vert +2\vert a_{2k}\vert \vert h_k\vert +h^2_k}\stackrel{\eqref{two-sim-ser}}{\sim}\\
&&\sum_{k\in \mathbb Z}\frac{
(\vert a_{1k}\vert+\vert a_{2k}\vert)^2}{
1+
2\vert a_{1k}\vert \vert h_k\vert +2\vert a_{2k}\vert \vert h_k\vert +h^2_k
}\stackrel{(*)}{>}
\sum_{k\in \mathbb Z_N}\frac{
(\vert a_{1k}\vert+\vert a_{2k}\vert)^2
}{
1+2\big(\vert a_{1k}\vert +
\vert a_{2k}\vert\big) \varepsilon+
\varepsilon^2}\stackrel{\eqref{two-ser-alpha} }{\sim}\\
&&
\sum_{k\in \mathbb Z_N}
(\vert a_{1k}\vert+\vert a_{2k}\vert)^2=\infty,
\end{eqnarray*}
where $\mathbb Z_N:=\{n\in \mathbb Z\mid \vert n\vert>N\}$, contradiction. The inequality (*) holds,
since $h\in l_2(\mathbb Z)$ and  we have $\sum_{k\in \mathbb Z_N}h^2_k<\varepsilon^2$
for sufficiently large $N\in \mathbb N$.
\qed\end{pf}

\begin{rem}
\label{r.Re-Im-(011)}
It is possible to prove an
analogue of
Lemma~\ref{l.Re-Im-exp.3} to approximate $D_{2n}$ with corresponding expressions $\Sigma_2(D,s),\,\Sigma^{\vee}_2(D,s)$ and $\Sigma_3(D),\,\Sigma^{\vee}_3(D)$. 
If one of the expressions $\Sigma_2(D,s),\,\Sigma^{\vee}_2(D,s),\,\Sigma_3(D,s)$ or $\Sigma^{\vee}_3(D,s)$ is divergent for some sequence $s=(s_k)_{k\in \mathbb Z}$, we can approximate $D_{2k}$ or $D_{3k}$ and the proof is finished when  $S=(0.1.1)$ in the cases (a) and  (b).
Suppose that for all sequence $s=(s_k)_{k\in \mathbb Z}$ we have
\begin{equation*}
\Sigma_2(D,s)+\Sigma^{\vee}_2(D,s)+ \Sigma_3(D,s)+\Sigma^{\vee}_3(D,s)<\infty.
\end{equation*}
Then, in particularly, we have for $s^{(r)}=(s_{rk})_{k\in \mathbb Z},\,\,2\leq r\leq 3$
with $\frac{s^2_{rk}}{b_{rk}}\equiv 1$
\begin{eqnarray}
\nonumber 
&& 
 \infty>\Sigma_2(D,s^{(2)})+\Sigma^{\vee}_2(D,s^{(2)})+ \Sigma_3(D,s^{(3)})+\Sigma^{\vee}_3(D,s^{(3)})\sim\\
\label{Sigma^_{23}(D)}
&&
\Sigma_2(D)+\Sigma^{\vee}_2(D)+ \Sigma_3(D)+\Sigma^{\vee}_3(D)=\\
\nonumber
&&
\sum_k\frac{
\frac{1}{2b_{2k}}+a^2_{2k}+\frac{1}{2b_{3k}}+a^2_{3k}}{C_k+a^2_{1k}+a^2_{2k}+a^2_{3k}}\stackrel{\eqref{two-sim-ser}}{\sim}
\sum_k\frac{
\frac{1}{2b_{2k}}+a^2_{2k}+\frac{1}{2b_{3k}}+a^2_{3k}}{\frac{1}{2b_{1k}}+a^2_{1k}}
=:\Sigma^{\vee}_{23}(D)\\
\nonumber
&&
\sim \sum_k\frac{
\frac{b_{1k}}{b_{2k}}+2b_{1k}a^2_{2k}+\frac{b_{1k}}{b_{3k}}+2b_{1k}a^2_{3k}}{1+2b_{1k}a^2_{1k}}\!\stackrel{\eqref{alpha}}{\sim}\!
\sum_k\frac{
a^2_{2k}+a^2_{3k}
}
{1+a^2_{1k}}
=:\Sigma^{a}_{23}(D).
\end{eqnarray}
\end{rem}
%
%
\begin{rem}
\label{r.Re-Im-10-30}
In the case (1.1.0) (resp. the case (1.2.0)) we have $\Vert Y_3 \Vert^2\sim\sum_{n\in \mathbb Z}a_{3n}^2<\infty$
(resp.  $\Vert Y_2 \Vert^2\sim\sum_{n\in \mathbb Z}a_{2n}^2<\infty$) and therefore,
$$
\Sigma^{a}_{23}(D)\sim
\sum_k\frac{a^2_{2k}}{1+a^2_{1k}}=\infty,\quad \text{resp.}\quad
\Sigma^{a}_{23}(D)\sim
\sum_k\frac{a^2_{3k}}{1+a^2_{1k}}=\infty,
$$
by Lemma~\ref{l.(011)-fin}, contradicting \eqref{Sigma^_{23}(D)}.
In the case (1.3.0) we have four cases:\\
(0) when $C_1C_2C_3\not=0$, $C_1a_1+C_2a_2+C_3a_3=h\in l_2(\mathbb Z)$,\\
(1) when $C_1=0$ hence, $C_2C_3\not=0$, $C_2a_2+C_3a_3=h\in l_2(\mathbb Z)$,\\
(2) when $C_2=0$ hence, $C_1C_3\not=0$, $C_1a_1+C_3a_3=h\in l_2(\mathbb Z)$,\\
(3) when $C_3=0$ hence, $C_1C_2\not=0$ $C_1a_1+C_2a_2=h\in l_2(\mathbb Z)$.
\end{rem}
In the case (0) we have $\Sigma^{a}_{23}(D)=\infty$ by Lemma~\ref{l.(011)-fin.3},
contradicting \eqref{Sigma^_{23}(D)}. In the cases (2) and (3) we get $\Sigma^{a}_{23}(D)=\infty$ by Lemma~\ref{l.(011)-fin}, contradicting \eqref{Sigma^_{23}(D)}.
Therefore, one of the expressions $\Sigma_2(D,s),\,\Sigma^{\vee}_2(D,s),\,$
$\Sigma_3(D,s)$ or $\Sigma^{\vee}_3(D,s)$ is convergent hence, we can approximate $D_{2n}$ or $D_{3n}$ and the proof is finished.
To study the case (1) we need the following statement.
\begin{lem}
\label{l.Delta(2,3)} 
 Let $ C_2Y_2+C_3Y_3=h_{23}\in l_2$ for some $(C_2,C_3)\in \big(\mathbb R \setminus \{0\}\big)^2$
and  $C_1Y_1\!+\!C_2Y_2\not\in l_2$ or $C_1Y_1\!+\!C_3Y_3\not\in l_2$ for all $ (C_1,C_r)\in \big(\mathbb R \setminus \{0\}\big)^2$,  then
\begin{equation}
\label{Delta(2,3)}
\Delta(Y_1,Y_2,Y_3)=
\infty.
\end{equation}
\end{lem} 
\begin{pf}
To prove \eqref{Delta(2,3)}
we have by \eqref{Delta(f,g,h)}
\begin{eqnarray}
\nonumber 
&&
\Delta(Y_1,Y_2,Y_3)=
\frac{\Gamma(Y_1)+\Gamma(Y_1,Y_2)+\Gamma(Y_1,Y_3)+\Gamma(Y_1,Y_2,Y_3)}{1+\Gamma(Y_2)+\Gamma(Y_3)+\Gamma(Y_2,Y_3)}\stackrel{(*)}{>}\\
\nonumber
&&
\frac{
\Gamma(Y_1,Y_2)\!+\!\Gamma(Y_1,Y_3)}{1+(1+c_2)\Gamma(Y_2)+\Gamma(Y_3)}\!\sim\!  \frac{
\Gamma(Y_1,Y_2)+\Gamma(Y_1,Y_3)}{\Gamma(Y_2)+\Gamma(Y_3)
}\stackrel{\eqref{sim(2,3)}}{\sim}\\
\label{(Y1,Y2)}
&&
 \frac{
\Gamma(Y_1,Y_2)+\Gamma(Y_1,Y_3)}{2\Gamma(Y_2)
}\stackrel{\eqref{sim(2,3)}}{\sim}
 \frac{
\Gamma(Y_1,Y_2)}{\Gamma(Y_2)}+\frac{\Gamma(Y_1,Y_3)}{\Gamma(Y_3)
}=\infty,
\\
\label{sim(2,3)}
&&
\Gamma(Y_2)\!\sim\! \Gamma(Y_3),\quad  \text{since}\quad C_2Y_2+C_3Y_3=h\in l_2.
\end{eqnarray}
The relation (*) holds by the inequality
$\Gamma(Y_2,Y_3)\leq c_2\Gamma(Y_2)$,  since $C_2Y_2\!+\!C_3Y_3\!\in\! l_2$ for some $ (C_1,C_3)\in \big(\mathbb R \setminus \{0\}\big)^2$, the relation  \eqref{(Y1,Y2)} holds by Lemma~\ref{l.min=proj}.   
To prove \eqref{sim(2,3)} we get since $Y_2\not\in l_2$ and $h\in l_2$, 
\begin{equation*}
\frac{\Gamma(Y_3)}{\Gamma(Y_2)}=
\frac{\Vert Y_3\Vert^2}{\Vert Y_2\Vert^2}=\frac{\Vert Y_2+h\Vert^2}{\Vert Y_2\Vert^2}\leq \Big(\frac{\Vert Y_2\Vert+\Vert h\Vert}{\Vert Y_2\Vert}\Big)^2=1, 
\end{equation*}
If $C_1Y_1\!+\!C_2Y_2\not\in l_2$ for all $ (C_1,C_2)\in \big(\mathbb R \setminus \{0\}\big)^2$, or $C_1Y_1\!+\!C_3Y_3\not\in l_2$ for all $ (C_1,C_3)\in \big(\mathbb R \setminus \{0\}\big)^2$, by Lemma~\ref{l.Delta(2,3)}  we get $\Delta(Y_1,Y_2,Y_3)=\infty$ hence, we can approximate $D_{1n}$ using Lemma~\ref{dx1.3} and the proof is finished. 
If $C_1Y_1\!+\!C_2Y_2=h_{12}\in l_2$ for some  for  $(C_1,C_2)\in \big(\mathbb R \setminus \{0\}\big)^2$ or 
$C_1Y_1\!+\!C_3Y_3=h_{13}\in l_2$ for some $ (C_1,C_3)\in \big(\mathbb R \setminus \{0\}\big)^2$,
then
we have $h_{12}+\alpha h_{23}=C_1Y_1\!+\!C_2Y_2+\!C_3Y_3\in l_2$ or $h_{12}+\beta h_{13}=C_1Y_1\!+\!C_2Y_2+\!C_3Y_3\in l_2$
with $C_1C_2C_3\not=0$ for an appropriate $\alpha\beta\not=0$,  and we are in the case (0).
\qed\end{pf}
%

\subsubsection{Case (c)}
In this case  both $\rho_2(C_2,-C_3)$ and $\rho_3(C_2,-C_3)$ are finite, i.e.,
we are in the case (2)
therefore, we can not approximate $x_{2n}x_{2t},\,x_{3n}x_{3t}$ by Lemmas~\ref{x2x2.3}--\ref{x3x3.3}.  By Lemma~\ref{l.light.123}  $\Sigma_{23}(s)<\infty$ 
and hence, 
$\Sigma_{23}(C_2,C_3)=\infty$.
Indeed, reasoning as in Remark~\ref{perp2-1} 
we see that
\begin{equation}
 \label{(c).m=3}
\mu^{L_{\tau_{23}(\phi,s)}}\perp\mu,\,\,\,\phi\in[0,2\pi),\,\,s>0\Leftrightarrow
\Sigma_{23}(s)+\Sigma_{23}(C_2,C_3)\!=\!\infty,\,\,s>0,
\end{equation}
for $(C_2,C_3)\in{\mathbb R}^2\setminus\{0\}$.
where $\tau_{23}(\phi,s),\,\,\Sigma_{23}(s)$ and $\Sigma_{23}(C_2,C_3)$
are defined as follows:
\begin{eqnarray}
\label{tau(ij)}
&&
\tau_{23}(\phi,s)= \left(
\begin{array}{ccc}
1&0&0\\
0&\cos\phi&s^2\sin\phi\\
0&s^{-2}\sin\phi&-\cos\phi
\end{array}
\right),\\
\label{sigma-ij(s)}
&&
\Sigma_{ij}(s)
=\sum_{n\in {\mathbb Z}}
\left(s^2\sqrt{\frac{b_{in}}{b_{jn}}}-s^{-2}\!\sqrt{\frac{b_{jn}}{b_{in}}}\right)^2,\quad s\in \mathbb R\setminus\{0\}, \\
\label{sigma_{ij}(C_i,C_j)}
&&
\Sigma_{ij}(C_i,C_j)
=\sum_{n\in {\mathbb Z}}(C^2_ib_{in}+C^2_jb_{jn})(C_ia_{in}+C_ja_{jn})^2.
\end{eqnarray}

In this case there are four  possibilities for the pair $(\Sigma^{12},\Sigma^{13})$:\\
(2.1)  $(\Sigma^{12},\Sigma^{13})=(0,0)$, i.e., $\Sigma^{12}<\infty$ and $\Sigma^{13}<\infty$, 
\\
(2.2)  $(\Sigma^{12},\Sigma^{13})=(0,1)$, i.e., $\Sigma^{12}<\infty$, but $\Sigma^{13}=\infty$,\\
(2.3)  $(\Sigma^{12},\Sigma^{13})=(1,0)$, i.e., $\Sigma^{12}=\infty$, but $\Sigma^{13}<\infty$,\\
(2.4)  $(\Sigma^{12},\Sigma^{13})=(1,1)$, i.e., $\Sigma^{12}=\infty$ and $\Sigma^{13}=\infty$.
\begin{lem}
\label{l.ex,011}
In the case {\rm (2.1)}, i.e.,  when $(\Sigma^{12},\Sigma^{13})=(0,0)$, we can approximate $D_{rn}$ for $1\leq r\leq 3$, hence the representation is irreducible. 
\end{lem}
\begin{pf}
Let $\Sigma^{12}<\infty$ and $\Sigma^{13}<\infty$ we have by \eqref{nu(CCC)}
\begin{eqnarray*}
&& \nu( C_1, C_2, C_3)
\sim \sum_{k\in \mathbb Z}b_{1k}(C_1a_{1k}+ C_2a_{2k}+C_3a_{3k})^2\\
&&\stackrel{(2.1)}{\sim} \!\sum_{k\in{\mathbb
Z}}\left[\frac{t^2}{4}\frac{b_{1k}}{b_{2k}}\!+\!\frac{s^2}{4}\frac{b_{1k}}{b_{3k}}\!+\!
\frac{b_{1k}}{2}(-2a_{1k}\!+\!ta_{2k}\!+\!sa_{3k})^2\right]
\stackrel{\eqref{S^{L}_{1,23}}}{=}
S^{L}_{1,23}(\mu,t,s)\!=\!\infty.
\end{eqnarray*}
Hence, $D_{1n},\,D_{2n},\,D_{3n}\,\,\eta\,\, {\mathfrak A}^3$ and the proof is finished.
\qed\end{pf}
\begin{rem}
\label{r.Sigma(123)=(0,1)}
The cases (2.2) and (2.3) do not  occur. 
\end{rem}
 Indeed, by Lemma~\ref{l.light.123} the three series $\Sigma_{23}(s)$ 
 (defined by \eqref{sigma-ij(s)}),  $\sum_{n\in \mathbb Z}c_n^2$ and 
$\sum_{n\in \mathbb Z}e_n^2$
are equivalent where $\frac{s^4b_{2n}}{b_{3n}}=(1+c_n)$, see Lemma~\ref{l.2-eq-ser}. In the case (c) we have $\sum_{n\in \mathbb Z}c_n^2<\infty$, therefore, $\lim_nc_n=0$ and hence, $\lim_nd_n^{-1}=\lim_n\frac{b_{2n}}{b_{3n}}=s^{-4}>0$. Recall that 
$d_n=\frac{d_{2n}}{d_{3n}}=\frac{b_{2n}}{b_{3n}}
$.
But this contradicts $(\Sigma^{12},\Sigma^{13})=(0,1)$, or $(\Sigma^{12},\Sigma^{13})=(1,0)$, since the two series
$$
\Sigma^{12}=\sum_nd_{2n}^{-1}\quad\text{and}\quad
\Sigma^{13}=\sum_nd_{3n}^{-1}
$$
should be equivalent by 
$\lim_n\frac{d_{2n}}{d_{3n}}=s^{-4}>0$.

In the {\bf case (2.4)} we have  
\begin{equation}
 \label{(2.4)}
\Sigma_{23}(s)\!<\!\infty, \quad
 \Sigma_{23}(C_2,C_3)=\infty,\quad 
\Sigma^{12}=\Sigma^{13}=\infty.
\end{equation}
To approximate $D_{rn}$ we need to estimate $\nu( C_1, C_2, C_3)$  defined by \eqref{nu(C,C,C)}.
By \eqref{nu(CCC)} we have
\begin{equation*}
\nu( C_1, C_2, C_3)\sim\sum_{n\in \mathbb Z}b_{1n}\Big(\sum_{r=1}^3C_ra_{rn}\Big)^2.
\end{equation*}
Since $\Vert Y_1\Vert^2=\sum_{n\in \mathbb Z}
b_{1n}a^2_{1n}\sim S^L_{11}(\mu)
=\infty$, in the case (2.4) we have  four 
possibilities
for $y_{23}:=(y_2,y_3)\in \{0,1\}^2$:
$$
\begin{array}{ccccc}
   &(2.4.1)&(2.4.2)&(2.4.3)&(2.4.4)\\
y_1&1&1&1&1\\
y_2&0&1&0&1\\
y_3&0&0&1&1\\
\end{array}
$$
\begin{rem}
The cases (2.4.1)--(2.4.3) are not compatible with the condition 
$ \Sigma_{23}(C_2,C_3)=\infty$ for all $(C_2,C_3)\in{\mathbb R}^2\setminus\{0\}$.
\end{rem}
So it suffices to consider only the   {\bf case (2.4.4)} when  $y_{123}=(1,1,1)$.

The case (2.4.4) splits into two subcases:\\
(2.4.4.1) when $\Sigma_{12}(s_{12})\!<\!\infty$ (resp. $\Sigma_{13}(s_{13})\!<\!\infty$) for some $s_{12},\,s_{13}>0$,\\
(2.4.4.2) when both  $\Sigma_{12}(s_{12})\!=\!\Sigma_{13}(s_{13})=\infty$ for all $s_{12},\,s_{13}>0$.

The case (2.4.4.1) does not occur.
 Indeed, we have  in this case  $\Sigma_{13}(s_{12}s_{23})\!<\!\infty$ (resp. $\Sigma_{12}(s_{13}s_{23}^{-1})\!<\!\infty$) since 
$$
\Sigma_{12}(s_{12})<\infty \Leftrightarrow
\mu_{(s_{12}^4b_1,0)}\sim \mu_{(b_2,0)},\quad
\Sigma_{23}(s_{23})<\infty \Leftrightarrow
\mu_{(s_{23}^4b_2,0)}\sim \mu_{(b_3,0)}
$$
where $\mu_{(b_r,0)}\!=\!\otimes_{n\in \mathbb Z}\mu _{(b_{rk},0)}$
for $1\!\leq\! r\leq 3$. Therefore,  
\begin{equation*}
 \mu_{((s_{12}s_{23})^4b_1,0)}\sim \mu_{(b_3,0)}\quad
\Leftrightarrow\quad\Sigma_{13}(s_{12}s_{23})\!<\!\infty.
\end{equation*}
Similarly, if $\Sigma_{13}(s_{13})\!<\!\infty$ and $\Sigma_{23}(s_{23})\!<\!\infty$ we have
\begin{equation*}
\mu_{(s_{13}^4b_1,0)}\sim \mu_{(b_3,0)},\quad
\mu_{(s_{23}^4b_2,0)}\sim \mu_{(b_3,0)} \Rightarrow \mu_{\big((s_{13}s_{23}^{-1})^4b_1,0\big)}\sim \mu_{(b_2,0)}
\end{equation*}
hence, $\Sigma_{12}(s_{13}s_{23}^{-1})<\infty$.
%
But condition $\Sigma_{13}(s_{12}s_{23})+\Sigma_{12}(s_{13}s_{23}^{-1})<\infty$ contradicts the first condition of \eqref{S_(123)}. 
%
Indeed, we have by Lemma~\ref{l.2-eq-ser}
\begin{eqnarray*}
&&
\Sigma_{12}(s)=\sum_{n\in \mathbb Z}\left(s^2\sqrt{\frac{b_{1n}}{b_{2n}}}-s^{-2}\sqrt{\frac{b_{2n}}{b_{1n}}}\right)^2
\sim
\sum_{n\in \mathbb Z}c_n^2<\infty,\quad s^2\sqrt{\frac{b_{1n}}{b_{2n}}}=1+c_n,\\
&&
\Sigma_{13}(s)=\sum_{n\in \mathbb Z}\left(s^2\sqrt{\frac{b_{1n}}{b_{3n}}}-s^{-2}\sqrt{\frac{b_{3n}}{b_{1n}}}\right)^2\sim
\sum_{n\in \mathbb Z}f_n^2<\infty,\quad s^2\sqrt{\frac{b_{1n}}{b_{3n}}}=1+f_n,\\
\end{eqnarray*}
and $\lim_nc_n=\lim_nf_n=0$. 
This contradicts
$S_1(3)\sim\sum_n\frac{b_{1n}^2}{b_{2n}b_{3n}}<\infty$. Indeed,
$$
\lim_{n\to\infty}\frac{b_{1n}^2}{b_{2n}b_{3n}}=s^{-4}\lim_{n\to\infty}(1+c_n)^2(1+f_n)^2=s^{-4}>0.
$$
Finally, to finish the case $S=(0,1,1)$, we need to consider only the case (2.4.4.2)
when $\Sigma_{12}(s_{12})\!=\!\Sigma_{13}(s_{13})=\infty$ for all $s_{12},\,s_{13}>0$.

By \eqref{S_(123)} and all the previous considerations we have the conditions:
\begin{eqnarray}
\nonumber
&&
S_1(3)\!\sim\!\sum_n\frac{b_{1n}^2}{b_{2n}b_{3n}}\!<\!\infty,\quad
S_2(3)\!\sim\!\sum_n\frac{b_{2n}}{b_{3n}}\!=\!\infty,
\quad 
S_3(3)\!\sim\!\sum_n\frac{b_{3n}}{b_{2n}}\!=\!\infty,\\
\label{S=(011)}
&&
 \Sigma_{23}(C_2,C_3)=\infty,\quad 
\Sigma^{12}=\sum_n\frac{b_{1n}}{b_{2n}}=\infty,\quad \Sigma^{13}=\sum_n\frac{b_{1n}}{b_{3n}}=\infty,\\
\nonumber
&&
\Sigma_{12}(s_{12})\!=\!\Sigma_{13}(s_{13})\!=\!\infty\,\,\text{for all}\,\, s_{12},\,s_{13}>0,\,\, \Sigma_{23}(s_{23})\!<\!\infty \,\,\text{for some}\,\, s_{23}>0.
\end{eqnarray}
\begin{rem}
  \label{r.(011).(c)}
 By \eqref{(b,a)-(b',a')} without loss of generality we can suppose that  $(b_{1n},b_{2n},b_{3n})$ is replaced with $ (1,d_{2n},d_{3n})$. Since $\Sigma_{23}(s)<\infty$, using notations
 \eqref{sigma(23)(s)} and \eqref{c_n,e_n} of Lemma~\ref{l.light.123} 
$$
\Sigma_{23}(s)\!=\!\sum_{n\in \mathbb Z}\Big(
\frac{s^2}{\sqrt{d_n}}\!-\!
\frac{\sqrt{d_n}}{s^2}\Big)^2\!=\!
\sum_{n\in \mathbb Z}\Big(
s^2\sqrt{\frac{d_{2n}}{d_{3n}}}\!-\!
s^{-2}\sqrt{\frac{d_{3n}}{d_{2n}}}\Big)^2,
$$
and taking into consideration \eqref{S=(011)}, we can chose $d_{2n}$ and     $d_{3n}$ as follows:
\begin{equation}
  \label{last(011)}
d_n\!=\!\frac{d_{3n}}{d_{2n}}\!=\!s^4(1+c_n),\,\,\,\sum_nc^2_{n}<\infty
,\,\,\,\sum_n\frac{1}{d^2_{2n}}<\infty,\,\,\,\,\,\sum_n\frac{1}{d_{n}}=\sum_n d_{n}\!=\!\infty.
\end{equation}
Since $\sum_nc^2_{n}<\infty$ 
we have 
$\sum_n\frac{b_{1n}^2}{b_{2n}b_{3n}}\sim \sum_n\frac{1}{d^2_{2n}}$ and
the measures 
$\mu_{(d^{c,s}_3,0)}$ and $\mu_{(d^s_3,0)}$ are equivalent, where
$$
\mu_{(d^{c,s}_3,0)}=\otimes_n\mu_{(s^4d_{2n}(1+c_n),0)},\quad
\mu_{(d^s_3,0)}=\otimes_n\mu_{(s^4d_{2n},0)},
$$
hence, we can choose $c_n\equiv0$ and $s=1$.
So, to finish the case $S=(0,1,1)$ we should prove the irreducibility for 
$b=(1,d_{2n},d_{2n})_{n\in \mathbb Z}$
with the only condition: 
\begin{equation}
\label{last(011)(2=3)}
\sum_nd^{-2}_{2n}<\infty. \quad\text{Since}\quad d_{n}\equiv 1,\quad \text{we have}\quad
\sum_nd^{-1}_{n}=\sum_n d_{n}\!=\!\infty.
\end{equation}
As usual, $a=(a_{1n},a_{2n},a_{3n})_{n\in \mathbb Z}$  should satisfy the orthogonality condition: 
\begin{equation*}
(\mu_{(b,a)}^3)^{L_t}\perp\mu_{(b,a)}^3\quad\text{for all}\quad  t\in {\rm GL}(3,{\mathbb R})\setminus\{e\}.                                                                                                                                          \end{equation*}
\end{rem}
%
%
\begin{ex}
\label{ex.a.1}
The pairwise conditions 
$$\Vert C_rY_r +C_sY_s\Vert^2=\infty\,\,\text{ for} \,\,1\leq r<s\leq 3\quad \text{do not imply}
\quad \Vert \sum_{r=1}^3C_rY_r \Vert^2=\infty.
$$
Let $a_{r,n}=a_{r,-n}$ for $n\in\mathbb N$ and $a_{1,0}=1,\,\,a_{2,0}=2,\,\,a_{3,0}=3$. 
We define  $a_{r,n}$  for $n\in \mathbb N$ as follows  
\begin{equation}
\label{a_n.c-ex}
a_{1n}=  \left\{\begin{array}{ll}
 2&n=2k+1\\
 1&n=2k
\end{array}\right.,\quad
a_{2n}=  \left\{\begin{array}{ll}
 1&n=2k+1\\
 2&n=2k
\end{array}\right.,\quad a_{3n}\equiv 3.
\end{equation}
Then we have clearly for arbitrary $(C_1,C_2,C_3)\in \mathbb R^3\setminus\{0\}$
\begin{eqnarray}
 \label{pairwise}
&&
\Vert C_1a_1\!+\!C_2a_2\Vert^2\!=\!\infty,\,\,\,\,
\Vert C_1a_1\!+\!C_3a_3\Vert^2\!=\!\infty,\,\,\Vert C_2a_2\!+\!C_3a_3\Vert^2\!=\!\infty,\qquad
\\
\label{triple}
&&\text{but}\quad a_1+a_2-a_3=0\quad\text{hence,}\quad \Vert a_1 +a_2-a_3\Vert^2=0.
\end{eqnarray}
\end{ex}
\begin{ex}
\label{ex.(011)1.1}
Let us consider the measure $\mu^3_{{(b,a)}}$ with 
$a=\big(a_{rn}\big)_{r,n}$ from Example~\ref{ex.a.1} and  $b=(b_{1n},b_{2n},b_{3n})$ defined as follows: 
\begin{equation}
\label{(011)1.1}
b_{1n}\equiv 1, \quad
d_{2n}=d_{3n}=\vert n\vert
\quad\text{for}\quad 
n\in \mathbb Z\setminus \{0\},\quad d_{20}=d_{30}=1.
\end{equation}
\end{ex}
\begin{lem}
\label{l.a.1}
 In  Example~\ref{ex.a.1} we have {\rm (}we consider only $n\in \mathbb N${\rm )}
\begin{equation}
\label{=2}
\Delta(a_1,a_2,a_3)=2,\quad \Delta(a_2,a_3,a_1)=2,\quad\Delta(a_3,a_1,a_2)=2,
\end{equation}
where $a_r=(a_{rn})_{n\in \mathbb N},\,\,1\leq r\leq 3$.
\end{lem}
\begin{pf}
Set $a_r(n)=(a_{rl})_{l=1}^n$ for $1\leq r\leq 3$ and $n\in \mathbb N$,  then for $1\leq k<r\leq 3$
$$
\Gamma(a_k(n))\!\sim\! \Gamma(a_1(n)+a_2(n))\!\sim\! n,\,\,\, \Gamma(a_k(n),a_r(n))\!\sim\! \frac{n(n-1)}{2},\,\,\,\Gamma(a_1,a_2,a_3)\!=\!0.
$$
We observe that $\Gamma(a_k,a_k+a_r)=\Gamma(a_k,a_r)$ for $k\not=r$.   Since $a_3=a_1+a_2$ we get
\begin{eqnarray*}
&& \Delta(a_1,a_2,a_3)=\frac{\Gamma(a_1)+\Gamma(a_1,a_2)+\Gamma(a_1,a_3)+\Gamma(a_1,a_2,a_3)}{1+\Gamma(a_2)+\Gamma(a_3)+\Gamma(a_2,a_3)}=\\
&&\frac{\Gamma(a_1)+\Gamma(a_1,a_2)+\Gamma(a_1,a_1+a_2)+\Gamma(a_1,a_2,a_1+a_2)}{1+\Gamma(a_2)+\Gamma(a_1+a_2)+\Gamma(a_2,a_1+a_2)}
=\\
&&\frac{\Gamma(a_1)+2\Gamma(a_1,a_2)}{1+\Gamma(a_2)+\Gamma(a_1+a_2)+\Gamma(a_1,a_2)}=2,\\
&&\Delta(a_2,a_3,a_1)=\frac{\Gamma(a_2)+\Gamma(a_2,a_3)+\Gamma(a_2,a_1)+\Gamma(a_2,a_3,a_1)}{1+\Gamma(a_3)+\Gamma(a_1)+\Gamma(a_3,a_1)}=\\
&&\frac{\Gamma(a_2)+\Gamma(a_2,a_1+a_2)+\Gamma(a_2,a_1)+\Gamma(a_2,a_1+a_2,a_1)}{1+\Gamma(a_1+a_2)+\Gamma(a_1)+\Gamma(a_1+a_2,a_1)}
=\\
&&\frac{\Gamma(a_2)+2\Gamma(a_2,a_1)}{1+\Gamma(a_1+a_2)+\Gamma(a_1)+\Gamma(a_2,a_1)}=2,
\end{eqnarray*}
\begin{eqnarray*}
&&\Delta(a_3,a_1,a_2)=\frac{\Gamma(a_3)+\Gamma(a_3,a_1)+\Gamma(a_3,a_2)+\Gamma(a_3,a_1,a_2)}{1+\Gamma(a_1)+\Gamma(a_2)+\Gamma(a_1,a_2)}=\\
&&\frac{\Gamma(a_1+a_2)+\Gamma(a_1+a_2,a_1)+\Gamma(a_1+a_2,a_2)+\Gamma(a_1+a_2,a_1,a_2)}{1+\Gamma(a_1)+\Gamma(a_2)+\Gamma(a_1,a_2)}
=\\
&&\frac{\Gamma(a_1+a_2)+2\Gamma(a_1,a_2)}{1+\Gamma(a_1)+\Gamma(a_2)+\Gamma(a_1,a_2)}=2.
\end{eqnarray*}
We use two facts for $1\leq r\leq 2$:
$$
\frac{\Gamma(a_1,a_2)}{\Gamma(a_r)}=\infty\quad\text{and}\quad \Gamma(a_1+a_2)\leq \Gamma(a_1)+\Gamma(a_2)+2\sqrt{\Gamma(a_1)\Gamma(a_2)}.
$$
The first relation follows from Lemma~\ref{l.min=proj} since 
$\Vert C_1a_1 +C_2a_2\Vert^2=\infty$. We get 
 $$\frac{\Gamma(a_1,a_2)}{\Gamma(a_r)}=\lim_{n\to\infty} \frac{\Gamma(a_1(n),a_2(n))}{\Gamma(a_r(n))}=\infty.$$
Recall that $\Gamma(a)=\Vert a\Vert^2$.  
The inequality follows from  $\Vert a_1+a_2\Vert\leq \Vert a_1\Vert+\Vert a_2\Vert$, i.e.,
$\sqrt{\Gamma(a_1+a_2)}\leq \sqrt{\Gamma(a_1)}+\sqrt{\Gamma(a_2)}$. 
\qed\end{pf}
By Lemma`\ref{l.mu^3(b,a)^t-perp} we have
$$
(\mu_{(b,a)}^3)^{L_t}\perp\mu_{(b,a)}^3\quad\Leftrightarrow\quad\Sigma^\pm(t):=\Sigma_1^\pm(t)+\Sigma_2(t)=\infty,
$$
where  
\begin{equation*}
\Sigma_1^+(t)\!=\!\sum_{n\in{\mathbb
Z}}
 \sum_{1\leq i<j\leq 3}\Big(t^i_j\sqrt{\frac{b_{in}}{b_{jn}}}\! -\!A^i_j(t)\sqrt{\frac{b_{jn}}{b_{in}}}\Big)^2,\quad 
\end{equation*}
\begin{equation*} 
\Sigma_1^{-}(t)\!=\!\sum_{n\in{\mathbb
Z}}
 \sum_{1\leq i<j\leq 3}\Big(t^i_j\sqrt{\frac{b_{in}}{b_{jn}}}\!+\!A^i_j(t)\sqrt{\frac{b_{jn}}{b_{in}}}\Big)^2,
\end{equation*}
\begin{eqnarray*}
&&\Sigma_2(t^{-1})=\sum_{n\in{\mathbb
Z}}\Big[b_{1n}\big((t_{11}-1)a_{1n}+t_{12}a_{2n}+t_{13}a_{3n}\big)^2+\\
\nonumber
&&
b_{2n}\big(t_{21}a_{1n}+(t_{22}\!-\!1)a_{2n}+t_{23}a_{3n}\big)^2
+b_{3n}\big(t_{31}a_{1n}+t_{32}a_{2n}+(t_{33}\!-\!1)a_{3n}\big)^2\Big].
\end{eqnarray*}
In Example~\ref{ex.(011)1.1} we can not approximate 
$x_{2n},\,x_{3n}$ since in this case we have
\begin{equation}
\label{(011),(c)-ex}
\Delta(Y_2^{(2)},Y_3^{(2)})=1,\quad \Delta(Y_3^{(3)},Y_2^{(3)})\!=1.
\end{equation}
Indeed, by \eqref{(011).m=3} we have
\begin{equation}
 \Delta(Y_2^{(2)},Y_3^{(2)})=
\!\frac{\Gamma(Y_2^{(2)})\!+\!\Gamma(Y_2^{(2)},Y_3^{(2)})}{1+\Gamma(Y_3^{(2)})},\,\, \Delta(Y_3^{(3)},Y_2^{(3)})=
\frac{\Gamma(Y_3^{(3)})\!+\!\Gamma(Y_3^{(3)},Y_2^{(3)})}{1+\Gamma(Y_2^{(3)})}.
\end{equation}

In Example~\ref{ex.(011)1.1} we have 
$d_n=\frac{d_{3n}}{d_{2n}}\equiv 1$
and hence, by  \eqref{S=(011)Y^r_s(alpha)} 
we have
\begin{eqnarray*}
&& \Vert Y_2^{(2)}\Vert^2\sim \sum_{n\in \mathbb Z}\frac{1}{1+2d_n}=\sum_{n\in \mathbb Z}\frac{1}{3},\quad
 \Vert Y_3^{(2)} \Vert^2\sim \sum_{n\in \mathbb Z}\frac{d_n^2}{1+2d_n}=\sum_{n\in \mathbb Z}\frac{1}{3},\\
 \nonumber
&&\Vert Y_2^{(3)} \Vert^2\sim \sum_{n\in \mathbb Z}\frac{1}{d_n^2+2d_n}=\sum_{n\in \mathbb Z}\frac{1}{3},\quad
\Vert Y_3^{(3)}\Vert^2 \sim \sum_{n\in \mathbb Z}\frac{d_n^2}{d_n^2+2d_n}=
\sum_{n\in \mathbb Z}\frac{1}{3}.
\end{eqnarray*}
Therefore,
$
\Gamma(Y_2^{(2)},Y_3^{(2)})
=\Gamma(Y_3^{(3)},Y_2^{(3)})=0,$
and
$$
 \Delta(Y_2^{(2)},Y_3^{(2)})=
\!\frac{\Gamma(Y_2^{(2)})
}{1+\Gamma(Y_3^{(2)})}=1,\quad  \Delta(Y_3^{(3)},Y_2^{(3)})=
\frac{\Gamma(Y_3^{(3)})}{1+\Gamma(Y_2^{(3)})}=1.
$$

Since $b_{1n}\equiv 1$, by \eqref{S=(011)Y_r(alpha)} we get
$$
\Vert Y_1\Vert^2\sim \sum_{n\in \mathbb Z}a_{1n}^2,\quad
\Vert Y_2\Vert^2\sim \sum_{n\in \mathbb Z}a_{2n}^2,\quad
\Vert Y_3\Vert^2\sim \sum_{n\in \mathbb Z}a_{3n}^2,
$$
so by \eqref{nu(CCC)} we have
\begin{equation*}
\nu( C_1, C_2, C_3)\sim\sum_{n\in \mathbb Z}b_{1n}\Big(\sum_{r=1}^3C_ra_{rn}\Big)^2=\sum_{n\in \mathbb Z}\Big(\sum_{r=1}^3C_ra_{rn}\Big)^2.
\end{equation*}
But in  Example~\ref{ex.a.1} there does not exist
$t\in \pm{\rm SL}(3,\mathbb R)\setminus \!\{e\}$ such that 
$\nu( C_1, C_2, C_3)\!=\!\infty$ for all
$(C_1,C_2,C_3)\!\in\!{\mathbb R}^3\setminus\!\{0\}$ to approximate some $D_{rn}$.

\subsubsection{Approximations of $x_{2k}x_{2r}+x_{3k}x_{3r}$ in the case (c)}
Since we can not approximate $x_{2n}x_{2t},\,x_{3n}x_{3t}$  using Lemmas~\ref{x2x2.3}--\ref{x3x3.3} in the case (c), we shall try to approximate 
$x_{2k}x_{2r}+s^4x_{3k}x_{3r}$
by an appropriate combinations of $A_{kn}A_{rn}$ for $n\in \mathbb Z$. Let $s=1$, the general case is similar.
\begin{lem}
\label{l.(xx+xx)(23)} 
For any  $k,r\in {\mathbb Z}$ one has
\begin{eqnarray}
 \label{(xx+xx)(23)}
 &&
(x_{2k}x_{2r}\!+\!x_{3k}x_{3r}){\bf 1}\in\langle A_{kn}A_{rn}{\bf 1}\mid n\in{\mathbb
Z}\rangle\,\,\Leftrightarrow\,\,
\Delta(Y^{(2)},Y^{(1)})\!=\!\infty,\\
\nonumber
&&
\text{where}\quad 
Y^{(r)}\!=\!\Big(\frac{b_{rn}}{\sqrt{\lambda_n}}\Big)_{n\in \mathbb Z},\,\,\,1\leq r\leq 2,\,\,\,\lambda_n\!=\!(b_{1n}+b_{2n}+b_{3n})^2-b^2_{1n}.
\end{eqnarray}
\end{lem}
\begin{pf}The proof of  Lemma~\ref{l.(xx+xx)(23)}  is based on Lemma~\ref{l.m=2}. 
We study when
$
(x_{2k}x_{2r}+x_{3k}x_{3r}){\bf 1}
\in\langle A_{kn}A_{rn}{\bf 1}\mid n\in{\mathbb
Z}\rangle.
$
Since
\begin{eqnarray*}
A_{kn}A_{rn}&&=(x_{1k}D_{1n}+x_{2k}D_{2n}+x_{3k}D_{3n})(x_{1r}D_{1n}+x_{2r}D_{2n}+x_{3r}D_{3n})\\
&&=x_{1k}x_{1r}D_{1n}^2+
x_{2k}x_{2r}D_{2n}^2+
x_{3k}x_{3r}D_{3n}^2+
(x_{1k}x_{2r}+x_{2k}x_{1r})D_{1n}D_{2n}\\
&&
+(x_{1k}x_{3r}+x_{3k}x_{1r})D_{1n}D_{3n}+
(x_{2k}x_{3r}+x_{3k}x_{2r})D_{2n}D_{3n},
\end{eqnarray*}
and $MD_{rn}^2{\bf 1}=-\frac{b_{rn}}{2}$, for $2\leq r\leq 3$ we take $t=(t_n)_{n=-m}^m$ as follows:
$$
(t,b_2)\!=\!(t,b_3)\!=\!1,\,\,
\text{where} \,\, t\!=\!(t_n)_{k=-m}^m,\,\, 
b_2\!=\!-(b_{2n}/2)_{n=-m}^m,\,\, 
b_3\!=\!-(b_{3n}/2)_{n=-m}^m.
$$ 
We have
\begin{eqnarray*}
&&\Vert \big[\sum_{n=-m}^mt_nA_{kn}A_{rn}-\big(x_{2k}x_{2r}+x_{3k}x_{3r}\big)\big]{\bf
1}\Vert^2=\\
&&\Vert \!\!\sum_{n=-m}^mt_n\big[
x_{1k}x_{1r}D_{1n}^2\!+\!
x_{2k}x_{2r}
\big(D_{2n}^2\!+\!\frac{b_{2n}}{2}\big)\!+\!
x_{3k}x_{3r}\big(D_{3n}^2\!+\!\frac{b_{3n}}{2}\big)+(x_{1k}x_{2r}+\\
&&
x_{2k}x_{1r})
D_{1n}D_{2n}\!+\!(x_{1k}x_{3r}+x_{3k}x_{1r})D_{1n}D_{3n}\!+
\!(x_{2k}x_{3r}+x_{3k}x_{2r})D_{2n}D_{3n}
\Big]{\bf 1}\Vert^2\\
&=&\sum_{-m\leq n,l\leq m}(f_n,f_l)t_nt_l=:(A_{2m+1}t,t),
\end{eqnarray*}
where $A_{2m+1}=(f_n,f_l)_{n,l=-m}^m$ and
\begin{eqnarray}
\label{f_n(r),f_n(ij)}
&&
 f_n= \sum_{i=1}^3f_n^i
 +\sum_{1\leq i<j\leq 3}f_n^{ij},\quad\text{with}\\
 \nonumber
 &&
f_n^i= x_{ik}x_{ir}\Big(D_{in}^2+\frac{b_{in}}{2}(1-\delta_{i1})
\Big){\bf 1},
\quad
f_n^{ij}=(x_{ik}x_{jr}+x_{jk}x_{ir})D_{in}D_{jn}{\bf 1} 
\end{eqnarray}
for $1\leq i\leq 3,\,\, 1\leq i<j\leq 3$.
Since 
$
 f_n^{i'}\perp f_n^{ij},\quad f_n^{ij}\perp f_n^{i'j'}
$ 
for different  $(ij),\,\,(i'j')$,
writing $c_{kn}=\Vert x_{kn}\Vert^2=\frac{1}{2b_{kn}}+a_{kn}^2$,  we get
\begin{eqnarray*}
&&(f_n,f_n)=
\sum_{i=1}^3\Vert f_n^i\Vert^2
+\sum_{1\leq i<j\leq 3}\Vert f_n^{ij}\Vert^2=
\\ 
&&c_{1k}c_{1r}3\Big(\frac{b_{1n}}{2}\Big)^2+
c_{2k}c_{2r}2\Big(\frac{b_{2n}}{2}\Big)^2+
c_{3k}c_{3r}2\Big(\frac{b_{3n}}{2}\Big)^2+
\end{eqnarray*}
\begin{eqnarray*}
&&
\big(c_{1k}c_{2r}+c_{2k}c_{1r}+2a_{1k}a_{2r}a_{2k}a_{1r}\big)\frac{b_{1n}}{2}\frac{b_{2n}}{2}+\big(c_{1k}c_{3r}
+c_{3k}c_{1r}+2a_{1k}a_{3r}a_{3k}a_{1r}\big) \\
&&
\times\frac{b_{1n}}{2}\frac{b_{3n}}{2}+
\big(c_{2k}c_{3r}+c_{3k}c_{2r}+2a_{2k}a_{3r}a_{3k}a_{2r}\big)\frac{b_{2n}}{2}\frac{b_{3n}}{2} \sim(b_{1n}+b_{2n}+b_{3n})^2,\\ 
&&(f_n,f_l)=(f_n^1,f_l^1)=
c_{1k}c_{1r}\frac{b_{1n}}{2}\frac{b_{1l}}{2}\sim b_{1n}b_{1l}.
\end{eqnarray*}
Finally, we get
\begin{equation}
\label{(A_{r,s}),xx+xx,m=3}
(f_n,f_n) \sim (b_{1n}+b_{2n}+b_{3n})^2, \quad
(f_n,f_l)\sim b_{1n}b_{1l}
,\quad n\not=l.
\end{equation}

Set
\begin{equation}
\lambda_n=(b_{1n}+b_{2n}+b_{3n})^2-b_{1n}^2,\quad g_n=(b_{1n}),
\end{equation}
then
\begin{equation}
(f_n,f_n)\sim \lambda_n+(g_n,g_n),\quad (f_n,f_l)\sim (g_n,g_l).
\end{equation}

For $A_{2m+1}=\big((f_n,f_l)\big)_{n,l=-m}^m$ and
$b_2=b_3=-(b_{2n}/2)_{n=-m}^m\in {\mathbb R}^{2m+1}$ we have
$$
A_{2m+1}= \sum_{n=-m}^m\lambda_nE_{nn}+\gamma(g_{-m},\dots,g_0,\dots,g_m).
$$
To finish the proof, it  suffices to use  Lemma~\ref{l.m=2}.
\qed\end{pf}
\begin{rem}
\label{r.(xx+xx)(23)}
In the case (c) we can approximate $x_{2k}x_{2r}+x_{3k}x_{3r}$ since $\Delta(Y^{(2)},Y^{(1)})\!=\!\infty$.
\end{rem}
Indeed, by \eqref{(xx+xx)(23)} we have
\begin{eqnarray*}
 &&
 \Delta(Y^{(2)},Y^{(1)})=\frac{\Gamma(Y^{(2)})\!+\!\Gamma(Y^{(2)},Y^{(1)})}{1+\Gamma(Y^{(1)})}>
 \frac{\Gamma(Y^{(2)})}{1+\Gamma(Y^{(1)})}=\\
 &&
 \frac{\sum_{n\in \mathbb Z}\frac{b^2_{2n}}{\lambda_n}}
 {1+\sum_{n\in \mathbb Z}\frac{b^2_{1n}}{\lambda_n}}=
 \frac{\sum_{n\in \mathbb Z}\frac{d^2_{2n}}{(1+2d_{2n})^2-1}}
 {1+\sum_{n\in \mathbb Z}\frac{1}{(1+2d_{2n})^2-1}}\sim
  \sum_{n\in \mathbb Z}\frac{d^2_{2n}}{d_n+d^2_{2n}}=\infty
\end{eqnarray*}
since  by \eqref{last(011)} we have   $\sum_n\frac{1}{d^2_{2n}}<\infty$. Therefore,
\begin{equation*}
\Gamma(Y^{(1)})=\sum_{n\in \mathbb Z}\frac{1}{(1+2d_{2n})^2-1}
 \sim
  \sum_{n\in \mathbb Z}\frac{1}{d_n+d^2_{2n}}<\infty. 
\end{equation*}

\begin{lem}
\label{l.x(2)} 
We have for all $k\in \mathbb Z$
\begin{equation}
 \label{x(2)}  
x_{2k}{\bf 1}\in\langle (x_{2k}x_{2n}+x_{3k}x_{3n}){\bf 1}\mid n\in{\mathbb
Z}\rangle \,\,\Leftrightarrow\,\,
\sigma_{2}(\mu)=\sum_{n\in\mathbb
Z}\frac{a_{2n}^2}{\frac{1}{2b_{2n}}+\frac{1}{2b_{3n}}+a_{3n}^2}=\infty.
\end{equation}
\end{lem}
\begin{pf}
 Recall the notation $c_{rn}=\frac{1}{2b_{rn}}+a_{rn}^2$.
Since $Mx_{2n}{\bf 1}=a_{2n}$  we take $t=(t_n)_{n=-m}^m$ as follows: $(t,a_2)=1$, where $a_2=(a_{2n})_{n=-m}^m.$ 
 We have
\begin{eqnarray*}
 &&\Vert \big[\sum_{n=-m}^mt_n\big(x_{2k}x_{2n}+x_{3k}x_{3n}\big)-x_{2k}x_{2n}\big]{\bf 1}\Vert^2=\\
&&
\Vert \big[\sum_{n=-m}^mt_n\big(x_{2k}(x_{2n}-a_{2n})+x_{3k}x_{3n}\big)\big]{\bf
1}\Vert^2= \Vert x_{2k}{\bf 1}\Vert^2 
 \Vert \sum_{n=-m}^mt_n (x_{2n}-a_{2n}){\bf 1}\Vert^2 \\
 &&
+\Vert x_{3k} {\bf 1}\Vert^2 
   \Vert \sum_{n=-m}^mt_n x_{3n} {\bf 1}\Vert^2=
c_{2k} \sum_{n=-m}^mt^2_n\frac{1}{2b_{2n}}  +
c_{3k} \sum_{n=-m}^mt^2_n\Big(\frac{1}{2b_{3n}}+ a^2_{3n}\Big)\\
&&
 \sim  \sum_{n=-m}^mt^2_n\Big(\frac{1}{2b_{2n}}+\frac{1}{2b_{3n}}+ a^2_{3n}\Big).
\end{eqnarray*}
By \eqref{A.min2} we get \eqref{x(2)}.
\qed\end{pf}
Similarly, we prove the following lemma.
\begin{lem}
\label{l.x(3)}
We have for all $k\in \mathbb Z$
\begin{equation}
 \label{x(3)}  
x_{3k}{\bf 1}\in\langle (x_{2k}x_{2n}+x_{3k}x_{3n}){\bf 1}\mid n\in{\mathbb
Z}\rangle \,\,\Leftrightarrow\,\,
\sigma_{3}(\mu)=\sum_{\in\mathbb
Z}\frac{a_{3n}^2}{\frac{1}{2b_{2n}}+\frac{1}{2b_{3n}}+a_{2n}^2}=\infty.
\end{equation}
\end{lem}
\begin{rem}
\label{r.S=(011)x(2)-x(3)} 
Suppose that $\sigma_{2}(\mu)+\sigma_{3}(\mu)<\infty$, this contradicts $\Sigma_{23}(C_2,C_3)=\infty$ for $(C_2,C_3)\in{\mathbb R}^2\setminus\{0\}$, where $\Sigma_{23}(C_2,C_3)$ is defined by 
 \eqref{sigma_{ij}(C_i,C_j)}: 
\begin{equation*}
\Sigma_{23}(C_2,C_3)
=\sum_{n\in {\mathbb Z}}(C^2_2b_{2n}+C^2_3b_{3n})(C_2a_{2n}+C_3a_{3n})^2.
\end{equation*}
\end{rem}
\begin{pf}
Indeed, we have
\begin{eqnarray*}
&&
\infty>\sigma_{2}(\mu)+\sigma_{3}(\mu)=
\sum_{\in\mathbb
Z}\frac{a_{2n}^2}{\frac{1}{2b_{2n}}+\frac{1}{2b_{3n}}+a_{3n}^2}+
\sum_{\in\mathbb
Z}\frac{a_{3n}^2}{\frac{1}{2b_{2n}}+\frac{1}{2b_{3n}}+a_{2n}^2}\sim\\
&&
\sum_{\in\mathbb
Z}\frac{a_{2n}^2+a_{3n}^2}{\frac{1}{2b_{2n}}\!+\!\frac{1}{2b_{3n}}+a_{2n}^2+a_{3n}^2}\!\sim\!
\sum_{\in\mathbb
Z}\frac{a_{2n}^2+a_{3n}^2}{\frac{1}{2b_{2n}}+\frac{1}{2b_{3n}}}\!
\stackrel{\eqref{last(011)}}{=}\!
\frac{2}{1\!+\!s^{-4}}\sum_{n\in\mathbb
Z}b_{2n}(a_{2n}^2+a_{3n}^2).
\end{eqnarray*}
This contradicts $\Sigma_{23}(C_2,C_3)\!=\!\infty$. Indeed,  by $b_{3n}\!=\!s^4b_{2n}$ (see \eqref{last(011)}) we have
$$
\hskip 1.6cm
\Sigma_{23}(C_2,C_3)=
\sum_{n\in {\mathbb Z}}\big(C^2_2+C^2_3s^4\big)b_{2n}\big(C_2a_{2n}+C_3a_{3n}\big)^2<\infty.
\hskip 1.6cm \Box
$$
\end{pf}
%
%
Finally, we have $\sigma_{2}(\mu)+\sigma_{3}(\mu)=\infty$, and therefore we have
$x_{rn}\,\eta\,\, {\mathfrak A}^3$ for some $2\leq r\leq 3$.
Let $x_{3n}\,\eta\,\, {\mathfrak A}^3$, then we can approximate $x_{2n}$ by combinations of $x_{2n}x_{2k},\,k\in \mathbb Z$ using an  analogue of Lemma~\ref{5.3}.
%
%
To approximate $D_{rn},\,\,1\leq r\leq 3$ we again follows Section~\ref{s.444}.
As in \eqref{Y(r)=} we get  
\begin{equation*}
\Vert Y_1 \Vert^2\sim
\sum_{n\in \mathbb Z}a_{1n}^2,\quad
\Vert Y_2 \Vert^2\sim
\sum_{n\in \mathbb Z}a_{2n}^2,\quad
\Vert Y_3\Vert^2\sim
\sum_{n\in \mathbb Z}a_{3n}^2. 
\end{equation*}
Indeed, for example, by \eqref{II}   we get
\begin{equation*}
\Vert Y_1 \Vert^2=
\sum_{n\in \mathbb Z}\frac{a_{1n}^2}
{\frac{1}{2b_{1n}}+\frac{1}{2b_{2n}}+\frac{1}{2b_{3n}}} =
\sum_{n\in \mathbb Z}\frac{a_{1n}^2}
{\frac{1}{2}+\frac{1}{d_{2n}}}\stackrel{\eqref{last(011)(2=3)}}{\sim}
\sum_{n\in \mathbb Z}a_{1n}^2.
\end{equation*}
Again, as in  \eqref{(y_1,y_2,y_3)(011)} we have four possibilities: (1.0), (1.1), (1.2) and (1.3). The corresponding expressions in \eqref{nu(12)}, \eqref{nu(13)}, \eqref{nu(123)} becomes as follows:
\begin{eqnarray*}
&& 
\nu_{12}( C_1, C_2):=\Vert C_1Y_1+C_2Y_2\Vert^2\sim \sum_{n\in \mathbb Z}\Big( C_1a_{1n}+C_2a_{2n}\Big)^2,\\
&& 
\nu_{13}( C_1, C_3):=\Vert C_1Y_1+C_3Y_3\Vert^2\sim \sum_{n\in \mathbb Z}\Big( C_1a_{1n}+C_3a_{3n}\Big)^2,\\
\nonumber
&& 
\nu( C_1, C_2, C_3)\!=
\sum_{n\in \mathbb Z}\Big( C_1a_{1n}\!+\!C_2a_{2n}
\!+\!C_3a_{3n}\Big)^2.
\end{eqnarray*}

To study the cases (1.1.1)--(1.3.1) we should use Remark~\ref{r.(011)-(1-3)(0-1)}. 
 We can approximate in these cases respectively   $D_{1n}$ and $ D_{2n}$ in \eqref{(1.1)},  $D_{1n}$ and $D_{3n}$ in \eqref{(1.2)} all $D_{1n},\,D_{2n},\,D_{3n}$ in \eqref{(1.3)}. The proof of irreducibilty is finished in these cases
because we have  $D_{rn},\, x_{2n},\,x_{3n}\,\eta\,\, {\mathfrak A}^3$ for some $1\leq r\leq 3$.
Following  Remark~\ref{r.Re-Im-(011)} 
we can use 
Lemma~\ref{l.Re-Im-exp.3}
and its analogue 
to approximate $D_{2n}$ and  $D_{3n}$  with corresponding  expressions $\Sigma_2(D,s),\,\Sigma^{\vee}_2(D,s)$ and $\Sigma_3(D),\,\Sigma^{\vee}_3(D)$. 
If one of the expressions $\Sigma_2(D,s),\,\Sigma^{\vee}_2(D,s),\,\Sigma_3(D,s)$ or $\Sigma^{\vee}_3(D,s)$ is divergent for some sequence $s=(s_k)_{k\in \mathbb Z}$, we can approximate $D_{2k}$ or $D_{3k}$ and the proof is finished. Suppose that for all sequence $s=(s_k)_{k\in \mathbb Z}$ we have
\begin{equation*}
\Sigma_2(D,s)+\Sigma^{\vee}_2(D,s)+ \Sigma_3(D,s)+\Sigma^{\vee}_3(D,s)<\infty.
\end{equation*}
Then, by \eqref{Sigma^_{23}(D)} we have
\begin{eqnarray*}
 &&
\infty>\Sigma^{\vee}_{23}(D)=
\sum_k\frac{
\frac{1}{2b_{2k}}+a^2_{2k}+\frac{1}{2b_{3k}}+a^2_{3k}}{\frac{1}{2b_{1k}}+a^2_{1k}}
=
\sum_k\frac{
\frac{1}{d_{2k}}+a^2_{2k}+a^2_{3k}}{\frac{1}{2}+a^2_{1k}}
\stackrel{\eqref{last(011)(2=3)}}{\sim}\\
&&
\sum_k\frac{a^2_{2k}+a^2_{3k}}{1+a^2_{1k}}=:
\Sigma^{a}_{23}(D).
\end{eqnarray*}
To study the cases (1.1.0)--(1.3.0) we should follow Remark~\ref{r.Re-Im-10-30}.

\subsection{Case $S=(1,1,1)$}
Denote by 
\begin{equation}
 \label{Sigma[123]}
 \Sigma_{123}(s)=(\Sigma_{12}(s_1), \Sigma_{23}(s_2),\Sigma_{13}(s_3)),
\end{equation}
where  $s=(s_1,s_2,s_3)$ and $\Sigma_{ij}(s)$ are defined by \eqref{sigma-ij(s)} for $1\leq i<j\leq 3$.
In terms of Remark~\ref{r.(0,1)}, we have   $2^3$ possibilities
for $\Sigma_{123}(s)
\in \{0,1\}^3$:
$$
\begin{array}{ccccccccc}
   &(0)&(1)&(2)&(3)&(4)&(5)&(6)&(7)\\
\Sigma_{12}(s_1)&0&0&0&0&1&1&1&1\\
\Sigma_{23}(s_2)&0&0&1&1&0&0&1&1\\
\Sigma_{13}(s_3)&0&1&0&1&0&1&0&1\\
\end{array}
$$
The cases (1), (2) and (4) and respectively the cases (3), (5) and (6) result from cyclic permutations of three  measures $\mu^{(1)},\mu^{(2)},\mu^{(3)}$ defined as follows:
\begin{equation}
\label{mu^(r)}
\mu^{(r)}=\otimes_{n\in \mathbb Z}\mu_{(b_{rn},a_{rn})},\,\,1\leq r\leq 3,\quad
\mu^{(r)}_0=\otimes_{n\in \mathbb Z}\mu_{(b_{rn},0)},\,\,1\leq r\leq 3.
\end{equation}
The case (1), (2) and (4) can not be realized. We prove this only in the case (1). By Lemma~\ref{l.eq-prod.m=1} we have
$\Sigma_{12}(s_1)<\infty \,\,\,
\Leftrightarrow \,\,\,\mu_0^{(1)}\sim \mu_0^{(2)}$ and $\Sigma_{23}(s_2)<\infty \Leftrightarrow \mu_0^{(2)}\sim \mu_0^{(3)}$ hence, $\mu_0^{(1)}\sim \mu_0^{(3)}$, that contradicts $\Sigma_{13}(s_2)=\infty \Leftrightarrow \mu_0^{(1)}\perp\mu_0^{(3)}.$
Finally, we are left with  the three cases
(0), (3) and (7):\\
the {\bf case (0)}, i.e., $\Sigma_{123}(s)=(0,0,0)$,\\
the {\bf case (3)}, i.e.,  $\Sigma_{123}(s)=(0,1,1)$,\\
the {\bf case (7)}, i.e., $\Sigma_{123}(s)=(1,1,1)$.

\subsubsection{Case $\Sigma_{123}(s)=(0,0,0)$}
In the case (0), we have  for some $s=(s_1,s_2,s_3)\in \big(\mathbb R_+\big)^3$
\begin{equation*}
\Sigma_{12}(s_1)<\infty,\quad \Sigma_{23}(s_2)<\infty,\quad\Sigma_{13}(s_3)<\infty. 
\end{equation*}
In this case we get $\mu_0^{(1)}\sim\mu_0^{(2)}\sim\mu_0^{(3)}.$
By \eqref{(b,a)-(b',a')} we can make the following change of the variables:
\begin{equation*}
\left(\begin{smallmatrix}
b_{1n}&b_{2n}&b_{3n}\\
a_{1n}&a_{2n}&a_{3n}
\end{smallmatrix}
\right)\to
\left(\begin{smallmatrix}
b_{1n}'&b_{2n}'&b_{3n}'\\
a_{1n}'&a_{2n}'&a_{3n}'
\end{smallmatrix}
\right)=
\left(\begin{smallmatrix}
1&\frac{b_{2n}}{b_{1n}}&\frac{b_{3n}}{b_{1n}}\\
a_{1n}\sqrt{b_{1n}}&a_{2n}\sqrt{b_{1n}}&a_{3n}\sqrt{b_{1n}}
\end{smallmatrix}
\right).
\end{equation*}
\begin{rem}
\label{r.Sigma-123=(000)} 
By Lemma~\ref{l.eq-prod}, we can suppose that 
\begin{equation}
 \label{b=1+c}
b\!=\!(b_{1n},b_{2n},b_{3n})_{n\in \mathbb Z}\!=\!(1,1+c_n,1+e
_n)_{n\in \mathbb Z},\,\,\,
\sum_nc_n^2<\infty,\,\,\,
\sum_ne_n^2<\infty.
\end{equation}
But the two measures 
$\mu_{(b,a)}$ and $\mu_{({\mathbb I},a)}$ are equivalent, where $b$ is defined by \eqref{b=1+c} and 
\begin{equation}
 \label{bf I}
{\mathbb I}:=(1,1,1)_{n\in \mathbb Z}.
\end{equation}
Finally, it is sufficient to consider the measure $\mu_{({\bf 1},a)}$. 
\end{rem}

\begin{ex}
 \label{ex.(111)b=1,a-alt}
Let $b_{1n}=b_{2n}=b_{3n}\equiv 1,\,\,n\in \mathbb Z$. \\
(a) Take  $a_n=(a_{1n},a_{2n},a_{3n}),\,\,n\in \mathbb Z$ as it was defined in Example~\ref{ex.a.1}: 
\begin{equation*}
a_{1n}=  \left\{\begin{array}{ll}
 2&n=2k+1\\
 1&n=2k
\end{array}\right.,\quad
a_{2n}=  \left\{\begin{array}{ll}
 1&n=2k+1\\
 2&n=2k
\end{array}\right.,\quad a_{3n}\equiv 3.
\end{equation*}
Then $a_1+a_2-a_3=0$, where $a_{r}=(a_{rn})_{n\in \mathbb Z}$.\\
(b) Take any $a_{r}=(a_{rn})_{n\in \mathbb Z}$ such that 
$a_1,a_2,a_3\not\in l_2$, but $C_1a_1+C_2a_2+C_3a_3\in l_2(\mathbb Z)$ for some $(C_1,C_2,C_3)\in  {\mathbb R}^3\setminus\{0\}$.
 \end{ex}
\begin{ex}
 \label{ex.(111)b=1,a-any}
Let  $b_{1n}=b_{2n}=b_{3n}\equiv 1,\,\,n\in \mathbb Z$ and $a=(a_{1n},a_{2n},a_{3n})_{n\in \mathbb Z}$  such that $a_1,a_2,a_3\not\in l_2$, 
but 
the measure $\mu_{(b,a)}^3$ satisfies 
the orthogonality conditions. 
The case 
$\Sigma_{123}(s)=(0,0,0)$ is reduced to this example.
 \end{ex}
\begin{rem}
\label{r.(111)-ex}
Since the measure $\mu_{(b,0)}^3$ is {\it standard} in Example~\ref{ex.(111)b=1,a-alt} and \ref{ex.(111)b=1,a-any}, i.e.,
it is invariant under rotations $\pm{\rm O}(3)$, we have 
\begin{equation}
 \label{r.O(3)-inv-ex}
(\mu_{(b,0)}^3)^{L_t}=\mu_{(b,0)}^3\quad \text{for all}\quad t\in \pm{\rm O}(3).
\end{equation}
By Lemma~\ref{l.SL(3)-perp}, the orthogonality condition 
 $(\mu_{(b,a)}^3)^{L_t}\perp\mu_{(b,a)}^3$ for 
$t\in \pm{\rm O}(3)
\setminus\{e\}$, 
is equivalent to  $$\Sigma^\pm_1(t)+\Sigma_2(t)=\infty,$$
where $\Sigma_1^+(t),\,\,\Sigma_1^{-}(t)$ are defined by \eqref{Sigma_1^+(t).m=3}, 
 \eqref{Sigma_1^{-}(t).m=3} and  $\Sigma_2(t)$  is defined by \eqref{Sigma_3}. 
By \eqref{r.O(3)-inv-ex} we get $\Sigma^\pm_1(t)<\infty$ in Example~\ref{ex.(111)b=1,a-alt} and \ref{ex.(111)b=1,a-any} hence the orthogonality condition 
 $(\mu_{(b,a)}^3)^{L_t}\perp\mu_{(b,a)}^3$ for 
$t\in \pm{\rm O}(3)
\setminus\{e\}$ is equivalent to $\Sigma_2(t)=\infty$.
Further, to prove the irreducibility in
Example~\ref{ex.(111)b=1,a-alt} and
\ref{ex.(111)b=1,a-any} we should show that $\Sigma_2(t)=\infty$ for all
$t\in \pm{\rm O}(3)\setminus\{e\}$
implies  $$\Vert C_1Y_1+ C_2Y_2+C_3Y_3\Vert^2=\infty\quad \text{for all} \quad(C_1,C_2,C_3)\in 
 {\mathbb R}^3\setminus\{0\}.$$
\end{rem}
\begin{lem}
(1) The representations corresponding to the measures in  Example~\ref{ex.(111)b=1,a-alt} (a) and (b) are reducible.\\
(2) The representations corresponding to the measures in  Example~\ref{ex.(111)b=1,a-any}  are irreducible.
\end{lem}
\begin{pf}
To prove the  part (1) of theorem, by Remark~\ref{r.(111)-ex} and \eqref{r.O(3)-inv-ex},
we should find for the measure 
in  Example~\ref{ex.(111)b=1,a-alt} an element
$t\!\in \pm{\rm O}(3)\setminus\{e\}$ 
such that $\Sigma_2(t)\!<\!\infty$. 
This will imply  $(\mu_{(b,a)}^3)^{L_t}\sim\mu_{(b,a)}^3$ hence, the {\it reducibility}.

Finally,
it is sufficient to find
$t\!\in \pm{\rm O}(3)\setminus\{e\}$ 
 such that 
\begin{equation}
\label{t-1}
t-1=\left(\begin{smallmatrix}
\lambda_1C_1 &\lambda_1C_2 &\lambda_1C_3 \\
\lambda_2C_1 &\lambda_2C_2 &\lambda_2C_3 \\
\lambda_3C_1 &\lambda_3C_2 &\lambda_3C_3
\end{smallmatrix}
\right),
\end{equation}
where $(C_1,C_2,C_3)=(1,1,-1)$, in  part (a), or for an arbitrary $(C_1,C_2,C_3)\in {\mathbb R}^3\setminus\{0\}$ in the part (b).
Such an element exists by Lemma~\ref{l.DtD=C.m=3} below. 
For such an element  $t$ we get 
respectively in the cases (a), (b) and  Example~\ref{ex.(111)b=1,a-any} 
(see \eqref{Sigma_3}):
\begin{eqnarray}
\nonumber
 &&
 \Sigma_2(t^{-1})=\sum_{n\in{\mathbb Z}}(b_{1n}\lambda_1^2+b_{2n}\lambda_2^2+b_{3n}\lambda_3^2)\big(a_{1n}+a_{2n}-a_{3n}\big)^2=0,\\
\nonumber
 &&
\Sigma_2(t^{-1})
=\sum_{n\in{\mathbb Z}}(b_{1n}\lambda_1^2+b_{2n}\lambda_2^2+b_{3n}\lambda_3^2)\big(C_1a_{1n}+C_2a_{2n}+C_3a_{3n}\big)^2<\infty,\\
\label{(111)(111)}
&&
\Sigma_2(t^{-1})\!=\!\sum_{n\in{\mathbb Z}}(b_{1n}\lambda_1^2\!+\!b_{2n}\lambda_2^2\!+\!b_{3n}\lambda_3^2)\big(C_1a_{1n}\!+\!C_2a_{2n}\!+\!C_3a_{3n}\big)^2\!=\!\infty.
\end{eqnarray}
{\it Note that the measure in Example~\ref{ex.(111)b=1,a-alt} does not satisfy the  orthogonality conditions}.

(2) {\it Irreducibility}.  In Example~\ref{ex.(111)b=1,a-any} we can not approximate  $x_{rn}$ by  Lemmas~\ref{x1x1.3}--\ref{x3x3.3}, since all the expressions
$$
\Delta(Y_1^{(1)},Y_2^{(1)},Y_3^{(1)}),\quad \Delta(Y_2^{(2)},Y_3^{(2)},Y_1^{(2)}),\quad \Delta(Y_3^{(3)},Y_1^{(3)},Y_2^{(3)})
$$
are bounded.
To approximate $D_{rn}$ using  Lemmas~\ref{d1.3}--\ref{d3.3},
we should estimate the following expressions:
 $$
 \Delta(Y_1,Y_2,Y_3),\quad \Delta(Y_2,Y_3,Y_1),\quad \Delta(Y_3,Y_1,Y_2).
 $$ 
By Lemma~\ref{l.min=proj.3},
 all these expressions  are infinite, if for all $(C_1,C_2,C_3)\in 
 {\mathbb R}^3\setminus\{0\}$  holds
\begin{equation*}
\nu( C_1, C_2, C_3):=\Vert C_1Y_1+C_2Y_2+C_3Y_3\Vert^2\!=\!\sum_{n\in \mathbb Z}\!\frac{\Big(\!C_1a_{1n}\!+\!C_2a_{2n}\!+\!C_3a_{3n}\!\Big)^2}
{\frac{1}{2b_{1n}}+\frac{1}{2b_{2n}}+\frac{1}{2b_{3n}}}=\infty.
\end{equation*}
In  Examples~\ref{ex.(111)b=1,a-any} we have 
\begin{eqnarray*}
&&
\nu( C_1, C_2, C_3)
\!=\!\Vert C_1Y_1+ C_2Y_2+C_3Y_3\Vert^2\!\sim\!
\sum_{k\in \mathbb Z}b_{1n}(C_1a_{1k}+ C_2a_{2k}+C_3a_{3k})^2\\
&&
\sim\sum_{n\in{\mathbb Z}}(b_{1n}\lambda_1^2+b_{2n}\lambda_2^2+b_{3n}\lambda_3^2)\big(C_1a_{1n}+C_2a_{2n}+C_3a_{3n}\big)^2\!=\!\Sigma_2(t^{-1})=\infty.\,\,\,\,\,\Box
\end{eqnarray*}
\end{pf}
\begin{lem}
\label{l.DtD=C.m=3}
For an arbitrary $(C_1,C_2,C_3)\in {\mathbb R}^3\setminus\{0\}$,  and  an arbitrary $D_3(s)={\rm diag}(s_1,s_2,s_3)$ with $(s_1,s_2,s_3)\in ({\mathbb R_+})^3$, 
there exists a unique  element $t\!\in \!\pm {\rm O}(3)\setminus\{e\}$  and $(\lambda_1,\lambda_2,\lambda_3)\in {\mathbb R}^3\setminus\{0\}$
such that 
\begin{equation}
\label{D_3tD^{-1}_3=(1,1,1)}
D_3(s)tD^{-1}_3(s)-{\rm I}=\left(\begin{smallmatrix}
\lambda_1C_1 &\lambda_1C_2 &\lambda_1C_3 \\
\lambda_2C_1 &\lambda_2C_2 &\lambda_2C_3 \\
\lambda_3C_1 &\lambda_3C_2 &\lambda_3C_3
\end{smallmatrix}
\right)\!=\!
{\small
\left(\begin{smallmatrix}
\lambda_1&0&0 \\
0&\lambda_2 &0\\ 0&0&\lambda_3
\end{smallmatrix}
\right)
\left(\begin{smallmatrix}
1&1&1 \\
1&1&1\\ 
1&1&1
\end{smallmatrix}
\right)
\left(\begin{smallmatrix}
C_1&0&0 \\
0&C_2 &0\\ 
0&0&C_3
\end{smallmatrix}
\right).
}
\end{equation}
\end{lem}
\begin{pf}
By \eqref{D_3tD^{-1}_3=(1,1,1)} 
we get
\begin{eqnarray}
\label{t(s)}
&& 
\left(
\begin{smallmatrix}
e_1\\
e_2\\
e_3
\end{smallmatrix}
\right):=
t=\left(
\begin{smallmatrix}
t_{11}&t_{12}&t_{13}\\
t_{21}&t_{22}&t_{23}\\
t_{31}&t_{32}&t_{33}
\end{smallmatrix}
\right)=
\left(\begin{smallmatrix}
C_1\lambda_1 +1&\frac{s_2}{s_1}C_2\lambda_1&\frac{s_3}{s_1}C_3\lambda_1\\
\frac{s_1}{s_2}C_1\lambda_2&C_2\lambda_2 +1&\frac{s_3}{s_2}C_3\lambda_2\\
\frac{s_1}{s_3}C_1\lambda_3&\frac{s_2}{s_3}C_2\lambda_3&C_3\lambda_3 +1
\end{smallmatrix}
\right),\\
\label{e(k)}
&& \text{where}\quad \Vert e_k\Vert^2=1\quad \text{and}\quad e_k\perp e_r,\quad 1\leq k<r\leq 3.
\end{eqnarray}

By  \eqref{t(s)} and the first relations in \eqref{e(k)} we get
\begin{equation}
 \label{lambda.m=3}
\lambda_k=-\frac{2s^2_kC_k}{s^2_1C^2_1+s^2_2C^2_2+s^2_3C^2_3},\quad 1\leq k\leq 3.
\end{equation}
Then the matrix elements $t=(t_{kr})_{k,r=1}^3$ are defined by \eqref{t(s)}.  To verify $e_k\perp e_r$ we need to show that
\begin{eqnarray*}
&& (e_1,e_2)=\frac{(s^2_1C^2_1+s^2_2C^2_2+s^2_3C^2_3)\lambda_1\lambda_2}{s_1s_2}+\frac{s^2_1C_1\lambda_2+s^2_2C_2\lambda_1}{s_1s_2}=0,\\
&& (e_1,e_3)=\frac{(s^2_1C^2_1+s^2_2C^2_2+s^2_3C^2_3)\lambda_1\lambda_3}{s_1s_3}+\frac{s^2_1C_1\lambda_3+s^2_3C_3\lambda_1}{s_1s_3}=0,\\
&& (e_2,e_3)=\frac{(s^2_1C^2_1+s^2_2C^2_2+s^2_3C^2_3)\lambda_2\lambda_3}{s_2s_3}+\frac{s^2_2C_2\lambda_3+s^2_3C_2\lambda_3}{s_2s_3}=0.
\end{eqnarray*}
Indeed, for example, for $(e_1,e_2)$ we have
\begin{eqnarray*}
&& (e_1,e_2)=\frac{(s^2_1C^2_1+s^2_2C^2_2+s^2_3C^2_3)\lambda_1\lambda_2}{s_1s_2}+\frac{s^2_1C_1\lambda_2+s^2_2C_2\lambda_1}{s_1s_2}=\\
&&\frac{1}{s_1s_2(s^2_1C_1^2+s^2_2C_2^2+s^2_3C_3^2)}\Big(4s^2_1s^2_2-(2s^2_1s^2_2+2s^2_1s^2_2)\Big)C_1C_2=0.
\end{eqnarray*}
The proofs of $e_1\perp e_3$ and $e_2\perp e_3$ are similar. 
\qed\end{pf}
Similarly, for  any $m\geq 2$ we can  prove the following lemma:
\begin{lem}
\label{l.DtD=C.m}
For an arbitrary $(C_k)_{k=1}^m\in {\mathbb R}^m\setminus\{0\}$,  and  $D_m(s)={\rm diag}(s_k)_{k=1}^m$ with 
$s_k\in \mathbb R_+,\,\,1\leq k\leq m$ 
there exists a unique  element $t\!\in \!\pm {\rm O}(m)\setminus\{e\}$  and $(\lambda_k)_{k=1}^m\in {\mathbb R}^m\setminus\{0\}$
such that 
\begin{equation}
\label{D_mtD^{-1}_m=(1,1,1)}
D_m(s)tD^{-1}_m(s)-{\rm I}=\left(\begin{smallmatrix}
\lambda_1C_1 &\lambda_1C_2 &\dots&\lambda_1C_m \\
\lambda_2C_1 &\lambda_2C_2 &\dots&\lambda_2C_m \\
&&\dots&&\\
\lambda_mC_1 &\lambda_mC_2 &\dots&\lambda_mC_m
\end{smallmatrix}
\right).
\end{equation}
The formulas for the corresponding $\lambda_k$ are as follows:
\begin{equation}
 \label{lambda.m}
\lambda_k=-\frac{2s^2_kC_k}{\sum_{r=1}^ms^2_rC^2_r},\quad 1\leq k\leq m.
\end{equation}
\end{lem}
\begin{lem}
\label{l.DtD=C.2}
For an arbitrary $(C_1,C_2)\in {\mathbb R}^2\setminus\{0\}$ and arbitrary $D_2(s)={\rm diag}(s_1,s_2)$ with $s_1,s_2\not=0$
there exists a unique  element $t\!\in \!\pm {\rm O}(2)\setminus\{e\}$  and $(\lambda_1,\lambda_2)\in {\mathbb R}^2\setminus\{0\}$
such that 
\begin{equation}
\label{D_2tD^{-1}_2=...}
D_2(s)tD^{-1}_2(s)-{\rm I}=\left(\begin{smallmatrix}
\lambda_1C_1 &\lambda_1C_2 \\
\lambda_2C_1 &\lambda_2C_2 
\end{smallmatrix}
\right).
\end{equation}
The formulas for the corresponding $\lambda_k$ are as follows:
\begin{equation}
 \label{lambda.2}
 \lambda_1=-\frac{2s_1^2C_1}{s_1^2C^2_1+s_2^2C^2_2},\quad \lambda_2=-\frac{2s_2C_2}{s_1^2C^2_1+s_2^2C^2_2}.
\end{equation}
In particular, we have
\begin{equation*}
 \left(\begin{smallmatrix}
 t_{11}&  \frac{s_1}{s_2}t_{12}\\  
\frac{s_2}{s_1}t_{21}&  t_{22}  
        \end{smallmatrix}\right)=
\left(\begin{smallmatrix}
C_1\lambda_1+1& 
C_2\lambda_1\\  
 C_1\lambda_2&C_2\lambda_2+1  
        \end{smallmatrix}\right)=
\left(\begin{smallmatrix} 
 \frac{s_2^2C^2_2-s^2C^2_1}{s_1^2C^2_1+s_2^2C^2_2}&  -\frac{2s_1^2C_1C_2}
 {s_1^2C^2_1+s_2^2C^2_2} \\    
 -\frac{2s_2^2C_1C_2}{s_1^2C^2_1+s_2^2C^2_2} & \frac{s_1^2C^2_1-s_2^2C^2_2}{s_1^2C^2_1+s_2^2C^2_2}   
 \end{smallmatrix}\right)=\tau_{12}(\phi,s_1,s_2).      
\end{equation*}
We can verify that 
$$\tau_{12}(\phi,s,s^{-1})=\tau_-(\phi,s)$$
where $\tau_-(\phi,s)$ is  defined by \eqref{m=2.crit-orthog.3}.
We just set $$\cos\phi=\frac{s^{-2}C^2_2-s^2C^2_1}{s^2C^2_1+s^{-2}C^2_2}\quad\text{and}\quad 
\sin\phi=-\frac{2C_1C_2}{s^2C^2_1+s^{-2}C^2_2}.$$
In addition ${\rm det}\,t=-1$.
\end{lem}
\subsubsection{Case $\Sigma_{123}(s)=(0,1,1)$}
We have for some $s_1\in  \mathbb R_+$ and all 
$(s_2,s_3)\in  \big(\mathbb R_+\big)^2$
$$,\quad \Sigma_{23}(s_2)=\infty,\quad\Sigma_{13}(s_3)=\infty.
$$

\begin{rem}
\label{r.Sigma-123=(011)} 
 Since $\Sigma_{12}(s_1)\!<\!\infty$, by \eqref{(b,a)-(b',a')} and  Lemma~\ref{l.eq-prod}, we can suppose that 
\begin{eqnarray*}
&&
 b\!=\!(b_{1n},b_{2n},b_{3n})_{n\in \mathbb Z}\!=\!(1,s_1^4(1\!+\!c_n),b_{3n})_{n\in \mathbb Z},\quad
\sum_nc_n^2<\infty,\\
&&
\text{therefore, we can take}\quad
b=(1,1,b_{3n})_{n\in \mathbb Z},\,\,s=1,\,\,c_n\equiv 0.
\end{eqnarray*}
\end{rem}
Since $\Sigma_{13}(s)=
\sum_{n\in \mathbb Z}\Big(
\frac{s^2}{\sqrt{b_{3n}}}-
\frac{\sqrt{b_{3n}}}{s^2}
\Big)^2=\infty$,
we have as in \eqref{d_n-gen} three cases:
\begin{equation}
\label{b_(3n)-gen}
\lim_nb_{3n}=\!\left\{\begin{array}{clcc}
(a)& \infty&&\\
(b)&b>0&\text{with}&
\sum_nb_n^2=\infty,\\
(c)&0&&
\end{array}\right.
\end{equation}
where $b_{3n}=b(1+b_n)$ with $\lim_nb_n=0$ in the case (b).
Note that condition $S_3(3)=\infty$, implies  $\sum_nb_{3n}^2=\infty$. 
Indeed, by \eqref{S_r(3).1} we have for $1\leq r\leq 3$
\begin{eqnarray}
\nonumber
&& 
S_r(3)=\sum_{n\in \mathbb Z}\frac{b_{rn}^2}{b_{1n}b_{2n}+b_{1n}b_{3n}+b_{2n}b_{3n}},\quad
S_1(3)=\sum_{n\in \mathbb Z}\frac{1}{1+2b_{3n}}=\infty,\\
\label{S_3(3)=}
&&
S_2(3)\!=\!\sum_{n\in \mathbb Z}\frac{1}{1\!+\!2b_{3n}}\!=\!\infty,\,\,\,\, 
\infty\!=\!S_3(3)\!=\!\sum_{n\in \mathbb Z}\frac{b_{3n}^2}{1\!+\!2b_{3n}}\!\stackrel{\eqref{two-ser-alpha}}{\sim} \!\sum_{n\in \mathbb Z}b_{3n}^2.\qquad\,\, 
\end{eqnarray}
By \eqref{Y_r^{(s)}} we have
\begin{eqnarray*}
&& \Vert Y_r^{(r)}\Vert^2 =\sum_{k\in \mathbb Z}\frac{b_{rk}^2}{b_{rk}^2+2(b_{1n}b_{2n}+b_{1n}b_{3n}+b_{2n}b_{3n})},\\
&& \Vert Y_r^{(s)}\Vert^2 =\sum_{k\in \mathbb Z}\frac{b_{rk}^2}{b_{sk}^2+2(b_{1n}b_{2n}+b_{1n}b_{3n}+b_{2n}b_{3n})},\,\,s\not=r.
\end{eqnarray*}
Let us denote
\begin{equation}
\label{Y^r_s(n)}
\left(\begin{array}{ccc}
Y_{1n}^{(1)}&Y_{2n}^{(1)}&Y_{3n}^{(1)}\\
Y_{1n}^{(2)}&Y_{2n}^{(2)}&Y_{3n}^{(2)}\\
Y_{1n}^{(3)}&Y_{2n}^{(3)}&Y_{3n}^{(3)}\\
 \end{array}\right)=
 \left(
\begin{array}{ccc}
\frac{1}{\sqrt{3+4b_{3n}}}&
\frac{1}{\sqrt{3+4b_{3n}}}&
\frac{b_{3n}}{\sqrt{3+4b_{3n}}}\\
\frac{1}{\sqrt{3+4b_{3n}}}&
\frac{1}{\sqrt{3+4b_{3n}}}&
\frac{b_{3n}}{\sqrt{3+4b_{3n}}}\\
\frac{1}{\sqrt{b^2_{3n}+4b_{3n}+2}}&
\frac{1}{\sqrt{b^2_{3n}+4b_{3n}+2}}&
\frac{b_{3n}}{\sqrt{b^2_{3n}+4b_{3n}+2}}
\end{array}
\right)
\end{equation}
We have $\Delta(Y_1^{(1)}, Y_2^{(1)},Y_3^{(1)})\!=\!\Delta(Y_2^{(2)}, Y_3^{(2)},Y_1^{(2)})\!<\!\infty$, indeed, since $Y_1^{(2)}\!=\!Y_2^{(2)}$ we
get for example
\begin{equation*}
 \Delta(Y_2^{(2)}, Y_3^{(2)},Y_1^{(2)})=
 \frac{\Gamma(Y_2^{(2)})+
 \Gamma(Y_2^{(2)},Y_3^{(2)})
 }{1+\Gamma(Y_3^{(2)})+\Gamma(Y_1^{(2)})
+\Gamma(Y_3^{(2)},Y_1^{(2)})
 }<1<\infty.
\end{equation*}
 \begin{lem}
\label{l.(111)-(011)}
In the cases (a), (b) and (c) given by \eqref{b_(3n)-gen} we have 
\begin{equation}
\label{(111)-(011)}
\Delta(Y_3^{(3)}, Y_1^{(3)},Y_2^{(3)})=\infty.
\end{equation}
\end{lem}
\begin{pf}
In all these cases we have $Y_1^{(3)}=Y_2^{(3)}$ hence,
$\Gamma(Y_3^{(3)},Y_1^{(3)},Y_2^{(3)})=0$ and 
$\Gamma(Y_1^{(3)},Y_2^{(3)})=0$. Therefore, 
by \eqref{Delta(f,g,h)}
\begin{eqnarray}
\nonumber
 &&
 \Delta(Y_3^{(3)},Y_1^{(3)},Y_2^{(3)})=
 \frac{\Gamma(Y_3^{(3)})+
 \Gamma(Y_3^{(3)},Y_1^{(3)})+
 \Gamma(Y_3^{(3)},Y_2^{(3)})
 }{1+\Gamma(Y_1^{(3)})+\Gamma(Y_2^{(3)})} \\
 \label{(011)(a-c)}
 &&
 =\frac{\Gamma(Y_3^{(3)})+2 \Gamma(Y_3^{(3)},Y_1^{(3)})}{1+2\Gamma(Y_1^{(3)})}
 \sim  
 \Delta(Y_3^{(3)},Y_1^{(3)}).
\end{eqnarray}
We have two cases:\\
(a.1) when $\Vert Y_1^{(3)}\Vert <\infty$, and (a.2) when 
$\Vert Y_1^{(3)}\Vert =\infty$.\\
In the case (a.1) we have 
$\Delta(Y_3^{(3)},Y_1^{(3)})\sim \Gamma(Y_3^{(3)})=\infty$.
Therefore,  \eqref{(111)-(011)} holds.
In the case (a.2) we should verify that
\begin{equation}
 \label{(111)(c)}
\Vert C_1Y_1^{(3)}+C_3Y_3^{(3)}\Vert^2=\infty\quad\text{ for all} \quad
(C_1,C_3)\in {\mathbb R}^2\setminus\{0\}.
\end{equation}
Then this will imply \eqref{(111)-(011)}. We have
\begin{equation*}
\Vert C_1Y_1^{(3)}+C_3Y_3^{(3)}\Vert^2=
\sum_{n\in \mathbb Z}\frac{\big(C_1+C_3b_{3n}\big)^2}{b^2_{3n}+4b_{3n}+2}=:\sum_{n\in \mathbb Z}g_n.
\end{equation*}
If $C_1=0$ or $C_3=0$ the later expression is divergent since $Y_1^{(3)}=Y_3^{(3)}=\infty$. Let $C_1C_3\!\not=\!0$. In this case $\lim_ng_n=C_3^2>0$ since $\lim_nb_{3n}=\infty$, case (a). Therefore, $\sum_{n\in \mathbb Z}g_n=\infty$. 
By Lemma~\ref{l.min=proj}  this implies 
$\Delta(Y_3^{(3)},Y_1^{(3)})=\infty$ therefore,
\eqref{(111)-(011)}.
In the case (b) we have by \eqref{(011)(a-c)}

\begin{equation*}
 \Delta(Y_3^{(3)},Y_1^{(3)},Y_2^{(3)})=
 \Delta(Y_3^{(3)},Y_1^{(3)}).
\end{equation*}
To prove that $ \Delta(Y_3^{(3)},Y_1^{(3)})=\infty$ using Lemma~
\ref{l.min=proj} we should verify \eqref{(111)(c)}.
We have  $\Vert Y_3^{(3)}\Vert^2=\infty$ since $S=(0,1,1)$. By \eqref{Y^r_s(n)} 
\begin{equation*}
\Vert Y_1^{(3)}\Vert^2\!=\!
\sum_{n\in \mathbb Z}\frac{1}{b^2_{3n}+4b_{3n}+2}\!\sim\!
\sum_{n\in \mathbb Z}\frac{1}{b^2+4b+2}\!=
\!\infty.
\end{equation*}
The expression $\Vert C_1Y_1^{(3)}+C_3Y_3^{(3)}\Vert^2$ can be finite only for  $(C_1,C_3)=\lambda(b,-1)$. Take 
$\lambda=1$, we get in the case (b)
\begin{eqnarray*}
 &&
 \Vert C_1Y_1^{(3)}+C_3Y_3^{(3)}\Vert^2=\sum_{n\in \mathbb Z}\frac{(b-b_{3n})^2}{b^2_{3n}+4b_{3n}+2}=
\sum_{n\in \mathbb Z}\frac{b^2b_n^2}{
b^2(1+b_n)^2+4b(1+b_n)+2
}
 \\
&& 
\stackrel{ \eqref{two-sim-ser}}{\sim} \sum_{n\in \mathbb Z}\frac{b_n^2}{
(4b+2b^2)b_{n}+b^2+4b+2
}\stackrel{ \eqref{two-ser-alpha}}{\sim}
\sum_{n\in \mathbb Z}b_n^2=\infty.
\end{eqnarray*}
In the case (c),  
we have by \eqref{(011)(a-c)} 
\begin{equation*} 
\Delta(Y_3^{(3)},Y_1^{(3)},Y_2^{(3)})
\sim
\Delta(Y_3^{(3)},Y_1^{(3)}).
\end{equation*}
To prove that $ \Delta(Y_3^{(3)},Y_1^{(3)})=\infty$ using Lemma~
\ref{l.min=proj} we should verify \eqref{(111)(c)}.
Again, we have  $\Vert Y_3^{(3)}\Vert^2=\infty$ since $S=(0,1,1)$.  Because of $\lim_nb_{3n}=0$, we have  by \eqref{Y^r_s(n)}  
\begin{equation*}
\Vert Y_1^{(3)} \Vert^2=\sum_{n\in \mathbb Z}\frac{1}{b^2_{3n}+4b_{3n}+2}\sim \sum_{n\in \mathbb Z}\frac{1}{2}
=\infty,
\end{equation*}
Let $C_1C_3\not=0$, then 
since $\lim_nb_{3n}=0$ we get
\begin{eqnarray*}
 &&
\Vert C_1 Y_1^{(3)}\!+\!C_3 Y_3^{(3)}\Vert^2\!=\!
 \sum_{n\in \mathbb Z}\frac{\big(C_1\!+\!C_3b_{3n}\big)^2}{b^2_{3n}\!+\!4b_{3n}\!+\!2}\!=\!
 \sum_{n\in \mathbb Z}\frac{C_3^2\big(b_{3n}\!+\!C_1C_3^{-1}\big)^2}{b^2_{3n}\!+\!4b_{3n}\!+\!2}=\infty. 
\hskip 0.3cm\Box
\end{eqnarray*}
\end{pf}
%

By  Lemma~\ref{l.(111)-(011)} we can approximate $x_{3n}$.
By \eqref{II}  we have
\begin{eqnarray*}
 &&\Vert Y_1\Vert^2=\sum_{n\in \mathbb Z}\frac{a_{1n}^2}{\frac{1}{2b_{1n}}+\frac{1}{2b_{2n}}+\frac{1}{2b_{3n}}}= 
\sum_{k\in \mathbb Z}\frac{a_{1n}^2}{1
+\frac{1}{2b_{3n}}
}=\sum_{k\in \mathbb Z}\frac{2b_{3n}a_{1n}^2}{1+2b_{3n}},
\\
&&
\Vert Y_2 \Vert^2
\!=\!\sum_{n\in \mathbb Z}\frac{a_{2n}^2}{
1+\frac{1}{2b_{3n}}}\!=\!\sum_{k\in \mathbb Z}\frac{2b_{3n}a_{2n}^2}{1+2b_{3n}},
\,\,\,\,
\Vert Y_3\Vert^2
\!=\! \sum_{n\in \mathbb Z}\frac{a_{3n}^2}{
1+
\frac{1}{2b_{3n}}} \!=\!
\sum_{k\in \mathbb Z}\frac{2b_{3n}a_{3n}^2}{1+2b_{3n}}.
\end{eqnarray*}
Therefore, in the case (a) and (b) we have
$$
\Vert Y_1 \Vert^2\sim \sum_{k\in \mathbb Z}a_{1n}^2,
\quad 
\Vert Y_2 \Vert^2\sim \sum_{k\in \mathbb Z}a_{2n}^2,
\quad 
\Vert Y_3 \Vert^2\sim \sum_{k\in \mathbb Z}a_{3n}^2,
$$
In the case (c) we get
$$
\Vert Y_1 \Vert^2\sim \sum_{k\in \mathbb Z}b_{3n}a_{1n}^2,
\quad 
\Vert Y_2 \Vert^2\sim \sum_{k\in \mathbb Z}b_{3n}a_{2n}^2,
\quad 
\Vert Y_3 \Vert^2\sim
\sum_{k\in \mathbb Z}b_{3n}a_{3n}^2.
$$
Since in the case (a), (b)
\begin{eqnarray*}
 &&
\Vert Y_1 \Vert^2\sim
\sum_{n\in \mathbb Z}a_{1n}^2=
\sum_{n\in \mathbb Z}b_{1n}a_{1n}^2\sim S^L_{11}(\mu)=\infty,\\
&&
\Vert Y_2 \Vert^2\sim
\sum_{n\in \mathbb Z}a_{2n}^2=\sum_{n\in \mathbb Z}b_{2n}a_{2n}^2=S^L_{22}(\mu)=\infty,
\end{eqnarray*}
we have two possibilities for $y_{23}:=(y_2,y_3)\in \{0,1\}^2$, 
see Section~\ref{s.444}:
\begin{equation}
\begin{array}{ccccc}
   &(1.1)&(1.3)\\
y_1&1&1\\
y_2&1&1\\
y_3&0&1\\
\end{array}
\end{equation}
In the case (c) we have 
$$
\Vert Y_3 \Vert^2\sim
\sum_{n\in \mathbb Z}b_{3n}a_{3n}^2\sim S^L_{33}(\mu)=\infty.
$$ 
Therefore,
we have four possibilities 
for $y_{12}:=(y_1,y_2)\in \{0,1\}^2$, see \eqref{(y_1,y_2,y_3)(011)},
\begin{equation}
\begin{array}{ccccc}
   &(1.0)&(1.1)&(1.2)&(1.3)\\
y_1&0&1&0&1\\
y_2&0&0&1&1\\
y_3&1&1&1&1\\
\end{array}
\end{equation}
Further, in the case (a), (b) we have four possibilities: (1.1.1), (1.3.1) and (1.1.0), (1.3.0), see Remark~\ref{r.(011)-(1-3)(0-1)}. In the case
(1.1.1) we can approximate 
$D_{1n}, \,\,D_{2n}$, in the case  (1.3.1) we can approximate all $D_{rn},\,1\leq r\leq 3$. In these cases the proof is finished, since we get 
respectively $D_{1n}, D_{2n},x_{3n}\,\eta\, {\mathfrak A}^3$.
The cases (a),(b) subcases 
(1.1.0), (1.3.0)  we consider below.

In the case (c) subcase (1.0)  we can approximate $D_{3n}$ using 
Lemma~\ref{d3.3}, since 
$\Delta(Y_3,Y_2,Y_1)\sim \Vert Y_3\Vert^2=\infty$, 
so we have $D_{3n},x_{3n}\,\eta\, {\mathfrak A}^3$,
and the proof is finished.

Further, in the case (c) we have six cases (1.1.1), (1.2.1), (1.3.1) and  
(1.1.0), (1.2.0), (1.3.0), 
according to whether 
corresponding expressions are divergent
(see analogue in Remark~\ref{r.(011)-(1-3)(0-1)}).
We can approximate in the three first cases by respectively
$D_{1n}$ and $ D_{3n}$ in 
the case (1.1.1),
$D_{2n}$ and $D_{3n}$ in
the case (1.1.2)
and all $D_{1n},\,D_{2n},\,D_{3n}$ in (1.1.3).
The proof of irreducibilty is finished in these cases
because we have respectively 
$D_{1n},\,D_{3n},\,x_{3n}\,\eta\,{\mathfrak A}^3$,
$D_{2n},\,D_{3n},\,x_{3n}\,\eta\,{\mathfrak A}^3$,
or $D_{1n},\,D_{2n},\,D_{3n},\,x_{3n}\,\eta\, {\mathfrak A}^3$.

If the opposite holds, in the cases (a), (b) or 
(c), i.e., we are in the cases (1.1.0), (1.2.0) and (1.3.0) respectively, we try to approximate $D_{3n}$  using Lemma~\ref{l.Re-Im-exp.3}. 
If one of the expressions  $\Sigma_3(D,s)$ or $\Sigma^{\vee}_3(D,s)$ is divergent, we can approximate $D_{3k}$ and the proof is finished,  since we have  $x_{3n},\,D_{3n}\,\,\eta\,\, {\mathfrak A}^3$. Let us suppose,  as in Remark~\ref{r.(001)},
that for every sequence $s=(s_k)_{k\in \mathbb Z}$ holds 
$$
\Sigma_3(D,s)+\Sigma^{\vee}_3(D,s)<\infty.
$$
Then, in particularly, we have for $s^{(3)}=(s_{k})_{k\in \mathbb Z}$
with $\frac{s^2_{k}}{b_{3k}}\equiv 1$
\begin{eqnarray}
 \nonumber
&& 
 \infty\!>\! \Sigma_3(D,s^{(3)})\!+\!\Sigma^{\vee}_3(D,s^{(3)})\!\sim\! \Sigma_3(D)+\Sigma^{\vee}_3(D)\!=\!\sum_k\frac{\frac{1}{2b_{3k}}+a^2_{3k}}{C_k\!+\!a^2_{1k}\!+\!a^2_{2k}\!+\!a^2_{3k}}\\
\label{Sigma_3(D)<} 
&&
\stackrel{\eqref{two-sim-ser}}{\sim}
\sum_k\frac{
\frac{1}{2b_{3k}}+a^2_{3k}}{\frac{1}{2b_{1k}}+a^2_{1k}+
\frac{1}{2b_{2k}}+a^2_{2k}}=
\sum_k\frac{\frac{1}{2b_{3k}}+a^2_{3k}}{1
+a^2_{1k}
+a^2_{2k}}=:\Sigma^{\vee,+}_3(D).
\end{eqnarray}
In the case (a), (b) and (c) we have respectively  $$\Sigma^{\vee,+}_3(D)\sim \Sigma^+_3(D)=\sum_k\frac{2a^2_{3k}}
{1+2a^2_{1k}+2a^2_{2k}},\quad 
\Sigma^{\vee,+}_3(D)=
\sum_k\frac{\frac{1}{2b_{3k}}+a^2_{3k}}{1
+a^2_{1k}
+a^2_{2k}}.
$$
In particular, in the case (c) we have by \eqref{Sigma_3(D)<} 
\begin{equation}
\label{Sigma-3(D)}
\infty>
\sum_k\frac{\frac{1}{2b_{3k}}+a^2_{3k}}{1
+a^2_{1k}
+a^2_{2k}}>
\sum_k\frac{a^2_{3k}}{1
+a^2_{1k}
+a^2_{2k}}\sim\Sigma^+_3(D).
\end{equation}
 The cases (a), subcase 
(1.1.0),  where $\Vert Y_3\Vert^2<\infty$ can not occur, because conditions
$\Sigma_{12}(s_1)<\infty$ and 
$\nu_{12}(C_1,C_2)<\infty$ defined by \eqref{nu(12)}, contradict the orthogonality condition for the matrix $\tau_{12}(\phi,s)$:
\begin{equation}
\tau_{12}(\phi,s)= \left(
\begin{array}{ccc}
\cos\phi&s^2\sin\phi&0\\
s^{-2}\sin\phi&-\cos\phi&0\\
0&0&1
\end{array}
\right), 
\end{equation}
Indeed, recall Remark~\ref{perp2-1} (instead of $\mu_{(b,a)}^2$ we can write $\mu_{(b,a)}^3$)
$$
\big(\mu_{(b,a)}^3\big)^{L_{\tau_{12}(\phi,s)}}\perp\mu_{(b,a)}^3
\Leftrightarrow 
\Sigma_{12}(s)+\Sigma_{12}(C_1,C_2)=\infty,
$$
where $\Sigma_{12}(s)=\Sigma_1(s)$ is defined by (\ref{sigma1(s)}) and
$\Sigma_{12}(C_1,C_2)=\Sigma_{2}(C_1,C_2)$ is defined by \eqref{sigma_2(C,C)} 
\begin{equation*} 
\Sigma_{12}(C_1,C_2):=
\sum_{n\in {\mathbb Z}}(C_1^2b_{1n}+C_2^2b_{2n})(C_1a_{1n}+C_2a_{2n})^2\sim\nu_{12}(C_1,C_2).
\end{equation*}
We get contradiction:
$$
\infty>\Sigma_{12}(s)+\nu_{12}(C_1,C_2)\sim
\Sigma_{12}(s)+\Sigma_{12}(C_1,C_2)=\infty.
$$
In the case (a), (b), subcase (1.3.0) we get
$\Sigma^+_3(D)=\infty$ by Lemma~\ref{l.(011)-fin.3},
contradiction with \eqref{Sigma_3(D)<}  hence, $D_{3n}\,\eta\, {\mathfrak A}^3$.

In the case (c), subcases (1.1.0) and (1.2.0) we have respectively $\Vert Y_2\Vert^2<\infty$ and $\Vert Y_1\Vert^2<\infty$ hence,
$$
\Sigma^+_3(D)\sim  \sum_k\frac{a^2_{3k}}{1
+a^2_{1k}}=\infty,\quad
\Sigma^+_3(D)\sim  \sum_k\frac{a^2_{3k}}{1
+a^2_{2k}}=\infty.
$$
by Lemma~\ref{l.(011)-fin},
contradiction with \eqref{Sigma_3(D)<}  hence, $D_{3n}\,\eta\, {\mathfrak A}^3$. In the case (c), subcase (1.3.0) we  get
\begin{equation*}
 \Sigma^+_3(D)=\sum_k\frac{a^2_{3k}}{1
+a^2_{1k}
+a^2_{2k}}=\infty
\end{equation*}
by Lemma~\ref{l.(011)-fin.3},
contradiction with \eqref{Sigma_3(D)<}  hence, $D_{3n}\,\eta\, {\mathfrak A}^3$.

\subsubsection{Case $\Sigma_{123}(s)=(1,1,1)$}
We have  for all $s=(s_{12},s_{23},s_{13})\in \mathbb R_+^3\setminus \{0\}$
\begin{eqnarray}
\label{Sigma(123)=(111)}
&&
\Sigma_{12}(s_{12})=\infty,\quad \Sigma_{23}(s_{23})=\infty,\quad\Sigma_{13}(s_{13})=\infty,\\
\nonumber
&&
  b\!=\!(b_{1n},b_{2n},b_{3n})_{n\in \mathbb Z}\stackrel{ \eqref{(b,a)-(b',a')}}{=}\!(1,d_{2n},d_{3n})_{n\in \mathbb Z}.
\end{eqnarray}
Recall \eqref{alpha}, that we denote $ D_n:=d^{-1}_{2n}+d^{-1}_{3n}+1$ and $d_n=\frac{d_{3n}}{d_{3n}}$. Set

\begin{equation}
\label{Y^r_s(n)(d,d)}
\left(\begin{array}{ccc}
Y_{1n}^{(1)}&Y_{2n}^{(1)}&Y_{3n}^{(1)}\\
Y_{1n}^{(2)}&Y_{2n}^{(2)}&Y_{3n}^{(2)}\\
Y_{1n}^{(3)}&Y_{2n}^{(3)}&Y_{3n}^{(3)}\\
 \end{array}\right)=
 \left(
\begin{array}{ccc}
\frac{1}{\sqrt{1\!+\!2D_nd_{2n}d_{3n}}}&
\frac{d_{2n}}{\sqrt{1\!+\!2D_nd_{2n}d_{3n}}}&
\frac{d_{3n}}{\sqrt{1\!+\!2D_nd_{2n}d_{3n}}}\\
\frac{1}{\sqrt{d_{2n}^2\!+\!2D_nd_{2n}d_{3n}}}&
\frac{d_{2n}}{\sqrt{d_{2n}^2\!+\!2D_nd_{2n}d_{3n}}}&
\frac{d_{3n}}{\sqrt{d_{2n}^2\!+\!2D_nd_{2n}d_{3n}}}\\
\frac{1}{\sqrt{d_{3n}^2\!+\!2D_nd_{2n}d_{3n}}}&
\frac{d_{2n}}{\sqrt{d_{3n}^2\!+\!2D_nd_{2n}d_{3n}}}&
\frac{d_{3n}}{\sqrt{d_{3n}^2\!+\!2D_nd_{2n}d_{3n}}}
\end{array}
\right).
\end{equation}
\begin{rem}
 \label{r.Sigma-123=(111)} 
For $(r,s)$ such that $1\leq r<s\leq 3$ the following equivalence hold 
%
%
\begin{eqnarray}
\label{Sigma(rs)(s)<}
 &&
 \Sigma_{rs}(s_{rs})\!<\!\infty\,\,\Leftrightarrow\,\,\sum_{n\in \mathbb Z}c^2_{rs,n}\!<\!\infty\,\,\Leftrightarrow\,\, \sum_{n\in \mathbb Z}c^2_{sr,n}\!<\!\infty,\,\, \text{where}\\
\label{c(rs,n)}
&&
\frac{b_{rn}}{b_{sn}}=:s^{-4}_{rs}(1+c_{rs,n}),\,\,\frac{b_{sn}}{b_{rn}}\!=\!s^4_{rs}(1\!+\!c_{sr,n}),\,\, \lim_n\frac{b_{rn}}{b_{sn}}\in (0,\infty).
\end{eqnarray}
\end{rem}
\begin{pf} 
By Lemma~\ref{l.eq-prod} we have
\begin{eqnarray}
 \nonumber
&&
\Sigma_{rs}(s_{rs})=
\sum_{n\in \mathbb Z}\Big(s_{rs}^2\sqrt{\frac{b_{rn}}{b_{sn}}}\!-\!
s_{rs}^{-2}\sqrt{\frac{b_{sn}}{b_{rn}}}\Big)^2=
\sum_{n\in \mathbb Z}\frac{c^2_{rs,n}}{1+c_{rs,n}}\sim \sum_{n\in \mathbb Z}c^2_{rs,n},\\
 \nonumber
&& 
\Sigma_{sr}(s^{-1}_{rs})=
\sum_{n\in \mathbb Z}\Big(s_{rs}^{-2}\sqrt{\frac{b_{sn}}{b_{rn}}}\!-\!
s_{rs}^2\sqrt{\frac{b_{rn}}{b_{sn}}}\Big)^2=
\sum_{n\in \mathbb Z}\frac{c^2_{sr,n}}{1+c_{sr,n}}\sim \sum_{n\in \mathbb Z}c^2_{sr,n},
\\
\label{(1+c)(1+e)=1}
&&
\text{note also that}\quad  1=\frac{b_{rn}}{b_{sn}}\frac{b_{sn}}{b_{rn}}=(1\!+\!c_{rs,n})(1\!+\!c_{sr,n}). 
\end{eqnarray}
\qed\end{pf}
By Remark~\ref{r.Sigma-123=(111)},  the condition $\Sigma_{rs}(s_{rs})=\infty$ means the following:
\begin{equation}
\label{l(sr)}
l_{sr}:=\lim_n\frac{b_{sn}}{b_{rn}}
=\!
\left\{\begin{array}{lccc}
(a)& \infty&&\\
(b)&s^{4}_{rs}>0&\text{with}&\sum_{n\in \mathbb Z}c^2_{sr,n}=\infty.\\
(c)&0&\\
(d)&\lim& \,\,\text{does not exist}
\end{array}\right.
\end{equation}
\begin{rem}
\label{r.b-arrange} 
In the case (d) we can use the fact that some {\it subsequence} of $\Big(\frac{b_{sn}}{b_{rn}}\Big)_{n\in \mathbb Z}$ has property (a), (b) or (c).
We can avoid the case (c). Namely, if $l_{sr}=0$ for some pair $(r,s)$ with $1\leq r<s\leq 3$, we can exchange the two lines $(b_{sn},a_{sn})$ and $(b_{rn},a_{rn})$ to obtain $l_{sr}=\infty$.
\end{rem}
Formally, we have $3^3=\#(A)^{\#(B)}$ possibilities where $A=\{(21),(32),(31)\}$ and $B=\{(a),(b),(d)\}$. Since $l_{32}l_{21}=l_{31}$  we get only the following cases: 
$$
\begin{array}{ccccc}
 e \setminus (rs) &(21)&(32)&(31)\\
(1)&b&b&b\\
(2)&a&a&a\\
(3)&a&b&a\\
(4)&b&a&b\\
\end{array}.
$$
To be able to approximate $x_{rn}$ for $1\leq r \leq 3$ we should study when the following expressions are infinite:
\begin{equation}
\label{rho(CCC)}
\rho_r(C_1,C_2,C_3)=
\Vert C_1Y_1^{(r)}+C_2Y_2^{(r)}+C_2Y_3^{(r)}\Vert ^2.
\end{equation}
By \eqref{Y^r_s(n)(d,d)} we have
\begin{eqnarray}
\label{rho_r(CCC)}
 &&
\rho_r(C_1,C_2,C_3)=
  :\sum_n\frac{\vert C_1+C_2d_{2n}+C_3d_{3n}\vert^2}
 {C_{rn}},\\
 \nonumber
  &&
\text{where}\,\,C_{1n}= 1\!+\!2D_nd_{2n}d_{3n},\,\,
C_{2n}= d_{2n}^2\!+\!2D_nd_{2n}d_{3n},\,\,C_{3n}= d_{3n}^2\!+\!2D_nd_{2n}d_{3n}.
 \end{eqnarray}
Consider the {\bf case 
 (1)=(bbb)}.
We prove the analogue of  Lemma~\ref{l.park} for the case $m=3$.
\begin{lem}
\label{l.park2}
Let 
for all $s=(s_{12},s_{23},s_{13})\in \big(\mathbb R_+\big)^3$ holds
\eqref{Sigma(123)=(111)}.
Then 
\begin{equation}
\label{3-Delta}
 \Delta(Y_3^{(3)},Y_1^{(3)},Y_2^{(3)})=
 \Delta(Y_2^{(2)},Y_3^{(2)},Y_1^{(2)})=
 \Delta(Y_1^{(1)},Y_2^{(1)},Y_3^{(1)})=\infty.
\end{equation}
\end{lem}
\begin{pf} 
 For $1\leq r<s\leq 3$  set 
 \begin{eqnarray*}
\frac{b_{sn}}{b_{rn}}\!=\!s^4_{rs}(1\!+\!c_{sr,n})\quad\text{with}\quad \sum_{n\in\mathbb Z}c^2_{sr,n}  =\infty,
\quad \lim_{n\to\infty}c_{sr,n}=0.
\end{eqnarray*}
For $b_{1n}\equiv 1$ we have
\begin{eqnarray*}
&&
b_{2n}=s^4_{12}(1+c_{21,n}),\,\,
b_{3n}=s^4_{13}(1+c_{31,n}),\,\,
\\
&&
\frac{b_{3n}}{b_{2n}}=\frac{s^4_{13}}{s^4_{12}}\frac{1+c_{31,n}}{1+c_{21,n}}=s^4_{23}(1+c_{32,n}),
\quad c_{32,n}=\frac{1+c_{31,n}}{1+c_{21,n}}-1,\,\,
s_{23}=\frac{s_{13}}{s_{12}},\\
&&
\sum_nc^2_{32,n}\! =\!\sum_n\Big(\frac{1+c_{31,n}}{1+c_{21,n}}\!-\!1\Big)^2\!\!
=\!\sum_n\Big(\frac{c_{31,n}\!-\!c_{21,n}}{1+c_{21,n}}\Big)^2\!\sim\!
\sum_n(c_{21,n}\!-\!c_{31,n})^2\!=\!\infty.
\end{eqnarray*}
Finally, we get
\begin{equation}
 \label{a,b,c-not-l_2}
\sum_nc_{21,n}^2=\infty,\quad \sum_nc_{31,n}^2=\infty,\quad \sum_n(c_{21,n}-c_{31,n})^2=\infty.
\end{equation}
By \eqref{rho(CCC)} and \eqref{rho_r(CCC)} we get
\begin{eqnarray*}
&&
\rho_r(C_1,C_2,C_3)\!=\!
\Vert C_1Y_1^{(r)}\!+\!C_2Y_2^{(r)}\!+\!C_2Y_3^{(r)}\Vert ^2\!=\!
\sum_n\frac{\vert C_1+C_2d_{2n}+C_3d_{3n}\vert^2}
 {C_{rn}} \\
 &&
 =\sum_n\frac{\vert C_1+C_2s^4_{12}(1+c_{21,n})+C_3s^4_{13}(1+c_{31,n})\vert^2}
 {C_{rn}}.
\end{eqnarray*}
The latter expression is divergent if $C_1+C_2s^4_{12}+C_3s^4_{13}\not=0$ since 
$\lim_{n\to\infty}c_{21,n}=\lim_{n\to\infty}c_{31,n}=0$ and $A_1\leq C_{rn}\leq A_2$.

In the case  when $C_1+C_2s^4_{12}+C_3s^4_{13}=0$ we get
\begin{equation}
 \label{rho(CC)}
 \rho_r(C_1,C_2,C_3)=\sum_{n\in \mathbb Z}
 \frac{\vert C_2s^4_{12}c_{21,n}+C_3s^4_{13}c_{31,n}\vert^2}
 {C_{rn}}=:\rho_r(C_2,C_3).
\end{equation}
The latter expression is divergent by the first two relations in \eqref{a,b,c-not-l_2} when 1) $C_2C_3>0$, 2)
$C_2=0$ and $C_3\not=0$, 3)  $C_2\not=0$ and $C_3=0$. 

If $C_2C_3<0$ we have  by the last relation in \eqref{a,b,c-not-l_2}
\begin{eqnarray*}
&&
\sum_{n\in \mathbb Z}
 \frac{\vert C_2s^4_{12}c_{21,n}-C_3s^4_{13}c_{31,n}\vert^2}
 {C_{rn}}\sim 
 \sum_{n\in \mathbb Z}\vert C_2s^4_{12}c_{21,n}-C_3s^4_{13}c_{31,n}\vert^2=\infty,
\end{eqnarray*}
since $(s_{12},s_{13})=\frac{1}{s_1}(s_{2},s_{3})\in(\mathbb R^*)^2$ are arbitrary.

Consider the {\bf case (2)=(aaa)}.
Now, see \eqref{l(sr)}, we have 
\begin{equation}
\label{(aaa)} 
l_{21}\!=\!\lim_n\frac{b_{2n}}{b_{1n}}\!=\!\infty,\,\,\, l_{32}\!=\!\lim_n\frac{b_{3n}}{b_{2n}}\!=\!\infty,\,\,\,\text{ therefore,}\,\,\,
l_{31}\!=\!\lim_n\frac{b_{3n}}{b_{1n}}\!=\!\infty.
\end{equation}
Since $b_{1n}\equiv 1$ we conclude that 
\begin{equation}
\label{d-infty}
 l_{21}\!=\!\lim_nd_{2n}=\infty\quad\text{and}  \quad
l_{31}\!=\!\lim_nd_{3n}=\infty. 
\end{equation}
Therefore, we get for some $C>0$ and all $n\in \mathbb Z$
\begin{equation}
\label{D_n<C}
1\leq D_n=1+\frac{1}{d_{2n}}+\frac{1}{d_{3n}}\leq C.
\end{equation}
By \eqref{rho(CCC)} and \eqref{rho_r(CCC)}
we obtain
\begin{eqnarray*}
&&
\rho_r(C_1,C_2,C_3)\!=\!
\Vert C_1Y_1^{(r)}\!+\!C_2Y_2^{(r)}\!+\!C_2Y_3^{(r)}\Vert ^2\!=\!
\sum_n\frac{\vert C_1+C_2d_{2n}+C_3d_{3n}\vert^2}
 {C_{rn}} \\
 &&
 \sim \sum_n\frac{\vert C_1+C_2d_{2n}+C_3d_{3n}\vert^2}
 {C'_{rn}}=:\rho'_r(C_1,C_2,C_3),\\
&&
\text{where}\quad
C'_{rn}=1+2d_{2n}d_{3n},\quad
C'_{rn}=d^2_{2n}+2d_{2n}d_{3n},\quad
C'_{rn}=d^2_{3n}+2d_{2n}d_{3n}.
\end{eqnarray*}
We should study when $\rho'_r(C_1,C_2,C_3)=\infty$ for some $1\leq r \leq 3$:
\begin{eqnarray*}
&&
\sum_n\frac{\vert C_1\!+\!C_2d_{2n}\!+\!C_3d_{3n}\vert^2}
 {1+2d_{2n}d_{3n}},\,\,
\sum_n\frac{\vert C_1\!+\!C_2d_{2n}\!+\!C_3d_{3n}\vert^2}
 {d^2_{2n}+2d_{2n}d_{3n}},\,\,
 \sum_n\frac{\vert C_1\!+\!C_2d_{2n}\!+\!C_3d_{3n}\vert^2}
 {d^2_{3n}+2d_{2n}d_{3n}}. 
\end{eqnarray*}
 Denoting  as before $d_n=\frac{d_{3n}}{d_{2n}}=:l^{-1}_n$, we get
\begin{eqnarray*}
&&
\rho'_1(C_1,C_2,C_3)=
\sum_n\frac{\vert C_1\!+\!C_2d_{2n}\!+\!C_3d_{3n}\vert^2}
 {1+2d_{2n}d_{3n}}=
 \sum_n\frac{\vert \frac{C_1}{d_{2n}}\!+\!C_2\!+\!C_3d_{n}\vert^2}
 {\frac{1}{d^2_{2n}}+2d_{n}}
\\
 &&
\sim \sum_n\frac{\vert C_2\!+\!C_3d_{n}\vert^2}{2d_{n}},\quad 
\rho'_2(C_1,C_2,C_3)=
\sum_n\frac{\vert C_1\!+\!C_2d_{2n}\!+\!C_3d_{3n}\vert^2}
 {d^2_{2n}+2d_{2n}d_{3n}}=
\\
 &&
\sum_n\frac{\vert \frac{C_1}{d_{2n}}\!+\!C_2\!+\!C_3d_{n}\vert^2}
 {1+2d_{n}}\sim \sum_n\frac{\vert C_2\!+\!C_3d_{n}\vert^2}{1+2d_{n}},\\
&&
\rho'_3(C_1,C_2,C_3)=
\sum_n\frac{\vert C_1\!+\!C_2d_{2n}\!+\!C_3d_{3n}\vert^2}
 {d^2_{3n}+2d_{2n}d_{3n}}=
\\
 &&
\sum_n\frac{\vert \frac{C_1}{d_{2n}}\!+\!C_2\!+\!C_3d_{n}\vert^2}
 {d_n^2+2d_{n}}\sim \sum_n\frac{\vert C_2\!+\!C_3d_{n}\vert^2}{d_n^2+2d_{n}}=\sum_n\frac{\vert C_2l_n\!+\!C_3\vert^2}{1+2l_{n}}.
\end{eqnarray*}
By Lemma~\ref{l.light.123}  we get when $C_2C_3^{-1}>0$
\begin{eqnarray}
\nonumber 
&&
\rho'_1(C_1,C_2,-C_3)\sim  \sum_n\frac{\vert C_2\!-\!C_3d_{n}\vert^2}{2d_{n}}>\sum_n\frac{\vert C_2\!-\!C_3d_{n}\vert^2}{1+2d_{n}}\sim \sum_nc^2_n\sim\Sigma_{23}(s),\\
\label{rho_r}
&&
\rho'_2(C_1,C_2,-C_3)\sim  \sum_n\frac{\vert C_2\!-\!C_3d_{n}\vert^2}{1+2d_{n}}\sim \sum_nc^2_n\sim\Sigma_{23}(s),\\
\nonumber
&&
\rho'_3(C_1,C_2,-C_3)\sim  \sum_n\frac{\vert C_2l_n\!-\!C_3\vert^2}{1+2l_{n}}\sim \sum_ne^2_n\sim\Sigma_{23}(s),\\
\nonumber
&&
\text{where}\quad d_n=
C_2C_3^{-1}(1+c_n),\,\,\,l_n=C_3C_2^{-1}(1+e_n),\quad\,s^4=C_2C_3^{-1}>0.
\end{eqnarray}
But $\Sigma_{23}(s)=\infty$ for all $s>0$ therefore, for $C_2C^{-1}_3>0$ we have
\begin{equation}
\label{(2)}
\rho_r(C_1,C_2,-C_3)\sim  \Sigma_{23}(s)=\infty.
\end{equation}
If $C_2C^{-1}_3>0$, by \eqref{rho_r} we get

\begin{eqnarray*}
\nonumber 
&&
\sum_n\frac{\vert C_2\!+\!C_3d_{n}\vert^2}{1+2d_{n}}>
\sum_n\frac{\vert C_2\!-\!C_3d_{n}\vert^2}{1+2d_{n}}\sim \sum_nc^2_n\sim\Sigma_{23}(s)=\infty,\\
&&
\sum_n\frac{\vert C_2\!+\!C_3d_{n}\vert^2}{1+2d_{n}}
\sum_n\frac{\vert C_2\!-\!C_3d_{n}\vert^2}{1+2d_{n}}\sim \sum_nc^2_n\sim\Sigma_{23}(s)=\infty,\\
\nonumber
&&
\sum_n\frac{\vert C_2l_n\!+\!C_3\vert^2}{1+2l_{n}}> 
\sum_n\frac{\vert C_2l_n\!-\!C_3\vert^2}{1+2l_{n}}\sim \sum_ne^2_n\sim\Sigma_{23}(s)=\infty.
\end{eqnarray*}
Therefore,
$\rho_r(C_1,C_2,C_3)\!=\!\infty$ for every $(C_1,C_2,C_3)\!\in\!\mathbb R^3\setminus \!\{0\}$.

Consider the {\bf case (3)=(aba)}.
In this case, see \eqref{l(sr)}, we have 
$$
l_{21}\!=\!\lim_n\frac{b_{2n}}{b_{1n}}\!=\!\infty,\quad l_{32}\!=\!\lim_n\frac{b_{3n}}{b_{2n}}<\infty,\quad\text{ therefore,}\quad 
l_{31}\!=\!\lim_n\frac{b_{3n}}{b_{1n}}\!=\!\infty.
$$
So, we have again, see \eqref{d-infty}
\begin{equation*}
 l_{21}\!=\!\lim_nd_{2n}=\infty,\quad\text{and}  \quad
l_{31}\!=\!\lim_nd_{3n}=\infty. 
\end{equation*}
We are reduced to the case (2).

Consider the {\bf case (4)=(baa).
}
Now, see \eqref{l(sr)}, we have 
\begin{equation*}
 l_{21}\!=\!\lim_nd_{2n}<\infty,\quad\text{and}  \quad
l_{31}\!=\!\lim_nd_{3n}=\infty. 
\end{equation*}
Hence \eqref{D_n<C}  holds too and we can use
all estimations of the case (1).
\qed\end{pf}
\begin{rem}
\label{r.(1)}
In the cases (1)--(4) by Lemma~\ref{l.park2} and Lemmas~\ref{x1x1.3}--  \ref{x3x3.3} we can approximate $x_{rn}$ for all $1\leq r\leq 3$ and $n\in \mathbb Z$ and the irreducibility is proved.
\end{rem}

\section{Approximation of $D_{kn}$ and $x_{kn}$}
 \subsection{Approximation of 
 $x_{kn}$ by $A_{nk}A_{tk}$  }
 \label{7.2}
For $m=3$, consider three rows as follows
$$
\left(
\begin{array}{cccccc}
...&b_{11}&b_{12}&...&b_{1n}&...\\
...&b_{21}&b_{22}&...&b_{2n}&...\\
...&b_{31}&b_{32}&...&b_{3n}&...
\end{array}
\right).
$$
Set
\begin{equation}
 \label{lambda(x,x)}
\lambda_k^{(r)}=(b_{1k}+b_{2k}+b_{3k})^2-(b_{1k}^2+b_{2k}^2+b_{3k}^2-b_{rk}^2),\,\,r=1,2,3,\,\,k\in \mathbb Z.
\end{equation}
Denote by $Y_r^{(s)}$  the following vectors:
\begin{equation}
 \label{y(1),y(3),y(3)}
 x_{rk}^{(s)}=
b_{rk}/ \sqrt{\lambda_k^{(s)}}
,\,\,k\in{\mathbb Z},\quad Y_r^{(s)}=(x_{rk}^{(s)})_{k\in {\mathbb Z}}.
\end{equation}
\begin{lem}
\label{x1x1.3}
For any  $n,t\in {\mathbb Z}$ one has
$$
x_{1n}x_{1t}{\bf 1}\in\langle A_{nk}A_{tk}{\bf 1}\mid k\in{\mathbb
Z}\rangle \,\,\Leftrightarrow\,\,
\Delta(Y_1^{(1)},Y_2^{(1)},Y_3^{(1)})=\infty.
$$
\end{lem}
Similarly, using the cyclic permutation of vectors, and changing $\lambda_k^{(r)}$ we arrive at the following lemma.
\begin{lem}
\label{x2x2.3}
  For any  $n,t\in {\mathbb Z}$ one has
$$
x_{2n}x_{2t}{\bf 1}\in\langle A_{nk}A_{tk}{\bf 1}\mid k\in{\mathbb
Z}\rangle \,\,\Leftrightarrow\,\,
\Delta(Y_2^{(2)},Y_3^{(2)},Y_1^{(2)}))=\infty.
$$
\end{lem}
\begin{lem}
\label{x3x3.3}
For any  $n,t\in {\mathbb Z}$ one has
$$
x_{3n}x_{3t}{\bf 1}\in\langle A_{nk}A_{tk}{\bf 1}\mid k\in{\mathbb
Z}\rangle \,\,\Leftrightarrow\,\,
\Delta(Y_3^{(3)},Y_1^{(3)},Y_2^{(3)})=\infty.
$$
\end{lem}
\begin{pf}The proof of  Lemma~\ref{x1x1.3} is also based on Lemma~\ref{l.m=3}. 
We study when
$
x_{1n}x_{1t}{\bf 1}\in\langle A_{nk}A_{tk}{\bf 1}\mid k\in{\mathbb
Z}\rangle.
$
Since
\begin{eqnarray*}
&&A_{nk}A_{tk}=(x_{1n}D_{1k}+x_{2n}D_{2k}+x_{3n}D_{3k})(x_{1t}D_{1k}+x_{2t}D_{2k}+x_{3t}D_{3k})\\
&&=x_{1n}x_{1t}D_{1k}^2+
x_{2n}x_{2t}D_{2k}^2+
x_{3n}x_{3t}D_{3k}^2+
(x_{1n}x_{2t}+x_{2n}x_{1t})D_{1k}D_{2k}+
\\
&&
(x_{1n}x_{3t}+x_{3n}x_{1t})D_{1k}D_{3k}+
(x_{2n}x_{3t}+x_{3n}x_{2t})D_{2k}D_{3k}
\end{eqnarray*}
and $MD_{1k}^2{\bf 1}=-\frac{b_{1k}}{2}$, we take $t=(t_k)$ as follows:
$-\sum_{k=-m}^mt_k\frac{b_{1k}}{2}=(t,b')=1$, where $t=(t_k)_{k=-m}^m$ and $b'=-(\frac{b_{1k}}{2})_{k=-m}^m\sim b=(b_{1k})_{k=-m}^m$. 
We have
\begin{eqnarray*}
&&\Vert \big[\sum_{k=-m}^mt_kA_{nk}A_{tk}-x_{1n}x_{1t}\big]{\bf
1}\Vert^2=\\
&&\Vert \sum_{k=-m}^mt_k\big[
x_{1n}x_{1t}\Big(D_{1k}^2\!+\!\frac{b_{1k}}{2}\Big)\!+\!
x_{2n}x_{2t}D_{2k}^2\!+\!
x_{3n}x_{3t}D_{3k}^2\!+\!
(x_{1n}x_{2t}+x_{2n}x_{1t})\times
\\
&&
 D_{1k}D_{2k}+(x_{1n}x_{3t}+x_{3n}x_{1t})D_{1k}D_{3k}+
(x_{2n}x_{3t}+x_{3n}x_{2t})D_{2k}D_{3k}
\Big]{\bf 1}\Vert^2\\
&=&\sum_{-m\leq k,r\leq m}(f_k,f_r)t_kt_r=:(A_{2m+1}t,t),
\end{eqnarray*}
where $A_{2m+1}=\big((f_k,f_r)\big)_{k,r=-m}^m$ and
$f_k=\sum_{r=1}^3f_k^r+\sum_{1\leq i<j\leq 3}f_k^{ij},$
\begin{equation}
\label{f_k(r),f_k(ij)}
f_k^r= x_{rn}x_{rt}\Big(D_{rk}^2+\frac{b_{rk}}{2}\delta_{1r}\Big){\bf 1},\quad
f_k^{ij}=(x_{in}x_{jt}+x_{jn}x_{it})D_{ik}D_{jk}{\bf 1}
\end{equation}
for $1\leq r\leq 3,\,\, 1\leq i<j\leq 3$.
Since 
$$
f_k^r\perp f_k^{ij},\quad 
f_k^{ij}\perp f_k^{i'j'}
$$
for different $(ij),\,\,(i'j')$,
writing $c_{kn}=\Vert x_{kn}\Vert^2=\frac{1}{2b_{kn}}+a_{kn}^2$,  we get
\begin{eqnarray*}
&&(f_k,f_k)=\sum_{r=1}^3\Vert f_k^r\Vert^2+\sum_{1\leq i<j\leq 3}\Vert f_k^{ij}\Vert^2=
\\
&&
c_{1n}c_{1t}2\Big(\frac{b_{1k}}{2}\Big)^2+
c_{2n}c_{2t}3\Big(\frac{b_{2k}}{2}\Big)^2+
c_{3n}c_{3t}3\Big(\frac{b_{3k}}{2}\Big)^2+
\\
&&
\big(c_{1n}c_{2t}+c_{1t}c_{2n}+2a_{1n}a_{2t}a_{1t}a_{2n}\big)\frac{b_{1k}}{2}\frac{b_{2k}}{2}+\big(c_{1n}c_{3t}
+c_{3t}c_{1n}+2a_{1n}a_{3t}a_{3t}a_{1n}\big)\times \\
&&
\frac{b_{1k}}{2}\frac{b_{3k}}{2}+
\big(c_{2n}c_{3t}+c_{3t}c_{2n}+2a_{2n}a_{3t}a_{3t}a_{2n}\big)\frac{b_{2k}}{2}\frac{b_{3k}}{2} \sim(b_{1k}+b_{2k}+b_{3k})^2,\\ 
&&(f_k,f_r)=(f_k^2,f_r^2)+(f_k^3,f_r^3)=
c_{2n}c_{2t}\frac{b_{2k}}{2}\frac{b_{2r}}{2}+c_{3n}c_{3t}\frac{b_{3k}}{2}\frac{b_{3r}}{2}\sim b_{2k}b_{2r}+b_{3k}b_{3r}.
\end{eqnarray*}
Finally, we have
\begin{equation}
\label{(A_{r,s}),xx1,m=2}
(f_k,f_k) \sim (b_{1k}+b_{2k}+b_{3k})^2, \quad(f_k,f_r)\sim b_{2k}b_{2r}+b_{3k}b_{3r},\quad k\not=r.
\end{equation}

Set
\begin{equation}
\lambda_k=(b_{1k}+b_{2k}+b_{3k})^2-(b_{2k}^2+b_{3k}^2),\quad g_k=(b_{2k},b_{3k}),
\end{equation}
then
\begin{equation}
(f_k,f_k)\sim \lambda_k+(g_k,g_k),\quad (f_k,f_r)\sim (g_k,r_r).
\end{equation}

For $A_{2m+1}=\big((f_k,f_r)\big)_{k,r=-m}^m,$ and
$b=-(b_{1k}/2)_{k=-m}^m\in {\mathbb R}^{2m+1}$ we have
$$
A_{2m+1}= \sum_{k=-m}^m\lambda_kE_{kk}+\gamma(g_{-m},\dots,g_0,\dots,g_m).
$$
To finish the proof, it  suffices to invoke Lemma~\ref{l.m=3}.
\qed\end{pf}

\subsection{Approximation of $D_{rn}$ by $A_{kn}$}
We will formulate several  lemmas, which will be  useful for
approximation of the independent variables $x_{kn}$ and operators
$D_{kn}$ by combinations of the generators $A_{kn}$. 
The generators $A_{kn}$ have the following form:
$$
A_{kn}=x_{1k}D_{1n}+x_{2k}D_{2n}+x_{3k}D_{3n},\,\,k,n\in \mathbb Z.
$$
For $m=3$, consider three rows as follows
$$
\left(
\begin{array}{cccccc}
...&a_{11}&a_{12}&...&a_{1n}&...\\
...&a_{21}&a_{22}&...&a_{2n}&...\\
...&a_{31}&a_{32}&...&a_{2n}&...
\end{array}
\right).
$$
Set
\begin{equation}
 \label{lambda(D)}
\lambda_k=\frac{1}{2b_{1k}}+\frac{1}{2b_{2k}}+\frac{1}{2b_{3k}}.
\end{equation}
Denote by $Y_1,Y_2$ and $Y_3$ the three following vectors:
\begin{equation}
 \label{Y(1),Y(3),Y(3)}
 x_{rk}=\frac{a_{rk}}{\sqrt{\lambda_k}},\,\,k\in{\mathbb Z}
 ,\quad Y_r=(x_{rk})_{k\in {\mathbb Z}}.
\end{equation}
The proofs of  Lemmas~\ref{d1.3}--\ref{d3.3} and  \ref{x1x1.3}--\ref{x3x3.3}
are based on Lemma~\ref{l.m=3}.
%
%
\begin{lem}
\label{d1.3}
For any  $l\in {\mathbb Z}$ we have
$$
D_{1l}{\bf 1}\in\langle A_{kl}{\bf 1}\mid k\in {\mathbb Z}\rangle
\quad\Leftrightarrow\quad \Delta(Y_1,Y_2,Y_3)=\infty.
$$
\end{lem}
Similarly, by cyclic permutation 
of the vectors, 
we obtain the following two lemmas.
\begin{lem}
\label{d2.3}
 For any  $l\in {\mathbb Z}$ we have
$$
D_{2l}{\bf 1}\in\langle A_{kl}{\bf 1}\mid k\in {\mathbb Z}\rangle
\quad\Leftrightarrow\quad \Delta(Y_2,Y_3,Y_1)=\infty.
$$
\end{lem}
\begin{lem}
\label{d3.3}
 For any  $l\in {\mathbb Z}$ we have
$$
D_{3l}{\bf 1}\in\langle A_{kl}{\bf 1}\mid k\in {\mathbb Z}\rangle
\quad\Leftrightarrow\quad \Delta(Y_3,Y_1,Y_2)=\infty.
$$
\end{lem}

\begin{pf} 
We determine when the inclusion
$$
D_{1n}{\bf 1}\in\langle A_{kn}{\bf
1}=(x_{1k}D_{1n}+x_{2k}D_{2n}+x_{3k}D_{3n}){\bf 1}\mid k\in {\mathbb Z}\rangle
$$
holds. Fix $m\in {\mathbb N}$, since $Mx_{1k}=a_{1k}$, we put 
$
\sum_{k=-m}^mt_ka_{1k}=(t,b)=1,
$
where $t=(t_k)_{k=-m}^m$ and $b=(a_{1k})_{k=-m}^m$. We have
\begin{eqnarray*}
 &&
\Vert \big[\sum_{k=-m}^mt_k(x_{1k}D_{1n}+x_{2k}D_{2n}+x_{3k}D_{3n})-D_{1n}\big]{\bf
1}\Vert^2\\
&=&
\Vert \sum_{k=-m}^mt_k[(x_{1k}-a_{1k})D_{1n}+x_{2k}D_{2n}+x_{3k}D_{3n}]{\bf
1}\Vert^2\\
&=&\sum_{-m\leq k,r\leq m}(f_k,f_r)t_kt_r=:(A_{2m+1}t,t),\,\,\text{where}\,\, A_{2m+1}=\big((f_k,f_r)\big)_{k,r=-m}^m,\\
&&\quad\text{ and}\quad
f_k=[(x_{1k}-a_{1k})D_{1n}+x_{2k}D_{2n}+x_{3k}D_{3n}]{\bf 1}.
\end{eqnarray*}  
We get
\begin{eqnarray*}
&&(f_k,f_k)=\Vert\left[(x_{1k}-a_{1k})D_{1n}+x_{2k}D_{2n}+x_{3k}D_{3n}\right]{\bf
1}\Vert^2\\
&&=\frac{1}{2b_{1k}}\frac{b_{1n}}{2}+\Big(\frac{1}{2b_{2k}}+a_{2k}^2\Big)
\frac{b_{2n}}{2} +\Big(\frac{1}{2b_{3k}}+a_{3k}^2\Big)
\frac{b_{3n}}{2} \\
&&\sim
\frac{1}{2b_{1k}}+\frac{1}{2b_{2k}}+\frac{1}{2b_{3k}}+a_{2k}^2+a_{3k}^2,\\
&&(f_k,f_r)=\Big(\big[(x_{1k}-a_{1k})D_{1n}+x_{2k}D_{2n}+x_{3k}D_{3n}\big]{\bf 1},\\
&&\big[(x_{1r}-a_{1r})D_{1n}+x_{2r}D_{2n}+x_{3r}D_{3n}\big]{\bf 1}\Big) \\
&&=(x_{2k}{\bf 1},x_{2r}{\bf 1})(D_{2n}{\bf 1},D_{2n}{\bf 1})+
(x_{3k}{\bf 1},x_{3r}{\bf 1})(D_{3n}{\bf 1},D_{3n}{\bf 1})\\
&&=a_{2k}a_{2r}\frac{b_{2n}}{2}+a_{3k}a_{3r}\frac{b_{3n}}{2}\sim a_{2k}a_{2r}+a_{3k}a_{3r}.
\end{eqnarray*}
Finally, we have
\begin{equation}
\label{(A_{r,s}),D1,m=3}
(f_k,f_k)\!\sim\!
\frac{1}{2b_{1k}}+\frac{1}{2b_{2k}}+\frac{1}{2b_{3k}}+a_{2k}^2+a_{3k}^2,\,\,\,\,\,
(f_k,f_r)\sim a_{2k}a_{2r}+a_{3k}a_{3r},\,\, k\not=r.
\end{equation}
If we denote
\begin{equation}
\lambda_k=\frac{1}{2b_{1k}}+\frac{1}{2b_{2k}}+\frac{1}{2b_{3k}},\quad g_k=(a_{2k},a_{3k}),
\end{equation}
then we have
\begin{equation}
(f_k,f_k)\sim \lambda_k+(g_k,g_k),\quad (f_k,f_r)\sim (g_k,g_r).
\end{equation}
For $A_{2m+1}=\big((f_k,f_r)\big)_{k,r=-m}^m,$ and
$b=(a_{1k})_{k=-m}^m\in {\mathbb R}^{2m+1}$ we have
$$
A_{2m+1}= \sum_{k=-m}^m\lambda_kE_{kk}+\gamma(g_{-m},\dots,g_0,\dots,g_m).
$$
To finish the proof, it  suffices to apply  Lemma~\ref{l.m=3}.
 \qed\end{pf}

The proofs of  Lemmas~\ref{d2.3} and \ref{d3.3} are exactly the same.

\subsection{Approximation of $x_{rk}$ by $D_{rn}A_{kn}$}
\begin{lem}
\label{x1,dA.m=3}
 For any  $k\in {\mathbb Z}$ we get
 $$
 x_{1k}{\bf 1}\in\langle D_{1n}A_{kn}{\bf 1}\mid n\in{\mathbb Z}\rangle
 \quad\Leftrightarrow\quad
\sum_{n\in{\mathbb Z}}\frac{b_{1n}}{b_{1n}+b_{2n}+b_{3n}}=\infty.
$$
\end{lem}
\begin{pf} Since
$$
D_{1n}A_{kn}=x_{1k}D_{1n}^2+x_{2k}D_{1n}D_{2n}+x_{3k}D_{1n}D_{3n}
$$
and $MD_{1k}^2{\bf 1}\!=\!-\!\frac{b_{1k}}{2}$, we take $t\!=\!(t_k)_{k=-m}^m$ as follows:
$
(t,b')=1$, where $t=(t_k)_{k=-m}^m$ and $b'=-(\frac{b_{1k}}{2})_{k=-m}^m\sim b=-(b_{1k})_{k=-m}^m$. 
We have
\begin{eqnarray*}
&&\Vert \big[\sum_{n=-m}^mt_nD_{1n}A_{kn}-x_{1k}\big]{\bf
1}\Vert^2\\
&=&\Vert \sum_{n=-m}^mt_n
\left[x_{1k}\Big(D_{1n}^2+\frac{b_{1n}}{2}\Big)+x_{2k}D_{1n}D_{2n}+x_{3k}D_{1n}D_{3n}\right]{\bf 1}
\Vert^2\\
&=&\sum_{-m\leq n,r\leq m}(f_n,f_r)t_kt_r=:(A_{2m+1}t,t),
\end{eqnarray*}
where $A_{2m+1}=\big((f_n,f_r)\big)_{n,r=-m}^m$
and
$$
f_n=\left[x_{1k}\Big(D_{1n}^2+\frac{b_{1n}}{2}\Big)+x_{2k}D_{1n}D_{2n}+x_{3k}D_{1n}D_{3n}\right]{\bf 1}.
$$
We have
\begin{equation*}
(f_n,f_n)\sim b_{1k}(b_{1n}+b_{2n}+b_{3n}) \quad\text{and}\quad (f_n,f_k)=0\quad\text{for}\quad n\not=k.
\end{equation*}
Therefore, by \eqref{A.min2}
$$
\min_{t\in \mathbb R^{2m+1}}\Vert \left[\sum_{n=-m}^mt_nD_{1n}A_{kn}-x_{1k}\right]{\bf
1}\Vert^2=\left(\sum_{n=-m}^m\frac{b_{1k}}{b_{1n}+b_{2n}+b_{3n}}\right)^{-1}.\hskip 1cm\Box
$$
\end{pf}
\begin{lem}
\label{x2,dA.m=3}
 For any  $k\in {\mathbb Z}$ we have
 $$
 x_{2k}{\bf 1}\in\langle D_{2n}A_{kn}{\bf 1}\mid n\in{\mathbb Z}\rangle
 \quad\Leftrightarrow\quad
\sum_{n\in{\mathbb Z}}\frac{b_{2n}}{b_{1n}+b_{2n}+b_{3n}}=\infty.
$$
\end{lem}
\begin{lem}
\label{x3,dA.m=3}
 For any  $k\in {\mathbb Z}$ we have
 $$
 x_{3k}{\bf 1}\in\langle D_{3n}A_{kn}{\bf 1}\mid n\in{\mathbb Z}\rangle
 \quad\Leftrightarrow\quad
\sum_{n\in{\mathbb Z}}\frac{b_{3n}}{b_{1n}+b_{2n}+b_{3n}}=\infty.
$$
\end{lem}

\subsection{Approximation of $D_{kn}$ by $x_{rk}A_{kn}$}
Set
\begin{eqnarray}
\label{lambda(1)}
&&\lambda_k^{(1)}\!=\!\Big(\frac{1}{2b_{1k}}\!+\!a_{1k}^2\Big)\Big(\frac{1}{2b_{1k}}\!+\!\frac{1}{2b_{2k}}\!+\!\frac{1}{2b_{3k}}
\!+\!a_{2k}^2+a_{3k}^2\Big)\!-\!
a_{1k}^2(a_{2k}^2\!+\!a_{3k}^2),\\
\label{y_1(123)}
&&Y_{11}\!=\!\left(
\frac{\frac{1}{2b_{1k}}\!+\!a_{1k}^2}
{\sqrt{\lambda_k^{(1)}}}
\right)_{k\in \mathbb Z},
Y_{12}\!=\!\left(\frac{a_{1k}a_{2k}}{\sqrt{\lambda_k^{(1)}}}\right)_{k\in \mathbb Z},
Y_{13}\!=\!\left(\frac{a_{1k}a_{3k}}{\sqrt{\lambda_k^{(1)}}}\right)_{k\in \mathbb Z}.
\end{eqnarray}
\begin{lem}
\label{dx1.3}
For any  $n\in {\mathbb Z}$ we have
$$
D_{1n}{\bf 1}\in\langle x_{1k}A_{kn}{\bf 1}\mid k\in {\mathbb Z}\rangle
\quad\Leftrightarrow\quad 
\Delta(Y_{11},Y_{12},Y_{13})=\infty.
$$
\end{lem}
\begin{pf} 
We determine  when the inclusion
$$
D_{1n}{\bf 1}\in\langle x_{1k}A_{kn}{\bf
1}=(x^2_{1k}D_{1n}+x_{1k}x_{2k}D_{2n}+x_{1k}x_{3k}D_{3n}){\bf 1}\mid k\in {\mathbb Z}\rangle
$$
holds. Fix $m\in {\mathbb N}$, since $Mx^2_{1k}=\frac{1}{2b_{1k}}+a^2_{1k}$, we put 
$$
\sum_{k=-m}^mt_k\big(\frac{1}{2b_{1k}}+a^2_{1k}\big)=(t,b)=1,
$$
where $t=(t_k)_{k=-m}^m$ and $b=(\frac{1}{2b_{1k}}+a^2_{1k})_{k=-m}^m$. We have
\begin{eqnarray*}
 &&
\Vert \big[\sum_{k=-m}^mt_k(x^2_{1k}D_{1n}+x_{1k}x_{2k}D_{2n}+x_{1k}x_{3k}D_{3n})-D_{1n}\big]{\bf
1}\Vert^2\\
&=&
\Vert 
\sum_{k=-m}^mt_k\Big[\Big(x^2_{1k}-\Big(\frac{1}{2b_{1k}}+a^2_{1k}\Big)\Big)D_{1n}+x_{1k}x_{2k}D_{2n}+x_{1k}x_{3k}D_{3n}\Big]{\bf1}
\Vert^2
\\
&=&\sum_{-m\leq k,r\leq m}(f_k,f_r)t_kt_r=:(A_{2m+1}t,t),\,\,\text{where}\,\, A_{2m+1}=\big((f_k,f_r)\big)_{k,r=-m}^m,\\
&&\quad\text{ and}\quad
f_k=\Big[\Big(x^2_{1k}-\Big(\frac{1}{2b_{1k}}+a^2_{1k}\Big)\Big)D_{1n}+x_{1k}x_{2k}D_{2n}+x_{1k}x_{3k}D_{3n}\Big]{\bf 1}.
\end{eqnarray*}  
Since $M\vert \psi-M\vert \psi\vert\vert^2=M \psi^2-\vert M \psi\vert^2$ we have
\begin{eqnarray*}
&& 
M\left|x^2_{1k}-\Big(\frac{1}{2b_{1k}}+a^2_{1k}\Big)\right|^2=M x^4_{1k}-\Big(\frac{1}{2b_{1k}}+a^2_{1k}\Big)^2\\
&&=\frac{3}{(2b_{1k})^2}+6\frac{1}{2b_{1k}}a_{1k}^2+a^4_{1k}-\Big(\frac{1}{2b_{1k}}+a^2_{1k}\Big)^2
=\frac{1}{2b_{1k}}\Big(\frac{2}{2b_{1k}}+4a^2_{1k}\Big),
\end{eqnarray*}
we get
\begin{eqnarray*}
&&(f_k,f_k)=\Vert\Big[\Big(x^2_{1k}-\Big(\frac{1}{2b_{1k}}+a^2_{1k}\Big)\Big)D_{1n}+x_{1k}x_{2k}D_{2n}+x_{1k}x_{3k}D_{3n}\Big]{\bf1}\Vert^2=\\
&&\frac{1}{2b_{1k}}\Big(\frac{2}{2b_{1k}}\!+\!4a^2_{1k}\Big)\frac{b_{1n}}{2}\!+\!\Big(\frac{1}{2b_{1k}}+a_{1k}^2\Big)
\Big(\frac{1}{2b_{2k}}+a_{2k}^2\Big)
\frac{b_{2n}}{2}\!+\!\Big(\frac{1}{2b_{1k}}+a_{1k}^2\Big)
\frac{b_{3n}}{2} \\
&&\times\Big(\frac{1}{2b_{3k}}+a_{3k}^2\Big)\sim
\Big(\frac{1}{2b_{1k}}+a_{1k}^2\Big)\Big(\frac{1}{2b_{1k}}+\frac{1}{2b_{2k}}+\frac{1}{2b_{3k}}+a_{2k}^2+a_{3k}^2\Big)\,\\
&&(f_k,f_r)=\Big(\Big[\Big(x^2_{1k}-\Big(\frac{1}{2b_{1k}}+a^2_{1k}\Big)\Big)D_{1n}+x_{1k}x_{2k}D_{2n}+x_{1k}x_{3k}D_{3n}\Big]{\bf 1},\\
&&\Big[\Big(x^2_{1r}-\Big(\frac{1}{2b_{1r}}+a^2_{1r}\Big)\Big)D_{1n}+x_{1r}x_{2r}D_{2n}+x_{1r}x_{3r}D_{3n}\Big]{\bf 1}\Big)= \\
&&(x_{1k}{\bf 1},x_{1r}{\bf 1})(x_{2k}{\bf 1},x_{2r}{\bf 1})(D_{2n}{\bf 1},D_{2n}{\bf 1})\!+\!
(x_{1k}{\bf 1},x_{1r}{\bf 1})(x_{3k}{\bf 1},x_{3r}{\bf 1})(D_{3n}{\bf 1},D_{3n}{\bf 1})\\
&&=a_{1k}a_{1r}a_{2k}a_{2r}\frac{b_{2n}}{2}+a_{1k}a_{1r}a_{3k}a_{3r}\frac{b_{3n}}{2}\simeq a_{1k}a_{1r}(a_{2k}a_{2r}+a_{3k}a_{3r}).
\end{eqnarray*}
Finally, we have
\begin{eqnarray}
 \label{(A_{r,s}),DX1,m=3}
&&(f_k,f_k)\!\sim\!
\Big(\frac{1}{2b_{1k}}+a_{1k}^2\Big)\Big(\frac{1}{2b_{1k}}+\frac{1}{2b_{2k}}+\frac{1}{2b_{3k}}+a_{2k}^2+a_{3k}^2\Big),\\
\nonumber
&&(f_k,f_r)\sim a_{1k}a_{1r}(a_{2k}a_{2r}+a_{3k}a_{3r}),\,\, k\not=r.
\end{eqnarray}
Set
\begin{eqnarray}
\nonumber
&&\lambda_k^{(1)}=\Big(\frac{1}{2b_{1k}}+a_{1k}^2\Big)\Big(\frac{1}{2b_{1k}}+\frac{1}{2b_{2k}}+\frac{1}{2b_{3k}}+a_{2k}^2+a_{3k}^2\Big)-
a_{1k}^2(a_{2k}^2+a_{3k}^2),\\
\label{l_k(xD)}
&& g_k=a_{1k}(a_{2k},a_{3k}),
\end{eqnarray}
then
\begin{equation}
(f_k,f_k)=\lambda_k^{(1)}+(g_k,g_k)
\quad (f_k,f_r)\sim (g_k,g_r),\,\,k\not=r.
\end{equation}
For $A_{2m+1}=\big((f_k,f_r)\big)_{k,r=-m}^m,$ and
$b=(a_{1k})_{k=-m}^m\in {\mathbb R}^{2m+1}$ we have
$$
A_{2m+1}= \sum_{k=-m}^m\lambda_kE_{kk}+\gamma(g_{-m},\dots,g_0,\dots,g_m).
$$
To finish the proof, it  suffices to apply  Lemma~\ref{l.m=3}. 
\qed\end{pf}

Similarly, by cyclic permutation of the indices, we obtain the following lemmas.
\begin{lem}
\label{dx2.3}
 For any  $l\in {\mathbb Z}$ we have
$$
D_{2l}{\bf 1}\in\langle x_{2k}A_{kl}{\bf 1}\mid k\in {\mathbb Z}\rangle
\quad\Leftrightarrow\quad 
\Delta(Y_{22},Y_{23},Y_{21})=\infty.
$$
\end{lem}
\begin{lem}
\label{dx3.3}
 For any  $l\in {\mathbb Z}$ we have
$$
D_{3l}{\bf 1}\in\langle x_{3k}A_{kl}{\bf 1}\mid k\in {\mathbb Z}\rangle
\quad\Leftrightarrow\quad 
\Delta(Y_{33},Y_{31},Y_{32})=\infty.
$$
\end{lem}
Here we set
\begin{eqnarray}
\label{lambda(2)}
&&\!\!\lambda_k^{(2)}\!=\!\Big(\frac{1}{2b_{2k}}\!+\!a_{2k}^2\Big)\Big(\frac{1}{2b_{1k}}\!+\!\frac{1}{2b_{2k}}\!+\!\frac{1}{2b_{3k}}+a_{1k}^2\!+\!
a_{3k}^2\Big)\!-\!a_{2k}^2(a_{1k}^2\!+\!a_{3k}^2),\,\,\\
\label{y_2(123)}
&&Y_{21}\!=\!\left(\frac{a_{1k}a_{2k}}{\sqrt{\lambda_k^{(2)}}}\right)_{k\in \mathbb Z}\!\!\!\!,
Y_{22}\!=\!\left(\frac{\frac{1}{2b_{2k}}+a_{2k}^2}{\sqrt{\lambda_k^{(2)}}}\right)_{k\in \mathbb Z}\!\!\!,
Y_{23}\!=\!\left(\frac{a_{2k}a_{3k}}{\sqrt{\lambda_k^{(2)}}}\right)_{k\in \mathbb Z}\!\!\!, \\
\label{lambda(3)}
&&\lambda_k^{(3)}\!=\!\Big(\frac{1}{2b_{3k}}\!+\!a_{3k}^2\Big)\Big(\frac{1}{2b_{1k}}\!+\!\frac{1}{2b_{2k}}\!+\!\frac{1}{2b_{3k}}\!+\!a_{1k}^2
\!+\!a_{2k}^2\Big)\!-\!
a_{3k}^2(a_{1k}^2\!+\!a_{2k}^2),\qquad\\
\label{y_3(123)}
&&Y_{31}=\left(\frac{a_{1k}a_{3k}}{\sqrt{\lambda_k^{(3)}}}\right)_{k\in \mathbb Z}, 
Y_{32}=\left(\frac{a_{2k}a_{3k}}{\sqrt{\lambda_k^{(3)}}}\right)_{k\in \mathbb Z},
Y_{33}=\left(\frac{\frac{1}{2b_{3k}}+a_{3k}^2}{\sqrt{\lambda_k^{(3)}}}\right)_{k\in \mathbb Z}.
\end{eqnarray}

\subsection{Approximation of $D_{rn}$ by  $(x_{3k}-a_{3k})A_{kn}$ and $\exp(i s_kx_{rk})A_{kn}$}
\begin{lem}
\label{l.x(3)A}
For any  $n\in {\mathbb Z}$ we have
\begin{eqnarray}
\nonumber
&&D_{3n}{\bf 1}\in\langle (x_{3k}-a_{3k})A_{kn}{\bf 1}\mid k\in {\mathbb Z}\rangle
\,\,\,\Leftrightarrow\,\,\,\\
\label{x(3)A}
&&
\Sigma_3(\mu):=
\sum_{k\in \mathbb Z}\frac{1}{2b_{3k}}\Big(\frac{1}{2b_{1k}}+\frac{1}{2b_{2k}}+\frac{1}{2b_{3k}}+a_{1k}^2+a_{2k}^2+a_{3k}^2\Big)^{-1}\!\!=\!\infty.\qquad
\end{eqnarray}
\end{lem}
\begin{pf} 
We determine  when the following inclusion holds
\begin{eqnarray*}
 &&
 D_{3n}{\bf 1}\in\langle (x_{3k}-a_{3k})A_{kn}{\bf 1}\\
&&=\big(x_{1k}(x_{3k}-a_{3k})D_{1n}+x_{2k}(x_{3k}-a_{3k})D_{2n}+x_{3k}(x_{3k}-a_{3k})D_{3n}\big){\bf 1}\mid k\in {\mathbb Z}\rangle.
\end{eqnarray*}
 Set $\xi_{3k}=x_{3k}(x_{3k}-a_{3k})$. Fix $m\in {\mathbb N}$. We have 
$$M\xi_{3k}{\bf 1}=M(x_{3k}-a_{3k})^2{\bf 1}=\frac{1}{2b_{3k}},\quad\text{chose}\quad(t_k)\,\,\text{as}\,\,
\sum_{k=-m}^mt_k\frac{1}{2b_{3k}}=(t,b)=1,
$$
where we denote $t=(t_k)_{k=-m}^m$ and $b=(\frac{1}{2b_{3k}})_{k=-m}^m$. We have
\begin{eqnarray}
\nonumber
 &&
\Vert \big[\sum_{k=-m}^mt_k\big(x_{1k}(x_{3k}-a_{3k})D_{1n}+x_{2k}(x_{3k}-a_{3k})D_{2n}+\xi_{3k}D_{3n}\big)-D_{3n}\big]{\bf
1}\Vert^2\\
\nonumber
&=&
\Vert 
\sum_{k=-m}^mt_k\Big[
x_{1k}(x_{3k}-a_{3k})D_{1n}+x_{2k}(x_{3k}-a_{3k})D_{2n}+\big(\xi_{3k}-M\xi_{3k}\big)D_{3n}
\Big]{\bf1}
\Vert^2
\\
\label{x(3)A-norm}
&=&\sum_{-m\leq k,r\leq m}(f_k,f_r)t_kt_r=\sum_{-m\leq k\leq m}(f_k,f_k)t_k^2,\,\,\,\text{since}\quad f_k\perp f_r,\quad\text{where} \\
\nonumber
&&
f_k=\Big[
x_{1k}(x_{3k}-a_{3k})D_{1n}+x_{2k}(x_{3k}-a_{3k})D_{2n}+\big(\xi_{3k}-M\xi_{3k}\big)D_{3n}
\Big]{\bf1}.
\end{eqnarray} 
To calculate $M\vert\xi_{3k}-M\xi_{3k}\vert^2$ set
\begin{equation}
 \label{mu(b,a),mu(b,0)}
d\mu_{(b,a)}(x)=\sqrt{\frac{b}{\pi}}
\exp(-b(x-a)^2)dx\quad\text{and}\quad
d\mu_{(b,0)}(x)=\sqrt{\frac{b}{\pi}}\exp(-bx^2)dx,
\end{equation}
then we get $M\xi_{3k}^2= \frac{3}{(2b_{3k})^2}+\frac{a_{3k}^2}{2b_{3k}}$. Indeed,
\begin{eqnarray*}
 &&Mx^2(x-a)^2=\int x^2(x-a)^2d\mu_{(b,a)}(x)=\int x^2(x+a)^2d\mu_{(b,0)}(x)\\
  &&=\int (x^4+2ax^3+a^2x^2)d\mu_{(b,0)}(x)=
  \frac{3}{(2b)^2}+\frac{a^2}{2b}.
\end{eqnarray*}
Since $M\vert \xi-M\vert \xi\vert\vert^2=M \xi^2-\vert M \xi\vert^2$ 
 we get
\begin{equation*}
 M \xi_{3k}^2-M\vert \xi_{3k}\vert^2=\frac{3}{(2b_{3k})^2}+\frac{a_{3k}^2}{2b_{3k}}-\frac{1}{(2b_{3k})^2}=
\frac{1}{2b_{3k}}\Big(\frac{2}{2b_{3k}}+a_{3k}^2\Big).
\end{equation*}
Set $f_{sk}=x_{sk}(x_{3k}\!-\!a_{3k})D_{1n}{\bf 1},\,\,1\leq s\leq 2$ and $f_{3k}=\big(\xi_{3k}\!-\!M\xi_{3k}\big)D_{3n}
{\bf1}$, then 
\begin{eqnarray*}
&&(f_k,f_k)\!=\!\Vert\Big[
x_{1k}(x_{3k}\!-\!a_{3k})D_{1n}\!+\!x_{2k}(x_{3k}\!-\!a_{3k})D_{2n}\!+\!\big(\xi_{3k}\!-\!M\xi_{3k}\big)D_{3n}
\Big]{\bf1}\Vert^2\\
&&=\Big(\frac{1}{2b_{1k}}+a^2_{1k}\Big)\frac{1}{2b_{3k}}\frac{b_{1n}}{2}\!+\!
\Big(\frac{1}{2b_{2k}}+a^2_{2k}\Big)\frac{1}{2b_{3k}}\frac{b_{2n}}{2}\!+\!
\Big(\frac{2}{2b_{3k}}+a^2_{3k}\Big)\frac{1}{2b_{3k}}\frac{b_{3n}}{2}
 \\
&&
\sim\Big(\frac{1}{2b_{1k}}+\frac{1}{2b_{2k}}+\frac{1}{2b_{3k}}+a_{1k}^2+a_{2k}^2+a_{3k}^2\Big)\frac{1}{2b_{3k}},\,\,
\text{since}\,\,f_{lk}\perp f_{sk},\,\,l\not=s,\\
&&(f_k,f_r)\!=\!\Big(
\Big[
x_{1k}(x_{3k}\!-\!a_{3k})D_{1n}\!+\!x_{2k}(x_{3k}\!-\!a_{3k})D_{2n}\!+\!\big(\xi_{3k}\!-\!M\xi_{3k}\big)D_{3n}
\Big]{\bf1},\\
&&\Big[
x_{1r}(x_{3r}-a_{3r})D_{1n}+x_{2r}(x_{3r}-a_{3r})D_{2n}+\big(\xi_{3r}-M\xi_{3r}\big)D_{3n}
\Big]{\bf1}\Big)=0.
\end{eqnarray*}
The previous equality holds since  $f_{lr}\perp f_{sk}$ for $1\leq l,s\leq 3$ and $r\!\not=\! k$. For 
$l\not= s$ this follows from 
$(D_{ln}{\bf 1},D_{sn}{\bf 1})\!=\!0$. For $ l\!=\!s$  it follows from the equalities: 
\begin{equation*}
 (x_{sk}{\bf 1},x_{sr}{\bf 1})=a_{sk}a_{sr},\,\,1\leq s\leq 2,\quad \Big((\xi_{3k}-M\xi_{3k}){\bf 1},(\xi_{3r}-M\xi_{3r}){\bf 1}\Big)=0.
\end{equation*}
and $\big((x_{3k}-a_{3k}){\bf 1},(x_{3r}-a_{3r}){\bf 1}\big)=0$. Finally, we have
\begin{equation*}
(f_k,f_k)\!\sim\!
\Big(\frac{1}{2b_{1k}}+\frac{1}{2b_{2k}}+\frac{1}{2b_{3k}}+a_{1k}^2+a_{2k}^2+a_{3k}^2\Big)\frac{1}{2b_{3k}},\,\,
(f_k,f_r)=0,\,\,k\not=r.
\end{equation*}
Set $a_k=(f_k,f_k)$ and $b_k=M\xi_k$, then by  Lemma~\ref{1.l.min}, \eqref{A.min2} and \eqref{x(3)A-norm} the proof is completed since
\begin{equation*}
\hskip 1.8cm 
\sum_{k}\frac{b_k^2}{a_k}=
\sum_{k\in \mathbb Z}\frac{1}{2b_{3k}}\Big(\frac{1}{2b_{1k}}+\frac{1}{2b_{2k}}+\frac{1}{2b_{3k}}+a_{1k}^2+a_{2k}^2+a_{3k}^2\Big)^{-1}.
\hskip 1.8cm \Box
\end{equation*}
\end{pf}
{\small
Now we would like to approximate $D_{3n}$ by combinations of $\exp\big(is_k(x_{3k}-a_{3k}\big)iA_{kn}$. Set $s=(s_k)_{k\in \mathbb Z}$
\begin{eqnarray}
\label{lambda(3)(s)}
\lambda_k^{(3)}(s_k)&=&\frac{1}{2b_{1k}}+\frac{1}{2b_{2k}}+\frac{1}{2b_{3k}}+a_{3k}^2-
\Big(\frac{s^2}{4b^2_{3k}}+a_{3k}^2\Big)\exp\big(-
\frac{s^2}{2b_{3k}}\big),\\
\label{Y_3(s)}
Y_{31}(s)&=&\left(
\frac{a_{1k}}{\sqrt{\lambda_k^{(3)}(s_k)}}\right)_{k\in \mathbb Z},\quad
Y_{32}(s)=\left(
\frac{a_{2k}}{\sqrt{\lambda_k^{(3)}(s_k)}}\right)_{k\in \mathbb Z},\\
\nonumber
Y_{32}(s)&=&\left(
\frac{\Big(-\frac{s_k}{2b_{3k}}+ia_{3k}\Big)\exp\big(-\frac{s_k^2}{4b_{3k}}\big)}{\sqrt{\lambda_k^{(3)}(s_k)}}\right)_{k\in \mathbb Z}.
\end{eqnarray}
In particular for $s_k=\sqrt{2b_{3k}}$ we get
\begin{equation}
\label{lambda(3)(s).min}
\lambda_k^{(3)}(s_k)=\frac{1}{2b_{1k}}+\frac{1}{2b_{2k}}+
\Big(\frac{1}{2b_{3k}}+a_{3k}^2\Big)\big(1-\frac{1}{e}
\big).
\end{equation}
}
\begin{lem}
\label{exp(x)A3}
For any  $n\in {\mathbb Z}$ we have
$$
D_{3n}{\bf 1}\in\langle \exp\big(is_k(x_{3k}-a_{3k}\big)iA_{kn}{\bf 1}\mid k\in {\mathbb Z}\rangle
\,\,\Leftrightarrow\,\, 
\Delta(Y_{33}(s),Y_{31}(s),Y_{32}(s))=\infty.
$$
\end{lem}
\begin{pf}
We determine when the inclusion
\begin{eqnarray*}
 &&
 D_{3n}{\bf 1}\in\langle \exp\big(is_k(x_{3k}-a_{3k}\big)iA_{kn}{\bf 1}=\Big(ix_{1k}\exp\big(is_k(x_{3k}-a_{3k})\big)D_{1n}\\
&& +ix_{2k}\exp\big(is_k(x_{3k}-a_{3k})\big)D_{2n}+ix_{3k}\exp\big(is_k(x_{3k}-a_{3k})\big)D_{3n}\Big){\bf 1}\mid k\in {\mathbb Z}\rangle
\end{eqnarray*}
holds. Set $\xi_{rk}(s_k)=ix_{rk}\exp\big(is_k(x_{3k}-a_{3k})\big)$ for $1\leq r \leq 3$ and
\begin{eqnarray}
\label{f_k(s)}
f_k(s_k)&=&\Big(\xi_{1k}(s_k)D_{1n}+\xi_{2k}(s_k)D_{2n}+\big[\xi_{3k}(s_k)-M\xi_{3k}(s_k)\big]D_{3n}\Big){\bf 1}\\
\nonumber
&=&\Big(ix_{1k}\exp\big(is_k(x_{3k}-a_{3k}\big)D_{1n}+ix_{2k}\exp\big(is_k(x_{3k}-a_{3k}\big)D_{2n}\\
\nonumber
&+&\big[ix_{3k}\exp\big(is_k(x_{3k}-a_{3k}\big)-M\xi_{3k}(s_k)\big]D_{3n}\Big){\bf 1}.
\end{eqnarray}
We show that 
\begin{eqnarray}
 \label{M(x.exp(x))}
&&M\xi_{3k}(s)=\Big(-\frac{s}{2b_{3k}}+ia_{3k}\Big)\exp\Big(-\frac{s^2}{4b_{3k}}\Big),\\
\label{(f(k),f(k))}
&& (f_k,f_k)\sim \lambda_k^{(3)}(s_k)+(g_k,g_k),\\
\label{(f(k),f(r))}
&& (f_k,f_r)\sim a_{1k}a_{1r}+a_{2k}a_{2r}=(g_k,g_r),
\end{eqnarray}
where $g_k=(a_{1k},a_{2k})\in \mathbb R^2$. Indeed, set
$
F_b(s) =\int_{\mathbb R}\exp\big(is(x-a)\big)d\mu_{(b,a)}(x),
$
then
\begin{equation}
\label{F_b(a)} 
F_b(s)=\int_{\mathbb R}\exp\big(isx\big)d\mu_{(b,0)}(x)=\exp\big(-\frac{s^2}{4b}\Big),
\end{equation}
where $d\mu_{(b,a)}(x)$ and $d\mu_{(b,0)}(x)$ are defined by \eqref{mu(b,a),mu(b,0)}. 
Therefore,
\begin{eqnarray}
&&
H_{a,b}(s)=\int_{\mathbb R}ix\exp\big(is(x-a)\big)d\mu_{(b,a)}(x)=\int_{\mathbb R}i(x+a)\exp(isx)d\mu_{(b,0)}(x)\\
\label{H(a,b,s)}
&&=\frac{dF_b(s)}{ds}+iaF_b(s)=\Big(-\frac{s}{2b}+ia\Big)\exp\Big(-\frac{s^2}{4b}\Big).
\end{eqnarray}

This implies \eqref{M(x.exp(x))}. Further, to obtain  \eqref{(f(k),f(k))} and \eqref{(f(k),f(r))} we write
\begin{eqnarray*}
 &&(f_k,f_r)=\sum_{1\leq t,l\leq 2}(x_{tk},x_{lr})(D_{tn}{\bf 1},D_{ln}{\bf 1})=
 \sum_{1\leq t\leq 2}(x_{tk},x_{tr})(D_{tn}{\bf 1},D_{tn}{\bf 1})
 \\
 &&=a_{1k}a_{1r}\frac{b_{1n}}{2}+a_{2k}a_{2r}\frac{b_{2n}}{2}\sim a_{1k}a_{1r}+a_{2k}a_{2r}=(g_k,g_r),\\
&& (f_k,f_k)=\sum_{1\leq t\leq 2}\Vert x_{tk}\Vert^2 \Vert D_{tn}{\bf 1}\Vert^2+
\Big(M\vert\xi_{3k}(s_k)\vert^2-\vert M\xi_{3k}(s_k)\vert^2\Big)\Vert D_{3n}{\bf 1}\Vert^2\\
&&=\Big(\frac{1}{2b_{1k}}+a_{1k}^2\Big)\frac{b_{1n}}{2}+
\Big(\frac{1}{2b_{2k}}+a_{2k}^2\Big)\frac{b_{2n}}{2}+\Big(M\vert\xi_{3k}(s_k)\vert^2-\vert M\xi_{3k}(s_k)\vert^2\Big)
\frac{b_{3n}}{2}
\end{eqnarray*}
\begin{eqnarray*}
&&\sim\Big(\frac{1}{2b_{1k}}+a_{1k}^2\Big)+\Big(\frac{1}{2b_{2k}}+a_{2k}^2\Big)+\Big(M \vert\xi_{3k}(s_k)\vert^2-\vert M\xi_{3k}(s_k)\vert^2\Big)\\
&&=\lambda_k^{(3)}(s_k)+(g_k,g_k).
\end{eqnarray*}
Set for brevity $\xi(s)=ix\exp\big(is(x-a)\big)$, then 
\begin{eqnarray}
\nonumber
 &&
 f(s):=M\vert \xi(s)\vert^2-\vert M\xi(s)\vert^2=\Big(\frac{1}{2b}+a^2\Big)-\Big(\frac{s^2}{4b^2}+a^2\Big)\exp\Big(-\frac{s^2}{2b}\Big),\\
 \label{min-s_k}
 &&
 \min_{s}f(s)=
  \left\{\begin{array}{ll}
\frac{1-e^{-(1-2ba^2)}}{2b}+a^2,&\,\text{if}\,\,1-2ba^2\geq 0\\
\frac{1}{2b},&\,\text{if}\,\,1-2ba^2<0
\end{array}\right..
\end{eqnarray}
Indeed, consider the function
$g_{a,b}(x)=\Big(\frac{x^2}{2b}+a^2\Big)\exp(-x^2)$ and set $x_0=\sqrt{1-2ba^2}$.
We have
\begin{equation}
\label{max-g(c,d)(x)}
\max_{x\in \mathbb R}g_{a,b}(x)\!=\!
 \left\{\begin{array}{ll}
g_{a,b}(x_0)= \frac{\exp(-(1-2ba^2))}{2b},
 &\text{if}\,\,1-2ba^2\geq 0,\\
g_{a,b}(0)= a^2,&\text{if}\,\,1-2ba^2<0.
\end{array}\right.
\end{equation}

To calculate $\lambda_k^{(3)}(s_k)$, we get finally
\begin{eqnarray*}
&&\lambda_k^{(3)}(s_k)\!=\!\frac{1}{2b_{1k}}\!+\!a_{1k}^2\!+\!\frac{1}{2b_{2k}}\!+\!a_{2k}^2\!+\!\Big(M\vert\xi_{3k}(s_k)\vert^2\!-\!\vert M\xi_{3k}(s_k)\vert^2\Big)\!-\!(g_k,g_k)\\
&&=\frac{1}{2b_{1k}}+\frac{1}{2b_{2k}}+\frac{1}{2b_{3k}}+a_{1k}^2+a_{2k}^2+a_{3k}^2-
\Big(\frac{s_k^2}{4b^2_{3k}}+a_{3k}^2\Big)\exp\big(-\frac{s_k^2}{2b_{3k}}\big)\\
&&-(g_k,g_k)=\frac{1}{2b_{1k}}+\frac{1}{2b_{2k}}+\frac{1}{2b_{3k}}+a_{3k}^2-
\Big(\frac{s_k^2}{4b^2_{3k}}+a_{3k}^2\Big)\exp\big(-
\frac{s_k^2}{2b_{3k}}\big)=\\
&&
\frac{1}{2b_{1k}}+\frac{1}{2b_{2k}}+\Big(\frac{1}{2b_{3k}}+a_{3k}^2\Big)\big(1-e^{-1}\big),\quad  \text{for}\quad 
\frac{s_k^2}{2b_{3k}}=1.
\end{eqnarray*}
Therefore, we get \eqref{lambda(3)(s)} and \eqref{lambda(3)(s).min}
$$
\lambda_k^{(3)}(s_k)=\frac{1}{2b_{1k}}+\frac{1}{2b_{2k}}+
\Big(\frac{1}{2b_{3k}}+a_{3k}^2\Big)\big(1-\frac{1}{e}\big).
$$
 For $A_{2m+1}=\big((f_k(s_k),f_r(s_r))\big)_{k,r=-m}^m,$ and
$b=(M\xi_{3k}(s_k))_{k=-m}^m\in {\mathbb R}^{2m+1}$ we have
$$
A_{2m+1}= \sum_{k=-m}^m\lambda_kE_{kk}+\gamma(g_{-m},\dots,g_0,\dots,g_m).
$$
The proof is now finished on invoking Lemma~\ref{l.m=3}.
\qed\end{pf}
\begin{lem}
\label{l.Re-Im-exp.3}
We have
\begin{eqnarray}
\label{Re-exp}
&&
D_{3k}{\bf 1}\in \langle 
\sin\big(s_k(x_{3k}-a_{3k})\big)A_{kn}{\bf 1}\mid k\in {\mathbb Z}\rangle
\,\,\Leftrightarrow\,\, 
\Sigma_3(D,s)=\infty,\\
\label{Im-exp}
&&D_{3k}{\bf 1}\in \langle 
\cos\big(s_k(x_{3k}-a_{3k})\big)A_{kn}{\bf 1}\mid k\in {\mathbb Z}\rangle
\,\,\Leftrightarrow\,\, 
\Sigma^{\vee}_3(D,s)=\infty,\\
\label{D(3.sin-cos)}
&&\text{where}\,\,\Sigma_3(D,s)\!=\!\sum_{k\in \mathbb Z}\frac{\vert M\eta_{3k}(s_k)\vert^2}{\Vert g_{k}(s_k)\Vert^2},
\,\, \Sigma^{\vee}_3(D,s)\!=\!\sum_{k\in \mathbb Z}\frac{\vert M\eta^{\vee}_{3k}(s_k)\vert^2}{\Vert g^{\vee}_{k}(s_k)\Vert^2},\qquad
\\
\label{sim-b.3}
&&
\text{moreover,}
\quad \Sigma_3(D,s^{(3)})
\sim \Sigma_3(D):=\sum_k\frac{\frac{1}{2b_{3k}}}{C_k+a^2_{1k}+a^2_{2k}+a^2_{3k}},\\
\label{sim-a.3}
&&
\text{and}\quad \Sigma^{\vee}_3(D,s^{(3)})\sim \Sigma^{\vee}_3(D):=\sum_k\frac{a^2_{3k}}{C_k+a^2_{1k}+a^2_{2k}+a^2_{3k}},
\end{eqnarray}
where $s^{(3)}=(s_{3k})_k$
with $\frac{s^2_{3k}}{b_{3k}}\equiv 1,\,\,k\in \mathbb Z$.
\end{lem}
\begin{pf}
We shall try 
to obtain separately the real part and imaginary part of $ M\xi_{3k}(s)$, where
$\xi_{3k}(s_k)=ix_{3k}\exp\big(is_k(x_{3k}-a_{3k})\big)$.
Using Lemma~\ref{exp(x)A3}
formulas \eqref{M(x.exp(x))} and \eqref{H(a,b,s)}  we get
\begin{eqnarray*}
\nonumber
&&
H_{a,b}(s)=\int_{\mathbb R}ix\exp\big(is(x-a)\big)d\mu_{(b,a)}(x)=\int_{\mathbb R}i(x+a)\exp(isx)d\mu_{(b,0)}(x)\\
&&=\frac{dF_b(s)}{ds}+iaF_b(s)=\Big(-\frac{s}{2b}+ia\Big)\exp\Big(-\frac{s^2}{4b}\Big)=:M\xi_{3k}(s),
\end{eqnarray*}
Recall the Euler formulas
\begin{eqnarray}
 \label{Euler}
&&
e^{it}=\cos t+i\sin t,\quad  e^{-it}=\cos t-i\sin t,\\
\nonumber
&&
\cos t=\frac{e^{it}+e^{-it}}{2},\quad
\sin t=\frac{e^{it}-e^{-it}}{2i}.
\end{eqnarray}

More precisely, we denote for $1\leq r\leq 3$
\begin{equation}
\label{eta-re.3}
\eta_{rk}(s)=x_{rk}\cos\big(s_k(x_{3k}-a_{3k})\big)
,\quad
\eta^{\vee}_{rk}(s)=x_{rk}\cos\big(s_k(x_{3k}-a_{3k})\big).
\end{equation}
We determine when the inclusion holds:
\begin{eqnarray*}
 &&
D_{3k}{\bf 1}\in \langle 
\sin\big(s_k(x_{3k}-a_{3k})\big)A_{kn}{\bf 1}
=\Big(x_{1k}
\sin\big(s_k(x_{3k}-a_{3k})\big)D_{1n}\\
&& +x_{2k}\sin\big(s_k(x_{3k}-a_{3k})\big)D_{2n}+x_{3k}\sin\big(s_k(x_{3k}-a_{3k})\big)D_{3n}\Big){\bf 1}\mid k\in {\mathbb Z}\rangle,\\
 &&
D_{3k}{\bf 1}\in \langle 
\cos\big(s_k(x_{3k}-a_{3k})\big)A_{kn}{\bf 1}
=\Big(x_{1k}
\cos\big(s_k(x_{3k}-a_{3k})\big)D_{1n}\\
&& +x_{2k}\cos\big(s_k(x_{3k}-a_{3k})\big)D_{2n}+x_{3k}\cos\big(s_k(x_{3k}-a_{3k})\big)D_{3n}\Big){\bf 1}\mid k\in {\mathbb Z}\rangle.
\end{eqnarray*} 
 Set 
\begin{eqnarray}
\label{g_k(s)}
&&
g_k(s_k)\!=\!\Big(\eta_{1k}(s_k)D_{1n}+\eta_{2k}(s_k)D_{2n}\!+\!\big[\eta_{3k}(s_k)\!-\!M\eta_{3k}(s_k)\big]D_{3n}\Big){\bf 1},\qquad \\
\label{g^v_k(s)}
&&
g^{\vee}_k(s_k)\!=\!\Big(\eta^{\vee}_{1k}(s_k)D_{1n}\!+\!\eta^{\vee}_{2k}(s_k)D_{2n}\!+\!\big[\eta^{\vee}_{3k}(s_k)\!-\!M\eta^{\vee}_{3k}(s_k)\big]D_{3n}\Big){\bf 1},
\end{eqnarray}
We show that (compare with \eqref{M(x.exp(x))}) 
\begin{eqnarray}
 \label{M-eta.3}
&&M\eta_{3k}(s)=
-\frac{1}{2}\Big(H_{a,b}(s)+\overline{H_{a,b}(s)}\Big)=
\frac{s}{2b_{3k}}\exp\Big(-\frac{s^2}{4b_{3k}}\Big),\\
 \label{M-eta^.3}
&&M\eta^{\vee}_{3k}(s)=
\frac{1}{2i}\Big(H_{a,b}(s)-\overline{H_{a,b}(s)}\Big)=
a_{3k}
\exp\Big(-\frac{s^2}{4b_{3k}}\Big).
\end{eqnarray}
Recall the definition of the  function $F_b(s)$
defined by \eqref{F_b(a)}:
\begin{equation}
\label{F_b(a).1}
F_b(s) =\int_{\mathbb R}\exp\big(is(x-a)\big)d\mu_{(b,a)}(x)=
\int_{\mathbb R}\exp\big(isx\big)d\mu_{(b,0)}(x)=\exp\Big(-\frac{s^2}{4b}\Big).
\end{equation}
We have 
\begin{eqnarray*}
&&
M\eta(s)=
\int_{\mathbb R}x\sin\big(s(x-a)\big)d\mu_{(b,a)}(x)=
\int_{\mathbb R}(x+a)\sin(sx)
d\mu_{(b,0)}(x)=\\
&&
\int_{\mathbb R}(x+a)\frac{e^{isx}-e^{-isx}}{2i}d\mu_{(b,0)}(x)=
-\frac{1}{2}\int_{\mathbb R}i(x+a)\big(e^{isx}-e^{-isx}\big)d\mu_{(b,0)}(x)=\\
&&
-\frac{1}{2}\Big(H_{a,b}(s)+\overline{H_{a,b}(s)}\Big)=\frac{s}{2b}\exp\Big(-\frac{s^2}{4b}\Big),
\end{eqnarray*}
this implies \eqref{M-eta.3}.
Similarly we get 
\begin{eqnarray*}
&&
M\eta^{\vee}(s)=
\int_{\mathbb R}x\cos\big(s(x-a)\big)d\mu_{(b,a)}(x)=
\int_{\mathbb R}(x+a)\cos(sx)
d\mu_{(b,0)}(x)=\\
&&
\frac{1}{2i}\int_{\mathbb R}i(x+a)\Big(e^{isx}+e^{-isx}\Big)d\mu_{(b,0)}(x)=
\frac{1}{2i}\Big(H_{a,b}(s)-\overline{H_{a,b}(s)}\Big)\\
&&
=a\exp\Big(-\frac{s^2}{4b}\Big)
\end{eqnarray*}
this implies \eqref{M-eta.3}.
Fix $m\in {\mathbb N}$, 
we put $\sum_{k=-m}^mt_kM\eta_{3k}(s_k)=(t,b)=1$,
where $t=(t_k)_{k=-m}^m$ and $b=(M\eta_{3k}(s_k))_{k=-m}^m$. We have
\begin{eqnarray}
\nonumber
&& \Vert \big[\sum_{k=-m}^mt_k
\sin\big(s_k(x_{3k}-a_{3k})\big)A_{kn}
-D_{3n}\big]{\bf 1}
\Vert^2\\
\nonumber
&=&
\Vert 
\sum_{k=-m}^mt_k
\Big(\eta_{1k}(s_k)D_{1n}+\eta_{2k}(s_k)D_{2n}+\big[\eta_{3k}(s_k)-M\eta_{3k}(s_k)\big]D_{3n}\Big){\bf 1}
\Vert^2
\\
\label{eta.3}
&&
=\sum_{k=-m}^mt^2_k\Vert g_{k}(s_k)\Vert^2,\,\,\,\text{since}\,\,\,
\Big(D_{rn}{\bf 1},D_{ln}{\bf 1}\Big)\!=\!0,\,\,\,1\leq r<l\leq 3,
\end{eqnarray} 
where the $g_k(s_k)$ are defined by \eqref{g_k(s)}.
To calculate $\Vert g_{k}(s_k)\Vert^2$ we have
\begin{eqnarray}
\nonumber
&&
\Vert g_{k}(s_k)\Vert^2=(g_{k}(s_k),g_{k}(s_k))=\\
\nonumber
&&\Big(\big(\eta_{1k}(s_k)D_{1n}+\eta_{2k}(s_k)D_{2n}+\big[\eta_{3k}(s_k)-M\eta_{3k}(s_k)\big]D_{3n}\big){\bf 1},\\
\nonumber
&&\big(\eta_{1k}(s_k)D_{1n}+\eta_{2k}(s_k)D_{2n}+\big[\eta_{3k}(s_k)-M\eta_{3k}(s_k)\big]D_{3n}\big){\bf 1}\Big)=\\
\nonumber
&&
\Vert x_{1k}{\bf 1}\Vert^2
\Vert 
\sin\big(s_k(x_{3k}-a_{3k})\big){\bf 1}\Vert^2
\Vert D_{1k}{\bf 1}\Vert^2+\\
\nonumber
&&
\Vert x_{2k}{\bf 1}\Vert^2
\Vert \sin\big(s_k(x_{3k}-a_{3k})\big){\bf 1}\Vert^2
\Vert D_{2k}{\bf 1}\Vert^2+\\
\nonumber
&&
 \Big(M\vert\eta_{kn}(s_k)\vert^2-\vert M\eta_{kn}(s_k)\vert^2\Big)
\Vert D_{3k}{\bf 1}\Vert^2=
\Big(\frac{1}{2b_{1k}}+a^2_{1k}\Big)I_3\frac{b_{1n}}{2}+
\\
\label{norm(g_k)}
&&
\Big(\frac{1}{2b_{2k}}+a^2_{2k}\Big)I_3\frac{b_{2n}}{2}+
 \Big(M\vert\eta_{kn}(s_k)\vert^2-\vert M\eta_{kn}(s_k)\vert^2\Big)\frac{b_{3n}}{2}.
\end{eqnarray}
We need to calculate $I_3=\Vert \sin\big(s_k(x_{3k}-a_{3k})\big){\bf 1}\Vert^2,\,\,M\vert\eta_{kn}(s_k)\vert^2$ and $\vert M\eta_{kn}(s_k)\vert^2$.
If we  set $a:=a_{3k},\,\,b:=b_{3k}$, we get
\begin{eqnarray}
\nonumber 
&&I_3=\Vert \sin\big(s_k(x_{3k}-a_{3k})\big){\bf 1}\Vert^2=
\int_{\mathbb R}
\frac{e^{isx}-e^{-isx}}{2i}\frac{e^{-isx}-e^{isx}}{-2i}d\mu_{(b,0)}(x)=
\\
\label{I_3}
&&
\frac{1}{2}\int_{\mathbb R}
\Big(1-\frac{e^{2isx}+e^{-2isx}}{2}\Big)
d\mu_{(b,0)}(x)\stackrel{\eqref{F_b(a).1}}{=}
\frac{1-e^{-\frac{s^2}{b}}}{2},
\\
\label{M-eta-3}
&&\vert M\eta_{kn}(s_k)\vert^2=\frac{s^2_k}{4b^2_{3k}}\exp\Big(-\frac{s_k^2}{2b_{3k}}\Big),\\
\nonumber 
&&M\vert\eta_{kn}(s_k)\vert^2=
\int_{\mathbb R}(x^2+2xa+a^2)
\frac{e^{isx}-e^{-isx}}{2i}\frac{e^{-isx}-e^{isx}}{-2i}d\mu_{(b,0)}(x)=\\
\nonumber
&&
\frac{1}{2}\int_{\mathbb R}(x^2+2xa+a^2)
\Big(1-\frac{e^{2isx}+e^{-2isx}}{2}\Big)
d\mu_{(b,0)}(x)=
\end{eqnarray}
\begin{eqnarray}
\nonumber
&&
\frac{1}{2}\Big[
\int_{\mathbb R}(x^2+a^2)d\mu_{(b,0)}(x)-
\int_{\mathbb R}(x^2+a^2)\frac{e^{2isx}+e^{-2isx}}{2}
d\mu_{(b,0)}(x)\Big]=\\
\nonumber
&&
\frac{1}{2}\Big[\frac{1}{2b}+a^2-\frac{d^2F_b(2s)}{ds^2}-a^2F_b(2s)\Big]\stackrel{\eqref{F_b(a).1}}{=}\frac{1}{2}\Big[\frac{1}{2b}+a^2-\frac{1}{(2i)^2}\Big[\Big(\frac{2s}{b}\Big)^2-\frac{2}{b}\Big]\times\\
\label{M-eta-fin}
&&
e^{-\frac{s^2}{b}}-a^2e^{-\frac{s^2}{b}}\Big]=
\frac{1}{2}\Big[\Big(\frac{1}{2b}+a^2\Big)(1-e^{-\frac{s^2}{b}})+
\frac{s^2}{b^2}e^{-\frac{s^2}{b}}\Big].
\end{eqnarray}
Finally, we get
\begin{equation}
\label{eta-eta}
M\vert\eta_{kn}(s_k)\vert^2- \vert M\eta_{kn}(s_k)\vert^2=
\frac{1}{2}\Big[\Big(\frac{1}{2b}+a^2\Big)(1-e^{-\frac{s^2}{b}})+
\frac{s^2}{b^2}e^{-\frac{s^2}{b}}\Big]-
\frac{s^2}{4b^2}e^{-\frac{s^2}{2b}}.
\end{equation}
By \eqref{eta^v.3.2},\eqref{norm(g_k)},  \eqref{I_3}, \eqref{eta-eta} and \eqref{A.min2} we prove \eqref{Re-exp}, where 
\begin{eqnarray}
\nonumber
&&
\Sigma_3(D,s)=\sum_{k\in \mathbb Z}\frac{\vert M\eta_{kn}(s_k)\vert^2}{\Vert g_{k}(s_k)\Vert^2}=\\
\nonumber
&&
\sum_{k\in \mathbb Z}\frac{
\frac{s^2_k}{4b^2_{3k}}
e^{-\frac{s_k^2}{2b_{3k}}}
}
{\Big(\frac{1}{2b_{1k}}+a^2_{1k}\Big)I_3\frac{b_{1n}}{2}+
\Big(\frac{1}{2b_{2k}}+a^2_{2k}\Big)I_3\frac{b_{2n}}{2}+
 \Big(M\vert\eta_{kn}(s_k)\vert^2-\vert M\eta_{kn}(s_k)\vert^2\Big)\frac{b_{3n}}{2}}\sim\\
 \nonumber
 &&
\sum_{k\in \mathbb Z}\frac{
\frac{s^2_k}{4b^2_{3k}}%
e^{-\frac{s_k^2}{2b_{3k}}}
}
{
\frac{1-e^{-\frac{s^2_k}{b_{3k}}}}{2}\Big(c_{1k}+c_{2k}\Big)+
\frac{1}{2}\Big[
c_{3k}(1-e^{-\frac{s^2_k}{b_{3k}}})+
\frac{s^2_k}{b_{3k}^2}e^{-\frac{s^2_k}{b_{3k}}}\Big]-
\frac{s^2_k}{4b_{3k}^2}e^{-\frac{s^2_k}{2b_{3k}}}
 }=
\end{eqnarray} 
\begin{equation}
 \label{D(3,x)}   
\sum_{k\in \mathbb Z}\frac{
\frac{x^2_k}{4b_{3k}}%
e^{-\frac{x_k^2}{2}}
}
{
\frac{1-e^{-x^2_k}
}{2}\Big(c_{1k}+c_{2k}\Big)\!+\!
\frac{1}{2}\Big[
c_{3k}(1\!-\!e^{-x^2_k})+
\frac{x^2_k}{b_{3k}}e^{-x_k^2
}\Big]\!-\!
\frac{x_k^2}{4b_{3k}}e^{-\frac{x^2_k}{2}}
 }=:\Sigma_3(D,x),
\end{equation}
where $x^2_k=\frac{s^2_k}{b_{3k}}$ and $c_{rk}=\frac{1}{2b_{rk}}+a^2_{rk}$. For $x^{(3)}=(x_k)_k$
with $x_k\equiv 1$ we get
\begin{eqnarray}
\nonumber
&&
\Sigma_3(D,x^{(3)})=
\sum_{k\in \mathbb Z}\frac{
\frac{1}{4b_{3k}}%
e^{-\frac{1}{2}}
}
{
\frac{1-e^{-1}
}{2}\Big(c_{1k}+c_{2k}\Big)+
\frac{1}{2}\Big[
c_{3k}(1-e^{-1})+
\frac{1}{b_{3k}}e^{-1
}\Big]-
\frac{1}{4b_{3k}}e^{-\frac{1}{2}}}=\\
\nonumber
&&
\sum_{k\in \mathbb Z}\frac{
\frac{1}{4b_{3k}}%
e^{-\frac{1}{2}}
}
{
\frac{1-e^{-1}
}{2}\Big(c_{1k}+c_{2k}+c_{3k}\Big)+
\frac{1}{2b_{3k}}\Big(e^{-1}-\frac{1}{2}e^{-\frac{1}{2}}\Big)}\stackrel{\eqref{r.two-sim-ser}}{\sim}\\
\label{D_3}
&&
\sum_{k\in \mathbb Z}\frac{\frac{1}{2b_{3k}}}
{c_{1k}+c_{2k}+c_{3k}}\!=\!
\sum_{k\in \mathbb Z}\frac{\frac{1}{2b_{3k}}}
{\frac{1}{2b_{1k}}+\frac{1}{2b_{2k}}+\frac{1}{2b_{3k}}+
a^2_{1k}+a^2_{2k}+a^2_{3k}}=\Sigma_3(D).
\end{eqnarray}
So, we have proved \eqref{sim-b.3}
for $x=(x_k)_k$
with $x_k\equiv 1$. 
To approximate $D_{3n}$ 
in terms of functions involving  the cosine,
fix $m\in {\mathbb N}$, 
and put $\sum_{k=-m}^mt_kM\eta^{\vee}_{3k}(s_k)=(t,b)=1$,
where $t=(t_k)_{k=-m}^m$ and $b=(M\eta^{\vee}_{3k}(s_k))_{k=-m}^m$. We have
\begin{eqnarray}
\nonumber
&& \Vert \big[\sum_{k=-m}^mt_k
\cos\big(s_k(x_{3k}-a_{3k})\big)A_{kn}
-D_{3n}\big]{\bf 1}
\Vert^2\\
\nonumber
&=&
\Vert 
\sum_{k=-m}^mt_k
\Big(\eta^{\vee}_{1k}(s_k)D_{1n}+\eta^{\vee}_{2k}(s_k)D_{2n}+\big[\eta^{\vee}_{3k}(s_k)-M\eta^{\vee}_{3k}(s_k)\big]D_{3n}\Big){\bf 1}
\Vert^2
\\
\label{eta^v.3.1}
&&
=\sum_{k=-m}^mt^2_k\Vert g^{\vee}_{k}(s_k)\Vert^2,\,\,\,\text{since}\,\,\,
\Big(D_{rn}{\bf 1},D_{ln}{\bf 1}\Big)\!=\!0,\,\,\,1\leq r<l\leq 3,
\end{eqnarray} 
where the $g^{\vee}_k(s_k)$ are defined by \eqref{g_k(s)}.
To calculate $\Vert g^{\vee}_{k}(s_k)\Vert^2$ we have
\begin{eqnarray}
\nonumber
&&
\Vert g^{\vee}_{k}(s_k)\Vert^2=(g^{\vee}_{k}(s_k),g^{\vee}_{k}(s_k))=\\
\nonumber
&&\Big(\big(\eta^{\vee}_{1k}(s_k)D_{1n}+\eta^{\vee}_{2k}(s_k)D_{2n}+\big[\eta^{\vee}_{3k}(s_k)-M\eta^{\vee}_{3k}(s_k)\big]D_{3n}\big){\bf 1},\\
\nonumber
&&\big(\eta^{\vee}_{1k}(s_k)D_{1n}+\eta^{\vee}_{2k}(s_k)D_{2n}+\big[\eta^{\vee}_{3k}(s_k)-M\eta^{\vee}_{3k}(s_k)\big]D_{3n}\big){\bf 1}\Big)=\\
\nonumber
&&
\Vert x_{1k}{\bf 1}\Vert^2
\Vert 
\cos\big(s_k(x_{3k}-a_{3k})\big){\bf 1}\Vert^2
\Vert D_{1k}{\bf 1}\Vert^2+\\
\nonumber
&&
\Vert x_{2k}{\bf 1}\Vert^2
\Vert \cos\big(s_k(x_{3k}-a_{3k})\big){\bf 1}\Vert^2
\Vert D_{2k}{\bf 1}\Vert^2+
\end{eqnarray}
\begin{eqnarray}
\nonumber
&&
 \Big(M\vert\eta^{\vee}_{kn}(s_k)\vert^2-\vert M\eta^{\vee}_{kn}(s_k)\vert^2\Big)
\Vert D_{3k}{\bf 1}\Vert^2=
\Big(\frac{1}{2b_{1k}}+a^2_{1k}\Big)I^{\vee}_3\frac{b_{1n}}{2}+
\\
\label{norm(g^v_k).1}
&&
\Big(\frac{1}{2b_{2k}}+a^2_{2k}\Big)I^{\vee}_3\frac{b_{2n}}{2}+
 \Big(M\vert\eta^{\vee}_{kn}(s_k)\vert^2-\vert M\eta^{\vee}_{kn}(s_k)\vert^2\Big)\frac{b_{3n}}{2}.
\end{eqnarray}
We need to calculate $I^{\vee}_3=\Vert \cos\big(is_k(x_{3k}-a_{3k})\big){\bf 1}\Vert^2,\,\,M\vert\eta^{\vee}_{kn}(s_k)\vert^2$ and $\vert M\eta^{\vee}_{kn}(s_k)\vert^2$.
Finally, we get (we set $b:=b_{3k}$)
To approximate $D_{3n}$ 
in terms of functions involving  $\cos$,
fix $m\in {\mathbb N}$, 
and put $\sum_{k=-m}^mt_kM\eta^{\vee}_{3k}(s_k)=(t,b)=1$,
where $t=(t_k)_{k=-m}^m$ and $b=(M\eta^{\vee}_{3k}(s_k))_{k=-m}^m$. We have
\begin{eqnarray}
\nonumber
&& \Vert \big[\sum_{k=-m}^mt_k
\cos\big(s_k(x_{3k}-a_{3k})\big)A_{kn}
-D_{3n}\big]{\bf 1}
\Vert^2\\
\nonumber
&=&
\Vert 
\sum_{k=-m}^mt_k
\Big(\eta^{\vee}_{1k}(s_k)D_{1n}+\eta^{\vee}_{2k}(s_k)D_{2n}+\big[\eta^{\vee}_{3k}(s_k)-M\eta^{\vee}_{3k}(s_k)\big]D_{3n}\Big){\bf 1}
\Vert^2
\\
\label{eta^v.3.2}
&&
=\sum_{k=-m}^mt^2_k\Vert g^{\vee}_{k}(s_k)\Vert^2,\,\,\,\text{since}\,\,\,
\Big(D_{rn}{\bf 1},D_{ln}{\bf 1}\Big)\!=\!0,\,\,\,1\leq r<l\leq 3,
\end{eqnarray} 
where the $g^{\vee}_k(s_k)$ are defined by \eqref{g_k(s)}.
To calculate $\Vert g^{\vee}_{k}(s_k)\Vert^2$ we have
\begin{eqnarray}
\nonumber
&&
\Vert g^{\vee}_{k}(s_k)\Vert^2=(g^{\vee}_{k}(s_k),g^{\vee}_{k}(s_k))=\\
\nonumber
&&\Big(\big(\eta^{\vee}_{1k}(s_k)D_{1n}+\eta^{\vee}_{2k}(s_k)D_{2n}+\big[\eta^{\vee}_{3k}(s_k)-M\eta^{\vee}_{3k}(s_k)\big]D_{3n}\big){\bf 1},\\
\nonumber
&&\big(\eta^{\vee}_{1k}(s_k)D_{1n}+\eta^{\vee}_{2k}(s_k)D_{2n}+\big[\eta^{\vee}_{3k}(s_k)-M\eta^{\vee}_{3k}(s_k)\big]D_{3n}\big){\bf 1}\Big)=\\
\nonumber
&&
\Vert x_{1k}{\bf 1}\Vert^2
\Vert 
\cos\big(s_k(x_{3k}-a_{3k})\big){\bf 1}\Vert^2
\Vert D_{1k}{\bf 1}\Vert^2+\\
\nonumber
&&
\Vert x_{2k}{\bf 1}\Vert^2
\Vert \cos\big(s_k(x_{3k}-a_{3k})\big){\bf 1}\Vert^2
\Vert D_{2k}{\bf 1}\Vert^2+\\
\nonumber
&&
 \Big(M\vert\eta^{\vee}_{kn}(s_k)\vert^2-\vert M\eta^{\vee}_{kn}(s_k)\vert^2\Big)
\Vert D_{3k}{\bf 1}\Vert^2=
\Big(\frac{1}{2b_{1k}}+a^2_{1k}\Big)I^{\vee}_3\frac{b_{1n}}{2}+
\\
\label{norm(g^v_k).2}
&&
\Big(\frac{1}{2b_{2k}}+a^2_{2k}\Big)I^{\vee}_3\frac{b_{2n}}{2}+
 \Big(M\vert\eta^{\vee}_{kn}(s_k)\vert^2-\vert M\eta^{\vee}_{kn}(s_k)\vert^2\Big)\frac{b_{3n}}{2}.
\end{eqnarray}
We need to calculate $I^{\vee}_3=\Vert \cos\big(s_k(x_{3k}-a_{3k})\big){\bf 1}\Vert^2,\,\,M\vert\eta^{\vee}_{kn}(s_k)\vert^2$ and $\vert M\eta^{\vee}_{kn}(s_k)\vert^2$.
If we  set $b:=b_{3k}$, we get
\begin{eqnarray}
\nonumber 
&&I^{\vee}_3=\Vert \cos\big(s_k(x_{3k}-a_{3k})\big){\bf 1}\Vert^2=
\int_{\mathbb R}
\frac{e^{isx}+e^{-isx}}{2}\frac{e^{-isx}+e^{isx}}{2}d\mu_{(b,0)}(x)=
\\
\label{I_3-vee}
&&
\frac{1}{2}\int_{\mathbb R}
\Big(1+\frac{e^{2isx}+e^{-2isx}}{2}\Big)
d\mu_{(b,0)}(x)\stackrel{\eqref{F_b(a).1}}{=}
\frac{1+e^{-\frac{s^2}{b}}}{2},
\\
\label{M-eta-vee}
&&\vert M\eta^{\vee}_{kn}(s_k)\vert^2=a^2_{3k}\exp\Big(-\frac{s_k^2}{2b_{3k}}\Big),\\
\nonumber 
&&M\vert\eta^{\vee}_{kn}(s_k)\vert^2=
\int_{\mathbb R}(x^2+2xa+a^2)
\frac{e^{isx}+e^{-isx}}{2}\frac{e^{-isx}+e^{isx}}{2}d\mu_{(b,0)}(x)=
\end{eqnarray}
\begin{eqnarray}
\nonumber
&&
\frac{1}{2}\int_{\mathbb R}(x^2+2xa+a^2)
\Big(1+\frac{e^{2isx}+e^{-2isx}}{2}\Big)
d\mu_{(b,0)}(x)=\\
\nonumber
&&
\frac{1}{2}\Big[
\int_{\mathbb R}(x^2+a^2)d\mu_{(b,0)}(x)+
\int_{\mathbb R}(x^2+a^2)\frac{e^{2isx}+e^{-2isx}}{2}
d\mu_{(b,0)}(x)\Big]=\\
\nonumber
&&
\frac{1}{2}\Big[\frac{1}{2b}+a^2+\frac{d^2F_b(2s)}{ds^2}+a^2F_b(2s)\Big]\stackrel{\eqref{F_b(a).1}}{=}\frac{1}{2}\Big[\frac{1}{2b}+a^2+\frac{1}{(2i)^2}\Big[\Big(\frac{2s}{b}\Big)^2-\frac{2}{b}\Big]\times\\
\label{M-eta-fin-v}
&&
e^{-\frac{s^2}{b}}+a^2e^{-\frac{s^2}{b}}\Big]=
\frac{1}{2}\Big[\Big(\frac{1}{2b}+a^2\Big)(1+e^{-\frac{s^2}{b}})-
\frac{s^2}{b^2}e^{-\frac{s^2}{b}}\Big].
\end{eqnarray}
Finally, we get
\begin{equation}
\label{eta-eta-v}
M\vert\eta^{\vee}_{kn}(s_k)\vert^2\!- \vert M\eta^{\vee}_{kn}(s_k)\vert^2\!=\!
\frac{1}{2}\Big[\Big(\frac{1}{2b}+a^2\Big)(1+e^{-\frac{s^2}{b}})-
\frac{s^2}{b^2}e^{-\frac{s^2}{b}}\Big]-
{\color{blue}a^2}
e^{-\frac{s^2}{2b}}.
\end{equation}
%
%
By \eqref{eta^v.3.1}, \eqref{norm(g^v_k).1},  \eqref{I_3-vee}, \eqref{eta-eta-v} and \eqref{A.min2} we prove \eqref{Im-exp}, where 
\begin{eqnarray}
\nonumber
&&
\Sigma^{\vee}_3(D,s)=\sum_{k\in \mathbb Z}\frac{\vert M\eta^{\vee}_{kn}(s_k)\vert^2}{\Vert g^{\vee}_{k}(s_k)\Vert^2}=\\
\nonumber
&&
\sum_{k\in \mathbb Z}\frac{
a^2_{3k}
e^{-\frac{s_k^2}{2b_{3k}}}
}
{\Big(\frac{1}{2b_{1k}}+a^2_{1k}\Big)I^{\vee}_3\frac{b_{1n}}{2}+
\Big(\frac{1}{2b_{2k}}+a^2_{2k}\Big)I^{\vee}_3\frac{b_{2n}}{2}+
 \Big(M\vert\eta^{\vee}_{kn}(s_k)\vert^2-\vert M\eta^{\vee}_{kn}(s_k)\vert^2\Big)\frac{b_{3n}}{2}}\\
 \nonumber
 &&
\sim\sum_{k\in \mathbb Z}\frac{
a^2_{3k}
e^{-\frac{s_k^2}{2b_{3k}}}
}
{
\frac{1+e^{-\frac{s^2_k}{b_{3k}}}}{2}\Big(c_{1k}+c_{2k}\Big)+
\frac{1}{2}\Big[
c_{3k}(1+e^{-\frac{s^2_k}{b_{3k}}})-
\frac{s^2_k}{b_{3k}^2}e^{-\frac{s^2_k}{b_{3k}}}\Big]-
a^2_{3k}
e^{-\frac{s^2_k}{2b_{3k}}}
 }=
\end{eqnarray}
\begin{equation} 
\label{D^v(3.x)}  
\sum_{k\in \mathbb Z}\frac{
a^2_{3k}
e^{-\frac{x_k^2}{2}}
}
{
\frac{1-e^{-x^2_k}
}{2}\Big(c_{1k}\!+\!c_{2k}\Big)\!+\!
\frac{1}{2}\Big[
c_{3k}(1+e^{-x^2_k})-
\frac{x^2_k}{b_{3k}}e^{-x_k^2
}\Big]\!-\!
a^2_{3k}
e^{-\frac{x^2_k}{2}}
 }\!=\!\Sigma_3(D,x),\quad
\end{equation}
where $x^2_k=\frac{s^2_k}{b_{3k}}$ and $c_{rk}=\frac{1}{2b_{rk}}+a^2_{rk}$. For $x^{(3)}=(x_k)_k$
with $x_k\equiv 1$ we get
\begin{eqnarray}
\nonumber 
&&
\Sigma^{\vee}_3(D,x^{(3)})=
\sum_{k\in \mathbb Z}\frac{
a^2_{3k}
e^{-\frac{1}{2}}
}
{
\frac{1+e^{-1}
}{2}\Big(c_{1k}+c_{2k}\Big)+
\frac{1}{2}\Big[
c_{3k}(1+e^{-1})-
\frac{1}{b_{3k}}e^{-1
}\Big]-
a^2_{3k}
e^{-\frac{1}{2}}}=\\
\nonumber
&&
\sum_{k\in \mathbb Z}\frac{
a^2_{3k}
e^{-\frac{1}{2}}
}
{
\frac{1+e^{-1}
}{2}\Big(c_{1k}+c_{2k}+c_{3k}\Big)-
\Big(
\frac{1}{2b_{3k}}e^{-1}+a^2_{3k}e^{-\frac{1}{2}}\Big)
}
\stackrel{\eqref{r.two-sim-ser}}{\sim}\\
\label{D^v_3} 
&&
\sum_{k\in \mathbb Z}\frac{
a^2_{3k}
}
{c_{1k}+c_{2k}+c_{3k}}\!=\!
\sum_{k\in \mathbb Z}\frac{
a^2_{3k}
}
{C_k+
a^2_{1k}+a^2_{2k}+a^2_{3k}}=
\Sigma^{\vee}_3(D).
\end{eqnarray}
So, we have proved \eqref{sim-a.3}
for $x=(x_k)_k$
with $x_k\equiv 1$. 
\qed\end{pf}

\section{How far is a vector from a hyperplane?}
\subsection{Some estimates}
\label{s.1.4.1}
We recall some material from \cite{Kos_B_09}, Section 1.4.1, pp.\,24--25.
\begin{lem} [\cite{Kos04}]
\label{1.l.min}
 For a strictly positive operator $A$ {\rm (}i.e., $(Af,f)>0$ for $f\not=0${\rm )} acting in ${\mathbb R}^n$ and a vector
$b\in{\mathbb R}^n\backslash\{0\}$, we have
\begin{equation}
\label{A.min3} \min_{x\in{\mathbb R}^n}\Big((Ax,x)\mid
(x,b)=1\Big)=
\frac{1}{(A^{-1}b,b)}.
\end{equation}
The minimum is assumed  for
$x=\frac{A^{-1}b}{(A^{-1}b,b)}.
$
\end{lem}
Lemma \ref{1.l.min} is a direct generalization of the well known
result (see, for example, \cite{BecBel61}, Chap. I, \S 52), 
stating that
for $a_k>0,\,\,1\leq k \leq n$ we have
\begin{equation}
\label{A.min1}
 \min_{x\in{\mathbb
R}^n}\Big(\sum_{k=1}^na_kx_k^2\mid \sum_{k=1}^n x_k=1\Big)=
\Big(\sum_{k=1}^n \frac{1}{a_k}\Big)^{-1}.
\end{equation}
We will also use  the same result   in a slightly different
form:
\begin{equation}
\label{A.min2}
 \min_{x\in{\mathbb
R}^n}\Big(\sum_{k=1}^na_kx_k^2\mid \sum_{k=1}^n x_kb_k=1\Big)=
\Big(\sum_{k=1}^n \frac{b_k^2}{a_k}\Big)^{-1},
\end{equation}
with the minimum being assumed 
for
$x_k=\frac{b_k}{a_k}\Big(\sum_{k=1}^n\frac{b_k^2}{a_k}\Big)^{-1}$.
\subsection{The distance of a vector from a hyperplane}
We follow closely the exposition  \cite{Kos-hpl-arx20}. We start with a classical result, see, e.g. 
\cite{Gan58}.
Consider the hyperplane $V_n$ generated by $n$ arbitrary vectors $f_1,\dots, f_{n}$ in some Hilbert space  $H$. 
\begin{lem} 
\label{l.d(f,v_n)}
The square of the distance $d(f_0,V_n)$ of a vector $f_0$ from the hyperplane $V_n$ is given by the ratio of two Gram determinants: 
\begin{equation}
\label{d(f,v_n)}
d^2(f_0,V_n)=
\frac{\Gamma(f_0,f_1,f_2,\dots, f_n)}{\Gamma(f_1,f_2,\dots, f_n)}.
\end{equation}
\end{lem}


\subsection{Gram determinants and Gram matrices}
\begin{df}
\label{d.Gram-det}
Let us recall the definition of
	the Gram  determinant and the Gram  matrix   (see
	\cite{Gan58}, Chap IX, \S 5). Given the  vectors $x_1,x_2,..., x_m$ in some
	Hilbert space $H$ the {\it Gram
		matrix} $\gamma(x_1,x_2,..., x_m)$ is defined  by the formula
$$
\gamma(x_1,x_2,..., x_m)=\big((x_k,x_n)\big)_{k,n=1}^m.
$$
The determinant of this matrix is called the {\it Gram  determinant} for
the vectors $x_1,x_2,..., x_m$ and is denoted by
$\Gamma(x_1,x_2,..., x_m)$.
Thus,
\begin{equation}
	\label{Gram-det}
	\Gamma(x_1,x_2,\dots, x_m):={\rm det}\,\gamma(x_1,x_2,\dots, x_m).
	\end{equation}
\end{df}
\subsection{The generalized characteristic polynomial
and its properties}
\label{sec.gen.har.pol}
\index{polynomial!characteristic!generalized}
{\bf Notations}. For a matrix $C\in {\rm Mat}(n,\mathbb R)$  and $1\leq i_1<i_2<\dots i_r\leq n,$ 
$1\leq j_1<j_2<\dots j_r\leq n,\,\,r\leq n$
denote by 
$$
M^{i_1i_2\dots i_r}_{j_1j_2\dots j_r}(C)\quad\text{ and} \quad
A^{i_1i_2\dots i_r}_{j_1j_2\dots j_r}(C)
$$
the corresponding {\it minors} and {\it cofactors} of the matrix $C$.
\begin{df}{\rm (\cite[Ch.1.4.3]{Kos_B_09})}
\label{d.G_k(lambda)} 
For the matrix $C\in {\rm
Mat}(m,{\mathbb C})$ and $\lambda =(\lambda_1,\dots,\lambda_m)\in {\mathbb C}^m$ define  the {\it
generalization of the characteristic polynomial}, $p_C(t)={\rm
det}\,(tI-C),\, t\in {\mathbb C}$  
as follows:
\begin{equation}
\label{P_C(lambda)}
P_C(\lambda)={\rm det}\,C(\lambda),\quad\text{where}\quad
C(\lambda)=
{\rm diag}(\lambda_1,\dots,\lambda_m)
+C.
\end{equation}
\end{df}
\begin{lem}
\label{l.detC-LI}
{\rm (\cite[Ch.1.4.3]{Kos_B_09})}
For the generalized characteristic polynomial $P_C(\lambda)$
of
$C\!\in\!{\rm Mat}(m,{\mathbb C})$ and
 $\lambda=(\lambda_1,\lambda_2,...,\lambda_m)\in {\mathbb
C}^m$ we have 
\begin{equation}
\label{detC-LI}
P_C(\lambda)=
{\rm det}\,C+
\sum_{r=1}^m\sum_{1\leq i_1<i_2<...<i_r\leq
m}\lambda_{i_1}\lambda_{i_2}...\lambda_{i_r}A^{i_1i_2...i_r}_{i_1i_2...i_r}(C).
\end{equation}
\end{lem}
\begin{rem}
\label{r.P_C(lam)}
If we set
$\lambda_\alpha=\lambda_{i_1}\lambda_{i_2}...\lambda_{i_r}$, where
$\alpha=(i_1,i_2,...,i_r)$ and
$A^\alpha_\alpha(C)=A^{i_1i_2...i_r}_{i_1i_2...i_r}(C),\,\,
M^\alpha_\alpha(C)=M^{i_1i_2...i_r}_{i_1i_2...i_r}(C),\,\,
\lambda_\emptyset=1,\,\,A^\emptyset_\emptyset(C)={\rm det}\,C$ we
may write (\ref{detC-LI}) as follows:
\begin{equation}
\label{A.detC-LI.2} 
P_C(\lambda)={\rm det}\,C(\lambda)=
\sum_{\emptyset\subseteq\alpha\subseteq\{1,2,...,m\}}\lambda_\alpha
A^\alpha_\alpha(C),
\end{equation}
\begin{equation}
\label{M.detC-LI.2}
P_C(\lambda)={\rm det}\,C(\lambda)=
\left(\prod_{k=1}^n\lambda_k\right)\sum_{\emptyset\subseteq\alpha\subseteq\{1,2,...,m\}}
\frac{M^\alpha_\alpha(C)}{\lambda_\alpha},
\end{equation}
\end{rem}

Let 
\begin{equation}
\label{X(mn)}
X=X_{mn} =\left(
\begin{array}{cccc}
x_{11}&x_{12}&...&x_{1n}\\
x_{21}&x_{22}&...&x_{2n}\\
...   &...   &...&...\\
x_{m1}&x_{m2}&...&x_{mn}
\end{array}
\right).
\end{equation}
Setting 
\begin{equation}
\label{x_k,y_r=}
x_k=(x_{1k},x_{2k},...,x_{mk})\in{\mathbb R}^m,
\quad
y_r=(x_{r1},x_{r2},...,x_{rn})\in{\mathbb R}^n,
\end{equation}
we get
\begin{equation}
\label{X^*X}
X^*X=
 \left(
\begin{array}{cccc}
(x_1,x_1)&(x_1,x_2)&...&(x_1,x_n)\\
(x_2,x_1)&(x_2,x_2)&...&(x_2,x_n)\\
...   &...   &...&...\\
(x_n,x_1)&(x_n,x_2)&...&(x_n,x_n)
\end{array}
\right)=\gamma(x_1,x_2,..., x_n),
\end{equation}
\begin{equation}
\label{XX^*}
XX^*= \left(
\begin{array}{cccc}
(y_1,y_1)&(y_1,y_2)&...&(y_1,y_m)\\
(y_2,y_1)&(y_2,y_2)&...&(y_2,y_m)\\
...   &...   &...&...\\
(y_m,y_1)&(y_m,y_2)&...&(y_m,y_m)
\end{array}
\right)=\gamma(y_1,y_2,..., y_m),
\end{equation}
therefore, we obtain
\begin{equation}
 \label{d(X^*X)=d(XX^*)}\
\Gamma(x_1,x_2,..., x_n)={\rm det}(X^*X)={\rm det}(XX^*)=\Gamma(y_1,y_2,..., y_m).
\end{equation}

\section{Explicit expressions for $C^{-1}(\lambda)$ and  $(C^{-1}(\lambda)a,a)$}
In this section we follow \cite{Kos-hpl-arx20}.
Fix $C\in {\rm Mat}(n,\mathbb R),$ $a\in \mathbb R^n$ and $\lambda\in \mathbb C^n$. 
Our aim is to find the explicit formulas for  
$C^{-1}(\lambda)$ and  $(C^{-1}(\lambda)a,a)$, where $C(\lambda)$ is defined by \eqref{P_C(lambda)}.  
Set $M(i_1i_2\dots i_r)(C)$
$=M^{i_1i_2\dots i_r}_{i_1i_2\dots i_r}(C)$ and 
$a_{i_1i_2\dots i_r}=(a_{i_1},a_{i_2},\dots,a_{i_r})$.
Let also $C_{i_1i_2\dots i_r}$ be the corresponding {\it submatrix} of the matrix $C$.
The elements of this matrix are on the intersection of $i_1,i_2,\dots, i_r$ rows and column of the matrix $C$.
Denote by $A(C_{i_1i_2\dots i_r})$ the matrix of the cofactors of the first order of the matrix $C_{i_1i_2\dots i_r}$, i.e. 
\begin{equation}
\label{A(C)}
A(C_{i_1i_2\dots i_r})=(A^{i}_j(C_{i_1i_2\dots i_r}))_{1\leq i,j\leq r} 
\end{equation}

Let $n=3$, then $A(C_{123})=A(C)$ is the following matrix: 
\begin{equation}
 \label{A(C)3}
A(C)=A(C_{123})=
\left(
\begin{array}{ccc}
A^1_1&A^1_2&A^1_3\\
A^2_1&A^2_2&A^2_3\\
A^3_1&A^3_2&A^3_3
\end{array}
\right)=
\left(
\begin{array}{ccc}
M^{23}_{23}&-M^{23}_{13}&M^{23}_{12}\\
-M^{13}_{23}&M^{13}_{13}&-M^{13}_{12}\\
M^{12}_{23}&-M^{12}_{13}&M^{12}_{12}
\end{array}
\right),
\end{equation}
where we write $M^{ij}_{rs}$ instead of $M^{ij}_{rs}(C)$ and $A^i_j$ instead of $A^i_j(C)$.
\begin{rem}
\label{r.C_{12}=}
Let $A^T$ be the transposed matrix
of $A$. Then
\begin{equation}
 \label{A(C_{12})=}
A^T(C_{i_1i_2\dots i_r})={\rm det}\,C^{-1}_{i_1i_2\dots i_r}\Big(C^{-1}_{i_1i_2\dots i_r}\Big),
\end{equation}

\end{rem}
In what follows we will  consider the submatrix $C_{i_1i_2\dots i_r}$ of the matrix $C\in {\rm Mat}(n,\mathbb R)$ as an appropriate element of ${\rm Mat}(n,\mathbb R)$.
\begin{thm}
\label{t.(C^{-1}a,a)}
 For the matrix $C(\lambda)$ defined by  \eqref{P_C(lambda)} 
 $a\in \mathbb R^n$  and $\lambda\in \mathbb C^n$ we have 
 {\small
   \begin{equation}
 \label{Del_n(lam,C)}
 P_C(\lambda)=
 \Big(\prod_{k=1}^n\lambda_k\Big)\sum_{r=1}^n
\sum_{1\leq i_1<i_2<\dots <i_r\leq n}
\frac{M(i_1i_2\dots i_r)}{\lambda_{i_1}\lambda_{i_2}\dots \lambda_{ i_r}},
\end{equation}
 \begin{equation}
 \label{C^{-1}(lam)}
C^{-1}(\lambda)=\frac{1}{P_C(\lambda)}\Big(\prod_{k=1}^n\lambda_k\Big)
\sum_{r=1}^n
\sum_{1\leq i_1<i_2<\dots i_r\leq n}
\frac{A(C_{i_1i_2\dots i_r})}{\lambda_{i_1}\lambda_{i_2}\dots \lambda_{ i_r}},
 \end{equation} 
 \begin{equation}
 \label{(C^{-1}(lam)a,a)}
(C^{-1}(\lambda)a,a)=\frac{1}{P_C(\lambda)}\Big(\prod_{k=1}^n\lambda_k\Big)
\sum_{r=1}^n
\sum_{1\leq i_1<i_2<\dots i_r\leq n}
\frac{(A(C_{i_1i_2\dots i_r})a_{i_1i_2\dots i_r},a_{i_1i_2\dots i_r})}{\lambda_{i_1}\lambda_{i_2}\dots \lambda_{ i_r}}.
 \end{equation}
 }
\end{thm}

\subsection{The case where $C$ is the Gram matrix}
Fix the matrix $X_{mn}$ defined by \eqref{X(mn)}. Denote by $C$ the Gram matrix 
$\gamma(x_1,x_2,..., x_n)$, 
i.e.,
\begin{equation}
\label{C=gamma}
C=\gamma(x_1,x_2,..., x_n),
\end{equation}
where $(x_1,x_2,..., x_n)$ are defined by \eqref{x_k,y_r=} 
and $\gamma(x_1,x_2,..., x_n)$  by \eqref{X^*X}.
In what follows we consider the operator $C(\lambda)$  defined  
by \eqref{P_C(lambda)}.
\begin{rem}
 \label{r.P_C(lam).7}
In this case  we have
\begin{eqnarray}
 \label{G.detC-LI.2} 
P_C(\lambda)&=&
{\rm det}\Big(\sum_{k=1}^n\lambda_kE_{kk}+\gamma(x_1,x_2,...,
x_n)\Big)\\
\nonumber
&=&\prod_{k=1}^n\lambda_k\Big(1+\sum_{r=1}^n\sum_{1\leq
	i_1<i_2<...<i_r\leq
	m}\Big(\lambda_{i_1}\lambda_{i_2}...\lambda_{i_r}\Big)^{-1}
\Gamma(x_{i_1},x_{i_2},...,x_{i_r})\Big)\\
\nonumber
&=&\prod_{k=1}^n\lambda_k\Big(1+\sum_{r=1}^{n} \sum_{
\substack{1\leq i_1<i_2<...<i_r\leq n;\\
1\leq j_1<j_2<...<j_r\leq n}
}
\Big(\lambda_{i_1}\lambda_{i_2}...\lambda_{i_r}\Big)^{-1}
\Big(M^{i_1i_2...i_r}_{j_1j_2...j_r}(X)\Big)^2\Big),
\end{eqnarray}
where we have used the following formula (see \cite{Gan58}, Chap IX, \S
5 formula (25)):
\begin{equation}
\label{Gramm(x,y)=M^2(X)}
\Gamma(x_{i_1},x_{i_2},...,x_{i_r})= \sum_{1\leq
	j_1<j_2<...<j_r\leq m}
\left(M^{i_1i_2...i_r}_{j_1j_2...j_r}(X)\right)^2.
\end{equation}
\end{rem}

\subsection{Case $m=2$}

Fix two natural numbers $n,m\in \mathbb N$ with $m\leq n$, two matrices $A_{mn}$ and $X_{mn}$, vectors $g_k\in\mathbb R^{m-1},\,\,1\leq k\leq n$
and $a\in \mathbb R^n$ as follows
\begin{equation} 
\label{A(mn)}
A_{mn}\!=\!\left(
	\begin{array}{cccc}
	a_{11}&a_{12}          &...&a_{1n}\\
	a_{21}&a_{22}          &...&a_{2n}\\
	&&...&\\
	a_{m1}&a_{m2}          &...&a_{mn}
	\end{array}
	\right),\,\, 
g_k=\left(
	\begin{array}{c}
	a_{2k}\\
	a_{3k}\\
	... \\
	a_{mk}
	\end{array}
	\right)\in \mathbb R^{m-1},\,\,\, a=(a_{1k})_{k=1}^n\in \mathbb R^n.	
\end{equation}
Set $C=\gamma(g_1,g_2,\dots,g_n)$. We calculate $\Delta_n(\lambda, C)$ and $(C^{-1}(\lambda)a,a)$ for an arbitrary $n$. 
Consider the matrix \eqref{X(mn)}
\begin{equation}
\label{X(mn).1}
X_{mn} =\left(
\begin{array}{cccc}
x_{11}&x_{12}&...&x_{1n}\\
x_{21}&x_{22}&...&x_{2n}\\
&&...&\\
x_{m1}&x_{m2}&...&x_{mn}
\end{array}
\right),\quad \text{where}\quad x_{kr}=\frac{a_{1k}}{\sqrt \lambda_k},\quad y_r=(x_{rk})_{k=1}^n\in \mathbb R^n.
\end{equation}
\begin{lem}[\cite{Kos-hpl-arx20}, Lemma~2.2]
\label{l.m=2}
For $m=2$ we have
\begin{equation}
\label{m=2}
(C^{-1}(\lambda)a,a)=\Delta(y_1,y_2)=\frac{\Gamma(y_1)+\Gamma(y_1,y_2)}{1+\Gamma(y_2)},
\end{equation}
where $y_1$ and $y_2$ are defined as follows
\begin{equation}
\label{y(1),y(2)=} 
y_1=\left(\frac{a_{1k}}{\sqrt{\lambda_k}}\right)_{k=1}^n\quad y_2=\left(\frac{a_{2k}}{\sqrt{\lambda_k}}\right)_{k=1}^n.
\end{equation} 
\end{lem}

\subsection{Case $m=3$}
By \eqref{G.detC-LI.2} we get
$$
\Delta(y_1,y_2,y_3)=\frac{\Gamma(y_1)+\Gamma(y_1,y_2)+\Gamma(y_1,y_3)+\Gamma(y_1,y_2,y_3)}{1+\Gamma(y_2)+\Gamma(y_3)+\Gamma(y_2,y_3)}.
$$
\begin{lem}
\label{l.m=3}
For $m=3$ we have
\begin{equation}
\label{m=3}
(C^{-1}(\lambda)a,a)=\Delta(y_1,y_2,y_3),
\end{equation}
where the $y_r$ are defined as follows:
\begin{equation}
\label{f,g,h=} 
y_r=
\left(\frac{a_{rk}}{\sqrt{\lambda_k}}\right)_{k=1}^n\in \mathbb R^n,\quad 1\leq r\leq 3.
\end{equation} 
\end{lem}
\section{Appendix}
\subsection{Comparison of two Gaussian measures}
\label{sec.Gauss-meas}
For two centered Gaussian measures $\mu_{(b,0)}$ and $\mu_{(b',0)}$ on the real line $\mathbb R$ defined by \eqref{mu(b,a)}  it is well known that  
\begin{equation}
\label{Hel-int-R}
H(\mu_{(b,0)},\mu_{(b',0)})=\left(\frac{4bb'}{(b+b')^2}\right)^{1/4}.
\end{equation}
By Kakutani's criterion for product measures on $\mathbb R^{\mathbb N}$ {\rm \cite {Kak48}}, and \eqref{Hel-int-R} we see that the following lemma holds true.
\begin{lem}
\label{l.eq-prod.m=1}
 Two Gaussian  measures $\mu_{(b,0)}=\otimes_{n\in \mathbb Z}\mu_{(b,0)}$ and $\mu_{(b',0)}=\otimes_{n\in \mathbb Z}\mu_{(b',0)}$ are equivalent if and only if the product
\begin{equation}
\label{eq-mu-prod.m=1} 
\prod_{n\in \mathbb Z}\frac{4b_{n}b'_{n}}{(b_{n}+b'_{n})^2}
\end{equation}
does not converge to $0$. The equivalent condition is 
\begin{equation}
 \label{eq-mu-sum.m=1} 
\sum_{n\in \mathbb Z}
\left(\sqrt{\frac{b_{n}}{b'_{n}}}- \sqrt{\frac{b'_{n}}{b_{n}}}\right)^2<\infty.
\end{equation}
\end{lem}
 Consider  
two measures: $\mu_{({\mathbb I},0)}=\otimes_{n\in \mathbb Z}\mu_{(1,0)}$ 
and  $\mu_{({\mathbb I}+c,0)}=\otimes_{n\in \mathbb Z}\mu_{(1+c_n,0)}$ on the space $X_1$,  where the measure $\mu_{(b,a)}$ on the real line $\mathbb R$ is defined by \eqref{mu(b,a)}.
\begin{lem}
\label{l.eq-prod} 
Two measures  $\mu_{({\mathbb I},0)}$ and  $\mu_{({\mathbb  I}+c,0)}$ are equivalent if and only if 
\begin{equation}
\label{eq-prod} 
\sum_{n\in \mathbb Z}c_n^2<\infty
\end{equation}
\end{lem}
\begin{pf}
By Lemma~\ref{l.eq-prod.m=1} and 
\eqref{eq-mu-sum.m=1},  the measures $\mu_{({\mathbb I},0)}$ and $\mu_{({\mathbb I}+c,0)}$
are equivalent if and only if 
\begin{equation*}
\sum_{n\in \mathbb Z}
\left(\frac{1}{\sqrt{1+c_n}}-\sqrt{1+c_n}\right)^2=\sum_{n\in \mathbb Z}
\frac{c_n^2}{1+c_n}
<\infty.
\end{equation*}
By Lemma~\ref{l.two-ser-alpha},
 two series $\sum_{n\in \mathbb Z}
\frac{c_n^2}{1+c_n}$ and $\sum_{n\in \mathbb Z}c_n^2$ are equivalent. 
\qed\end{pf}

The next lemma is also a consequence of Kakutani's criterion {\rm \cite {Kak48}}.
\begin{lem}
 \label{l.eq-prod.m=3}
 Two Gaussian  measures $\mu_{(b,0)}^m$ and $\mu_{(b',0)}^m$ are equivalent if and only if the product
\begin{equation}
\label{eq-mu-prod.m=3} 
\prod_{r=1}^m\prod_{n\in \mathbb Z}\frac{4b_{rn}b'_{rn}}{(b_{rn}+b'_{rn})^2}
\end{equation}
does not converge to $0$. The equivalent condition is 
\begin{equation}
 \label{eq-mu-sum.m=3} 
\sum_{r=1}^m\sum_{n\in \mathbb Z}
\left(\sqrt{\frac{b_{rn}}{b'_{rn}}}- \sqrt{\frac{b'_{rn}}{b_{rn}}}\right)^2<\infty.
\end{equation}
\end{lem}

Lemma \ref{perp1} follows from Lemmas~\ref{perp3} -- \ref{ldet*}.
\begin{lem}
\label{perp3} For $t\!\in\!{\rm GL}(m,{\mathbb R})\backslash\{e\}$ we have
$(\mu_{(b,a)}^m)^{L_t}\!\perp\!\mu_{(b,a)}^m$ if and only if
 \begin{equation}
  \label{7.1}
(\mu_{(b,0)}^m)^{L_t}\perp\mu_{(b,0)}^m\quad\text{or}\quad
\mu_{(b,{L_t}a)}^m\perp\mu_{(b,a)}^m.
 \end{equation}
\end{lem}
Let us define the following measures on the spaces ${\mathbb R}^m$
and $X_m$:
$$
\mu_m^{(B_n,0)}=\otimes_{k=1}^m\mu_{(b_{kn},0)},\quad
\mu_m^{(B_n,a_n)}=\otimes_{k=1}^m\mu_{(b_{kn},a_{kn})},
$$
where $a_n\!=\!(a_{1n},...,a_{mn})\in{\mathbb R}^m$ and
$B_n\!=\!{\rm diag}(b_{1n},...,b_{mn})\in {\rm Mat}(m,{\mathbb R})$. Since
$$
\mu_{(b,a)}^m=\otimes_{n\in{\mathbb Z}}\mu_m^{(B_n,a_n)},\quad
\mu_{(b,0)}^m=\otimes_{n\in{\mathbb Z}}\mu_m^{(B_n,0)},
$$
$$
\left(\mu_{(b,a)}^m\right)^{L_t}=\otimes_{n\in{\mathbb
Z}}\left(\mu_m^{(B_n,a_n)}\right)^{L_t},\quad
\left(\mu_{(b,0)}^m\right)^{L_t}=\otimes_{n\in{\mathbb
Z}}\left(\mu_m^{(B_n,0)}\right)^{L_t},
$$
and
$$
\mu_{(b,{L_t}a)}^m=\otimes_{n\in{\mathbb
Z}}\mu_m^{(B_n,{L_t}a_n)},
$$
by the Kakutani criterion \cite{Kak48}, we derive the following  two lemmas:

\begin{lem}
\label{l7.2}
 For the measures $\mu_{(b,0)}^m,\,m\in{\mathbb N}$ and $t\in
 {\rm GL}(m,{\mathbb R})\backslash\{e\}$, we obtain
$$
(\mu_{(b,0)}^m)^{L_t}\perp\mu_{(b,0)}^m\,\, \Leftrightarrow
\prod_{n\in{\mathbb
Z}}H\left(\left(\mu_m^{(B_n,0)}\right)^{L_t},\mu_m^{(B_n,0)}\right)=0.
$$
\end{lem}
\begin{lem}
\label{l7.3}
For the  measures $\mu_{(b,0)}^m,\,m\in{\mathbb N}$ and $t\in
 {\rm GL}(m,{\mathbb R})\backslash\{e\}$, we get
$$
\mu_{(b,{L_t}a)}^m\perp\mu_{(b,a)}^m \Leftrightarrow
\prod_{n\in{\mathbb
Z}}H\left(\mu_m^{(B_n,{L_t}a_n)},\mu_m^{(B_n,a_n)}\right)=0.
$$
\end{lem}
To prove Lemma \ref{perp1} it is sufficient to show, by Lemma
\ref{perp3}, that
\begin{equation}
\label{Hel1-sl}
H_{m,n}(t)\!=\!H\Big(\Big(\mu_m^{(B_n,0)}\Big)^{L_t}\!\!\!\!,\mu_m^{(B_n,0)}\Big)\!=\!
\Big(\frac{1}{2^m\vert{\rm det}\,\,t\vert} {\rm
det}\left(I\!+\!X_n^*(t)X_n(t)\right)\Big)^{-1/2}\!\!\!\!\!,\quad
\end{equation}
to  prove the  equivalence
\begin{equation}
\label{Hel2-sl}
 \prod_{n\in{\mathbb
Z}}H\left(\mu_m^{(B_n,{L_t}a_n)},\mu_m^{(B_n,a_n)}\right)=0
\Leftrightarrow \sum_{n\in {\mathbb Z}}\sum_{r=1}^m
b_{rn}\Big(\sum_{s=1}^m(t_{rs}-\delta_{rs})a_{sn}\Big)^2 =\infty,
\end{equation}
and to  apply  the following lemma.
\begin{lem}
\label{ldet*}
 For $X\in {\rm Mat}(m,{\mathbb R})$ we have
\begin{equation}
\label{det*} {\rm det}\left(I+X^*X\right)= 1+\sum_{r=1}^{m}
\sum_{1\leq i_1<i_2<...<i_r\leq m;1\leq j_1<j_2<...<j_r\leq m}
\left(M^{i_1i_2...i_r}_{j_1j_2...j_r}(X)\right)^2.
\end{equation}
\end{lem}

The proof of  the equivalence (\ref{Hel2-sl}) is based on the following theorem
that one can find, e.g., in \cite[ Ch.\,III, \S 16, Theorem 2]{Sko74}.
\begin{thm}
\label{t.(B,a)sim(B,0)}
 Two Gaussian measures $\mu_{B,a}$ and $\mu_{B,b}$ in a Hilbert space $H$ are
equivalent if and only of $B^{-1/2}(a-b)\in H$.
\end{thm}

{\small Indeed, we have
$$
\Vert C^{-1/2}(ta-a)\Vert^2_H\!=\!\!\sum_{n\in {\mathbb Z}}
\Vert C_n^{-1/2}(t-I)a_n\Vert^2_{H_n}\!=\!2\!\sum_{n\in
{\mathbb Z}}\sum_{r=1}^m
\frac{b_{kn}}{d_{kn}}\Big(\sum_{s=1}^m(t_{rs}\!-\delta_{rs})a_{sn}\Big)^2\!d_{kn}.
$$
To explain the latter equality let us describe $H$ and $C$. To
find an operator $C$ we present the measure $\mu^m_{(b,a)}$ in the
canonical form $\mu_{C,a}$ defined by its Fourier transform:
\begin{equation}
\label{A.F(H,mu_B),}
 \int_{H}\exp i(y,x)d\mu_{C,a}(x)=\exp\left(i(a,y)-\frac{1}{2}(Cy,y)\right),\,\,
y\in H,
\end{equation}
where $C$ is a positive {\it nuclear operator} (called the {\it
covariance operator}) on the Hilbert space $H$, and $a\in H$ is
the {\it mathematical expectation} or {\it mean}.
}

Recall the {\it Kolmogorov zero-one law}. Let us consider in the space ${\mathbb R}^\infty={\mathbb
R}\times {\mathbb R}\times\cdots$ the infinite tensor product  $\mu_b=\otimes_{n\in
{\mathbb N}}\mu_{b_k}$  of one-dimensional Gaussian measures $\mu_{b_k}$ on $\mathbb R$ defined as follows:
\begin{equation}
\label{A.mu-(b,a)}
d\mu_{b}(x)=\sqrt{\frac{b}{\pi}}\exp(-bx^2)dx.
\end{equation}
Consider a Hilbert space $l_2(a)$ defined by
\begin{equation*}
l_2(a)=\big\{x\in {\mathbb R}^\infty\,:\, \Vert x
\Vert^2_{l_2(a)}=\sum_{k\in {\mathbb N}}x_k^2a_k<\infty\big\},
\end{equation*}
where $a=(a_k)_{k\in {\mathbb N}}$ is an infinite sequence of
positive numbers.
\begin{thm}[Kolmogorov's zero-one law, \cite{ShFDT67}]
\label{A.Kol-0-1}
We have
$$
\mu_b(l_2(a))= \left\{\begin{array}{ccc} 0,&\text{if}&\sum_{k\in
{\mathbb N}}\frac{a_k}{b_k}=\infty;\\
1,&&\text{otherwise}.
\end{array}\right.
$$
\end{thm}
\subsection{Properties of  two vectors $f,g\not\in l_2$}
In what follows we will use systematically the following notation. 
For $k$ vectors $f_1,f_2,\dots,f_k\in \mathbb R^n$ with $k\leq n$ we set
\begin{equation}
\label{Delta(f(1),...,f(k)).1}
\Delta(f_1,f_2,\dots,f_k)=
\frac{{\rm det}(I+\gamma(f_1,f_2,\dots,f_k))}{{\rm det}(I+\gamma(f_2,\dots,f_k))}-1.
\end{equation}
For $k=2$ and $k=3$ we get respectively:
\begin{equation}
\label{Delta(f,g)}
\Delta(f_1,f_2)=\frac{{\rm det}(I+\gamma(f_1,f_2))}{{\rm det}(I+\gamma(f_2))}-1
=\frac{I+\Gamma(f_1)+\Gamma(f_2)+\Gamma(f_1,f_2)}{I+\Gamma(f_2)}.
\end{equation}
\begin{eqnarray}
 \label{Delta(f,g,h)}
&&\Delta(f_1,f_2,f_3)=\frac{{\rm det}(I+\gamma(f_1,f_2,f_3))}{{\rm det}(I+\gamma(f_2,f_3))}-1
=\\
\nonumber
&&\frac{\Gamma(f_1)+\Gamma(f_1,f_2)+\Gamma(f_1,f_3)+\Gamma(f_1,f_2,f_3)}{1+\Gamma(f_2)+\Gamma(f_3)+\Gamma(f_2,f_3)}.
\end{eqnarray}
\begin{lem} [\cite{KosJFA17}, and  \cite{Kos_B_09}, Ch.10]
\label{l.min=proj0}
Let $f=(f_k)_{k\in{\mathbb N}}$ and $g=(g_k)_{k\in{\mathbb N}}$ be two real vectors  such that
$\Vert f\Vert^2=\infty$,  where $\Vert f\Vert^2=\sum_kf_k^2$.  Denote by $f_{(n)}$,
$g_{(n)}\in {\mathbb R}^n$ their projections to the subspace
${\mathbb R}^n$, i.e., $f_{(n)}=(f_k)_{k=1}^n,\quad  g_{(n)}=(g_k)_{k=1}^n$. Then
\begin{equation}
\label{Delta-to-infty}
\Delta(f,g):=
\lim_{n\to\infty}\Delta(f_{(n)},g_{(n)})=\infty,
\end{equation}
where $\Delta(f_{(n)},g_{(n)})$ is defined by \eqref{Delta(f,g)}, 
in the following cases:
\begin{eqnarray*}
(a)&\Vert g\Vert^2<\infty,\\
(b)&\Vert g\Vert^2=\infty,\quad\text{and}\quad  \lim_{n\to\infty}\frac{\Vert f_{(n)}\Vert}{\Vert g_{(n)}\Vert}=\infty,\\
(c)&\Vert f\Vert^2=\Vert g\Vert^2=\Vert f+s g\Vert^2=\infty,\quad\text{for all}\quad
     s\in {\mathbb R}\setminus\{0\}.
\end{eqnarray*}
\end{lem}
\begin{pf} Obviously $\lim_{n\to\infty}\Delta(f_{(n)},g_{(n)})=\infty$ if conditions (a)  or (b) hold.
The implication $ (c)\Rightarrow (\ref{Delta-to-infty})$ is based on the following lemma.
\qed\end{pf}
\begin{lem}
\label{l.min=proj}
Let $f=(f_k)_{k\in{\mathbb N}}$ and $g=(g_k)_{k\in{\mathbb N}}$ be two real vectors  such that
\begin{equation}
\label{norm=infty}
\Vert f\Vert^2=\Vert g\Vert^2=\Vert C_1 f+C_2 g\Vert^2=\infty\quad\text{for all}\quad
(C_1,C_2)\in {\mathbb R}^2\setminus\{0\},
\end{equation}
\begin{equation}
\label{final.2}
\text{then}\quad
\lim_{n\to\infty}\frac{\Gamma(f_{(n)},g_{(n)})}{\Gamma(g_{(n)})}=\infty\quad\text{and}\quad
\lim_{n\to\infty}\frac{\Gamma(f_{(n)},g_{(n)})}{\Gamma(f_{(n)})}=\infty.
\end{equation}
\end{lem}
\begin{lem}
\label{l.opposite-l_2.2}
Let $f_1,f_2\not\in l_2$ and 
$\Delta(f_1,f_2)<\infty$, then for some
$(C_1,C_2)\in {\mathbb R}^2\setminus\{0\}$ we have
$C_1f_1+C_2f_2\in l_2$.
\end{lem}
\begin{pf}
Let assume the opposite, i.e., $C_1f_1+C_2f_2\not\in l_2$ for all $(C_1,C_2)\in {\mathbb R}^2\setminus\{0\}$. 
Then by Lemma~\ref{l.min=proj}
\begin{equation*}
\hskip 1.3cm
\Delta(f_1,f_2)= \frac{\Gamma(f_1)+\Gamma(f_1,f_2)}{1+\Gamma(f_2)}>
\frac{\Gamma(f_1,f_2)}
{1+\Gamma(f_2)}\sim\frac{\Gamma(f_1,f_2)}
{\Gamma(f_2)}=\infty.\hskip 1.3cm \Box
\end{equation*}
\end{pf}

\begin{lem}
\label{l.opposite-l_2.3}
Let $f_1,f_2,f_3\not\in l_2$ and 
$\Delta(f_1,f_2,f_3)<\infty$, then for some
$(C_1,C_2,C_3)\in {\mathbb R}^3\setminus\{0\}$ we have
$C_1f_1+C_2f_2+C_3f_3\in l_2$.
\end{lem}
\begin{pf}
Let assume the opposite, i.e., that  $C_1f_1+C_2f_2+C_3f_3\not\in l_2$ for all $(C_1,C_2,C_3)\in {\mathbb R}^3\setminus\{0\}$.
Then by  Lemmas~\ref{l.min=proj.3}--\ref{l.min=proj}
and \eqref{Delta(f,g,h)}  we get
 \begin{eqnarray*}
&&
\Delta(f_1,f_2,f_3)=
\frac{\Gamma(f_1)+\Gamma(f_1,f_2)+\Gamma(f_1,f_3)+\Gamma(f_1,f_2,f_3)}{1+\Gamma(f_2)+\Gamma(f_3)+\Gamma(f_2,f_3)}>\\
&&
\hskip 1.9cm 
\frac{\Gamma(f_1,f_2,f_3)}{1+\Gamma(f_2)+\Gamma(f_3)+\Gamma(f_2,f_3)}\sim
\frac{\Gamma(f_1,f_2,f_3)}{
\Gamma(f_2,f_3)}
=\infty.\hskip 1.9cm 
\Box
 \end{eqnarray*}
\end{pf}
%
\subsection{Properties of three vectors $f,g,h\not\in l_2$
}

\label{6.6}
\begin{lem}
\label{l.min=proj0.3}
Let $f=(f_k)_{k\in{\mathbb N}},\,\,g=(g_k)_{k\in{\mathbb N}}$ and $h=(h_k)_{k\in{\mathbb N}}$ be three real vectors  such that
$\Vert f\Vert^2=\infty$  where $\Vert f\Vert^2=\sum_kf_k^2$.  Denote by $f_{(n)},\,\,g_{(n)},\,\,h_{(n)}\in {\mathbb R}^n$ 
their projections to the subspace ${\mathbb R}^n$, i.e., $f_{(n)}=(f_k)_{k=1}^n,\,  g_{(n)}=(g_k)_{k=1}^n,\,\,$ 
$ h_{(n)}=(h_k)_{k=1}^n$. Then
\begin{equation}
\label{Delta-to-infty.3}
\Delta(f,g,h):=
\lim_{n\to\infty}\Delta(f_{(n)},g_{(n)},h_{(n)})=\infty,
\end{equation}
where $\Delta(f,g,h)$ is defined by \eqref{Delta(f,g,h)},  in the following cases:
\begin{eqnarray*}
&(a)&\Vert g\Vert^2<\infty\quad\text{and}\quad \Vert h\Vert^2<\infty,\\
&(b)&\Vert g\Vert^2<\infty\quad\text{and}\quad \Vert h\Vert^2=\infty,\,\,\text{or}\,\,
\Vert g\Vert^2=\infty\quad\text{and}\quad \Vert h\Vert^2<\infty,\\
&(c)&
\Vert C_1 f+C_2 g+C_3 h\Vert^2=\infty,\,\,\,\text{for all}\,\,\,(C_1,C_2,C_3)\in {\mathbb R}^3\setminus\{0\}.
\end{eqnarray*}
\end{lem}
\begin{pf} 
(a) In this case $\Vert g_{(n)}\Vert^2\leq C,\,\,\,\Vert h_{(n)}\Vert^2\leq C$ and $\Gamma(g_{(n)},h_{(n)})\leq C$
and therefore, 
$$\lim_{n\to\infty}\Delta(f_{(n)},g_{(n)},h_{(n)})\geq \lim_{n\to\infty}\frac{\Vert f_{(n)}\Vert^2}{1+3C}=\infty.$$ 

(b) Let $\Vert g\Vert^2<\infty,\,\,\Vert h\Vert^2=\infty$.  In this case we have 
$$
\Gamma(g_{(n)},h_{(n)})\leq \Vert g_{(n)}\Vert^2\Vert h_{(n)}\Vert^2\sin^2(\alpha_n)\leq C_1\Gamma(h_{(n)}),
$$
where $\alpha_n$ is the angle betwen two vectors $g_{(n)}$ and $h_{(n)}$.
Therefore,
$$
1+\Gamma(g_{(n)})+\Gamma(h_{(n)})+\Gamma(g_{(n)},h_{(n)})\leq (1+C)
(1+\Gamma(h_{(n)}),
$$
$$
\Delta(f_{(n)},g_{(n)},h_{(n)})\geq \frac{\Gamma(f_{(n)})+\Gamma(f_{(n)},h_{(n)})}
{(1+C)
(1+\Gamma(h_{(n)})
}\sim
\Delta(f_{(n)},h_{(n)}).
$$
So, this case is reduced to the case $m=2$, see Lemma~\ref{l.min=proj}.

 The implication $ (c)\Rightarrow (\ref{Delta-to-infty.3})$ is based on the following lemma
(the proof of which will appear  in \cite{Kos-m-Arx23}. Compare with Lemma~\ref{l.min=proj0}.
\qed\end{pf}

\begin{lem}
[see \cite{Kos-m-Arx23}]
\label{l.min=proj.3}
Let $f_0,f_1,f_2$ be three infinite real vectors\\ $f_r=(f_{rk})_{k\in \mathbb N},$
$0\leq r\leq 2$
such that for all $(C_0,C_1,C_2)\in {\mathbb R}^{3}\setminus\{0\}$  holds
\begin{equation}
\label{norm=infty.3}
\sum_{r=0}^2C_rf_r\not\in l_2,\quad\text{i.e.,} \quad
\sum_{k\in \mathbb N}\vert C_0f_{0k}+C_1f_{1k}+C_2f_{2k}\vert^2=\infty.
\end{equation}
Denote by $f_r(n)=(f_{rk})_{k=1}^n\in \mathbb R^n$ the projections of the vectors $f_r$ on the subspace $\mathbb R^n$. Then
\begin{equation}
\label{final.3}
\frac{\Gamma(f_0,f_1,f_2)}{\Gamma(f_1,f_2)}:=
\lim_{n\to\infty}\frac{\Gamma(f_0(n),f_1(n),f_2(n))}{\Gamma(f_1(n),f_2(n))}
=\infty.
\end{equation}
\end{lem}


{\bf The idea of the proof} (for details see \cite{Kos-m-Arx23}).
We assume there exists a real number $C$ such that 
\begin{equation}
\label{<C.m=2}
\quad\quad\frac{\Gamma(f_0(n),f_1(n),f_2(n))}{\Gamma(f_1(n),f_2(n))}\leq C,
\end{equation}
for every integer $n$ and will show that this leads to a contradiction.

For $t\in \mathbb R^2$ and $f_0,f_1,f_2\in H$ we define the function
\begin{eqnarray*}
F_{20}(t)&=&\Vert\sum_{k=1}^2t_kf_k-f_0\Vert ^2=
\sum_{k,r=1}^2t_kt_r(f_k,f_r)-2\sum_{k=1}^2t_k(f_k,f_0)+(f_0,f_0)\\
\nonumber
&=&(At,t)-2(t,b)+(f_0,f_0),
\end{eqnarray*}
where $b=(f_k,f_0)_{k=1}^2\in \mathbb R^2$ and $A$ is the Gram matrix 
$$A=\gamma(f_1,f_2)=\big((f_k,f_r)\big)_{k,r=1}^2.
$$
Suppose that $At_0=b$, then we have
\begin{equation*}
(At,t)-2(t,b)+(f_0,f_0)
=(A(t-t_0),(t-t_0))+
\frac{\Gamma(f_0,f_1,f_2)}{\Gamma(f_1,f_2)},
\end{equation*}  
and therefore
\begin{eqnarray}
\nonumber
&& 
 F_{20}(t)=
 \big(A(t-t_0),(t-t_0)\big)+
\frac{\Gamma(f_0,f_1,f_2)}{\Gamma(f_1,f_2)},\\
 \label{F_3(t)=(At,t).1}
&& 
 F_{20}(t_0)=\min_{(t_1,t_2)\in \mathbb R^2}\Vert\sum_{k=1}^2t_kf_k-f_0\Vert ^2=
\frac{\Gamma(f_0,f_1,f_2)}{\Gamma(f_1,f_2)},
\end{eqnarray}

Consider the matrix
\begin{equation}
\label{X(mn).3}
X_{3n}=
\left(\begin{array}{cccc}
f_{11}&f_{12}&\dots&f_{1n}\\
f_{21}&f_{22}&\dots&f_{2n}\\
f_{01}&f_{02}&\dots&f_{0n}\\
\end{array}\right)
\end{equation}
and its minors
$$
M_{krs}^{123}=\left|\begin{array}{ccc}
f_{1k}&f_{1r}&f_{1s}\\
f_{2k}&f_{2r}&f_{2s}\\
f_{0k}&f_{0r}&f_{0s}\\
\end{array}\right|,\quad 
M_{kr}^{12}=
\left|\begin{array}{cc}
f_{1k}&f_{1r}\\
f_{2k}&f_{2r}
\end{array}\right|.
$$
Then by 
\cite{Gan58}
we have
\begin{eqnarray*}
&&\Gamma(f_1(n),f_2(n))=\sum_{1\leq r<s\leq n}\vert M_{rs}^{12}\vert^2,\\
 &&\Gamma(f_0(n),f_1(n),f_2(n))=\sum_{1\leq k<r<s\leq n}\vert M_{krs}^{123}\vert^2.
\end{eqnarray*}
Therefore, the inequality 
\eqref{<C.m=2}
will have the following form
\begin{equation}
\label{G(02):G(12)<C.2}
\frac{\Gamma(f_0(n),f_1(n),f_2(n))}{\Gamma(f_1(n),f_2(n))}
=\frac{\sum_{1\leq k<r<s\leq n}\vert M_{krs}^{123}\vert^2}
{\sum_{1\leq r<s\leq n}\vert M_{rs}^{12}\vert^2}\leq C
\end{equation}
for all $n\in \mathbb N$. Set now
\begin{eqnarray}
\nonumber
&&a_{2n}=\gamma(f_1(n),f_2(n),f_0(n)),\\
\label{A_3n.3}
&&A_{2n}=\gamma(f_1(n),f_2(n)),\quad A_{2n}t_0^{(n)}=b_{2n},\\
\label{b(n),t_0(n).3}
&&b_{2n}=(f_k(n),f_0(n))_{k=1}^2\in \mathbb R^2,\quad
t_0^{(n)}=(t_{0r}^{(n)})_{r=1}^2.
\end{eqnarray}
More explicitly
\begin{equation}
\label{a(mn).3}
a_{2n}\!=\!
\left(\!\!\begin{array}{ccc}
(f_1(n),f_1(n))&(f_1(n),f_2(n))&(f_1(n),f_0(n))\\
(f_2(n),f_1(n))&(f_2(n),f_2(n))&(f_2(n),f_0(n))\\
(f_0(n),f_1(n))&(f_0(n),f_2(n))&(f_0(n),f_0(n))
\end{array}\!\!\right),
\end{equation}
\begin{equation}
\label{A(3n).3}
A_{2n}\!=\!
\left(\!\!\begin{array}{cc}
(f_1(n),f_1(n))&(f_1(n),f_2(n))\\
(f_2(n),f_1(n))&(f_2(n),f_2(n))\\
\end{array}\!\!\right).
\end{equation}
If we replace the vectors $f_0,f_1,f_2$ with $f_0(n),f_1(n),f_2(n)$, 
the equality  \eqref{F_3(t)=(At,t).1} then becomes 
\begin{equation}
 \label{F_3n(t)=(At,t).3}
 F^{(n)}_{20}(t^{(n)}_0)=\min_{(t_1,t_2)\in \mathbb R^2}
 \Vert\sum_{k=1}^2t_kf_k(n)-f_0(n)\Vert ^2=
\frac{\Gamma(f_0(n),f_1(n),f_2(n))}{\Gamma(f_1(n),f_2(n))}.
\end{equation}
For $t\in \mathbb R^2$ and $f_0(n),f_1(n),f_2(n)\in \mathbb R^n$
and $0\leq s\leq 2$,  define the functions
\begin{equation}
 \label{F(2s)(t)}
F^{(n)}_{2s}(t)=\Vert \sum_{0\leq r\leq 2,r\not=s}t_rf_r(n)-f_s(n)\Vert ^2.
\end{equation}
The minimum of the corresponding 
expressions $F^{(n)}_{2s}(t)$ for $0\leq s \leq 2$ is attained respectively  at $t^{(n)}_0,\,t^{(n)}_1,\,t^{(n)}_2$. 
The proof of the fact that one of the sequences 
$
t^{(n)}_0,\,t^{(n)}_1,\,t^{(n)}_2
$ 
is bounded, is based on the positive definiteness 
of the matrices 
$\gamma(f_1(n),f_2(n),f_0(n))$, for the details see \cite{Kos-m-Arx23}. Let for example, the sequence $t^{(n)}_0\in \mathbb R^2$ is bounded. Therefore,  there exists a subsequence $(t^{(n_k)})_{k\in \mathbb N}$ that converges to some $t\in \mathbb R^2$. This contradics \eqref{<C.m=2}. Indeed
\begin{equation}
\label{Contr.m=2}
\lim_{n\to \infty}F^{(n)}_{20}(t)=\infty,\quad  F^{(n)}_{20}(t_0^{(n_k)})\leq C,\quad
\lim_{k\to \infty}t_0^{(n_k)}=t.
\end{equation}

\noindent{\it Acknowledgement.} 
\vskip 0.3 cm
A.~Kosyak is very grateful to Prof. K.-H. Neeb, Prof. M.~Smirnov and 
Dr Moree
for their personal efforts  to make academic  stays possible at their respective institutes. A.~Kosyak  visited: MPIM from March to April 2022 and from January to April 2023,
University of Augsburg from June to July 2022, and  University of Erlangen-Nuremberg 
from August to December 2022, all during the Russian invasion in Ukraine.
Also, Prof. R.~Kashaev  kindly invited A.~Kosyak to Geneva.

Further, both authors would like  to pay their respect to Prof. P. Teichner at MPIM, for his
immediate efforts started to help mathematicians
in Ukraine  after the Russian invasion.

Since the spring of 2023 A.~Kosyak is an Arnold Fellow at the
London Institute for Mathematical Sciences, and he would like to express
his gratitude to  Mrs Myers Cornaby   
and to Miss Ker Mercer 
and especially to the Director of LIMS Dr T.~Fink and Prof. Y.-H.~He.

P.~Moree, not an expert in this area, thanks 
A.~Kosyak for his patient clarifications of the big picture,  thus making him see the forest for the trees.

\end{document}